\documentclass[11pt,a4paper]{amsart}
\usepackage{amsfonts}
\usepackage{amsthm}
\usepackage{amsmath}
\usepackage{amscd}
\usepackage{amssymb}
\usepackage{subfigure}
\usepackage[latin2]{inputenc}
\usepackage{t1enc}
\usepackage[mathscr]{eucal}
\usepackage{indentfirst}
\usepackage{graphicx}
\usepackage{graphics}
\usepackage{pict2e}
\usepackage{color}
\usepackage{epic}
\usepackage{bbm}
\usepackage{enumitem}
\usepackage{tikz}
\usepackage{tkz-graph}
\usepackage{xcolor}
\usetikzlibrary{shapes,snakes}
\usetikzlibrary{decorations.pathreplacing}
\usepackage{bm}
\usetikzlibrary{arrows}
\numberwithin{equation}{section}
\usepackage[margin=2.6cm]{geometry}
\usepackage{epstopdf} 
\usepackage{textcomp} %fur mu in PM10
 \usepackage{autobreak}
 \usepackage{bigstrut}
\usepackage{colortbl}
\usepackage{rotating}
\usepackage{multirow}
\usepackage{natbib}
\bibpunct{(}{)}{,}{a}{}{;}
\renewcommand{\citep}[1]{[\citealp{#1}]}

\newtheorem{defi}{Definition}[section]

\newtheorem{ann}[defi]{Assumption}
\newtheorem{thm}[defi]{Theorem}
\newtheorem{rem}[defi]{Remark}

\newtheorem{cor}[defi]{Corollary}

\newcommand{\F}{F_{i,j}}

\newcommand{\C}{C_{i,j}}

\newcommand{\f}{f_j}
\newcommand{\fh}{\widehat{f}_j}
\newcommand{\s}{\sigma_j^2}

\numberwithin{equation}{section}

\begin{document}
%\begin{frontmatter}
\title{BOOTSTRAP CONSISTENCY FOR THE MACK BOOTSTRAP}

%\runtitle{Asymptotic Theory and Bootstrap Consistency of the CDR}

\author[Julia Steinmetz \& Carsten Jentsch]{Julia Steinmetz$^*$   \\  \texttt{Department of Statistics,  TU Dortmund University, D-44221 Dortmund, Germany; steinmetz@statistik.tu-dortmund.de}\vspace*{0.5cm} \\  \vspace*{0.5cm}
\and  \\ 
Carsten Jentsch \\ \texttt{ Department of Statistics,  TU Dortmund University, D-44221 Dortmund, Germany; jentsch@statistik.tu-dortmund.de \\  
}}

\begin{abstract}
Mack's distribution-free chain ladder reserving model belongs to the most popular approaches in non-life insurance mathematics. %While p
Proposed to determine the first two moments of the reserve, it does not allow to identify the whole distribution of the reserve. For this purpose, Mack's model is usually equipped with a tailor-made bootstrap procedure.
%, which combines a residual-based non-parametric resampling step together with a parametric bootstrap. 
Although
 %it is 
widely used in practice
%applications 
to estimate the reserve risk, no theoretical bootstrap consistency results exist that justify this approach.

%In this paper, t
To fill this gap in the literature, we adopt the %theoretical 
framework proposed by \cite{jentschasympThInf} to derive asymptotic theory in Mack's model. By splitting the reserve into two 
%additive 
parts corresponding to process and estimation uncertainty, %it 
this enables - for the first time - a rigorous investigation also of the validity of the Mack bootstrap. We prove that the (conditional) distribution of the asymptotically dominating process uncertainty part is correctly mimicked by Mack's bootstrap if the parametric family of distributions of the individual development factors is correctly specified.
 %in Mack's bootstrap proposal. 
Otherwise, this is
%will be generally 
not the case. In contrast, the 
%corresponding 
(conditional) distribution of the estimation uncertainty part is generally not correctly captured by Mack's bootstrap. To tackle this, we propose an alternative Mack-type bootstrap, which is designed to capture also the distribution of the estimation uncertainty part.

We illustrate our findings by simulations and show that the newly proposed alternative Mack
%-type 
bootstrap performs superior to the 
%original 
Mack bootstrap.
 %in finite samples.
\end{abstract}


\maketitle

\section{Introduction}\label{sec:intro}

In a non-life insurance business an insurer needs to build up a reserve to be able to meet future obligations arising from incurred claims. The actual sizes of the claims are unknown at the time the reserves have to be built, since the claims are incurred, but either not been reported yet or they have been reported, but not settled yet. This process of forecasting of outstanding claims is called \textit{reserving}. An accurate estimation of the outstanding claims is crucial for pricing future policies and for the assessment of the solvency of the insurer. A popular and	widely used technique in practice to forecast future claims is the Chain Ladder Model (CLM), which provides an algorithm to predict future claims.
 %and to construct the (best) estimate of the reserve using a set of development factors and variance parameters. 
In this respect, the most popular model is the recursive model proposed by \cite{mack1993distribution}, which extends the CLM by allowing also the calculation of the standard deviation of the reserve.
	
%Also, 
Alternatively, frameworks based on general linear models (GLMs) considered e.g.~in \cite{renshaw1998stochastic} %, which 
make %e.g.~
use of over-dispersed Poisson and Log-normal distributions for modeling mean and variance of the reserve.
%, have been proposed for the calculation of the first two moments of the reserve. 
%%%For incremental claims data, \cite{harnau2018over} and \cite{kuang2019generalized}  developed an asymptotic theory for an over-dispersed Poisson and a Log-normal framework, respectively.Both models are based on the assumption that the incremental claims belong to infinitely divisible distributions. They derive central limit theorems (CLTs) for the estimators and proposed a t-distribution for the centered and standardized reserve. 
However, such parametric assumptions are often restrictive and the knowledge of the first two moments of the reserve is not satisfactory for actuaries
%, since it does not allow 
to draw sufficient conclusions about the reserve risk and the solvency of the insurance company. The \emph{reserve risk} is defined as the risk that the economic-valued reserve 
%(best estimate) of the reserve 
does not suffice to pay for all outstanding claims, which inevitably requires
%. To get such insights, 
the knowledge %about the \emph{whole} distribution or at least 
of high quantiles of the reserve.
 %is inevitable. 
For this purpose, \cite{england2006predictive} proposed the Mack bootstrap which equips Mack's model with a tailor-made bootstrap procedure.
 %that makes use of the recursive structure of CLMs and relies on additional parametric assumptions. 
Alternative bootstrap procedures for GLM-based setups have been addressed also in \cite{england1999analytic, england2006predictive, england2002addendum} and \cite{pinheiro2003bootstrap}. Without providing any consistency results, \cite{Bjoerkwall2009} review these bootstrap techniques and suggest alternative non-parametric and parametric bootstrap procedures.
 %without providing any consistency results. 
Similarly, \cite{bjorkwall2010bootstrapping} suggest %also a 
bootstrap techniques for the separation method, %which 
that takes calendar year effects into account. In recent years, bootstrap-based approaches have been favored by many actuaries, because such methods usually produce plausible distributions 
%that appear to be plausible 
in practice. However, as demonstrated by \cite{GRIRO2007} and \citet{GRIRO2008}, Mack's model and GLM-type models in combination with the bootstrap do not produce satisfactory results in certain situations. 
%However, the existing literature lacks a theoretical framework and rigorous asymptotic results that would allow to identify such scenarios a priori. 
In this regard, %more 
refined approaches have been proposed to improve the finite sample performance. For example, \cite{verdonck2011influence} investigate the influence of outliers for the parameter estimation in the GLM framework and calculate its leverage on the CLM. %Alternatively, 
\cite{hartl2010bootstrapping} propose to use deviance residuals instead of Pearson residuals for the GLM framework. \cite{tee2017comparison} provide an extensive case study for bootstrapping the GLM using a (over-dispersed) Poisson model, the Gamma model and the Log-normal model in combination with different residual types.  
\cite{peremans2017robust} propose a %more 
robust bootstrap procedure in a GLM setting based on M-estimators using influence functions. % and suggest to use weighted bootstrap resampling. 
\cite{peters2010chain} compare the Mack bootstrap with a Bayesian bootstrap. % and show that their Bayesian approach gives the same results.

Nevertheless, already for the original Mack bootstrap method, the existing literature lacks a deeper and mathematically rigorous understanding. For this purpose, it is desirable to provide a suitable theoretical framework to be able to justify the application of the Mack bootstrap.
%Nevertheless, a deeper and mathematically rigorous understanding of the original Mack bootstrap method and its underlying stochastic model is desirable to be able to justify the application of the Mack bootstrap.
Only recently, 
%to enable a rigorous asymptotic treatment of Mack's model, 
\citet{jentschasympThInf} proposed a suitable theoretical (stochastic and asymptotic) framework, which allows the derivation of conditional %as well as
and unconditional asymptotic 
%distributions of 
theory for the reserve in Mack's model. %Precisely, 
They split the reserve (centered around its best estimate) into two %random additive 
parts, that carry the \emph{process uncertainty} and the \emph{estimation uncertainty}, respectively. This allows to derive unconditional limiting distributions for both parts of the reserve, and when conditioning on the latest observed cumulative claims %, but also unconditionally. 
%In this regard, when addressing the question of bootstrap consistency for 
As risk reserving
%, which 
is generally a prediction task, these \emph{conditional} limiting distributions %are crucial and 
serve well as benchmarks for the corresponding Mack bootstrap distributions, when addressing the question of bootstrap consistency. %Whereas 
While the conditional limiting distribution of 
%the second part, which corresponds to 
the estimation uncertainty part %will 
turns out to be Gaussian under mild regularity conditions and when properly inflated, the conditional limiting distribution of the 
%first part corresponding to the 
process uncertainty part will be generally non-Gaussian. Considering both parts jointly, the process uncertainty part dominates asymptotically, which leads to a non-Gaussian limiting distribution of the reserve in total.

In this paper, we adopt the %stochastic and asymptotic 
theoretical framework introduced in \citet{jentschasympThInf} to %rigorously
investigate the long-standing question of Mack bootstrap consistency. Our contributions are twofold. First,
%to fill this gap in the literature, 
we derive bootstrap asymptotic theory for both parts of the (centered) Mack bootstrap reserve corresponding to process uncertainty and estimation uncertainty, respectively. We prove that the (conditional) bootstrap distribution of the asymptotically dominating process uncertainty part is correctly mimicked if the parametric family of distributions of the Mack bootstrap individual development factors is correctly specified.
 %in Mack's bootstrap proposal. 
Otherwise, this will be generally not the case. In contrast, the %corresponding 
(conditional) 
%bootstrap 
distribution of the estimation uncertainty part is generally not correctly captured. Second, %motivated 
inspired from our asymptotic findings, we propose an alternative Mack-type bootstrap, which is designed to capture also the distribution of the estimation uncertainty part.

The paper is organized as follows. Section \ref{CLM_section} introduces the required notation and assumptions for the CLM, discusses parameter estimation in Mack's model, and provides the asymptotic and stochastic framework %as introduced 
of \cite{jentschasympThInf}. In Section \ref{sec_boot}, we %introduce
discuss the Mack bootstrap approach as %originally 
proposed by \cite{england2006predictive}. % and discuss its construction.
In Section \ref{sec_boot_consistency}, we summarize the (conditional) asymptotic results from \cite{jentschasympThInf} for the process uncertainty and estimation uncertainty terms in Section \ref{sec_asymp_summary}, which will serve as benchmarks %results 
for the Mack bootstrap results. Then, in Section \ref{sec_asymp_boot}, we derive bootstrap asymptotic theory for both parts of the (centered) Mack bootstrap reserve corresponding to process uncertainty and estimation uncertainty, respectively. 
%For this purpose, asymptotic bootstrap theory for (smooth functions of) bootstrap development factor estimators has to be established, which might be of independent interest and can be found in the appendix. 
Based on these results, 
%based on the findings from the previous section, 
we propose an alternative Mack-type bootstrap in Section \ref{sec_alternative_boot} and derive its asymptotic properties in Section \ref{sec_alternative_boot_consistency}. We illustrate our findings in simulations in Section \ref{sec_simulation} and show that the newly proposed alternative Mack-type bootstrap performs superior to the original Mack bootstrap in finite samples. Section \ref{sec_conclusion} concludes. All proofs, auxiliary results and additional simulations are deferred to %supplementary material \cite{Mack_boot_supplement}. 
the appendix.

%\bigskip

\section{The Chain Ladder Model}\label{CLM_section}\

Reserves are the major part of the balance sheet for non-life insurance companies such that their accurate prediction is crucial. For this purpose, insurers summarize all observed claims of a business line in a loss triangle (upper-left triangle in Table \ref{loss_triangle}). Its entries, the cumulative amount of claims $\C$, are sorted by their years of accident $i$ (vertical axis) and their years of occurrence $j$ (horizontal axis), where $i,j=0,\ldots,I$ with $i+j\leq I$.	Hence, the (observed) loss triangle contains all cumulative claims $\C$ that have already been observed up to calendar year $I$. It constitutes the available data basis and is denoted by
\begin{align}\label{def_D}
\mathcal{D}_I=\{C_{i,j}|i,j=0, \dots, I,\ 0\leq i+j\leq I\}.
\end{align}
The total aggregated amount of claims of the same calendar year $k$ with $k=0,\ldots,I$ are lying on the same diagonal (from lower-left to upper-right corner) of the loss triangle. We denote these diagonals by
%\begin{align}\label{def_Q}
$\mathcal{Q}_k=\{C_{k-i,i}|i=0,\dots, k\}$.
%\end{align}
In this setup, $I$ is the current calendar year corresponding to the most recent accident year and development period such that the diagonal $\mathcal{Q}_I$ (orange diagonal in Table \ref{loss_triangle}) summarizes the latest cumulative claim amounts collected in year $I$.

\begin{table}
	\centering
	%\vspace{1cm}
	%\hspace{-0.5cm}
	%\begin{minipage}[b]{\textwidth}\centering
		\begin{tiny}
			%	\begin{table}[t]
			%		\centering
			\begin{tabular}{|p{0.1cm}|p{0.5cm}p{0.5cm}p{0.5cm}p{0.5cm}p{0.5cm}p{0.5cm}p{0.5cm}p{0.5cm}p{0.5cm}}
		\hline
		& \multicolumn{1}{r}{} & \multicolumn{8}{c|}{Development Year $j$}\\ % \bigstrut[t]\\
		\cline{2-10} 	& \multicolumn{1}{c}{} & \multicolumn{1}{c}{$0\vphantom{I^2}$} & \multicolumn{1}{c}{$1$} & $2$ &  & $\cdots$ &  &\multicolumn{1}{c}{ $I-1$} &\multicolumn{1}{c|}{$\quad \;I\;\quad $}  \\ %\bigstrut[b]\\
		\cline{3-10}    \multirow{8}[20]{*}{\begin{sideways}\hspace{2cm } 
				Accident Year $i$ \end{sideways}}          &  \multicolumn{1}{c|}{$0$} & \multicolumn{1}{c}{$C_{0,0}\vphantom{I^2}$} & \multicolumn{1}{c}{$C_{0,1}$} & &  &  &\multicolumn{1}{c}{ } & \multicolumn{1}{c}{$C_{0,I-1}$} & 
		\multicolumn{1}{c|}{\cellcolor[rgb]{.957,  .69,  .518}$C_{0,I}$}   \\
		\cline{10-10}          &  \multicolumn{1}{c|}{ $1$}& 
		\multicolumn{1}{c}{$C_{1,0}$} & \multicolumn{1}{c}{$\cdots$ }&  &  &  &  & \multicolumn{1}{c|}{\cellcolor[rgb]{ .957,  .69,  
				.518}$C_{1,I-1}$}  &   	\multicolumn{1}{c|}{\cellcolor[rgb]{ .663,  .816,  .557}$C_{1,I}$}     \\
		\cline{9-9}          &  \multicolumn{1}{c|}{} & 
		\multicolumn{1}{c}{$\cdots$} & \multicolumn{1}{c}{} & $\cdots$ &  & & 
		\multicolumn{1}{c|}{\cellcolor[rgb]{ .957,  .69,  .518}} &   	\multicolumn{1}{c}{\cellcolor[rgb]{ .663,  .816,  .557}}     &    	\multicolumn{1}{c|}{\cellcolor[rgb]{ .663,  .816,  .557}}    
		\\
		\cline{8-8}          &  \multicolumn{1}{c|}{ $\cdots$} & 
		\multicolumn{1}{c}{$\cdots$} & \multicolumn{1}{c}{} &  & \multicolumn{1}{c}{} & 
		\multicolumn{1}{c|}{\cellcolor[rgb]{ .957,  .69,  .518}} &   	\multicolumn{1}{c}{\cellcolor[rgb]{ .663,  .816,  .557}}     &     	\multicolumn{1}{c}{\cellcolor[rgb]{ .663,  .816,  .557}}   
		&   	\multicolumn{1}{c|}{\cellcolor[rgb]{ .663,  .816,  .557}$\cdots$}       \\
		\cline{7-7}          &  \multicolumn{1}{c|}{} & 
		\multicolumn{1}{c}{} & \multicolumn{1}{c}{} & \multicolumn{1}{c}{} & 
		\multicolumn{1}{c|}{\cellcolor[rgb]{ .957,  .69,  .518}} &  	\multicolumn{1}{c}{\cellcolor[rgb]{ .663,  .816,  .557}$\cdots$}      &       	\multicolumn{1}{c}{\cellcolor[rgb]{ .663,  .816,  .557}} 
		&  	\multicolumn{1}{c}{\cellcolor[rgb]{ .663,  .816,  .557}}      &    	\multicolumn{1}{c|}{\cellcolor[rgb]{ .663,  .816,  .557}}     \\
		\cline{6-6}          &  \multicolumn{1}{c|}{} & 
		\multicolumn{1}{c}{} &\multicolumn{1}{c}{} & 
		\multicolumn{1}{c|}{\cellcolor[rgb]{ .957,  .69,  .518}} &  	\multicolumn{1}{c}{\cellcolor[rgb]{ .663,  .816,  .557}}      &	\multicolumn{1}{c}{\cellcolor[rgb]{ .663,  .816,  .557}}        
		&     	\multicolumn{1}{c}{\cellcolor[rgb]{ .663,  .816,  .557}}   & 	\multicolumn{1}{c}{\cellcolor[rgb]{ .663,  .816,  .557}}       & 	\multicolumn{1}{c|}{\cellcolor[rgb]{ .663,  .816,  .557}}          \\
		\cline{5-5}          &  \multicolumn{1}{c|}{$I-1$}& 
		\multicolumn{1}{c}{$C_{I-1,0}$} & 
		\multicolumn{1}{c|}{\cellcolor[rgb]{ .957,  .69,  .518} $C_{I-1,1}$} &   	\multicolumn{1}{c}{\cellcolor[rgb]{ .663,  .816,  .557}} 
		&    	\multicolumn{1}{c}{\cellcolor[rgb]{ .663,  .816,  .557}}    & 	\multicolumn{1}{c}{\cellcolor[rgb]{ .663,  .816,  .557}$\cdots$}       & 	\multicolumn{1}{c}{\cellcolor[rgb]{ .663,  .816,  .557}}       & 	\multicolumn{1}{c}{\cellcolor[rgb]{ .663,  .816,  .557}}       & 	\multicolumn{1}{c|}{\cellcolor[rgb]{ .663,  .816,  .557}}          \\	
		\cline{4-4}          & \multicolumn{1}{c|}{$I$} & 
		\multicolumn{1}{c|}{\cellcolor[rgb]{ .957,  .69,  .518}$C_{I,0}$} &   	\multicolumn{1}{c}{\cellcolor[rgb]{ .663,  .816,  .557}$C_{I,1}$}    
		&    	\multicolumn{1}{c}{\cellcolor[rgb]{ .663,  .816,  .557}}    & 	\multicolumn{1}{c}{\cellcolor[rgb]{ .663,  .816,  .557}}       & 	\multicolumn{1}{c}{\cellcolor[rgb]{ .663,  .816,  .557}}       &  	\multicolumn{1}{c}{\cellcolor[rgb]{ .663,  .816,  .557}}      & 	\multicolumn{1}{c}{\cellcolor[rgb]{ .663,  .816,  .557}}       &    	\multicolumn{1}{c|}{\cellcolor[rgb]{ .663,  .816,  .557} $C_{I,I}$}      \\
		\cline{1-3}
		\hline
		\multicolumn{5}{c}{}	&       &       &       &             &  \\
	\end{tabular}	
			%	\end{table}
		\end{tiny}
	%\end{minipage}
	%	}
	%}
	\caption{Observed upper loss triangle $\mathcal{D}_I$ (upper-left triangle; white and orange) with accident years (vertical axis), development years (horizontal axis), diagonal $\mathcal{Q}_I$ (orange), and unobserved lower loss triangle $\mathcal{D}_I^c$ (lower-right triangle; green).}
	\label{loss_triangle}
\end{table}

For the theoretical analysis of the prediction of the outstanding (unobserved) claims, it is useful to augment the (observed) \textit{upper} loss triangle $\mathcal{D}_I$ by an unobserved \textit{lower} triangle
\begin{align*}
\mathcal{D}_I^c=\{C_{i,j}|i,j=0, \dots, I,\ i+j>I\}
\end{align*}
that contains all future claims that have \emph{not} been observed (yet) up to time $I$ (green triangle in Table \ref{loss_triangle}). The resulting cumulative claim matrix is denoted by $\mathcal{C}_I=(C_{i,j})_{i,j=0,...,I}=\mathcal{D}_I\cup\mathcal{D}_I^c$. For each accident year $i$, the main interest lies in the \textit{reserves} at terminal time $I$, denoted by $R_{i,I}$, which is computed by taking the difference of the ultimate claim $C_{i,I}$ (last column), which is not observed (for $i>0$) at time $I$, minus the latest observed claim $C_{i,I-i}$ (on the diagonal) at time $I$. Precisely, we define the reserve for accident year $i$ by $R_{i,I}=C_{i,I}-C_{i,I-i}$ for $i=0,\dots, I$ and the aggregated total amount of the reserve $R_I$ by
\begin{equation}\label{eqtotal_R}
R_I=\sum^{I}_{i=0}R_{i,I},
\end{equation}
noting that $R_{0,I}=C_{0,I}-C_{0,I}=0$ by construction. Hence, for each accident year $i$ and being in calendar year $I$, to get an estimate of $R_{i,I}$, we have to predict the unobserved ultimate claim $C_{i,I}$. % based on the observed upper triangle $\mathcal{D}_I$. 
Starting from %the last observed claim 
$C_{i,I-i}$, this is done by predicting sequentially all future, yet (at time $I$) unobserved claims $\{C_{i,j}|j=I-i+1, \dots, I\}$. %That is, 
By doing this for all %accident years 
$i=0, \dots, I$, the whole unobserved lower loss triangle $\mathcal{D}_I^c$ has to be predicted, and by summing-up all predictions for $R_{i,I}$, we get a prediction also for $R_I$.

However, to make the CLM setup above accessible for the derivation of %meaningful limiting 
asymptotic theory for predictive inference, \citet{jentschasympThInf} introduced a suitable stochastic and asymptotic framework for Mack's model, which is adopted here as well and will be described in the following.

\begin{table}
	\centering
	%\begin{minipage}[b]{\textwidth}\centering
		\begin{tiny}
			%\begin{table}[t]			
			\centering
			\begin{tabular}
			%{|p{0.1cm}|p{0.7cm}p{0.7cm}p{0.7cm}p{0.7cm}p{0.7cm}p{0.7cm}p{0.7cm}p{0.7cm}p{0.7cm}p{0.7cm}p{0.7cm}}
			{|p{0.1cm}|p{1cm}p{1cm}p{1cm}p{1cm}p{1cm}p{1cm}p{1cm}p{1cm}p{1cm}p{1cm}p{1cm}}
				\hline
				&  & \multicolumn{10}{c|}{Development Year $j$}\\ % \bigstrut[t]\\
					\cline{2-12} 
				& \multicolumn{1}{c}{} & \multicolumn{1}{c}{$0\vphantom{I^2}$} & \multicolumn{1}{c}{$1$} & $2$ &  & $\cdots$ &  &  & $\quad \;I\;\quad $ &\multicolumn{1}{c}{$I+1$} & \multicolumn{1}{c|}{$I+2$}\\ %\bigstrut[b]\\
			\cline{3-12}    \multirow{10}[20]{*}{\begin{sideways}\hspace{2cm } Accident Year $i$ \end{sideways}} & \multicolumn{1}{c|}{$0$} & \multicolumn{1}{c}{$C_{0,0}$} & \multicolumn{1}{c}{$C_{0,1}\vphantom{I^2}$} &  &  & $\cdots$ &  &  &  $\;\;C_{0,I}\;\;\quad $&\multicolumn{1}{c}{ \cellcolor[rgb]{ .957,  .69,  .518} $C_{0,I+1}$ }& \multicolumn{1}{c|}{\cellcolor[rgb]{ .663,  .816,  .557}$C_{0,I+2}$} \\
				\cline{12-12}          & \multicolumn{1}{c|}{ $1$ }& \multicolumn{1}{c}{$C_{1,0}$} & &  & & $\cdots$ &  &  & \cellcolor[rgb]{ .957,  .69,  .518} & \multicolumn{1}{c|}{\cellcolor[rgb]{ .663,  .816,  .557} $C_{1,I+1}$} & \\
				\cline{11-11}          &\multicolumn{1}{c|}{ $2$} & \multicolumn{1}{c}{} &  &  &  &  &  & \cellcolor[rgb]{ .957,  .69,  .518} & \multicolumn{1}{c|}{\cellcolor[rgb]{ .663,  .816,  .557}} &       &  \\
				\cline{10-10}          & \multicolumn{1}{c|}{} & \multicolumn{1}{c}{} & $\cdots$ &  &  &  & \cellcolor[rgb]{ .957,  .69,  .518} & \multicolumn{1}{c|}{\cellcolor[rgb]{ .663,  .816,  .557}} &       &       &  \\
				\cline{9-9}          &  \multicolumn{1}{c|}{} & \multicolumn{1}{c}{$\cdots$} &  & & \multicolumn{1}{c}{$\cdots$}  & \cellcolor[rgb]{ .957,  .69,  .518} & \multicolumn{1}{c|}{\cellcolor[rgb]{ .663,  .816,  .557}} &       &       &       & \\
				\cline{8-8}          &  \multicolumn{1}{c|}{$\cdots$}  & \multicolumn{1}{c}{$\cdots$} &  &  & \cellcolor[rgb]{ .957,  .69,  .518} & \multicolumn{1}{c|}{\cellcolor[rgb]{ .663,  .816,  .557}} &       &       &       &       &  \\
				\cline{7-7}          &   \multicolumn{1}{c|}{ $I-1$} & $C_{I-1,0}$& $C_{I-1,1}$ & \cellcolor[rgb]{ .957,  .69,  .518}$C_{I-1,2}$ & \multicolumn{1}{c|}{\cellcolor[rgb]{ .663,  .816,  .557} } &       &       &       &       &       & \\
				\cline{6-6}          &  \multicolumn{1}{c|}{$I$} & \multicolumn{1}{c}{$C_{I,0}$} & \cellcolor[rgb]{ .957,  .69,  .518}$C_{I,1}$ & \multicolumn{1}{c|}{\cellcolor[rgb]{ .663,  .816,  .557}$C_{I,2}$} &       &       &       &       &       &       &  \\
				\cline{5-5}          &  \multicolumn{1}{c|}{$I+1$} & \multicolumn{1}{c}{\cellcolor[rgb]{ .957,  .69,  .518}$C_{I+1,0}$} & \multicolumn{1}{c|}{\cellcolor[rgb]{ .663,  .816,  .557}$C_{I+1,1}$} &       &       &       &       &       &       &       &  \\
				\cline{4-4}          &   \multicolumn{1}{c|}{$I+2$} & \multicolumn{1}{c
					|}{\cellcolor[rgb]{ .663,  .816,  .557}$C_{I+2,0}$} &       &       &       &       &       &       &       &       &  \\
				\cline{1-3} 
			\end{tabular}	
		\end{tiny}
	%\end{minipage}
	
	\vspace{0.5cm}
	%\hspace{-0.5cm}
	%\begin{minipage}[b]{\textwidth}\centering
		\begin{tiny}
			%	\begin{table}[t]
			%		\centering
			\begin{tabular}%{|p{0.1cm}|p{0.8cm}p{0.8cm}p{0.8cm}p{0.7cm}p{0.7cm}p{0.7cm}p{0.7cm}p{0.7cm}p{0.7cm}p{0.7cm}p{0.7cm}}
			{|p{0.1cm}|p{1cm}p{1cm}p{1cm}p{1cm}p{1cm}p{1cm}p{1cm}p{1cm}p{1cm}p{1cm}p{1cm}}
				\hline
				& \multicolumn{1}{r}{} & \multicolumn{10}{c|}{Development Year $j$}\\ % \bigstrut[t]\\
%		&  & \multicolumn{1}{c}{$0$} & $1$ & $2$ &  & $\cdots$ & & & $\quad \;I\;\quad $ &\multicolumn{1}{c}{\cellcolor[rgb]{  .957,  .69,  .518} $I+1$} & \multicolumn{1}{c|}{\cellcolor[rgb]{ .663,  .816,  .557} $I+2$}}\\
	\cline{2-12} &&\multicolumn{1}{c}{$0\vphantom{I^2}$} &\multicolumn{1}{c}{$1$} &\multicolumn{1}{c}{$2$} &&&$\cdots$&&\multicolumn{1}{c}{$I$} &\multicolumn{1}{c}{$I+1$} & \multicolumn{1}{c|}{$I+2$} \\
%				%\hline  %\bigstrut[b]\\
			%	\cline{2-12} 
			\multirow{10}[20]{*}{\begin{sideways}\hspace{2cm } Accident Year $i$ \end{sideways}} & \multicolumn{1}{c|}{$-2$}  &  
  \multicolumn{1}{c}{\cellcolor[rgb]{ .663,  .816,  .557}$C_{-2,0}\vphantom{I^2}$} & \multicolumn{1}{c}{\cellcolor[rgb]{ .663,  .816,  .557}} & \cellcolor[rgb]{ .663,  .816,  .557} & \cellcolor[rgb]{ .663,  .816,  .557} & \cellcolor[rgb]{ .663,  .816,  .557}$\cdots$ & \cellcolor[rgb]{ .663,  .816,  .557} &  \multicolumn{1}{c}{\cellcolor[rgb]{ .663,  .816,  .557}} & \multicolumn{1}{c}{ \cellcolor[rgb]{ .663,  .816,  .557}} &\multicolumn{1}{c}{ \cellcolor[rgb]{ .663,  .816,  .557} $C_{-2,I+1}$ }& \multicolumn{1}{c|}{\cellcolor[rgb]{ .663,  .816,  .557}$C_{-2,I+2}$} \\
				\cline{12-12}          &  \multicolumn{1}{c|}{ $-1$} & \multicolumn{1}{c}{\cellcolor[rgb]{  .957,  .69,  .518}$C_{-1,0}$} &\multicolumn{1}{c}{ \cellcolor[rgb]{ .957,  .69,  .518}} & \cellcolor[rgb]{  .957,  .69,  .518} & \cellcolor[rgb]{  .957,  .69,  .518} & \cellcolor[rgb]{  .957,  .69,  .518}$\cdots$ & \cellcolor[rgb]{  .957,  .69,  .518} &  \multicolumn{1}{c}{\cellcolor[rgb]{  .957,  .69,  .518}} &  \multicolumn{1}{c}{\cellcolor[rgb]{ .957,  .69,  .518}} & \multicolumn{1}{c|}{\cellcolor[rgb]{ .957,  .69,  .518} $C_{-1,I+1}$} & \\
				\cline{11-11}          &  \multicolumn{1}{c|}{$0$}& \multicolumn{1}{c}{$C_{0,0}$} & \multicolumn{1}{c}{$C_{0,1}$} &  & &$\cdots$  &  & & \multicolumn{1}{c|}{$C_{0,I}$} &       &  \\
				\cline{10-10}          &  \multicolumn{1}{c|}{ $1$} & \multicolumn{1}{c}{$C_{1,0}$} & $\cdots$ &  &  &  &  & \multicolumn{1}{c|}{$C_{1,I-1}$} &       &       &  \\
				\cline{9-9}          &  \multicolumn{1}{c|}{} & \multicolumn{1}{c}{$\cdots$} &  & $\cdots$ &  &  & \multicolumn{1}{c|}{$\cdots$} &       &       &       & \\
				\cline{8-8}          &  \multicolumn{1}{c|}{$\cdots$} & \multicolumn{1}{c}{$\cdots$} &  &  &  & \multicolumn{1}{c|}{} &       &       &       &       &  \\
				\cline{7-7}          & \multicolumn{1}{c|}{} & \multicolumn{1}{c}{} & & & \multicolumn{1}{c|}{} &       &       &       &       &       & \\
				\cline{6-6}          & \multicolumn{1}{c|}{} & \multicolumn{1}{c}{} &  & \multicolumn{1}{c|}{} &       &       &       &       &       &       &  \\
				\cline{5-5}          & \multicolumn{1}{c|}{$I-1$} & \multicolumn{1}{c}{$C_{I-1,0}$} & \multicolumn{1}{c|}{$C_{I-1,1}$} &       &       &       &       &       &       &       &  \\
				\cline{4-4}          & \multicolumn{1}{c|}{$I$} & \multicolumn{1}{c|}{$C_{I,0}$} &       &       &       &       &       &       &       &       &  \\
				\cline{1-3} 
				\multicolumn{5}{c}{}	&       &       &       &             &  \\
			\end{tabular}		
			%	\end{table}
		\end{tiny}
	%\end{minipage}
	%	}
	%}
	\caption{Two asymptotic frameworks of growing loss triangles based on adding diagonals (upper panel) and by adding rows (lower panel). Both approaches lead to loss triangles that are equal in distributions (adapted from \cite{jentschasympThInf}).}
	\label{tabelI_I1}
\end{table}

\subsection{Asymptotic framework for reserve prediction}\label{sec_asymptotic_framework}\

With the loss triangle $\mathcal{D}_I$ at hand, %a conditional 
an asymptotic analysis %given 
conditional on the diagonal $\mathcal{Q}_I$, which contains the most up-to-date information in the loss triangle, is of much interest for insurers. However, for this purpose, we will not rely on a seemingly ''natural`` asymptotic frameworkbased on $I\rightarrow \infty$, where increasing $I$ means adding new \emph{diagonals} $\mathcal{Q}_{I+h}=\{C_{I-i,i}|i=0, \dots, I+h\}$, $h\geq 1$ to the loss triangle $\mathcal{D}_I$ (see Table \ref{tabelI_I1}, upper panel). Instead, as common in predictive inference (see e.g.~\citet{PaparoditisShang2021}), we employ a different asymptotic framework throughout this paper. That is, we keep the \emph{latest} cumulative claims in $\mathcal{D}_I$, that is, $\mathcal{Q}_I$, fixed and let $\mathcal{D}_I$ grow by adding new \emph{rows} of cumulative claims $\{C_{-h,i}|i=0, \dots, I+h\}$, $h\geq 1$ (see Table \ref{tabelI_I1}, lower panel). Nevertheless, both versions of differently growing loss triangles %indicated 
displayed in Table \ref{tabelI_I1} are equal in distribution.
%Hence, 
In what follows, all asymptotic results are derived under the framework that %we observe 
a sequence of (upper) loss triangles
\begin{align}\label{def_D_n}
	\mathcal{D}_{I,n}=\left\{C_{i,j}|i=-n,\ldots,I,~j=0,\ldots,I+n,~-n\leq i+j\leq I\right\},	\quad n\in\mathbb{N}_0=\{0,1,2,\ldots\},
\end{align}
is observed, where
%where $\mathbb{N}_0=\{0,1,2,\ldots\}$, with (main) diagonals 
\begin{align}\label{def_Q_n}
\mathcal{Q}_{I,n}=\left\{C_{I-i,i}|i=0,\ldots,I+n\right\}, \quad	n\in\mathbb{N}_0,
\end{align}
denote the corresponding diagonals. Note that $\mathcal{D}_{I,0}=\mathcal{D}_I$, $\mathcal{Q}_{I,0}=\mathcal{Q}_I$ and that $\mathcal{D}_{I,n}$ (and $\mathcal{Q}_{I,n}$) is obtained by sequentially adding $n$ rows of lengths $I+2,I+3,\ldots, I+n+1$, respectively, on top to $\mathcal{D}_I$ (see Table \ref{tabelI_I1}, lower panel). As before, for all $n\in\mathbb{N}_0$, we augment the (observed) \textit{upper} loss triangle $\mathcal{D}_{I,n}$ by an unobserved \textit{lower} triangle $\mathcal{D}_{I,n}^c=\{C_{i,j}|i=-n,\ldots,I,~j=0,\ldots,I+n,~i+j> I\}$ that contains all future claims that have not been observed (yet) up to time $I$. Further, according to \eqref{eqtotal_R}, the aggregated total amount of the reserve is denoted by
\begin{equation}\label{eqtotal_Rn}
R_{I,n}=\sum^{I}_{i=-n}R_{i,I+n},	\quad	n\in\mathbb{N}_0,
\end{equation}
where $R_{i,I+n}=C_{i,I+n}-C_{i,I-i}$, $n\in\mathbb{N}_0$ and $R_{-n,I+n}=C_{-n,I+n}-C_{-n,I+n}=0$ by construction. 

While we keep $I$ and $n$ fixed in the expositions of the remainder of this section and of Section \ref{sec_boot}, we let $n\rightarrow \infty$ to derive the limiting distribution of the reserve in Section \ref{sec_boot_consistency}. 
%According to \eqref{def_D_n} and \eqref{def_Q_n}, the limiting upper loss triangle $\mathcal{D}_{I,\infty}$ and the limiting diagonal $\mathcal{Q}_{I,\infty}$ are defined by
%\begin{align}\label{def_D_Q_infty}
%\mathcal{D}_{I,\infty}=\left\{C_{i,j}|i\in\mathbb{Z},i\leq I,~j\in\mathbb{N}_0,~i+j\leq I\right\}	\quad	\text{and}	\quad	\mathcal{Q}_{I,\infty}=\left\{C_{I-i,i}|i\in\mathbb{N}_0\right\},
%\end{align}
%respectively.

\subsection{Mack's distribution-free chain ladder reserving}\label{paraestMack_subsec}\

%The distribution-free chain ladder model proposed by \citet{mack1993distribution}, often denoted simply as Mack's Model, is widely used in practice to determine both the mean and the variance of the reserve. 
By adopting the notion of the asymptotic framework described in Section \ref{sec_asymptotic_framework}, the conditions of Mack's Model originally proposed in \citet{mack1993distribution} can be summarized as follows.
	
\begin{ann}[Mack's Model]\label{MackM}
For any $n\in\mathbb{N}_0$, let $\mathcal{C}_{I,n}=(C_{i,j},i=-n,\ldots,I,~j=0,\ldots,I+n)$ denote random variables on some probability space $(\Omega, \mathcal{A}, P)$ and suppose the following holds:
\begin{itemize}
	\item[(i)] There exist so-called development factors $f_0,\dots, f_{I+n-1}$ such that
\begin{align}\label{MackM_mean}
E(C_{i,j+1}| C_{i,j})=f_{j}C_{i,j},	\quad	i=-n, \ldots, I,\ j=0,\ldots, I+n-1.
\end{align}
%for $ i=-n, \ldots, I$ and $j=0,\ldots, I+n-1$.
	\item[(ii)]  There exist variance parameters $\sigma^2_0,\dots, \sigma^2_{I+n-1}$ such that
\begin{align}\label{MackM_var}
Var(C_{i,j+1}|C_{i,j} )=\sigma_{j}^2C_{i,j},	\quad	i=-n, \ldots, I,\ j=0,\ldots, I+n-1.
\end{align}
%for $i=-n, \ldots, I$ and $j=0, \ldots, I+n-1$.
	\item[(iii)] The cumulative claims are stochastically independent over the accident years $i=-n, \dots, I$, that is, the cumulative claim matrix $\mathcal{C}_{I,n}$ consists of independent rows $C_{i,\bullet}=(C_{i,0},\dots, C_{i,I+n})$, $i=-n,\dots, I$.
	\end{itemize}
\end{ann}	

For any $n\in\mathbb{N}_0$, based on the available data $\mathcal{D}_{I,n}$, all development factors $\f$ and variance parameters $\s$ for $j=0, \dots, I+n-1$ are unknown and have to be estimated from $\mathcal{D}_{I,n}$. %As proposed by the CLM, 
The development factors $f_0, \ldots,f_{I+n-1}$ can be (consistently) estimated by $\widehat{f}_{0,n},\dots, \widehat{f}_{I+n-1,n}$, where
\begin{align}\label{eqf}
\widehat{f}_{j,n}=\frac{\sum\limits_{i=-n}^{I-j-1}C_{i,j+1}}{\sum\limits_{i=-n}^{I-j-1}C_{i,j}},	\quad	j=0, \dots, I+n-1.
\end{align}
According to \cite{mack1993distribution}, these estimators are unbiased, i.e.~$E(\widehat f_{j,n})=\f$, and pairwise uncorrelated, i.e. $Cov(\widehat f_{j,n}, \widehat{f}_{k,n})=0$ for all $j\neq k$. By plugging-in the $\widehat f_{j,n}$'s, the best estimate of the ultimate claim $\widehat{C}_{i,I+n}$ (point predictor) of the ultimate claim $C_{i,I+n}$ is calculated by
\begin{align*}
\widehat{C}_{i,I+n}=C_{i,I-i}\prod^{I+n-1}_{j=I-i}\widehat{f}_{j,n},	\quad	i=-n, \dots, I.	
\end{align*}
Consequently, given $C_{i,I-i}$, the best estimate $\widehat R_{i,I+n}$ of the reserve $R_{i,I+n}$ is given by 
\begin{equation}\label{defpointRi}
\widehat{R}_{i,I+n}=\widehat C_{i,I+n}- C_{i,I-i}=C_{i,I-i}\left(\prod^{I+n-1}_{j=I-i}\widehat{f}_{j,n}-1\right), \quad  i=-n,\dots, I, 
\end{equation}
and the best estimate $\widehat{R}_{I,n}$ of the total reserve $R_{I,n}$ defined in \eqref{eqtotal_Rn} computes to
\begin{equation}\label{defpointRtotal}
\widehat R_{I,n}=\sum_{i=-n}^{I}\widehat{R}_{i,I+n}
\end{equation}
noting that $\widehat{R}_{-n,I+n}=0$ due to $\prod^{I+n-1}_{j=I+n}\widehat{f}_{j,n}:=1$. Furthermore, \cite{mack1993distribution} proposed to estimate the variance parameters $\sigma_{0}^2,\ldots,\sigma_{I+n-1}^2$ by
\begin{align}\label{def_sigma}
\widehat{\sigma}_{j,n}^2=\frac{1}{I+n-j-1}\sum\limits_{i=-n}^{I-j-1}C_{i,j}\left(\frac{C_{i,j+1}}{C_{i,j}}-\widehat{f}_{j,n}\right)^2,	\quad	j=0, \dots, I+n-2,
\end{align}
which are unbiased estimators, i.e.~$E(\widehat{\sigma}_{j,n}^2)=\sigma_j^2$, and by setting $\widehat \sigma_{I+n-1,n}^2=0$.

Of particular interest is in the distribution of the difference of the stochastic (unobserved) reserve $R_{I,n}$ and its best estimate $\widehat R_{I,n}$ (based on the observed data $\mathcal{D}_{I,n}$), which is denoted as the \emph{predictive root of the reserve} in the following. That is, by combining \eqref{eqtotal_Rn} and \eqref{defpointRtotal}, it computes to
\begin{equation}\label{eqpredroot_firstappearance}
R_{I,n}-\widehat{R}_{I,n} = \sum_{i=-n}^{I}\left({R}_{i,I+n} - \widehat{R}_{i,I+n}\right) %= \sum_{i=-n}^{I}\left(\widehat C_{i,I+n}-C_{i,I-i} - \left(C_{i,I+n}- C_{i,I-i}\right)\right) 
= \sum_{i=-n}^{I}\left( C_{i,I+n}-\widehat C_{i,I+n}\right).
\end{equation}

While a common approach to approximate an unknown (finite sample) distributions is the derivation of asymptotic theory,
%While a common solution to such problems is the derivation of asymptotic theory to approximate the unknown (finite sample) distributions, 
Mack's conditions summarized in Assumption \ref{MackM} are not (yet) sufficient to establish limiting distributions for the predictive root of the reserve $R_{I,n}-\widehat{R}_{I,n}$.
 %For this purpose, as summarized in the next section, \citet{jentschasympThInf} proposed a suitable stochastic model which slightly strengthens Mack's model assumptions and enables the derivation of (conditional and unconditional) limiting distributions for the predictive root of the reserve $R_{I,n}-\widehat{R}_{I,n}$ defined in \eqref{eqpredroot_firstappearance}.

\subsection{A fully-described stochastic framework of Mack's Model}\label{StochasticModel}\

Following \citet[Section 2.2]{jentschasympThInf}, to establish a theoretical framework sufficient to %be able to 
derive asymptotic theory for parameter estimators $\widehat f_{j,n}$ and $\widehat \sigma_{j,n}^2$, which finally also enables the derivation of the limiting distributions of the predictive root of the reserve $R_{I,n}-\widehat{R}_{I,n}$, we introduce Assumptions \ref{initial_lem}, \ref{annF}, and \ref{ann_highermoments_F} on the stochastic mechanism that generates the cumulative claim matrix $\mathcal{C}_{I,n}$, and Assumption \ref{parameterAss} on the sequences of development factors and of variance parameters. They resemble the Assumptions 2.2, 2.3 and 3.3 %in \cite{jentschasympThInf} %, when adopted to the asymptotic framework of Section \ref{sec_asymptotic_framework}, 
as well as Assumption 4.1 in \cite{jentschasympThInf}, respectively. This framework will also allow to rigorously investigate %bootstrap 
consistency properties of the Mack bootstrap in Section \ref{sec_boot_consistency}. 

The first assumption addresses the initial claims, i.e.~the first column %$C_{\bullet,0}=(C_{-n,0},\ldots,C_{I,0})'$ 
of $\mathcal{C}_{I,n}$ (and of $\mathcal{D}_{I,n}$).

\begin{ann}[Initial claims]\label{initial_lem}
Let the initial claims $(C_{I-n,0},n\in\mathbb{N}_0)$ 
%%For any $n\in\mathbb{N}_0$, 
%Let the initial claims %$C_{\bullet,0}=
%$(C_{-n,0},\ldots,C_{I,0})'$, $n\in\mathbb{N}$, 
be independent and identically distributed (i.i.d.) random variables with support $[1,\infty)$, i.e.~$C_{i,0}\geq 1$ for all $i$. Further, let 
$\mu_0:=E(C_{i,0})\in [1,\infty)$ and $\tau^2_0:=Var(C_{i,0})\in (0,\infty)$.
\end{ann}

%Note that 
The independence %between 
of the initial claims is a %common assumption which is also a 
direct consequence of Assumption \ref{MackM} (iii). In addition, Assumption \ref{initial_lem} %also 
imposes an identical distribution for the initial claims. In practice, the condition on the support $[1,\infty)$ of $C_{i,0}$ is not restrictive and can be relaxed to %the condition that 
$C_{i,0}$ %is 
being bounded away from zero.

In view of %the multiplicative structure of 
$E(C_{i,j+1}|C_{i,j})$ in \eqref{MackM_mean}, suppose that the %(random)
cumulative claims $C_{i,j+1}$
%, $i=-n,\ldots,I$, $j=0,\ldots,I+n-1$, 
are recursively defined by 
\begin{align}\label{eq:C_recursive}
C_{i,j+1}=C_{i,j}F_{i,j}=C_{i,0}\prod_{k=0}^{j}F_{i,k},	\quad	i=-n,\ldots,I,\ j=0,\ldots,I+n-1,
\end{align}
where the \textit{individual development factors} $F_{i,j}$ %, which satisfy $\F=\frac{C_{i,j+1}}{C_{i,j}}$ by construction, 
are assumed to fulfill the following condition.

\begin{ann}[Conditional distribution of the individual development factors]\label{annF}
Let the individual development factors $(F_{I-i,j}, i\in\mathbb{N}_0,\ j\in\mathbb{N}_0)$
%For any $n\in\mathbb{N}_0$, let the individual development factors $F_{i,j}$, $i=-n,\dots, I$, $j=0, \dots, I+n-1$ 
be random variables with support $(\epsilon,\infty)$ for some $\epsilon\geq 0$ such that $F_{i,j}$ and $F_{k,l}$ are independent given $(C_{i,j},C_{k,l})$ for all $(i,j)\neq (k,l)$ with %conditional mean and conditional variance 
\begin{equation}\label{annFstMack}
	E(\F|\C)=\f	\quad \text{ and }\quad 	Var(\F|\C)=\frac{\s}{\C}.
\end{equation}
\end{ann}

Note that Mack's original model setup in Assumption \ref{MackM} is implied by Assumptions \ref{initial_lem} and \ref{annF} together. Also note that the stochastic mechanism determined by \eqref{eq:C_recursive} and Assumption \ref{annF} are assumed for the whole cumulative claim matrix $\mathcal{C}_{I,n}$. However, recall that only those $\C$ in $\mathcal{C}_{I,n}$ are observed %up to year $I$ 
that are contained in the upper loss triangle $\mathcal{D}_{I,n}$. Hence, by using the multiplicative relationship in
 %the first identity of 
\eqref{eq:C_recursive}, we have also \emph{perfect} knowledge of %the individual development factors 
$F_{i,j}$, $i=-n,\dots, I-1$, $j=0, \dots, I-i-1$.

According to Lemma 2.4 in \citet{jentschasympThInf}, %the stochastic framework determined by 
Assumptions \ref{initial_lem} and \ref{annF} allow to derive formulas for the (unconditional) means and variances of $C_{i,j}$, $i=-n,\dots, I$, $j=0, \dots, I+n$ leading to
\begin{align*}
E(\C) = \mu_0\prod^{j-1}_{k=0}f_k=:\mu_j	\quad	\text{and}	\quad
Var(\C) = \tau^2_0\prod^{j-1}_{k=0}f_k^2+\mu_0\sum_{l=0}^{j-1}\sigma^2_l\prod^{j-1}_{n=l+1}f_n^2\prod^{l-1}_{m=0}f_m=:\tau_j^2,
\end{align*}
%\begin{align}
%E(\C) &= \mu_0\prod^{j-1}_{k=0}f_k=:\mu_j\label{meanC2},\\
%Var(\C) &= \tau^2_0\prod^{j-1}_{k=0}f_k^2+\mu_0\sum_{l=0}^{j-1}\sigma^2_l\prod^{j-1}_{n=l+1}f_n^2\prod^{l-1}_{m=0}f_m=:\tau_j^2,\label{varC2}
%\end{align}
where $\mu_0$ and $\tau^2_0$ are defined in Assumption \ref{initial_lem}. Together with Assumption \ref{parameterAss} below, according to Lemma 4.2 in \citet{jentschasympThInf}, both sequences $(\mu_j,j\in\mathbb{N}_0)$ and $(\tau_j^2,j\in\mathbb{N}_0)$ are non-negative, monotonically non-decreasing, and converging with $\mu_j\rightarrow \mu_\infty$ and $\tau_j^2\rightarrow \tau_\infty^2$ as $j\rightarrow \infty$, where $\mu_\infty:=\mu_0\prod^\infty_{j=0}f_j$ and $\tau_\infty^2:=\tau^2_0\prod^\infty_{k=0}f_k^2+\mu_0\sum_{l=0}^\infty\left(\prod^{l-1}_{m=0}f_m\right)\sigma^2_l\left(\prod^\infty_{n=l+1}f_n^2\right)$.

\begin{ann}[Development Factors and Variance Parameters]\label{parameterAss}\
Letting $n\rightarrow\infty$ in the setup of Assumptions \ref{initial_lem} and \ref{annF} leads to
	\begin{itemize}
		\item[(i)] a sequence of development factors $(\f,j\in\mathbb{N}_0)$ with $f_j\geq 1$ for all $j\in \mathbb{N}_0$ and $\f\rightarrow 1$ as $j\rightarrow \infty$ such that $\prod\limits^\infty_{j=0}f_j<\infty$, which is equivalent to $\sum\limits^\infty_{j=0}(f_j-1)<\infty$.
		\item[(ii)] a sequence of variance parameters $(\s,j\in\mathbb{N}_0)$ with $\sigma_0^2>0$ and $\sigma_j^2\geq 0$ for all $j\in \mathbb{N}$ with $\s\rightarrow 0$ as $j\rightarrow \infty$ such that $\sum\limits^\infty_{j=0}(j+1)^2\s<\infty$.
	\end{itemize}
\end{ann}

The conditions imposed on the sequences of development factors $(\f,j\in\mathbb{N}_0)$ and variance parameters $(\s,j\in\mathbb{N}_0)$ in Assumption \ref{parameterAss} are rather mild. In practice, each claim has a finite, but possibly unknown horizon until it is finally settled, %. The time horizons of claim developments vary by the insurance lines, which are usually categorized in short-term and long-term.
which varies by the insurance lines. Altogether, as done in \citet[Section 3.1]{jentschasympThInf}, this setup allows to derive central limit theorems (CLTs) for (smooth functions of) the parameter estimators $\widehat f_{j,n}$ for $n\rightarrow\infty$.

According to \citet[Section 3.2]{jentschasympThInf}, the following additional assumption has to be imposed to derive a CLT also for $\widehat \sigma^2_{j,n}$. Although the distributional properties of $\widehat \sigma^2_{j,n}$ do not show asymptotically in the distribution of the reserve, $\sqrt{I+n}$-consistency of $\widehat \sigma^2_{j,n}$ as obtained in \citet[Theorem 3.5]{jentschasympThInf} is required for establishing the bootstrap asymptotic theory in Section \ref{sec_boot_consistency}.

%For the derivation of a similar CLT result for $\widehat \sigma^2_{j,n}$ for fixed $j$, according to \citet[Section 3.2]{jentschasympThInf}, the following additional assumption is imposed. However, although the distributional properties of $\widehat \sigma^2_{j,n}$ do not show asymptotically in the distribution of the reserve, $\sqrt{I+n}$-consistency of $\widehat \sigma^2_{j,n}$ as obtained in \citet[Theorem 3.5]{jentschasympThInf} is required for establishing the bootstrap asymptotic theory in Section \ref{sec_boot_consistency}.

\begin{ann}[Higher-order conditional moments of individual development factors]\label{ann_highermoments_F}
For all $i\in\mathbb{Z}$, $i\leq I$, $j\in\mathbb{N}_0$, suppose that 
%For any $n\in\mathbb{N}_0$, suppose that, 
conditional on $C_{i,j}$, the third and fourth (central) moments of the individual development factors $F_{i,j}$, that is, $E((F_{i,j}-f_j)^3|\C)$ and $E((F_{i,j}-f_j)^4|\C)$ exist 
%for $j=0,\dots, I+n-1$ 
such that both
\begin{align}\label{kappa4}
\kappa_j^{(3)}=E(\C^2E((\F-f_j)^3|\C)) \quad \text{and}\quad \kappa_j^{(4)}=E(\C^2E((\F-f_j)^4|\C))
\end{align}
exist and are finite, respectively.
\end{ann}

Using %the recursive stochastic model for the claims in
\eqref{eq:C_recursive}, conditional on $\mathcal{Q}_{I,n}$, the reserve $R_{I,n}$ can be written as
\begin{align}\label{eq_stochreserve}
R_{I,n}=\sum^I_{i=-n}C_{i,I-i}\left(\prod^{I+n-1}_{j=I-i}F_{i,j}-1\right).
\end{align}
Hence, by plugging-in \eqref{defpointRi} and \eqref{eq_stochreserve}, the predictive root of the reserve from \eqref{eqpredroot_firstappearance} becomes 
\begin{equation}\label{eqpredroot}
R_{I,n}-\widehat{R}_{I,n}=\sum^I_{i=-n}C_{i,I-i}\left(\prod^{I+n-1}_{j=I-i}\F-\prod^{I+n-1}_{j=I-i}\widehat{f}_{j,n}\right)
	=\sum^{I+n}_{i=0}C_{I-i,i}\left(\prod^{I+n-1}_{j=i}F_{I-i,i}-\prod^{I+n-1}_{j=i}\widehat{f}_{j,n}\right),
\end{equation}
where we flipped the index $i$ to $I-i$ in the last step.

%\bigskip

\section{Mack's Bootstrap Scheme}\label{sec_boot}

The Mack bootstrap, %proposal, 
introduced by \cite{england2006predictive}, equips Mack's Model with a resampling procedure to estimate the whole distribution of the (predicted) reserve. %The resulting Mack bootstrap
It is very popular and widely used in practice as it describes a rather simple to implement algorithm to estimate the \emph{reserve risk} %such as e.g.~the value-at-risk, 
by estimating high quantiles of the reserve distribution.% of the . 

As proposed by \cite{england2006predictive}, to mimic the distribution of the predictive root of the reserve $R_{I,n}-\widehat{R}_{I,n}$, the Mack bootstrap constructs a certain bootstrap version $R_{I,n}^*-\widehat{R}_{I,n}$ of it. On the one hand, this bootstrap predictive root relies on the \emph{same} best estimate of the reserve and centers $R_{I,n}^*$ also around $\widehat R_{I,n}$. %This is motivated by the definition of the reserve risk, which captures the risk that the best estimate of the reserve is not sufficient to pay for all outstanding claims. 
On the other hand, it constructs a certain \emph{double}-bootstrap version of the reserve $R_{I,n}$, that is $R_{I,n}^*$, by combining \emph{two} complementing (non-parametric and parametric) bootstrap approaches for resampling the individual development factors in the upper triangle and in the lower triangle:

\begin{itemize}
\item[(i)] First, a \emph{non-parametric} residual-based bootstrap (see Step 4 below) is applied to construct bootstrap individual development factors $F_{i,j}^*$, $j=0,\ldots,I+n-1$, $i=-n,\ldots,I-j-1$, that is, for the \emph{upper triangle}, in order to get bootstrap development factor estimators $\widehat f_{j,n}^*$, $j=0,\ldots,I+n-1$.
\item[(ii)] Second, the bootstrap development factor estimators $\widehat f_{j,n}^*$ from (i) together with a \emph{parametric} bootstrap (see Step 5 below) are used to construct also bootstrap individual development factors $F_{i,j}^*$, $i=-n,\ldots,I$, $j=0,\ldots,I+n-1$ and $i+j\geq I$, that is for the \emph{lower triangle}. For this purpose, a parametric family of (conditional) bootstrap distributions %(such as e.g.~a gamma or a log-normal distribution) 
has to be chosen. %In contrast to the non-parametric approach in (i), the parametric approach is favored here to assure that $F_{i,j}^*>0$ (a.s.). 
\end{itemize}

Finally, as we are dealing with a prediction problem when estimating the reserve risk, the limiting properties of the predictive root of the reserve \emph{conditional} on the latest observed cumulative claims are relevant and have to be mimicked by a suitable resampling procedure. For this purpose, the Mack bootstrap is employed to estimate the conditional distribution of $R_{I,n}-\widehat R_{I,n}$ given $\mathcal{Q}_{I,n}$ by the conditional bootstrap distribution of $R_{I,n}^*-\widehat R_{I,n}$ given $\mathcal{Q}_{I,n}^*=\mathcal{Q}_{I,n}$ and $\mathcal{D}_{I,n}$.

\subsection{Mack's Bootstrap Algorithm}\label{subsec_boot}\

With the upper triangle $\mathcal{D}_{I,n}$ at hand, Mack's bootstrap algorithm is defined as follows:

\begin{itemize}
	\item[Step 1.] Estimate the development factors $f_j$ and the variance parameters $\sigma_j^2$ from $\mathcal{D}_{I,n}$ by computing $\widehat f_{j,n}$ and $\widehat \sigma^2_{j,n}$ for $j=0, \dots, I+n-1$ as defined in \eqref{eqf} and \eqref{def_sigma}, respectively.
	
	\item[Step 2.] For all $j=0, \dots, I+n-1$ with $\widehat\sigma_{j,n}^2>0$, compute 'residuals'
\begin{align}\label{residuals}
\widehat r_{i,j}=\sqrt{\C}(\F-\widehat f_{j,n})/\widehat \sigma_{j,n},	\quad	i=-n, \dots, I-j-1,
\end{align}
and %for $i=-n, \dots, I-j-1$. 
re-center and re-scale them %se $\widehat r_{i,j}$'s 
to get
%\begin{align*}
$\widetilde r_{i,j} = \frac{1}{s}\left(\widehat r_{i,j}-\overline{r}\right)$,
%\end{align*}
where\footnote{Note that $\widehat\sigma_{I+n-1,n}^2=0$ by construction such that $\widehat r_{-n,I+n-1}$ is excluded in \eqref{residuals} such that (at most) $(I+n+1)(I+n)/2-1=((I+n+1)(I+n)-2)/2$ residuals can be computed. If $\widehat\sigma_{j,n}^2=0$ holds also for other $j$, the corresponding residuals are excluded in \eqref{residuals} as well and the formulas for $\overline{r}$ and $s$ %in \eqref{r_bar} and \eqref{s}, respectively,	
have to be adjusted accordingly. In the following, for notational convenience, we assume that only $\widehat\sigma_{I+n-1,n}^2=0$ and $\widehat\sigma_{j,n}^2>0$ holds for all $j=0,\ldots,I+n-2$ and all $n\in\mathbb{N}_0$.}
\begin{align*}%\label{r_bar_s}
\overline{r} = \frac{2}{(I+n+1)(I+n)-2}\sum_{k=0}^{I+n-2} \sum_{l=-n}^{I-k-1}\widehat r_{l,k},%,	\label{r_bar}	\\
\quad	%\text{and}	\quad	
s^2 = \frac{2}{(I+n+1)(I+n)-2}\sum_{k=0}^{I+n-2} \sum_{l=-n}^{I-k-1} \left(\widehat r_{l,k}-\overline{r}\right)^2.	%\label{s}
\end{align*}
%\begin{align}
%\overline{r} &= \frac{2}{(I+n+1)(I+n)-2}\sum_{k=0}^{I+n-2} \sum_{l=-n}^{I-k-1}\widehat r_{l,k},	\label{r_bar}	\\
%s^2 &= \frac{2}{(I+n+1)(I+n)-2}\sum_{k=0}^{I+n-2} \sum_{l=-n}^{I-k-1} \left(\widehat r_{l,k}-\overline{r}\right)^2.	\label{s}
%\end{align}
	
	\item[Step 3.] Draw randomly with replacement from the re-centered and re-scaled residuals $\widetilde{r}_{i,j}$, $j=0,\dots,I+n-2$, $i=-n,\dots,I-j-1$ to get 'bootstrap errors' $r_{i,j}^*$, $j=0,\dots,I+n-1$, $i=-n,\dots,I-j-1$.
	
	\item[Step 4.] Define the bootstrap individual development factors 
\begin{align}	
\F^*=\widehat f_{j,n}+\frac{\widehat \sigma_{j,n}}{\sqrt{\C}} r_{i,j}^*,	\quad	j=0, \dots, I+n-1,\ 	i=-n, \dots, I-j-1,
\end{align}
%for $j=0, \dots, I+n-1$ and $i=-n, \dots, I-j-1$, 
let $\mathcal{F}_{I,n}^*=\{F_{i,j}^*|j=0,\dots,I+n-1,~i=-n,\dots,I-j-1\}$, and compute the Mack bootstrap development factor estimators
\begin{align}\label{eqf_boot}
\widehat f_{j,n}^*=\frac{\sum^{I-j-1}_{i=-n}\C\F^*}{\sum^{I-j-1}_{i=-n}\C}=\widehat f_{j,n}+\frac{\widehat \sigma_{j,n}\sum^{I-j-1}_{i=-n}\sqrt{\C} r_{i,j}^*}{\sum^{I-j-1}_{i=-n}\C},	\quad	j=0, \dots, I+n-1.
\end{align}
%for $j=0, \dots, I+n-1$.

	\item[Step 5.] Choose a parametric family for the (conditional) bootstrap distributions of $\F^*$ given $\C^*,\mathcal{D}_{I,n}$ and $\mathcal{F}_{I,n}^*$ such that $\F^*>0$ a.s.~with 
\begin{align*}
E^*(F_{i,j}^*|C_{i,j}^*, \mathcal{F}_{I,n}^*)=\widehat{f}_{j,n}^*,	\quad	Var^*(F_{i,j}^*|C_{i,j}^*, \mathcal{F}_{I,n}^*)=\frac{\widehat{\sigma}^2_{j,n}}{C_{i,j}^*},	\quad	i=-n, \dots, I,\ j=I-i,\dots, I+n-1,%, i+j\geq I
\end{align*}
%mean $E^*(F_{i,j}^*|C_{i,j}^*, \mathcal{F}_{I,n}^*)=\widehat{f}_{j,n}^*$ and variance $Var^*(F_{i,j}^*|C_{i,j}^*, \mathcal{F}_{I,n}^*)=\frac{\widehat{\sigma}^2_{j,n}}{C_{i,j}^*}$ for $i=-n, \dots, I$, $j=0,\dots, I+n-1$ and $i+j\geq I$
where $E^*(\cdot):=E^*(\cdot|\mathcal{D}_{I,n})$, $Var^*(\cdot):=Var^*(\cdot|\mathcal{D}_{I,n})$, etc.~denote the Mack bootstrap mean, %Mack bootstrap
variance, etc., respectively, that is, conditional on the data $\mathcal{D}_{I,n}$. Then, given $\mathcal{Q}_{I,n}^*=\mathcal{Q}_{I,n}$, generate the bootstrap ultimate claims $C_{i,I+n}^*$ and the reserves $R_{i,I+n}^*=C_{i,I+n}^*-C_{i,I-i}^*$ for $i=-n,\dots,I$ using the recursion
\begin{align}
C_{i,j+1}^*=C_{i,j}^*F_{i,j}^*,	\quad	j=I-i,\ldots,I+n-1.
\end{align}
	
	\item[Step 6.] Compute the bootstrap total reserve $R_{I,n}^*=\sum_{i=0}^{I+n} R_{I-i,I+n}^*$ and %the 
	its bootstrap predictive root %of the reserve
\begin{align}\label{eqpredroot_boot}
	R^*_{I,n}-\widehat{R}_{I,n} = \sum^{I+n}_{i=0}C_{I-i,i}^*\left(\prod^{I+n-1}_{j=i}F_{I-i,i}^*-\prod^{I+n-1}_{j=i}\widehat{f}_{j,n}\right).
\end{align}	
	\item[Step 7.] Repeat Steps	3 - 6 above $B$ times, where $B$ is large, to get bootstrap predictive roots $(R^*_{I,n}-\widehat{R}_{I,n})^{(b)}$, $b=1,\ldots,B$, and denote by $q^*(\alpha)$ the $\alpha$-quantile of their empirical distribution. 
	
	\item[Step 8.] Construct the $(1-\alpha)$ equal-tailed prediction interval for $R_{I,n}$ as
\begin{align*}
\left[\widehat{R}_{I,n}+q^*(\alpha/2), \widehat{R}_{I,n}+q^*(1-\alpha/2)\right].
\end{align*}
\end{itemize}

\begin{rem}[On Mack's bootstrap proposal]\ 
\begin{itemize}
	\item[(i)] While the Mack bootstrap predictive root of the reserve $R^*_{I,n}-\widehat{R}_{I,n}$ uses the same best estimate $\widehat{R}_{I,n}$ for centering (as in $R_{I,n}-\widehat{R}_{I,n}$), it relies on a certain type of double-bootstrap version $R^*_{I,n}$ of the total reserve $R_{I,n}$, which employs $\widehat f_{j,n}^*$ instead of just $\widehat f_{j,n}$, but uses $\widehat \sigma_{j,n}^2$. However,
%although $\widehat f_{j,n}^*$ instead of $\widehat f_{j,n}$ is used, 
although $E^*(F_{i,j}^*|C_{i,j}^*,\mathcal{D}_{I,n},\mathcal{F}_{I,n}^*)=\widehat f_{j,n}^*$ holds, we still have $E^*(F_{i,j}^*|C_{i,j}^*)=\widehat f_{j,n}$ for the lower triangle individual development factors. %, but $E^*(F_{i,j}^*|C_{i,j}^*,\mathcal{D}_{I,n},\mathcal{F}_{I,n}^*)=\widehat f_{j,n}^*$. 
In contrast, for the variances, we have 
$Var^*(F_{i,j}^*|C_{i,j}^*,\mathcal{F}_{I,n}^*)=\frac{\widehat \sigma_{j,n}^2}{C_{i,j}^*}$, but% and %, using the law of total variance, we get 
\begin{align*}
Var^*(F_{i,j}^*|C_{i,j}^*) =
%& E^*\left(Var^*(F_{i,j}^*|C_{i,j}^*,\mathcal{F}_{I,n}^*)|C^*_{i,j}\right)+Var^*\left(E^*(F_{i,j}^*|C_{i,j}^*,\mathcal{F}_{I,n}^*)|C^*_{i,j}\right)	\\
%=& 
\frac{\widehat \sigma_{j,n}^2}{C^*_{i,j}}+\frac{\widehat \sigma_{j,n}^2}{\sum^{I-j-1}_{k=-n}C_{k,j}}.
\end{align*}
	\item[(ii)] Due to the fixed-design bootstrap in Step 4, which does not generate bootstrap cumulative claims $C_{i,j}^*$ (and consequently no bootstrap upper loss triangle $\mathcal{D}_{I,n}^*$), but only $F_{i,j}^*$'s, the bootstrap development factor estimators $\widehat f_{j,n}^*$ and $\widehat{f}_{k,n}^*$ defined in \eqref{eqf_boot} are independent for $j\not= k$ conditional on $\mathcal{D}_{I,n}$. This is on contrast to the development factor estimators $\widehat{f}_{j,n}$ and $\widehat{f}_{k,n}$, which are asymptotically independent for $j\not= k$, but only uncorrelated in finite samples such that $E(\widehat{f}_{j,n}^2\widehat{f}_{k,n}^2)<0$ for $j\neq k$.
%Whether the development factors are independent is also reflected in the formula for the MSEP of the reserve. Mack's formula takes into account the uncorrelatedness of $\widehat{f}_{j,n}$ and $\widehat{f}_{k,n}$, whereas in the formula of the MSEP by \citet{buchwalder2006mean} the estimates of the development factors are independent.
	\item[(iii)] The non-parametric bootstrap used to construct the $\widehat f_{j,n}^*$'s in Step 4 uses residuals, but according to Assumption \ref{initial_lem} and \ref{annF}, %the claims $C_{i,j}$ are defined recursively such that 
there are no errors in Mack's model that are approximated by these residuals. In fact, each (possibly parametric) bootstrap proposal that successfully mimics the first and second conditional moments of $C_{i,j+1}$ given $C_{i,j}$ will correctly mimic the limiting distribution of the $\widehat f_{j,n}$'s.
	%\item[(iv)] While $\widehat f_{j,n}^*>0$ is not guaranteed in finite samples by the non-parametric bootstrap proposal in Step 4, the parametric bootstrap used to construct the $F_{i,j}^*$'s in Step 5 assures that all individual development factors $F_{i,j}^*$ are a.s.~positive. However, the generation of $F_{i,j}^*$ requires $\widehat f_{j,n}^*>0$ to hold. 
	\item[(iv)] In view of the discussion above, a fully parametric implementation that uses the same parametric family from Step 5 also in Step 4 to get bootstrap development factors $F_{i,j}^*$'s %and, consequently, the $\widehat f_{j,n}^*$'s 
	can be used. 
	\item[(v)] A fully non-parametric approach that uses the non-parametric bootstrap from Step 4 also in Step 5 is thinkable, but would suffer from issues arising from potentially negative $F_{i,j}^*$'s leading to a reduced finite sample performance. 
\end{itemize}
\end{rem}

%\bigskip

\section{Asymptotic Theory for the  Mack Bootstrap}\label{sec_boot_consistency}

Although the Mack bootstrap as proposed by \cite{england2006predictive} and described in Section \ref{sec_boot} is widely used in practice for reserve risk estimation, limiting results that confirm its consistency are still missing in the literature. In this section, based on the asymptotic and stochastic framework described in Section \ref{CLM_section}, we derive asymptotic theory for the Mack bootstrap, which enables a rigorous investigation of its consistency properties. 

The Mack bootstrap is designed to mimic the distribution of the predictive root of the reserve $R_{I,n}-\widehat R_{I,n}$ conditional on $\mathcal{Q}_{I,n}$ based on the bootstrap distribution of the corresponding Mack bootstrap predictive root of the reserve $R_{I,n}^*-\widehat R_{I,n}$ conditional on $\mathcal{Q}_{I,n}^*=\mathcal{Q}_{I,n}$ and $\mathcal{D}_{I,n}$. Hence, a closer inspection of both expressions is advisable. Picking-up the representation of the predictive root of the reserve $R_{I,n}-\widehat R_{I,n}$ in \eqref{eqpredroot}, it can be decomposed into two additive parts that account for the prediction error and the estimation error, respectively. Precisely, by subtracting and adding $\sum^{I+n}_{i=0}C_{I-i,i}\prod^{I+n-1}_{j=i}f_j$, we get
\begin{align}
R_{I,n}-\widehat{R}_{I,n}
& =\sum^{I+n}_{i=0}C_{I-i,i}\left(\prod^{I+n-1}_{j=i}F_{I-i,j}-\prod_{j=i}^{I+n-1}{f}_j\right) + \sum^{I+n}_{i=0}C_{I-i,i}\left(\prod_{j=i}^{I+n-1}{f}_j-\prod_{j=i}^{I+n-1}\widehat{f}_{j,n}\right)	\nonumber	\\
&=:\left(R_{I,n}-\widehat{R}_{I,n}\right)_1+\left(R_{I,n}-\widehat{R}_{I,n}\right)_2,	\label{eqpredroot_decomp_first_appearance}
\end{align}
where $(R_{I,n}-\widehat{R}_{I,n})_1$ represents the process uncertainty (that carries the process variance) and $(R_{I,n}-\widehat{R}_{I,n})_2$ the estimation uncertainty (that carries the estimation variance).

Similarly, for the Mack bootstrap predictive root of the reserve $R_{I,n}^*-\widehat R_{I,n}$ from \eqref{eqpredroot_boot}, by subtracting and adding $\sum^{I+n}_{i=0}C_{I-i,i}^*\prod^{I+n-1}_{j=i}\widehat f_{j,n}^*$, we get
\begin{align}
R^*_{I,n}-\widehat{R}_{I,n}
&=\sum^{I+n}_{i=0}C_{I-i,i}^*\left(\prod^{I+n-1}_{j=i}F_{I-i,j}^*-\prod_{j=i}^{I+n-1}\widehat{f}_{j,n}^*\right)  +\sum^{I+n}_{i=0}C_{I-i,i}^*\left(\prod_{j=i}^{I+n-1}\widehat{f}_{j,n}^*-\prod_{j=i}^{I+n-1}\widehat{f}_{j,n}\right)	\nonumber	\\
&=:\left(R^*_{I,n}-\widehat{R}_{I,n}\right)_1+\left(R^*_{I,n}-\widehat{R}_{I,n}\right)_2,	\label{eqpredroot_boot_decomp_appearance}
\end{align} 
where $(R^*_{I,n}-\widehat{R}_{I,n})_1$ and $(R^*_{I,n}-\widehat{R}_{I,n})_2$ are the Mack bootstrap versions of $(R_{I,n}-\widehat{R}_{I,n})_1$ and $(R_{I,n}-\widehat{R}_{I,n})_2$, respectively.

As main interest is in the distribution of the predictive root of the reserve $R_{I,n}-\widehat R_{I,n}$ conditional on $\mathcal{Q}_{I,n}$, in view of the decompositions \eqref{eqpredroot_decomp_first_appearance} and \eqref{eqpredroot_boot_decomp_appearance}, it is instructive to first consider separately the (limiting) distributions of $(R_{I,n}-\widehat{R}_{I,n})_1$ and $(R_{I,n}-\widehat{R}_{I,n})_2$ conditional on $\mathcal{Q}_{I,n}$, respectively. %, as well as jointly. 
They will serve as valuable benchmark distributions %to enable a meaningful 
for the investigation of consistency properties of the Mack bootstrap in Section \ref{sec_asymp_boot}. Such asymptotic results have been established in \citet[Section 4]{jentschasympThInf}. We will briefly summarize the relevant conditional limiting distributions below in Section \ref{sec_asymp_summary}.

\subsection{Conditional asymptotics for the predictive root of the reserve}\label{sec_asymp_summary}\

%Based on the same stochastic and asymptotic framework used throughout this paper, \citet{jentschasympThInf} established asymptotic theory for both parts of the predictive root of the reserve $R_{I,n}-\widehat{R}_{I,n}$, that is, for $(R_{I,n}-\widehat{R}_{I,n})_1$ and $(R_{I,n}-\widehat{R}_{I,n})_2$ separately, as well as jointly for $R_{I,n}-\widehat{R}_{I,n}$. As risk reserving relies on the \emph{prediction} of $R_{I,n}$ using $\widehat{R}_{I,n}$, which is computed from the observed upper loss triangle $\mathcal{D}_{I,n}$ which consists of all cumulative claims up to the diagonal $\mathcal{Q}_{I,n}$, the main interest is in asymptotic theory for $(R_{I,n}-\widehat{R}_{I,n})_1$, $(R_{I,n}-\widehat{R}_{I,n})_2$, and $R_{I,n}-\widehat{R}_{I,n}$ conditional on $\mathcal{Q}_{I,n}$, respectively.

In the following, we review the conditional asymptotic results established in \citet[Theorems 4.3, 4.10, 4.12, and Corollary 4.13]{jentschasympThInf} separately for the process uncertainty term $(R_{I,n}-\widehat{R}_{I,n})_1$ in Section \ref{sec_asymp_summary_process}, for the estimation uncertainty term $(R_{I,n}-\widehat{R}_{I,n})_2$ in Section \ref{sec_asymp_summary_estimation}, as well as jointly for $R_{I,n}-\widehat{R}_{I,n}$ in Section \ref{sec_asymp_summary_joint}, respectively.

\subsubsection{Conditional asymptotics for reserve prediction: process uncertainty}\label{sec_asymp_summary_process}\

Based on Theorem 4.3 from \citet{jentschasympThInf}, the following theorem provides the limiting distribution of the process uncertainty term $(R_{I,n}-\widehat{R}_{I,n})_1$ conditional on $\mathcal{Q}_{I,n}$.

\begin{thm}[Asymptotics for $(R_{I,n}-\widehat{R}_{I,n})_1$ conditional on $\mathcal{Q}_{I,n}$]\label{Dist_R1}
Suppose Assumptions \ref{initial_lem}, \ref{annF} and \ref{parameterAss} hold. Then, as $n\rightarrow \infty$, conditionally on $\mathcal{Q}_{I,n}$, $(R_{I,n}-\widehat{R}_{I,n})_1$ converges in $L_2$-sense to the non-degenerate random variable $(R_{I,\infty}-\widehat{R}_{I,\infty})_1$. That is, we have
\begin{align}
E\left(\left((R_{I,n}-\widehat{R}_{I,n})_1-(R_{I,\infty}-\widehat{R}_{I,\infty})_1\right)^2|\mathcal{Q}_{I,n}\right) &\overset{p}{\rightarrow} 0,	\label{cond_L_2_convergence}	
\end{align}
where
\begin{align}\label{R_1_diff_infty}
(R_{I,\infty}-\widehat{R}_{I,\infty})_1 := \sum_{i=0}^\infty C_{I-i,i}\left(\prod^\infty_{j=i}F_{I-i,j}-\prod_{j=i}^\infty{f}_j\right)\sim \mathcal{G}_1.
\end{align}
Conditional on $\mathcal{Q}_{I,\infty}=(C_{I-i,i}|i\in\mathbb{N}_0)$, $\mathcal{G}_1$ has mean zero, $E((R_{I,\infty}-\widehat{R}_{I,\infty})_1|\mathcal{Q}_{I,\infty})=0$, and variance
\begin{align}
Var\left((R_{I,\infty}-\widehat{R}_{I,\infty})_1|\mathcal{Q}_{I,\infty}\right) &= \sum^\infty_{i=0}C_{I-i,i}\sum^\infty_{j=i}\left(\prod^{j-1}_{k=i}f_k\right)\sigma^2_j\left(\prod^\infty_{l=j+1}f_l^2\right)=O_P(1).	
\label{var_R1_cond}
\end{align}
\end{thm}

The (conditional) $L_2$-convergence result in Theorem \ref{Dist_R1} immediately implies also (conditional) convergence in distribution. That is, for $n\rightarrow \infty$, we have
\begin{align}\label{eq_convergence_in_distribution}
(R_{I,n}-\widehat{R}_{I,n})_1|\mathcal{Q}_{I,n} \overset{d}{\rightarrow} (R_{I,\infty}-\widehat{R}_{I,\infty})_1|\mathcal{Q}_{I,\infty} \sim \mathcal{G}_1|\mathcal{Q}_{I,\infty}.
\end{align}
Moreover, according to Theorem \ref{Dist_R1} (see also the discussion in \cite[Remark 4.4]{jentschasympThInf}), the conditional limiting distribution $\mathcal{G}_1|\mathcal{Q}_{I,\infty}$ will be typically \emph{non}-Gaussian and depending on the (conditional) distribution of the individual development factors $F_{i,j}|C_{i,j}$. 

\subsubsection{Conditional asymptotics for reserve prediction: estimation uncertainty}\label{sec_asymp_summary_estimation}\

In comparison to the conditional limiting result for $(R_{I,n}-\widehat{R}_{I,n})_1$ displayed in Theorem \ref{Dist_R1}, the derivation of asymptotic results for $(R_{I,n}-\widehat{R}_{I,n})_2$ is rather different and also much more cumbersome. In particular, to obtain non-degenerate limiting distributions, we have to inflate $(R_{I,n}-\widehat{R}_{I,n})_2$ by $\sqrt{I+n+1}$ and the obtained (Gaussian) distribution relies on CLTs for (smooth functions of) development factor estimators $\widehat f_{j,n}$ established in \cite[Section 3 and Appendix C]{jentschasympThInf}. 
For the derivation of asymptotic theory, conditional on $\mathcal{Q}_{I,n}$,
%\footnote{Note that conditioning on $\mathcal{Q}_{I,n}$ or $\mathcal{Q}_{I,\infty}$ is equivalent.}, 
it is instructive to further decompose $(R_{I,n}-\widehat{R}_{I,n})_2$ to get
\begin{align}
& (R_{I,n}-\widehat{R}_{I,n})_2	\nonumber	\\
&= \sum^{I+n}_{i=0}C_{I-i,i}\left(\prod_{j=i}^{I+n-1}f_j-\prod_{j=i}^{I+n-1}f_{j,n}(\mathcal{Q}_{I,n})\right) + \sum^{I+n}_{i=0}C_{I-i,i}\left(\prod_{j=i}^{I+n-1}f_{j,n}(\mathcal{Q}_{I,n})-\prod_{j=i}^{I+n-1}\widehat{f}_{j,n}\right)	\nonumber	\\
&= (R_{I,n}-\widehat{R}_{I,n})_2^{(1)}+(R_{I,n}-\widehat{R}_{I,n})_2^{(2)},	\label{eqpredroot_decomp}
\end{align}
where $(R_{I,n}-\widehat{R}_{I,n})_2^{(1)}$ is measurable wrt $\mathcal{Q}_{I,n}$ and $f_{j,n}(\mathcal{Q}_{I,n}):=\mu_{j+1,n}^{(1)}(\mathcal{Q}_{I,n})/\mu_{j,n}^{(2)}(\mathcal{Q}_{I,n})$ with $\mu_{j+1,n}^{(1)}(\mathcal{Q}_{I,n}) = E(\frac{1}{I+n-j}\sum_{i=-n}^{I-j-1} C_{i,j+1}|\mathcal{Q}_{I,n})$ and $\mu_{j,n}^{(2)}(\mathcal{Q}_{I,n}) = E(\frac{1}{I+n-j}\sum_{i=-n}^{I-j-1} C_{i,j}|\mathcal{Q}_{I,n})$.
%\begin{align*}
%\mu_{j+1,n}^{(1)}(\mathcal{Q}_{I,n}) = E\left(\frac{1}{I+n-j}\sum_{i=-n}^{I-j-1} C_{i,j+1}|\mathcal{Q}_{I,n}\right),	\quad
%\mu_{j,n}^{(2)}(\mathcal{Q}_{I,n}) = E\left(\frac{1}{I+n-j}\sum_{i=-n}^{I-j-1} C_{i,j}|\mathcal{Q}_{I,n}\right).
%\end{align*}
%
%
%\begin{align*}
%\mu_{j+1,n}^{(1)}(\mathcal{Q}_{I,n}) &:= E\left(\frac{1}{I+n-j}\sum_{i=-n}^{I-j-1} C_{i,j+1}|\mathcal{Q}_{I,n}\right) = \frac{1}{I+n-j}\sum_{i=-n}^{I-j-1} E(C_{i,j+1}|C_{i,I-i}),	\\
%\mu_{j,n}^{(2)}(\mathcal{Q}_{I,n}) &:= E\left(\frac{1}{I+n-j}\sum_{i=-n}^{I-j-1} C_{i,j}|\mathcal{Q}_{I,n}\right) =\frac{1}{I+n-j}\sum_{i=-n}^{I-j-1} E(C_{i,j}|C_{i,I-i}).
%\end{align*}
The derivation of (conditional) asymptotic theory for $(R_{I,n}-\widehat{R}_{I,n})_2$ requires additional assumptions on the stochastic properties of the individual development factors $F_{i,j}$ summarized in Assumptions \ref{support_condition} and \ref{reverse_moments_condition} below, which resemble Assumptions 4.6 and 4.8 in \citet{jentschasympThInf}. 

\begin{ann}[Support condition and variance parameters]\label{support_condition}
The individual development factors $(F_{i,j}$, $i\in\mathbb{Z}$, $i\leq I$, $j\in\mathbb{N}_0)$ are random variables with support $(\epsilon,\infty)$ for some $\epsilon>0$ and the sequence of variance parameters $(\sigma_j^2,j\in\mathbb{N}_0)$ converges to $0$ as $j\rightarrow \infty$ such that $\sum\limits^\infty_{j=0}(j+1)^2\frac{\sigma_j^2}{\epsilon^j}<\infty$.
\end{ann}

In addition to the condition on the support and the variance parameters in Assumption \ref{support_condition}, a regularity condition for the \emph{backward} conditional distribution of cumulative claim $C_{i,j}$ given $C_{i,j+1}$ is required.

\begin{ann}[Backward conditional moments]\label{reverse_moments_condition}
Assumptions \ref{initial_lem}, \ref{annF}, \ref{parameterAss} and \ref{support_condition} are fulfilled such that, for all $K\in\mathbb{N}_0$, $k\geq 0$ and $j,j_1,j_2\in\{0,\ldots,K\}$, $j_1\leq j_2$, we have
\begin{align*}
|E(C_{i,j}|C_{i,j+k})-E(C_{i,j}|C_{i,j+k+1})| &\leq a_k X_i,	\\
|Cov(C_{i,j_1},C_{i,j_2}|C_{i,j_2+k})-Cov(C_{i,j_1},C_{i,j_2}|C_{i,j_2+k+1})| &\leq b_k Y_i,
\end{align*}
where $(X_i,i\in\mathbb{Z},i\leq I)$, $(Y_i,i\in\mathbb{Z},i\leq I)$ are sequences of non-negative i.i.d.~random variables with $E(X_i^{2+\delta})<\infty$ for some $\delta>0$ and $E(Y_i^2)<\infty$, and $(a_j,j\in\mathbb{N}_0)$ and $(b_j,j\in\mathbb{N}_0)$ are non-negative real-valued sequences with $\sum_{j=0}^\infty (j+1)^2a_j<\infty$ and $\sum_{j=0}^\infty (j+1)^2b_j<\infty$.
\end{ann}

While Mack's model is designed to generate loss triangles in a rather simple \emph{forward} way according to the recursion \eqref{eq:C_recursive}, which allows to easily calculate \emph{forward} conditional means $E(C_{i,j+1}|C_{i,j})$ and variances $Var(C_{i,j+1}|C_{i,j})$, it is \emph{not} straightforward to calculate \emph{backward} conditional means $E(C_{i,j}|C_{i,j+1})$ and variances $Var(C_{i,j}|C_{i,j+1})$; see Example 4.9 in \citet{jentschasympThInf}. %Hence, Assumption \ref{reverse_moments_condition} is required to control the backward conditional distributions of cumulative claims.

Based on Theorem 4.10 in \citet{jentschasympThInf}, which relies on conditional CLTs for (smooth functions of) development factor estimators $\widehat f_{j,n}$ given $\mathcal{Q}_{I,n}$ stated in \cite[Appendix C]{jentschasympThInf}, the following theorem provides the limiting distribution of the estimation uncertainty term $(R_{I,n}-\widehat{R}_{I,n})_2$ conditional on $\mathcal{Q}_{I,n}$. While $(R_{I,n}-\widehat{R}_{I,n})_2^{(1)}$ is measurable with respect to $\mathcal{Q}_{I,n}$, Assumptions \ref{support_condition} and \ref{reverse_moments_condition} %together 
allow to %establish 
prove asymptotic normality of $\sqrt{I+n+1}(R_{I,n}-\widehat{R}_{I,n})_2^{(2)}$ conditional on $\mathcal{Q}_{I,n}$.

\begin{thm}[Asymptotics for $(R_{I,n}-\widehat{R}_{I,n})_2$ conditional on $\mathcal{Q}_{I,n}$]\label{Dist_R2}
Suppose Assumptions \ref{initial_lem}, \ref{annF}, \ref{parameterAss}, \ref{support_condition} and \ref{reverse_moments_condition} hold. Then, as $n\rightarrow \infty$, the following holds:
\begin{itemize}
	\item[(i)] Unconditionally, $\sqrt{I+n+1}(R_{I,n}-\widehat{R}_{I,n})_2^{(1)}$ converges in distribution to a non-degenerate limiting distribution $\mathcal{G}_2^{(1)}$. That is, we have
	\begin{align}
	\sqrt{I+n+1}(R_{I,n}-\widehat{R}_{I,n})_2^{(1)} \overset{d}{\longrightarrow} \left\langle \mathcal{Q}_{I,\infty},\mathbf{Y}_\infty^{(1)}\right\rangle\sim\mathcal{G}_2^{(1)},
	\end{align}
	where $\mathbf{Y}_\infty^{(1)}=(Y_i^{(1)},i\in\mathbb{N}_0)$ denotes a centered Gaussian process with covariances
	\begin{align*}%\label{var_R_2_uncon}
	Cov(Y_{i_1}^{(1)},Y_{i_2}^{(1)}) = \lim_{K\rightarrow \infty}\bm{\Sigma}_{K,\prod \f}^{(1)}(i_1,i_2),	\quad	i_1,i_2\in\mathbb{N}_0,
	\end{align*}
	%for $i_1,i_2\in\mathbb{N}_0$, 
	where $\bm{\Sigma}_{K,\prod \f}^{(1)}(i_1,i_2)$ is defined in Corollary C.2 in \cite{jentschasympThInf}. Here, the two random sequences $\mathcal{Q}_{I,\infty}$ and $\mathbf{Y}_\infty^{(1)}$ are stochastically independent. 
	\item[(ii)] Conditionally on $\mathcal{Q}_{I,n}$, $\sqrt{I+n+1}(R_{I,n}-\widehat{R}_{I,n})_2^{(2)}$ converges in distribution to a centered normal distribution. That is, we have%\footnote{CJ: Muss hier noch (steht im ersten Paper auch nicht) noch die Bedingung $|\mathcal{Q}_{I,\infty}$ hinter $\left\langle \mathcal{Q}_{I,\infty},\mathbf{Y}_\infty^{(2)}\right\rangle$!? \tr{wäre konsistenter}}
	\begin{align}
	\sqrt{I+n+1}(R_{I,n}-\widehat{R}_{I,n})_2^{(2)}|\mathcal{Q}_{I,n}\overset{d}{\longrightarrow} \left\langle \mathcal{Q}_{I,\infty},\mathbf{Y}_\infty^{(2)}\right\rangle|\mathcal{Q}_{I,\infty}\sim\mathcal{G}_2^{(2)}|\mathcal{Q}_{I,\infty},
	\end{align}
where $\mathcal{G}_2^{(2)}|\mathcal{Q}_{I,\infty}\sim\mathcal{N}(0,\Xi(\mathcal{Q}_{I,\infty}))|\mathcal{Q}_{I,\infty}$ is Gaussian with mean zero and variance 
	\begin{align}
	\Xi(\mathcal{Q}_{I,\infty}) = \lim_{K\rightarrow\infty}
	\mathcal{Q}_{I,K-I}\bm{\Sigma}_{K,\prod \f}^{(2)}\mathcal{Q}_{I,K-I}^\prime = \lim_{K\rightarrow\infty} \mathcal{Q}_{I,K-I}(\bm{\Sigma}_{K,\prod \f}-\bm{\Sigma}_{K,\prod \f}^{(1)})\mathcal{Q}_{I,K-I}^\prime,
	\end{align}
	where $\bm{\Sigma}_{K,\prod \f}$ and $\bm{\Sigma}_{K,\prod \f}^{(1)}$ and $\bm{\Sigma}_{K,\prod \f}^{(2)}$ are defined in Corollary C.2 in \cite{jentschasympThInf}. %\ref{CLTProdf_cond}. 
\end{itemize}
\end{thm}

According to Theorem \ref{Dist_R2}(ii), in contrast to $\mathcal{G}_1|\mathcal{Q}_{I,\infty}$ in Theorem \ref{Dist_R1}, the %conditional 
limiting distribution $\mathcal{G}_2^{(2)}|\mathcal{Q}_{I,\infty}$ is Gaussian. Together with Theorem \ref{Dist_R2}(i), conditional on $\mathcal{Q}_{I,n}$, % and inflated with $\sqrt{I+n+1}$, 
the estimation uncertain term $\sqrt{I+n+1}(R_{I,n}-\widehat{R}_{I,n})_2$ is Gaussian with mean $\langle \mathcal{Q}_{I,\infty},\mathbf{Y}_\infty^{(2)}\rangle$ and variance $\Xi(\mathcal{Q}_{I,\infty})$. 

\subsubsection{Conditional asymptotics for the whole predictive root of the reserve}\label{sec_asymp_summary_joint}\

By combining the results
%limiting conditional distributions 
derived for $(R_{I,n}-\widehat R_{I,n})_1$ and $(R_{I,n}-\widehat R_{I,n})_2$ in Theorems \ref{Dist_R1} and \ref{Dist_R2}, respectively, joint asymptotic results for
 %the whole predictive root of the reserve 
$R_{I,n}-\widehat{R}_{I,n}$ conditional on $\mathcal{Q}_{I,n}$ can also be established. 
%Based on Theorem 4.12 and Corollary 4.13 from \citet{jentschasympThInf}, the following theorem provides the limiting distribution of the whole predictive root of the reserve $R_{I,n}-\widehat{R}_{I,n}$ conditional on $\mathcal{Q}_{I,n}$.

\begin{thm}[Asymptotics for $R_{I,n}-\widehat R_{I,n}$ conditional on $\mathcal{Q}_{I,n}$]\label{Dist_R1_R2_joint}	\
Suppose the assumptions of Theorems \ref{Dist_R1} and \ref{Dist_R2} hold. Then, conditional on $\mathcal{Q}_{I,n}$, $(R_{I,n}-\widehat R_{I,n})_1$ and $(R_{I,n}-\widehat R_{I,n})_2$ are stochastically independent, and $R_{I,n}-\widehat R_{I,n}|\mathcal{Q}_{I,n}$ converges in distribution to $\mathcal{G}_1|\mathcal{Q}_{I,\infty}$. That is, we have
\begin{align}	
	R_{I,n}-\widehat{R}_{I,n}|\mathcal{Q}_{I,n}=(R_{I,n}-\widehat R_{I,n})_1+(R_{I,n}-\widehat R_{I,n})_2|\mathcal{Q}_{I,n}\overset{d}{\longrightarrow}\mathcal{G}_1|\mathcal{Q}_{I,\infty}.
\end{align}
\end{thm}

According to Theorem \ref{Dist_R2}, $(R_{I,n}-\widehat{R}_{I,n})_2$ requires an inflation factor $\sqrt{I+n+1}$ to get convergence to a non-degenerate limiting distribution. As this is not the case for $(R_{I,n}-\widehat{R}_{I,n})_1$ in Theorem \ref{Dist_R1}, the %latter 
process uncertainty term $(R_{I,n}-\widehat{R}_{I,n})_1$ %corresponding to the process uncertainty will 
asymptotically dominates the predictive root of the reserve $R_{I,n}-\widehat{R}_{I,n}$.

Hence, we can conclude that asymptotic normality of the (predictive root of the) reserve does generally \emph{not} hold, which casts the common practice to use a normal approximation for the reserve in Mack's model into doubt. Moreover, the shape of $\mathcal{G}_1|\mathcal{Q}_{I,\infty}$ does depend on the true (conditional) distribution family of the individual development factors $F_{i,j}|C_{i,j}$.

\subsection{Conditional bootstrap asymptotics for the Mack bootstrap predictive root of the reserve}\label{sec_asymp_boot}\

In view of the decomposition $R_{I,n}-\widehat{R}_{I,n}=(R_{I,n}-\widehat{R}_{I,n})_1+(R_{I,n}-\widehat{R}_{I,n})_2$ in \eqref{eqpredroot_decomp_first_appearance} and the conditional limiting distributions of $(R_{I,n}-\widehat{R}_{I,n})_1$ and $(R_{I,n}-\widehat{R}_{I,n})_2$ gathered in Section \ref{sec_asymp_summary}, it is instructive to consider the corresponding Mack bootstrap quantities $(R_{I,n}^*-\widehat{R}_{I,n})_1$ and $(R_{I,n}^*-\widehat{R}_{I,n})_2$ from \eqref{eqpredroot_boot_decomp_appearance} and check whether they are correctly mimicking such limiting distributions. While $(R_{I,n}-\widehat{R}_{I,n})_1$ and $(R_{I,n}-\widehat{R}_{I,n})_2$ are analyzed conditional on $\mathcal{Q}_{I,n}$, the bootstrap quantities $(R_{I,n}^*-\widehat{R}_{I,n})_1$ and $(R_{I,n}^*-\widehat{R}_{I,n})_2$ have to be considered conditional on $\mathcal{Q}_{I,n}^*=\mathcal{Q}_{I,n}$, but also %(as common in the bootstrap literature) conditional 
on $\mathcal{D}_{I,n}$. %, that is, on the available cumulative claim data.

\subsubsection{Conditional bootstrap asymptotics for reserve prediction: process uncertainty}\label{sec_asymp_boot_process}\

For the derivation of bootstrap asymptotics, we have to impose additional smoothness properties of the parametric family of (conditional) distributions of the individual development factors to assure that consistent estimation of development factors and variance parameters implies also consistent estimation of the whole distribution. 

\begin{ann}[Parametric family of (conditional) distributions of $F_{i,j}$]\label{parametric_family}\
The (conditional) distribution $F_{i,j}|C_{i,j}$, $i\in\mathbb{Z}$, $i\leq I$, $j\in\mathbb{N}_0$, belongs to a parametric family of distributions $\mathbf{H}$ %, that fulfills the following properties:
such that:
\begin{itemize}
	\item[(i)] A distribution $\mathcal{H}\in\mathbf{H}$ is uniquely specified by its first two (conditional) moments. That is, for all %$j=0,\ldots,I+n-1$
$i\in\mathbb{Z}$, $i\leq I$, $j\in\mathbb{N}_0$	and all $c\in(0,\infty)$, the conditional distribution of $F_{i,j}|C_{i,j}=c$ is %(almost surely) 
uniquely determined by %its conditional mean 
$E(F_{i,j}|C_{i,j}=c)=f_j$ and %its conditional variance 
$Var(F_{i,j}|C_{i,j}=c)=\frac{\sigma_j^2}{c}$ according to \eqref{annFstMack}.
	\item[(ii)] The distributions $\mathcal{H}\in\mathbf{H}$ are continuous in %the parameters 
$f_j$ and $\sigma_j^2$. That is, for all %$j=0,\ldots,I+n-1$ 
$i\in\mathbb{Z}$, $i\leq I$, $j\in\mathbb{N}_0$ and for all $c\in(0,\infty)$, the conditional distribution of $F_{i,j}|C_{i,j}=c$ is continuous in a neighborhood of $(f_j,\sigma_j^2)$.
\end{itemize}
\end{ann}

As the limiting distribution derived in Theorem \ref{Dist_R1} is generally non-Gaussian and depends on the (conditional) distribution (family) of the individual development factors, we require also that the bootstrap individual development factors $F_{i,j}^*$, $i=-n, \dots, I$, $j=I-i,\dots, I+n-1$, that is, for the \emph{lower} triangle, follow the true parametric family of (conditional) distributions as the $F_{i,j}$'s according to Assumption \ref{parametric_family}.

\begin{ann}[(Conditional) distributions of $F_{i,j}^*$ in lower triangle]\label{parametric_family_bootstrap}
For any $n\in\mathbb{N}_0$, the (conditional) distribution of $F_{i,j}^*|C_{i,j}^*,\mathcal{D}_{I,n},\mathcal{F}_{I,n}^*$, $j=0,\ldots,I+n-1$ and $i=-n,\ldots,I-j-1$ in Step 5 of the Mack Bootstrap in Section \ref{subsec_boot} belongs to the true parametric family of (conditional) distributions $\mathbf{H}$ used to generate $F_{i,j}|C_{i,j}$ according to Assumption \ref{parametric_family}. That is, we have
\begin{align*}
F_{i,j}^*|(C_{i,j}^*=x,\widehat f_{j,n}^*=y,\widehat \sigma_{j,n}^2=z)\overset{d}{=}F_{i,j}|(C_{i,j}=x,f_j=y,\sigma_j^2=z)\quad	\text{for all}	\quad	(x,y,z)'\in(0,\infty).
\end{align*}
\end{ann}

Together with the setup of %the assumptions imposed in 
Theorem \ref{Dist_R1}, the Assumptions \ref{parametric_family} and \ref{parametric_family_bootstrap} allow to prove the following result. % theorem.

\begin{thm}[Bootstrap asymptotics for $(R_{I,n}^*-\widehat{R}_{I,n})_1$ conditional on $\mathcal{Q}_{I,n}^*=\mathcal{Q}_{I,n}$ and $\mathcal{D}_{I,n}$]\label{Dist_R1_boot_new}
Suppose Assumptions \ref{initial_lem}, \ref{annF}, \ref{parameterAss}, \ref{ann_highermoments_F}, \ref{parametric_family}, and \ref{parametric_family_bootstrap} hold. Then, as $n\rightarrow \infty$, conditionally on $\mathcal{Q}_{I,n}^*=\mathcal{Q}_{I,n}$ and $\mathcal{D}_{I,n}$, $(R_{I,n}^*-\widehat{R}_{I,n})_1$ converges in distribution to $\mathcal{G}_1|\mathcal{Q}_{I,\infty}$ in probability, which is the (limiting) distribution of $(R_{I,\infty}-\widehat{R}_{I,\infty})_1|\mathcal{Q}_{I,\infty}$ according to \eqref{eq_convergence_in_distribution} described in Theorem \ref{Dist_R1}. Moreover, for all $n\in\mathbb{N}_0$, it holds $E^*((R_{I,n}^*-\widehat{R}_{I,n})_1|\mathcal{Q}_{I,n}^*=\mathcal{Q}_{I,n})=0$ and, for $n\rightarrow \infty$, we have
\begin{align}
Var^*\left((R_{I,n}^*-\widehat{R}_{I,n})_1|\mathcal{Q}_{I,n}^*=\mathcal{Q}_{I,n}\right)\longrightarrow Var\left((R_{I,\infty}-\widehat{R}_{I,\infty})_1|\mathcal{Q}_{I,\infty}\right)
\end{align}
in probability, where $Var((R_{I,\infty}-\widehat{R}_{I,\infty})_1|\mathcal{Q}_{I,\infty})=O_P(1)$ as given in \eqref{var_R1_cond}. Consequently, as $n\rightarrow\infty$, we have
\begin{align*}
d_2\left(\mathcal{L}\left((R_{I,n}-\widehat{R}_{I,n})_1|\mathcal{Q}_{I,n}\right), \mathcal{L}^*\left((R^*_{I,n}-\widehat{R}_{I,n})_1|\mathcal{Q}_{I,n}^*=\mathcal{Q}_{I,n}\right)\right)\longrightarrow 0
\end{align*}
in probability, where $\mathcal{L}^*(\cdot)$ denotes a bootstrap distribution conditional on $\mathcal{D}_{I,n}$, and $d_2$ is the Mallows metric, that is defined for two %random variables $X$ and $Y$ with 
distributions $G$ and $H$ %, respectively, 
as $d_2(G,H) = inf(E||X-Y||^2)^{\frac{1}{2}}$, where the infimum is taken over all joint distributions of $(X,Y)$ with marginals $X\sim G$ and $Y\sim H$.
\end{thm}

\subsubsection{Conditional bootstrap asymptotics for reserve prediction: estimation uncertainty}\label{subsec_asymp_boot_process}\

In view of the decomposition $(R_{I,n}-\widehat{R}_{I,n})_2=(R_{I,n}-\widehat{R}_{I,n})_2^{(1)}+(R_{I,n}-\widehat{R}_{I,n})_2^{(2)}$ in \eqref{eqpredroot_decomp}, for the derivation of corresponding bootstrap asymptotic theory, it is \emph{seemingly} instructive to further decompose also its bootstrap counterpart $(R_{I,n}^*-\widehat{R}_{I,n})_2$ in the same way conditional on $\mathcal{Q}_{I,n}^*=\mathcal{Q}_{I,n}$ and $\mathcal{D}_{I,n}$. That is, by taking into account the specific definition of $\widehat f_{j,n}^*$ in \eqref{eqf_boot}, %which relies on the original $C_{i,j}$'s and adds just one bootstrap individual development factor in each summand, 
we get
\begin{align}
& (R_{I,n}^*-\widehat{R}_{I,n})_2	\nonumber	\\
=& \sum^{I+n}_{i=0}C_{I-i,i}^*\left(\prod_{j=i}^{I+n-1}\widehat f_{j,n}^*-\prod_{j=i}^{I+n-1}f_{j,n}^*(\mathcal{Q}_{I,n})\right) + \sum^{I+n}_{i=0}C_{I-i,i}^*\left(\prod_{j=i}^{I+n-1}f_{j,n}^*(\mathcal{Q}_{I,n})-\prod_{j=i}^{I+n-1}\widehat{f}_{j,n}\right)	\nonumber	\\
=& (R_{I,n}^*-\widehat{R}_{I,n})_2^{(1)}+(R_{I,n}^*-\widehat{R}_{I,n})_2^{(2)},	\label{eqpredroot_boot_decomp}
\end{align}
where $f_{j,n}^*(\mathcal{Q}_{I,n}):=\mu_{j+1,n}^{*(1)}(\mathcal{Q}_{I,n})/\mu_{j,n}^{*(2)}(\mathcal{Q}_{I,n})$ with $\mu_{j+1,n}^{*(1)}(\mathcal{Q}_{I,n}) := E^*(\frac{1}{I+n-j}\sum_{i=-n}^{I-j-1}C_{i,j}F_{i,j}^*|\mathcal{Q}_{I,n}^*=\mathcal{Q}_{I,n})$ and $\mu_{j,n}^{*(2)}(\mathcal{Q}_{I,n}) := E^*(\frac{1}{I+n-j}\sum_{i=-n}^{I-j-1}C_{i,j}|\mathcal{Q}_{I,n}^*=\mathcal{Q}_{I,n})$.
%\begin{align*}
%\mu_{j+1,n}^{*(1)}(\mathcal{Q}_{I,n}) &:= E^*\left(\frac{1}{I+n-j}\sum_{i=-n}^{I-j-1}C_{i,j}F_{i,j}^*|\mathcal{Q}_{I,n}^*=\mathcal{Q}_{I,n}\right)	\\%=\frac{1}{I+n-j}\sum_{i=-n}^{I-j-1} E^*(C_{i,j}F_{i,j}^*|C_{i,I-i}^*=C_{i,I-i}),	\\
%\mu_{j,n}^{*(2)}(\mathcal{Q}_{I,n}) &:= E^*\left(\frac{1}{I+n-j}\sum_{i=-n}^{I-j-1}C_{i,j}|\mathcal{Q}_{I,n}^*=\mathcal{Q}_{I,n}\right) %= \frac{1}{I+n-j}\sum_{i=-n}^{I-j-1} C_{i,j}.
%\end{align*}
%\begin{align*}
%\mu_{j+1,n}^{*(1)}(\mathcal{Q}_{I,n}) &:= E^*\left(\frac{1}{I+n-j}\sum_{i=-n}^{I-j-1}C_{i,j}F_{i,j}^*|\mathcal{Q}_{I,n}^*=\mathcal{Q}_{I,n}\right)=\frac{1}{I+n-j}\sum_{i=-n}^{I-j-1} E^*(C_{i,j}F_{i,j}^*|C_{i,I-i}^*=C_{i,I-i}),	\\
%\mu_{j,n}^{*(2)}(\mathcal{Q}_{I,n}) &:= E^*\left(\frac{1}{I+n-j}\sum_{i=-n}^{I-j-1}C_{i,j}|\mathcal{Q}_{I,n}^*=\mathcal{Q}_{I,n}\right) = \frac{1}{I+n-j}\sum_{i=-n}^{I-j-1} C_{i,j}.
%\end{align*}
Now, for $\mu_{j+1,n}^{*(1)}(\mathcal{Q}_{I,n})$, we get
\begin{align*}
\mu_{j+1,n}^{*(1)}(\mathcal{Q}_{I,n}) = \frac{1}{I+n-j}\sum_{i=-n}^{I-j-1} C_{i,j}E^*(F_{i,j}^*|C_{i,I-i}^*=C_{i,I-i})=\frac{1}{I+n-j}\sum_{i=-n}^{I-j-1} C_{i,j}E^*(F_{i,j}^*)
%&= \left(\frac{1}{I+n-j}\sum_{i=-n}^{I-j-1} C_{i,j}\right)\widehat f_{j,n}, %= \widehat f_{j,n}\left(\frac{1}{I+n-j}\sum_{i=-n}^{I-j-1} C_{i,j}\right),
\end{align*}
where we used that $C_{i,j}$ is measurable with respect to $\mathcal{D}_{I,n}$ and that $F_{i,j}^*$ is stochastically independent of the condition $C_{i,I-i}^*=C_{i,I-i}$ given $\mathcal{D}_{I,n}$. %that $F_{i,j}^*$ is generated independently from the diagonal. 
This is because the Mack bootstrap relies on a fixed-design approach based on the $C_{i,j}$'s instead of recursively generating $C_{i,j}^*$ to get a whole bootstrap loss triangle $\mathcal{D}_{I,n}^*$. 
%Hence, $F_{i,j}^*$ is stochastically independent of the condition $C_{i,I-i}^*=C_{i,I-i}$ given $\mathcal{D}_{I,n}$ for all $i=-n,\ldots,I-j-1$. 
Altogether, using $E^*(F_{i,j}^*)=\widehat f_{j,n}$,  %conditional on $\mathcal{Q}_{I,n}^*=\mathcal{Q}_{I,n}$ and $\mathcal{D}_{I,n}$, 
we get
\begin{align*}
f_{j,n}^*(\mathcal{Q}_{I,n}) = \frac{\widehat f_{j,n}\left(\frac{1}{I+n-j}\sum_{i=-n}^{I-j-1} C_{i,j}\right)}{\frac{1}{I+n-j}\sum_{i=-n}^{I-j-1} C_{i,j}} = \widehat f_{j,n}
\end{align*}
leading to $(R_{I,n}^*-\widehat{R}_{I,n})_2^{(2)}=0$ such that $(R_{I,n}^*-\widehat{R}_{I,n})_2=(R_{I,n}^*-\widehat{R}_{I,n})_2^{(1)}$. Hence, in comparison to $(R_{I,n}-\widehat{R}_{I,n})_2$, which was decomposed into two %additive 
parts $(R_{I,n}-\widehat{R}_{I,n})_2^{(1)}$ and $(R_{I,n}-\widehat{R}_{I,n})_2^{(2)}$, such an analogous decomposition of $(R_{I,n}^*-\widehat{R}_{I,n})_2$ does \emph{not} exist.
 %Hence, in contrast to the limiting behavior of $\sqrt{I+n+1}(R_{I,n}-\widehat{R}_{I,n})_2$, which was derived after further decomposing $(R_{I,n}-\widehat{R}_{I,n})_2$ into two %additive 
%parts $(R_{I,n}-\widehat{R}_{I,n})_2^{(1)}$ and $(R_{I,n}-\widehat{R}_{I,n})_2^{(2)}$ %as stated 
%in Theorem \ref{Dist_R2}, %such 
%an analogous decomposition of $(R_{I,n}^*-\widehat{R}_{I,n})_2$ does \emph{not} exist. 
However, for $n\rightarrow \infty$, it remains to check the limiting properties of $(R_{I,n}^*-\widehat{R}_{I,n})_2$ %, which are presented in the following theorem.
in the following. %For this purpose, it is important to note that, 
In contrast to the derivation of the \emph{conditional} limiting result obtained in Theorem \ref{Dist_R2}(ii), which relies on conditional CLTs for the development factor estimators $\widehat f_{j,n}$ as stated in \cite[Appendix C]{jentschasympThInf}, the derivation of the limiting properties of $(R_{I,n}^*-\widehat{R}_{I,n})_2$ rely on \emph{unconditional} bootstrap CLTs for the Mack bootstrap development factor estimators $\widehat f_{j,n}^*$, that is, \emph{without} conditioning on $\mathcal{Q}_{I,n}^*=\mathcal{Q}_{I,n}$. For this purpose, %Nevertheless, 
to prove asymptotic normality for the $\widehat f_{j,n}^*$'s by justifying a Lyapunov condition, we have to impose additional regularity conditions on the estimators for the development factors and variance parameters. %This is required to control the behavior of the non-parametric bootstrap in Step 3 of the Mack bootstrap in Section \ref{subsec_boot}, which leads to bootstrap errors $r_{i,j}^*$ that are drawn randomly from \emph{all} (re-centered and re-scaled) residuals $\widetilde r_{i,j}$, that is, for all $j=0,\ldots,I+n-2$ and $i=-n,\ldots,I-j-1$.

\begin{ann}[Uniform boundedness condition]\label{uniform_boundedness_condition}
%For $n\rightarrow \infty$, 
Suppose that the development factor estimators $\widehat f_{j,n}$, $j=0,\ldots,I+n-1$ and the variance parameter estimators $\widehat \sigma_{j,n}^2$, $j=0,\ldots,I+n-2$ fulfill %the uniform boundedness conditions
\begin{align*}
\sup_{j=0,\ldots,I+n-1}\frac{\widehat f_{j,n}}{f_j}=O_P(1)	\quad	\text{and}	\quad	\sup_{j=0,\ldots,I+n-2}\frac{\sigma_j^2}{\widehat \sigma_{j,n}^2}=O_P(1)
\end{align*}
for $n\rightarrow \infty$. Moreover, for $\kappa_j^{(4)}$ defined in \eqref{kappa4}, suppose that $((\kappa_j^{(4)}/\sigma_j^4), j\in\mathbb{N}_0)$ is a bounded sequence.
\end{ann}

This allows for the following asymptotic result for the Mack bootstrap estimation uncertainty part.

\begin{thm}[Bootstrap asymptotics for $(R_{I,n}^*-\widehat{R}_{I,n})_2$ conditional on $\mathcal{Q}_{I,n}^*=\mathcal{Q}_{I,n}$ and $\mathcal{D}_{I,n}$]\label{Dist_R2_boot_new}
Suppose Assumptions \ref{initial_lem}, \ref{annF}, \ref{ann_highermoments_F}, \ref{parameterAss}, \ref{support_condition} and \ref{uniform_boundedness_condition} hold. 
Then, as $n\rightarrow \infty$, conditionally on $\mathcal{Q}_{I,n}^*=\mathcal{Q}_{I,n}$ and $\mathcal{D}_{I,n}$, $\sqrt{I+n+1}(R_{I,n}^*-\widehat{R}_{I,n})_2$ converges in distribution to $\mathcal{\widetilde G}_2|\mathcal{Q}_{I,\infty}$ in probability, where $\mathcal{\widetilde G}_2|\mathcal{Q}_{I,\infty}\sim \mathcal{N}(0, \widetilde \Xi(\mathcal{Q}_{I,\infty}))|\mathcal{Q}_{I,\infty}$ is a conditional Gaussian distribution with conditional mean zero and conditional variance
\begin{align}
\widetilde \Xi(\mathcal{Q}_{I,\infty}) = \lim_{K\rightarrow\infty} \mathcal{Q}_{I,K-I}\bm{\Sigma}_{K,\prod \f}\mathcal{Q}_{I,K-I}^\prime,
\end{align}
where $\bm{\Sigma}_{K,\prod \f}=\bm{\Sigma}_{K,\prod \f}^{(1)}+\bm{\Sigma}_{K,\prod \f}^{(2)}$ is defined in Corollary C.2 in \cite{jentschasympThInf}.
%\ref{CLTProdf_cond} and
\begin{align*}%\label{eq_var_R_2_con}
\widetilde \Xi(\mathcal{Q}_{I,\infty})=\sum_{i_1,i_2=0}^\infty C_{I-i_1,i_1}C_{I-i_2,i_2}\sum^\infty_{j=\max(i_1,i_2)}\frac{\sigma^2_j}{\mu_j}\prod^{\infty}_{l=\max(i_1,i_2), l \neq j}f^2_l\prod^{\max(i_1,i_2)-1}_{m=\min(i_1,i_2)}f_m=O_P(1).
\end{align*}
Consequently, as $n\rightarrow\infty$, we have
\begin{align*}
d_2\left(\mathcal{L}\left(\sqrt{I+n+1}(R_{I,n}-\widehat{R}_{I,n})_2|\mathcal{Q}_{I,n}\right), \mathcal{L}^*\left(\sqrt{I+n+1}(R^*_{I,n}-\widehat{R}_{I,n})_2|\mathcal{Q}_{I,n}^*=\mathcal{Q}_{I,n}\right)\right)\not\rightarrow 0	\quad	\text{in prob.,}
\end{align*}
%in probability. 
%This is 
because the limiting normal distribution of $\sqrt{I+n+1}(R^*_{I,n}-\widehat{R}_{I,n})_2$ conditional on $\mathcal{Q}_{I,n}^*=\mathcal{Q}_{I,n}$ and $\mathcal{D}_{I,n}$ deviates in its (zero) mean and %its 
variance $\widetilde \Xi(\mathcal{Q}_{I,\infty})$ from %the limiting distribution 
that of $\sqrt{I+n+1}(R_{I,n}-\widehat{R}_{I,n})_2$ conditional on $\mathcal{Q}_{I,n}$, which has mean $\langle \mathcal{Q}_{I,\infty},\mathbf{Y}_\infty^{(2)}\rangle$ and variance $\Xi(\mathcal{Q}_{I,\infty})$ according to Theorem \ref{Dist_R2}.
\end{thm}

\subsubsection{Conditional bootstrap asymptotics for the whole predictive root of the reserve}\label{sec_asymp_boot_joint}\

As in Section \ref{sec_asymp_summary_joint}, combining the results for $(R_{I,n}^*-\widehat R_{I,n})_1$ and $(R_{I,n}^*-\widehat R_{I,n})_2$ from Theorems \ref{Dist_R1_boot_new} and \ref{Dist_R2_boot_new}, respectively, joint asymptotics for $R_{I,n}^*-\widehat{R}_{I,n}$ conditional on $\mathcal{Q}_{I,n}^*=\mathcal{Q}_{I,n}$ and $\mathcal{D}_{I,n}$ can be obtained. 
%As in Section \ref{sec_asymp_summary_joint}, we can combine the results on the limiting distributions for $(R_{I,n}^*-\widehat R_{I,n})_1$ and $(R_{I,n}^*-\widehat R_{I,n})_2$ from Theorems \ref{Dist_R1_boot_new} and \ref{Dist_R2_boot_new}, respectively, to get the limiting bootstrap distribution of the whole bootstrap predictive root of the reserve $R_{I,n}^*-\widehat{R}_{I,n}$ conditional on $\mathcal{Q}_{I,n}^*=\mathcal{Q}_{I,n}$ and $\mathcal{D}_{I,n}$.

\begin{thm}[Bootstrap asymptotics for $R_{I,n}^*-\widehat R_{I,n}$ conditional on $\mathcal{Q}_{I,n}^*=\mathcal{Q}_{I,n}$ and $\mathcal{D}_{I,n}$]\label{Dist_R1_R2_joint_boot}
Suppose the assumptions of Theorems \ref{Dist_R1_boot_new} and \ref{Dist_R2_boot_new} hold. Then, conditional on $\mathcal{Q}_{I,n}^*=\mathcal{Q}_{I,n}$ and $\mathcal{D}_{I,n}$, $(R_{I,n}^*-\widehat R_{I,n})_1$ and $(R_{I,n}^*-\widehat R_{I,n})_2$ are uncorrelated, and 
 $R_{I,n}^*-\widehat R_{I,n}|(\mathcal{Q}_{I,n}^*=\mathcal{Q}_{I,n},\mathcal{D}_{I,n})$ converges in distribution to $\mathcal{G}_1|\mathcal{Q}_{I,\infty}$. That is, we have
\begin{align*}	
R_{I,n}^*-\widehat{R}_{I,n}|\left(\mathcal{Q}_{I,n}^*=\mathcal{Q}_{I,n},\mathcal{D}_{I,n}\right)=(R_{I,n}^*-\widehat R_{I,n})_1+(R_{I,n}^*-\widehat R_{I,n})_2|\left(\mathcal{Q}_{I,n}^*=\mathcal{Q}_{I,n},\mathcal{D}_{I,n}\right)\overset{d}{\longrightarrow}\mathcal{G}_1|\mathcal{Q}_{I,\infty}
\end{align*}
in probability.
\end{thm}

As already observed in Theorem \ref{Dist_R2} for the estimation uncertainty term $(R_{I,n}-\widehat{R}_{I,n})_2$, its Mack bootstrap version requires also an inflation factor $\sqrt{I+n+1}$ to establish convergence towards a non-degenerate limiting distribution. As this is not the case for the process uncertainty term $(R_{I,n}-\widehat{R}_{I,n})_1$ in Theorem \ref{Dist_R1} and its bootstrap version in Theorem \ref{Dist_R1_boot_new}, the process uncertainty terms will asymptotically dominate the predictive roots $R_{I,n}-\widehat{R}_{I,n}$ and $R_{I,n}^*-\widehat{R}_{I,n}$. Hence, although the limiting bootstrap distribution of $\sqrt{I+n+1}(R_{I,n}^*-\widehat{R}_{I,n})_2$ conditional on $\mathcal{Q}_{I,n}^*=\mathcal{Q}_{I,n}$ and $\mathcal{D}_{I,n}$ in Theorem \ref{Dist_R2_boot_new} does \emph{not} correctly mimic the corresponding limiting behavior of $\sqrt{I+n+1}(R_{I,n}-\widehat{R}_{I,n})_2$ conditional on $\mathcal{Q}_{I,n}$ in Theorem \ref{Dist_R2}, the whole bootstrap predictive root $R_{I,n}^*-\widehat{R}_{I,n}$ still mimics the limiting distribution of the predictive root $R_{I,n}-\widehat{R}_{I,n}$ correctly. 

Hence, in view of the concepts of \emph{asymptotic validity} and \emph{asymptotic pertinence} of a bootstrap prediction approach discussed in \cite{pan2016bootstrap}, the Mack bootstrap can be regarded as asymptotically valid, but \emph{not} as asymptotically pertinent under the stated conditions.

%This motivates the construction of an alternative Mack-type bootstrap proposed in the following section. We conclude this section with some remarks.

\begin{rem}[On the asymptotic results for the Mack bootstrap]\ 
\begin{itemize}
	\item[(i)] A closer inspection of the decompositions in \eqref{eqpredroot_decomp_first_appearance} and \eqref{eqpredroot_boot_decomp_appearance} reveals some inconsistencies:
	\begin{itemize}
		\item While a term based on products of %development parameters 
		$f_j$'s is added to and subtracted from $R_{I,n}-\widehat R_{I,n}$ to get $(R_{I,n}-\widehat R_{I,n})_1$ and $(R_{I,n}-\widehat R_{I,n})_2$, a term using products of %bootstrap development factor estimators 
$\widehat f_{j,n}^*$'s instead of % the more natural choice of development factor estimators 
$\widehat f_{j,n}$'s, which would be the natural choice, is added to and subtracted from $R_{I,n}^*-\widehat R_{I,n}$ to get $(R_{I,n}^*-\widehat R_{I,n})_1$ and $(R_{I,n}^*-\widehat R_{I,n})_2$.   
		\item Consequently, while $(R_{I,n}-\widehat R_{I,n})_1$ relies on products of %individual development factors 
$F_{i,j}$'s centered around products of %development parameters 
$f_j$'s, its Mack bootstrap version $(R_{I,n}^*-\widehat R_{I,n})_1$ relies on products of % bootstrap individual development factors 
$F_{i,j}^*$'s, which are not (naturally) centered around %development factor estimators 
products of $\widehat f_{j,n}$'s, but %around bootstrap quantities 
around products of $\widehat f_{j,n}^*$'s.
		\item Moreover, while $(R_{I,n}-\widehat R_{I,n})_2$ relies on differences between products of development parameters $f_j$ and products of their estimators $\widehat f_{j,n}$, its Mack bootstrap version $(R_{I,n}^*-\widehat R_{I,n})_2$ relies on differences between products of bootstrap development factor estimators $\widehat f_{j,n}^*$ and products of estimators $\widehat f_{j,n}$. Hence, the sign of $(R_{I,n}^*-\widehat R_{I,n})_2$ is flipped in comparison to $(R_{I,n}-\widehat R_{I,n})_2$. This may have a negative effect in finite samples, but as the limiting conditional distribution is Gaussian and hence symmetric, this will not be an issue asymptotically.
		\item According to the latter observation, also the terms $(R_{I,n}^*-\widehat R_{I,n})_2^{(1)}$ and $(R_{I,n}^*-\widehat R_{I,n})_2^{(2)}$ in the seemingly natural decomposition of the bootstrap estimation uncertainty term in \eqref{eqpredroot_boot_decomp} are switched in comparison to $(R_{I,n}-\widehat R_{I,n})_2^{(1)}$ and $(R_{I,n}-\widehat R_{I,n})_2^{(2)}$.
	\end{itemize}

	\item[(ii)] The bootstrap consistency result for the Mack bootstrap process uncertainty part conditional on $\mathcal{Q}_{I,n}$ in Theorem \ref{Dist_R1_boot_new} requires the correct choice of the true family of (conditional) distributions of the $F_{i,j}$'s also for the $F_{i,j}^*$'s in Step 5 of Section \ref{subsec_boot}. Otherwise, only the first and second moments of the conditional distribution will be correctly mimicked asymptotically, but %in general not 
not necessarily the whole distribution. % and, consequently, also not the quantiles.

	\item[(iii)] The uniform boundedness conditions in Assumption \ref{uniform_boundedness_condition} are required to establish a Lyapunov condition for bootstrap CLTs for (smooth functions of) $\widehat f_{j,n}^*$ %in Theorem \ref{CLT_f_boot}, 
because the Mack bootstrap draws bootstrap errors $r_{i,j}^*$ from residuals computed from all columns in $\mathcal{D}_{I,n}$. 
	
	%\item[(iv)] In contrast to $(R_{I,n}-\widehat R_{I,n})_1$ and $(R_{I,n}-\widehat R_{I,n})_2$, which are independent conditional on $\mathcal{Q}_{I,n}$, both parts $(R_{I,n}^*-\widehat R_{I,n})_1$ and $(R_{I,n}^*-\widehat R_{I,n})_2$ of the bootstrap predictive root are in general only uncorrelated conditional on $\mathcal{Q}_{I,n}^*=\mathcal{Q}_{I,n}$ and $\mathcal{D}_{I,n}$.
	
	\item[(iv)] The bootstrap inconsistency result for the Mack bootstrap estimation uncertainty part conditional on $\mathcal{Q}_{I,n}^*=\mathcal{Q}_{I,n}$ and $\mathcal{D}_{I,n}$ %established 
in Theorem \ref{Dist_R2_boot_new} is %due to the fact that 
because the bootstrap approach in Step 4 is not taking the condition $\mathcal{Q}_{I,n}^*=\mathcal{Q}_{I,n}$ into account. Hence, the (always larger!) variance-covariance matrix $\bm{\Sigma}_{K,\prod \f}$ shows in the conditional limiting distribution instead of $\bm{\Sigma}_{K,\prod \f}^{(2)}$ obtained in Theorem \ref{Dist_R2}. Moreover, a decomposition of $(R_{I,n}^*-\widehat R_{I,n})_2$ resembling the decomposition of $(R_{I,n}-\widehat R_{I,n})_2$ in \eqref{eqpredroot_decomp} does not exist. %Consequently, the behavior of $(R_{I,n}-\widehat R_{I,n})_2^{(1)}$ is also not correctly mimicked.
	
	\item[(v)] The requirement of a bootstrap procedure to not only mimic the asymptotically dominating part of the (conditional) predictive distribution that captures the prediction (i.e.~process) uncertainty (asymptotic validity), but also the asymptotically negligible part capturing the uncertainty due to model parameter estimation is closely related to the concept coined asymptotic pertinence in \cite{pan2016bootstrap} for time series prediction, which is also discussed by \cite{beutner2021} from a slightly different perspective. %\cite{pan2016bootstrap} argue that asymptotic validity of predictive inference is a fundamental property, but capturing the uncertainty due to model estimation is beneficial in finite samples.

	%\item[(vi)] The discussion above motivates an alternative notion of a Mack-type bootstrap to be introduced in Section \ref{sec_alternative_boot} that will be designed to eliminate the raised issues. In particular, it respects the conditioning on $\mathcal{Q}_{I,n}^*=\mathcal{Q}_{I,n}$ and generates a bootstrap loss triangle $\mathcal{D}_{I,n}^*$ in a backward manner starting from the diagonal $\mathcal{Q}_{I,n}^*=\mathcal{Q}_{I,n}$. See e.g.~\cite{PaparoditisShang2021} for bootstrap predictive inference in a functional time series setup. 
\end{itemize}
\end{rem}

The discussion above motivates an alternative notion of a Mack-type bootstrap to be introduced in the following section that is designed to eliminate the raised issues. In particular, it should respect the \emph{conditioning on $\mathcal{Q}_{I,n}^*=\mathcal{Q}_{I,n}$} and it should generate a whole \emph{bootstrap loss triangle $\mathcal{D}_{I,n}^*$} in a \emph{backward} manner starting from the diagonal $\mathcal{Q}_{I,n}^*=\mathcal{Q}_{I,n}$. See e.g.~\cite{PaparoditisShang2021} for bootstrap predictive inference in a functional time series setup.

%\bigskip

\section{An alternative Mack-type Bootstrap Scheme}\label{sec_alternative_boot}

According to the findings and the discussion in Section \ref{sec_boot_consistency}, the original Mack bootstrap proposal is not capable of mimicking the conditional distribution of the estimation uncertainty part correctly. Although it is asymptotically dominated by the process uncertainty part, it is generally desirable to construct a Mack-type bootstrap that addresses this issue to enable a better finite sample performance.

For this purpose, we propose an alternative Mack-type bootstrap in this section to mimic the distribution of the predictive root of the reserve $R_{I,n}-\widehat{R}_{I,n}$ using an \emph{alternative} bootstrap predictive root of the reserve $R_{I,n}^+-\widehat{R}_{I,n}^+$ to be defined below. To distinguish it from the original Mack bootstrap proposal in Section \ref{sec_boot}, we denote all related bootstrap quantities and operations with a ``$+$'' instead of a ``$*$''. This novel approach deviates from the original Mack bootstrap scheme from Section \ref{sec_boot} in several ways:

\begin{itemize}
	\item[(i)] First, given the loss triangle $\mathcal{D}_{I,n}$ and conditional on %the bootstrap diagonal 
$\mathcal{Q}_{I,n}^+=\mathcal{Q}_{I,n}$, where $\mathcal{Q}_{I,n}^+=\{C_{I-i,i}^+|i=0,\ldots,I+n\}$, a recursive \emph{backward} bootstrap approach is employed to generate a whole bootstrap \emph{upper triangle}
\begin{align*}	
\mathcal{D}_{I,n}^+=\left\{C_{i,j}^+|i=-n,\ldots,I,~j=0,\ldots,I+n,~-n\leq i+j\leq I\right\}.
\end{align*}
%Then, bootstrap estimators for the development factors $\widehat f_{j,n}^+$, $j=0,\ldots,I+n-1$ are computed according to formula \eqref{eqf}, but based on the bootstrap upper loss triangle $\mathcal{D}_{I,n}^+$.
	\item[(ii)] Second, %instead of the bootstrap development factor estimators $\widehat f_{j,n}^+$ computed from $\mathcal{D}_{I,n}^+$, 
	the development factor estimators $\widehat f_{j,n}$ computed from $\mathcal{D}_{I,n}$ are used for a \emph{parametric} bootstrap to construct bootstrap individual development factors $F_{i,j}^+$, $i=-n,\ldots,I$, $j=0,\ldots,I+n-1$ and $i+j\geq I$, that is for the \emph{lower triangle}, which also allows to construct $R_{I,n}^+$.
	\item[(iii)] Third, for the construction of the bootstrap predictive root of the reserve $R_{I,n}^+-\widehat{R}_{I,n}^+$, the bootstrap reserve $R_{I,n}^+$ is \emph{not} centered around its best estimate $\widehat{R}_{I,n}$, but around a suitable bootstrap version $\widehat{R}_{I,n}^+$.
\end{itemize}

%Finally, 
Analogous to the original Mack bootstrap, the alternative Mack bootstrap is employed to estimate the conditional distribution of $R_{I,n}-\widehat{R}_{I,n}$ given $\mathcal{Q}_{I,n}$ %is estimated 
by the conditional bootstrap distribution of $R_{I,n}^+-\widehat{R}_{I,n}^+$ given $\mathcal{Q}_{I,n}^+=\mathcal{Q}_{I,n}$ and $\mathcal{D}_{I,n}$. 

\subsection{An alternative Mack-type Bootstrap Algorithm}\label{subsec_alternative_boot}\

With the upper triangle $\mathcal{D}_{I,n}$ at hand, the alternative Mack-type bootstrap algorithm is defined as follows:

\begin{itemize}
	\item[Step 1.] Estimate the development factors $f_j$ and the variance parameters $\sigma_j^2$ from the data by computing $\widehat f_{j,n}$ and $\widehat \sigma^2_{j,n}$ for $j=0, \dots, I+n-1$ as defined in \eqref{eqf} and \eqref{def_sigma}, respectively.
	
	\item[Step 2.] Choose a parametric family for the (conditional) bootstrap distributions of the \emph{backward} individual development factors $G_{i,j}^+$ given $C_{i,j+1}^+$ and $\mathcal{D}_{I,n}$ such that $G_{i,j}^+>0$ a.s.~with \begin{align*}
E^+(G_{i,j}^+|C_{i,j+1}^+)=\widehat{f}_{j,n}^{-1},	\quad	Var^+(G_{i,j}^+|C_{i,j+1}^+)=\frac{\widehat{\sigma}^2_{j,n}}{C_{i,j+1}^+},	\quad	j=0, \dots, I+n-1,\ i=-n, \dots, I-j-1.
\end{align*}
Then, given $\mathcal{Q}_{I,n}^+=\mathcal{Q}_{I,n}$, generate \emph{backwards} a bootstrap loss triangle $\mathcal{D}_{I,n}^+$ using the recursion
\begin{align*}
C_{i,j}^+=C_{i,j+1}^+G_{i,j}^+,	\quad	j=0, \dots, I+n-1,\ i=-n, \dots, I-j-1.
\end{align*}	

	\item[Step 3.] Compute bootstrap development factor estimators % $f_j$ by computing 
$\widehat f_{j,n}^+$ for $j=0, \dots, I+n-1$, which are defined %analogously to 
as $\widehat f_{j,n}$ %defined 
in \eqref{eqf}, but %$\widehat f_{j,n}^+$ is 
are calculated from the bootstrap loss triangle $\mathcal{D}_{I,n}^+$. That is, we compute
\begin{align}\label{eq_def_f+}
\widehat f_{j,n}^+=\frac{\sum^{I-j-1}_{i=-n}C_{i,j+1}^+}{\sum^{I-j-1}_{i=-n}C_{i,j}^+}=\frac{\sum^{I-j-1}_{i=-n}C_{i,j+1}^+}{\sum^{I-j-1}_{i=-n}C_{i,j+1}^+G_{i,j}^+},	\quad	j=0, \dots, I+n-1.
\end{align}
%for $j=0, \dots, I+n-1$. 
	\item[Step 4.] Choose a parametric family for the (conditional) bootstrap distributions of $\F^+$ given $\C^+$ and $\mathcal{D}_{I,n}$ such that $\F^+>0$ a.s.~with 
\begin{align*}
E^+(F_{i,j}^+|C_{i,j}^+)=\widehat{f}_{j,n},	\quad	Var^+(F_{i,j}^+|C_{i,j}^+)=\frac{\widehat{\sigma}^2_{j,n}}{C_{i,j}^+},	\quad	i=-n, \dots, I,\ j=I-i,\dots, I+n-1.
\end{align*}
Then, given $\mathcal{Q}_{I,n}^+=\mathcal{Q}_{I,n}$, generate the bootstrap ultimate claims $C_{i,I+n}^+$ and the reserves $R_{i,I+n}^+=C_{i,I+n}^+-C_{i,I-i}^+$ for $i=-n,\dots,I$ using the recursion
\begin{align}
C_{i,j+1}^+=C_{i,j}^+F_{i,j}^+,	\quad	j=I-i,\ldots,I+n-1.
\end{align}
	
	\item[Step 5.] Compute the bootstrap total reserve $R_{I,n}^+=\sum_{i=-n}^I R_{i,I+n}^+$ and the alternative Mack bootstrap predictive root of the reserve
	\begin{align*}	
	R^+_{I,n}-\widehat{R}^+_{I,n} = \sum^{I+n}_{i=0}C_{I-i,i}^+\left(\prod^{I+n-1}_{j=i}F_{I-i,i}^+-\prod^{I+n-1}_{j=i}\widehat{f}_{j,n}^+\right),
\end{align*}
where the centering term $\widehat{R}^+_{I,n}$ is a bootstrap version of the best estimate $\widehat{R}_{I,n}$, that is defined by
\begin{align}\label{centering_term_alternative}
\widehat{R}^+_{I,n}=\sum^I_{i=-n}C_{i,I-i}^+\prod^{I+n-1}_{j=I-i}\widehat{f}_{j,n}^+ = \sum^{I+n}_{i=0}C_{I-i,i}^+\prod^{I+n-1}_{j=i}\widehat{f}_{j,n}^+.
\end{align}

	\item[Step 6.] Repeat Steps	2 - 5 above $B$ times, where $B$ is large, to get $(R^+_{I,n}-\widehat{R}_{I,n}^+)^{(b)}$, $b=1,\ldots,B$ bootstrap predictive roots, and denote by $q^+(\alpha)$ the $\alpha$-quantile of their empirical distribution. 
	
	\item[Step 7.] Construct the $(1-\alpha)$ equal-tailed prediction interval for $R_{I,n}$ as
\begin{align*}
\left[\widehat{R}_{I,n}+q^+(\alpha/2), \widehat{R}_{I,n}+q^+(1-\alpha/2)\right].
\end{align*}
\end{itemize}

\begin{rem}[On the alternative Mack-type bootstrap]\ 
\begin{itemize}
	\item[(i)] In comparison to the Mack bootstrap from Section \ref{subsec_boot}, the bootstrap reserve is not a double-bootstrap quantity anymore, the centering is based on a bootstrap version of the best estimate, and the bootstrap for the upper loss triangle is backwards starting in the diagonal.
	\item[(ii)] The conditional distribution for the $G_{i,j}^+|C_{i,j+1}^+$ can be chosen in different ways. For instance, this can be done %in a
non-parametrically %way 
similar to Steps 2 - 4 in Section \ref{subsec_boot} or using the parametric family of distributions used in Step 5 in Section \ref{subsec_boot}. However, it is crucial to mimic sufficiently well the first and second backward conditional moments, that is, $E(C_{ij}|C_{i,j+1})$ and $Var(C_{i,j}|C_{i,j+1})$, respectively. %This will be reflected by Assumption \ref{consistent_estimation_backward_moments} to be imposed in Section \ref{sec_alternative_boot_consistency}.
\end{itemize}
\end{rem}

%\bigskip

\section{Asymptotic Theory for the alternative Mack Bootstrap}\label{sec_alternative_boot_consistency}

By adopting the general strategy of Section \ref{sec_boot_consistency} to investigate the consistency properties of the original Mack bootstrap, the alternative Mack predictive root of the reserve $R^+_{I,n}-\widehat{R}^+_{I,n}$ can be decomposed also into %two additive parts that account for the 
a prediction error part and an estimation error part, respectively. That is, by adding and subtracting $\sum^{I+n}_{i=0}C_{I-i,i}^+\prod^{I+n-1}_{j=i}\widehat{f}_{j,n}$, we get
\begin{align}
R^+_{I,n}-\widehat{R}^+_{I,n}
&=\sum^{I+n}_{i=0}C_{I-i,i}^+\left(\prod^{I+n-1}_{j=i}F_{I-i,j}^+-\prod_{j=i}^{I+n-1}\widehat{f}_{j,n}\right)+\sum^{I+n}_{i=0}C_{I-i,i}^+\left(\prod_{j=i}^{I+n-1}\widehat{f}_{j,n}-\prod_{j=i}^{I+n-1}\widehat{f}_{j,n}^+\right)	\nonumber	\\
&=:\left(R^+_{I,n}-\widehat{R}^+_{I,n}\right)_1+\left(R^+_{I,n}-\widehat{R}^+_{I,n}\right)_2,	\label{eqpredroot_alternative_boot_decomp_appearance}
\end{align} 
where $(R^+_{I,n}-\widehat{R}^+_{I,n})_1$ and $(R^+_{I,n}-\widehat{R}^+_{I,n})_2$ are the alternative Mack bootstrap versions of $(R_{I,n}-\widehat{R}_{I,n})_1$ and $(R_{I,n}-\widehat{R}_{I,n})_2$, respectively.

\subsection{Conditional bootstrap asymptotics for the alternative Mack bootstrap predictive root %of the reserve
}\label{sec_asymp_alternative_boot}\

As in Section \ref{sec_asymp_boot} for the %original 
Mack bootstrap, %conditional on $\mathcal{Q}_{I,n}^+=\mathcal{Q}_{I,n}$ and $\mathcal{D}_{I,n}$, 
we have to check whether, conditional on $\mathcal{Q}_{I,n}^+=\mathcal{Q}_{I,n}$ and $\mathcal{D}_{I,n}$, the alternative Mack bootstrap quantities $(R_{I,n}^+-\widehat{R}_{I,n}^+)_1$ and $(R_{I,n}^+-\widehat{R}_{I,n}^+)_2$ in \eqref{eqpredroot_boot_decomp_appearance} are correctly mimicking the limiting distributions of $(R_{I,n}-\widehat{R}_{I,n})_1$ and $(R_{I,n}-\widehat{R}_{I,n})_2$ given $\mathcal{Q}_{I,n}$, respectively. %, as summarized in Section \ref{sec_asymp_summary}.

\subsubsection{Conditional bootstrap asymptotics for reserve prediction: process uncertainty}\label{sec_asymp_alternative_boot_process}\

The process uncertainty part $(R_{I,n}^+-\widehat{R}_{I,n}^+)_1$ of the alternative Mack bootstrap differs from the $(R_{I,n}^*-\widehat{R}_{I,n})_1$ % in two aspects. On the one hand, 
as the $F_{i,j}^+$'s in $(R_{I,n}^+-\widehat{R}_{I,n}^+)_1$ use $\widehat f_{j,n}$ instead of $\widehat f_{j,n}^*$ and %, on the other hand, 
as $\prod^{I+n-1}_{j=i}F_{I-i,j}^+$ is centered around $\prod_{j=i}^{I+n-1}\widehat{f}_{j,n}$ instead of $\prod_{j=i}^{I+n-1}\widehat{f}_{j,n}^*$ accordingly. However, %as the proof to establish the asymptotic distribution in Theorem \ref{Dist_R1_boot_new} exclusively relies on $\widehat f_{j,n}-f_j=O_P((I+n-1)^{-1/2})$, $\widehat f_{j,n}^*-\widehat f_{j,n}=O_{P^*}((I+n-1)^{-1/2})$ and $\widehat \sigma_{j,n}^2-\sigma_j^2=O_P((I+n-1)^{-1/2})$ for all fixed $j\in\mathbb{N}_0$, 
by using %the same 
very similar arguments, we get %immediately 
the same limiting result also for the process uncertainty part $(R_{I,n}^+-\widehat{R}_{I,n}^+)_1$ of the alternative Mack bootstrap.
 
\begin{thm}[Bootstrap asymptotics for $(R_{I,n}^+-\widehat{R}_{I,n}^+)_1$ conditional on $\mathcal{Q}_{I,n}^+=\mathcal{Q}_{I,n}$ and $\mathcal{D}_{I,n}$]\label{Dist_R1_alternative_boot_new}
Suppose Assumptions \ref{initial_lem}, \ref{annF}, \ref{parameterAss}, \ref{ann_highermoments_F}, \ref{parametric_family} and \ref{parametric_family_bootstrap} (for $F_{i,j}^+$ instead of $F_{i,j}^*$) hold.
Then, as $n\rightarrow \infty$, conditionally on $\mathcal{Q}_{I,n}^+=\mathcal{Q}_{I,n}$ and $\mathcal{D}_{I,n}$, $(R_{I,n}^+-\widehat{R}_{I,n}^+)_1$ converges in distribution to $\mathcal{G}_1|\mathcal{Q}_{I,\infty}$ in probability, which is the (limiting) distribution of $(R_{I,\infty}-\widehat{R}_{I,\infty})_1|\mathcal{Q}_{I,\infty}$ according to \eqref{eq_convergence_in_distribution} described in Theorem \ref{Dist_R1}. Moreover, for all $n\in\mathbb{N}_0$, it holds $E^+((R_{I,n}^+-\widehat{R}_{I,n}^+)_1|\mathcal{Q}_{I,n}^+=\mathcal{Q}_{I,n})=0$ and, for $n\rightarrow \infty$, we have
\begin{align}
Var^+\left((R_{I,n}^+-\widehat{R}_{I,n}^+)_1|\mathcal{Q}_{I,n}^+=\mathcal{Q}_{I,n}\right)\longrightarrow Var\left((R_{I,\infty}-\widehat{R}_{I,\infty})_1|\mathcal{Q}_{I,\infty}\right)	\quad	\text{in prob.,}
\end{align}
%in probability, 
where $Var((R_{I,\infty}-\widehat{R}_{I,\infty})_1|\mathcal{Q}_{I,\infty})=O_P(1)$ as given in \eqref{var_R1_cond}. Consequently, as $n\rightarrow\infty$, we have
\begin{align*}
d_2\left(\mathcal{L}\left((R_{I,n}-\widehat{R}_{I,n})_1|\mathcal{Q}_{I,n}\right), \mathcal{L}^+\left((R^+_{I,n}-\widehat{R}_{I,n}^+)_1|\mathcal{Q}_{I,n}^+=\mathcal{Q}_{I,n}\right)\right)\longrightarrow 0	\quad	\text{in prob..}
\end{align*}
%in probability. 
\end{thm}

\subsubsection{Conditional bootstrap asymptotics for reserve prediction: estimation uncertainty}\label{subsec_asymp_alternative_boot_process}\

In view of the decomposition of $(R_{I,n}-\widehat{R}_{I,n})_2$ %=(R_{I,n}-\widehat{R}_{I,n})_2^{(1)}+(R_{I,n}-\widehat{R}_{I,n})_2^{(2)}$ 
in \eqref{eqpredroot_decomp}, conditional on $\mathcal{Q}_{I,n}^+=\mathcal{Q}_{I,n}$ and $\mathcal{D}_{I,n}$, its alternative Mack bootstrap counterpart $(R_{I,n}^+-\widehat{R}_{I,n}^+)_2$ can be also decomposed further. That is, we have
\begin{align}
& (R_{I,n}^+-\widehat{R}_{I,n}^+)_2	\nonumber	\\
=& \sum^{I+n}_{i=0}C_{I-i,i}^+\left(\prod_{j=i}^{I+n-1}\widehat f_{j,n}-\prod_{j=i}^{I+n-1}f_{j,n}^+(\mathcal{Q}_{I,n})\right) + \sum^{I+n}_{i=0}C_{I-i,i}^+\left(\prod_{j=i}^{I+n-1}f_{j,n}^+(\mathcal{Q}_{I,n})-\prod_{j=i}^{I+n-1}\widehat{f}_{j,n}^+\right)	\nonumber	\\
=& (R_{I,n}^+-\widehat{R}_{I,n}^+)_2^{(1)}+(R_{I,n}^+-\widehat{R}_{I,n}^+)_2^{(2)},	\label{eqpredroot_alternative_boot_decomp}
\end{align}
where $(R_{I,n}^+-\widehat{R}_{I,n}^+)_2^{(1)}$ is measurable with respect to $\mathcal{Q}_{I,n}^+=\mathcal{Q}_{I,n}$ and $\mathcal{D}_{I,n}$ and $f_{j,n}^+(\mathcal{Q}_{I,n}):=\mu_{j+1,n}^{+(1)}(\mathcal{Q}_{I,n})/\mu_{j,n}^{+(2)}(\mathcal{Q}_{I,n})$ with $\mu_{j+1,n}^{+(1)}(\mathcal{Q}_{I,n}) := E^+(\frac{1}{I+n-j}\sum_{i=-n}^{I-j-1} C_{i,j+1}^+|\mathcal{Q}_{I,n}^+=\mathcal{Q}_{I,n})$ as well as $\mu_{j,n}^{+(2)}(\mathcal{Q}_{I,n}) := E^+(\frac{1}{I+n-j}\sum_{i=-n}^{I-j-1} C_{i,j}^+|\mathcal{Q}_{I,n}^+=\mathcal{Q}_{I,n})$.
%\begin{align*}
%\mu_{j+1,n}^{+(1)}(\mathcal{Q}_{I,n}) &:= E^+\left(\frac{1}{I+n-j}\sum_{i=-n}^{I-j-1} C_{i,j+1}^+|\mathcal{Q}_{I,n}^+=\mathcal{Q}_{I,n}\right) = \frac{1}{I+n-j}\sum_{i=-n}^{I-j-1} E^+(C_{i,j+1}^+|C_{i,I-i}^+=C_{i,I-i}),	\\
%\mu_{j,n}^{+(2)}(\mathcal{Q}_{I,n}) &:= E^+\left(\frac{1}{I+n-j}\sum_{i=-n}^{I-j-1} C_{i,j}^+|\mathcal{Q}_{I,n}^+=\mathcal{Q}_{I,n}\right) =\frac{1}{I+n-j}\sum_{i=-n}^{I-j-1} E^+(C_{i,j}^+|C_{i,I-i}^+=C_{i,I-i}).
%\end{align*}

In comparison to Theorem \ref{Dist_R2_boot_new}, the uniform boundedness condition in Assumption \ref{uniform_boundedness_condition} can be dropped, but  
%Besides a correctly chosen and sufficiently smooth parametric family of (conditional) bootstrap distributions of the individual development factors for the \emph{lower triangle}, 
the derivation of (conditional) bootstrap asymptotic theory and consistency results for $(R_{I,n}^+-\widehat{R}_{I,n}^+)_2$ requires additional assumptions on the \emph{backward} individual development factors $G_{i,j}^+$ from Step 3 in Section \ref{subsec_alternative_boot}. Precisely, it has to be guaranteed that the backward conditional mean $E(C_{i,j}|C_{i,j+1})$ and the backward conditional variance $Var(C_{i,j}|C_{i,j+1})$ are consistently mimicked by their alternative Mack bootstrap counterparts $E^+(C_{i,j}^+|C_{i,j+1}^+)$ and $Var^+(C_{i,j}^+|C_{i,j+1}^+)$, respectively, such that the corresponding limiting distributions obtained in \citet[Theorem C.1]{jentschasympThInf} %(see also Theorem \ref{CLT_f_cond} in the appendix) 
are correctly mimicked.

\begin{ann}[Consistent estimation of backward moments]\label{consistent_estimation_backward_moments}
For $n\rightarrow \infty$, suppose that the (conditional) bootstrap distributions of the backward individual development factors $G_{i,j}^+$, $j=0, \dots, I+n-1$ and $i=-n, \dots, I-j-1$ given $C_{i,j+1}^+$ and $\mathcal{D}_{I,n}$ are chosen in Step 2 in Section \ref{subsec_alternative_boot} such that:
 %the following holds:
\begin{itemize}
	\item[(i)] For each fixed $K\in\mathbb{N}_0$, let $\underline{f}_{K,n}^+(\mathcal{Q}_{I,n})=(f_{0,n}^+(\mathcal{Q}_{I,n}),f_{1,n}^+(\mathcal{Q}_{I,n}),\ldots,f_{K,n}^+(\mathcal{Q}_{I,n}))^\prime$ and	define $\underline{\widehat f}_{K,n}=(\widehat f_{0,n},\widehat f_{1,n},\ldots,\widehat f_{K,n})^\prime$. Then, conditional on $\mathcal{Q}_{I,n}^+=\mathcal{Q}_{I,n}$, we have
\begin{align*}
J_n^{1/2}\left(\underline{f}_{K,n}^+(\mathcal{Q}_{I,n})-\underline{\widehat f}_{K,n}\right)|\left(\mathcal{Q}_{I,n}^+=\mathcal{Q}_{I,n}\right) \overset{d}{\longrightarrow} \mathcal{N}\left(0, \bm{\Sigma}_{K,\underline{f}}^{(1)}\right),
\end{align*}
where $J_n^{1/2}=diag\left(\sqrt{I+n-j\vphantom{I^2}},j=0,\ldots,K\right)$ is a diagonal $(K+1)\times (K+1)$ matrix of inflation factors and the variance-covariance matrix $\bm{\Sigma}_{K,\underline{f}}^{(1)}$ is defined in Theorem C.1 in \cite{jentschasympThInf}.
%\ref{CLT_f_cond}.
	\item[(ii)] For each fixed $K\in\mathbb{N}_0$, let $\underline{\widehat f}_{K,n}^+=(\widehat f_{0,n}^+,\widehat f_{1,n}^+,\ldots,\widehat f_{K,n}^+)^\prime$. Then, conditional on $\mathcal{Q}_{I,n}^+=\mathcal{Q}_{I,n}$ and $\mathcal{D}_{I,n}$, we have
\begin{align*}
J_n^{1/2}\left(\underline{\widehat f}_{K,n}^+-\underline{f}_{K,n}^+(\mathcal{Q}_{I,n})\right)|\left(\mathcal{Q}_{I,n}^+=\mathcal{Q}_{I,n},\mathcal{D}_{I,n}\right) \overset{d}{\longrightarrow} \mathcal{N}\left(0,\bm{\Sigma}_{K,\underline{f}}^{(2)}\right)	\quad	\text{in prob.}
\end{align*}
%in probability, 
where the variance-covariance matrix $\bm{\Sigma}_{K,\underline{f}}^{(2)}$ is defined in Theorem C.1 in \cite{jentschasympThInf}.
%\ref{CLT_f_cond}.
\end{itemize}
\end{ann}

In concordance to the derivation of the conditional limiting result obtained in Theorem \ref{Dist_R2}(ii), which relies on conditional CLTs for the development factor estimators $\widehat f_{j,n}$ given in \cite[Appendix C]{jentschasympThInf}, the conditional bootstrap CLTs in Assumption \ref{consistent_estimation_backward_moments} allow to state the following theorem, which provides the limiting distribution of the alternative Mack bootstrap estimation uncertainty term $(R_{I,n}^+-\widehat{R}_{I,n}^+)_2$ conditional on $\mathcal{Q}_{I,n}^+=\mathcal{Q}_{I,n}$ and $\mathcal{D}_{I,n}$. Precisely, while $(R_{I,n}^+-\widehat{R}_{I,n}^+)_2^{(1)}$ is measurable with respect to $\mathcal{D}_{I,n}$, Assumption \ref{consistent_estimation_backward_moments} allows to establish asymptotic normality of $\sqrt{I+n+1}(R_{I,n}^+-\widehat{R}_{I,n}^+)_2^{(2)}$ conditional on $\mathcal{Q}_{I,n}^+=\mathcal{Q}_{I,n}$ and $\mathcal{D}_{I,n}$.

\begin{thm}[Bootstrap asymptotics for $(R_{I,n}^+-\widehat{R}_{I,n}^+)_2$ conditional on $\mathcal{Q}_{I,n}^+=\mathcal{Q}_{I,n}$ and $\mathcal{D}_{I,n}$]\label{Dist_R2_alternative_boot_new}
Suppose Assumptions \ref{initial_lem}, \ref{annF}, \ref{ann_highermoments_F}, \ref{parameterAss}, \ref{support_condition} and \ref{consistent_estimation_backward_moments}
hold. Then, as $n\rightarrow \infty$, the following holds:
\begin{itemize}
	\item[(i)] Conditional on $\mathcal{Q}_{I,n}^+=\mathcal{Q}_{I,n}$, $\sqrt{I+n+1}(R_{I,n}^+-\widehat{R}_{I,n}^+)_2^{(1)}$ converges in distribution to the non-degenerate limiting distribution $\mathcal{G}_2^{(1)}$. That is, we have
	\begin{align}
	\sqrt{I+n+1}(R_{I,n}^+-\widehat{R}_{I,n}^+)_2^{(1)}|\left(\mathcal{Q}_{I,n}^+=\mathcal{Q}_{I,n}\right) \overset{d}{\longrightarrow} \left\langle \mathcal{Q}_{I,\infty},\mathbf{Y}_\infty^{(1)}\right\rangle\sim\mathcal{G}_2^{(1)},
	\end{align}
	where $\mathbf{Y}_\infty^{(1)}=(Y_i^{(1)},i\in\mathbb{N}_0)$ denotes a centered Gaussian process with covariances
	\begin{align*}%\label{var_R_2_uncon}
	Cov(Y_{i_1}^{(1)},Y_{i_2}^{(1)}) = \lim_{K\rightarrow \infty}\bm{\Sigma}_{K,\prod \f}^{(1)}(i_1,i_2),	\quad	i_1,i_2\in\mathbb{N}_0,
	\end{align*}
%for $i_1,i_2\in\mathbb{N}_0$, 
where $\bm{\Sigma}_{K,\prod \f}^{(1)}(i_1,i_2)$ is defined in Corollary C.2 in \cite{jentschasympThInf}.
%\ref{CLTProdf_cond}. 
Here, the %two random 
sequences 
$\mathcal{Q}_{I,\infty}$ and $\mathbf{Y}_\infty^{(1)}$ are %stochastically 
independent. 
	\item[(ii)] Conditionally on $\mathcal{Q}_{I,n}^+=\mathcal{Q}_{I,n}$ and $\mathcal{D}_{I,n}$, $\sqrt{I+n+1}(R_{I,n}^+-\widehat{R}_{I,n}^+)_2^{(2)}$ converges in distribution to	$\mathcal{G}_2|\mathcal{Q}_{I,\infty}$ in probability, where $\mathcal{G}_2|\mathcal{Q}_{I,\infty}\sim \mathcal{N}(0,\Xi(\mathcal{Q}_{I,\infty}))|\mathcal{Q}_{I,\infty}$ is the (conditional) limiting distribution obtained in Theorem \ref{Dist_R2}(ii).
\end{itemize}
Consequently, as $n\rightarrow\infty$, we have
\begin{align*}
d_K\left(\mathcal{L}\left((R_{I,n}-\widehat{R}_{I,n})_2|\mathcal{Q}_{I,n}\right), \mathcal{L}^+\left((R^+_{I,n}-\widehat{R}_{I,n}^+)_2|\mathcal{Q}_{I,n}^+=\mathcal{Q}_{I,n}\right)\right)\longrightarrow 0	\quad	\text{in prob..}
\end{align*}
%in probability.
\end{thm}

\subsubsection{Conditional bootstrap asymptotics for the whole predictive root of the reserve}\label{sec_asymp_boot_joint_alternative}\

As in Sections \ref{sec_asymp_summary_joint} and \ref{sec_asymp_boot_joint}, combining the results for $(R_{I,n}^+-\widehat R_{I,n}^*)_1$ and $(R_{I,n}^+-\widehat R_{I,n}^*)_2$ from Theorems \ref{Dist_R1_alternative_boot_new} and \ref{Dist_R2_alternative_boot_new}, we get joint asymptotics for $R_{I,n}^+-\widehat{R}_{I,n}$ conditional on $\mathcal{Q}_{I,n}^+=\mathcal{Q}_{I,n}$ and $\mathcal{D}_{I,n}$. % can be obtained.

%As in Section \ref{sec_asymp_summary_joint}, we can combine the results on the limiting distributions for $(R_{I,n}^+-\widehat R_{I,n}^+)_1$ and $(R_{I,n}^+-\widehat R_{I,n}^+)_2$ from Theorems \ref{Dist_R1_alternative_boot_new} and \ref{Dist_R2_alternative_boot_new}, respectively, to get the limiting bootstrap distribution for the whole bootstrap predictive root of the reserve $R_{I,n}^+-\widehat{R}_{I,n}^+$ conditional on $\mathcal{Q}_{I,n}^+=\mathcal{Q}_{I,n}$ and $\mathcal{D}_{I,n}$.

\begin{thm}[Bootstrap asymptotics 
for $R_{I,n}^+-\widehat R_{I,n}^+$ conditional on $\mathcal{Q}_{I,n}^+=\mathcal{Q}_{I,n}$ and $\mathcal{D}_{I,n}$]\label{Dist_R1_R2_joint_alternative_boot}
Suppose the assumptions of Theorems \ref{Dist_R1_alternative_boot_new} and \ref{Dist_R2_alternative_boot_new} hold. Then, conditional on $\mathcal{Q}_{I,n}^+=\mathcal{Q}_{I,n}$ and $\mathcal{D}_{I,n}$, $(R_{I,n}^+-\widehat R_{I,n}^+)_1$ and $(R_{I,n}^+-\widehat R_{I,n}^+)_2$ are stochastically independent, and $R_{I,n}^+-\widehat R_{I,n}^+$ converges in distribution to $\mathcal{G}_1|\mathcal{Q}_{I,\infty}$ in probability. That is, we have
\begin{align*}	
R_{I,n}^+-\widehat{R}_{I,n}^+|\left(\mathcal{Q}_{I,n}^+=\mathcal{Q}_{I,n},\mathcal{D}_{I,n}\right)=(R_{I,n}^+-\widehat R_{I,n}^+)_1+(R_{I,n}^+-\widehat R_{I,n}^+)_2|\left(\mathcal{Q}_{I,n}^+=\mathcal{Q}_{I,n},\mathcal{D}_{I,n}\right)\overset{d}{\longrightarrow}\mathcal{G}_1|\mathcal{Q}_{I,\infty}
\end{align*}
in probability.
%\end{itemize}
\end{thm}

According to the discussion below Theorem \ref{Dist_R1_R2_joint_boot} and in view of the concepts of \emph{asymptotic validity} and \emph{asymptotic pertinence} of %a 
bootstrap %prediction approach as discussed 
predictive inference in \cite{pan2016bootstrap}, the alternative Mack bootstrap can be regarded as asymptotically valid \emph{and} asymptotically pertinent under the stated conditions.

\begin{rem}[Backward vs.~forward bootstrapping]
While a backward bootstrap %approach 
appears to be natural in time series setups addressed in \cite{pan2016bootstrap}, they also propagate a %somewhat
simpler forward bootstrap %approach 
to capture the estimation uncertainty in bootstrap prediction. Asymptotically, in their setup, both approaches are indeed equivalent due to the intrinsic stationarity assumption. However, in Mack's Model setup considered here, this is not the case and %a 
the (fixed-design) forward bootstrap %as proposed by 
of \cite{england2006predictive} does not correctly capture the conditional limiting distribution of the estimation uncertainty part. %, while a backward bootstrap is capable of doing this.
\end{rem}

%\bigskip

\section{Simulation Study}\label{sec_simulation}

In this section, we compare the original Mack bootstrap from Section \ref{sec_boot} and the alternative Mack bootstrap from Section \ref{sec_alternative_boot} to illustrate our theoretical findings from Sections \ref{sec_boot_consistency} and \ref{sec_alternative_boot_consistency} by means of simulations of several parameter scenarios. Additionally, we simulate a Mack-type bootstrap, which uses a \emph{forward} bootstrap approach in Step 2 of Section \ref{subsec_alternative_boot}, but coincides otherwise with the alternative Mack bootstrap. %The inclusion of this intermediate Mack-type bootstrap allows to disentangle the effects caused by the backward resampling proposed in Step 2 and by the different centering used in Step 5 of Section \ref{subsec_alternative_boot} on the finite sample performance. Both aspects constitute the deviance of the alternative Mack bootstrap from the original Mack bootstrap.

\subsection{Simulation setup}\

%To assure comparability, 
We pick up the simulation setup employed in \cite{jentschasympThInf}. That is, in the notion of the asymptotic framework introduced in Section \ref{sec_asymptotic_framework}, let $I=10$ and choose $n\in\{0,10,20,30,40\}$ leading to effective number of accident years $I+n+1\in\{11,21,31,41,51\}$. For each $n$ and for different parameter scenarios %to be 
specified below, we generate $M=500$ loss triangles $\mathcal{D}_{I,n}^{(m)}=\{C_{i,j}^{(m)}|i=-n,\ldots,I,~j=0,\ldots,I-i\}$, $m=1,\ldots,500$, %having 
with diagonals $\mathcal{Q}_{I,n}^{(m)}$ by generating the entries in their first columns $C_{\bullet,0}^{(m)}$ (independently) from a uniform distribution and the individual developments factors $F_{i,j}$ given $C_{i,j}$ from a
\begin{itemize}
	\item[(DGP1)] conditional gamma distribution,
	\item[(DGP2)] conditional log-normal distribution,
	\item[(DGP3)] conditional left-tail truncated normal distribution (truncated at 0.1).
\end{itemize}
In all scenarios, the development factors and the variance parameters 
%are specified to 
fulfill $\f>1$ and $\sigma^2_j>0$ %for all $j=0, \dots, I+n-1$ such that 
with $\f$ and $\s$ decreasing to $1$ and $0$, respectively. Precisely, we use exponentially decreasing sequences $(f_j)_{j\in\mathbb{N}_0}$ and $(\sigma_j^2)_{j\in\mathbb{N}_0}$ with $f_j=1+e^{-1-0.2j}$ and $\sigma^2_j=509,518\cdot e^{-1-0.7j}$. %for $j=0, \dots, I+n-1$ for  $I=10$ and $n\in\{0,10,20,30,40\}$.
%Further, 
We distinguish between two %different 
Setups a) and b), where the parameter %settings 
are exactly the same in both cases, but the first column $C_{\bullet,0}=(C_{-n,0},\ldots,C_{I,0})'$ of the (upper) loss triangle is uniformly distributed on $[120\times 10^6,350\times 10^6]$ in Setup a) and on $[120\times 10^4,350\times 10^4]$ in b). The results for both setups are similar and, we show only the those for Setup a) here and report the results for Setup b) in the appendix. % supplementary material \cite{Mack_boot_supplement}.

In the following, to evaluate the performance of all bootstrap procedures under study, for each diagonal $\mathcal{Q}_{I,n}^{(m)}$, $m=1,\ldots,500$, we would like to know the exact distribution $R_{I,n}^{(m)}-\widehat{R}_{I,n}^{(m)}$ conditional on $\mathcal{Q}_{I,n}^{(m)}$. However, although knowing exactly the stochastic mechanism to generate a loss triangle $\mathcal{D}_{I,n}$, it is not straightforward %at all 
to simulate $R_{I,n}^{(m)}-\widehat{R}_{I,n}^{(m)}|\mathcal{Q}_{I,n}^{(m)}$. This is because $\widehat{R}_{I,n}^{(m)}$ requires a \emph{backward} generation of a loss triangle $\mathcal{D}_{I,n}$ starting with $\mathcal{Q}_{I,n}^{(m)}$. Hence, as a workaround, 
%\footnote{CJ: Vielleicht dieses Problem dann so zu offensiv formuliert? Julia?}
we simulate %instead 
the distribution of the ``true'' predictive root $R_{I,n}^{(m)}-\widehat{R}_{I,n}^{(m)}$ conditional on $\mathcal{Q}_{I,n}$ by a Monte Carlo simulation with $B=10,000$, since we know the true underlying parametric family of distributions of the individual development factors for each observed triangle $\mathcal{D}_{I,n}^{(m)}$ and the true parameters for the simulation of $R_{I,n}^{(m)}$ for each setup (DGP1)-(DGP3) such that 
\begin{align}
    \F|\C \sim \left(f_j, \frac{\sigma^2_j}{\C}\right) \quad \text{ for } j=I-i, \dots, I+n-1 \text{ and } i=-n, \dots I.
\end{align}
Next, for each setup (DGP1)-(DGP3) above and for each loss triangle $\mathcal{D}_{I,n}^{(m)}$, $m=1,\ldots,500$, we perform three different Mack-type bootstraps based on 10,000 bootstrap replications each to estimate the conditional distributions of the predictive roots of the reserve. 
That is, we apply the following three bootstrap approaches:
\begin{itemize}
	\item[(oMB)] original Mack bootstrap (from Section \ref{sec_boot}), 
	\item[(aMB)] alternative Mack-type bootstrap (from Section \ref{sec_alternative_boot}), 
	\item[(iMB)] intermediate Mack-type bootstrap (using a \emph{forward} bootstrap in Step 2 of Section \ref{sec_alternative_boot}).
\end{itemize}
The third intermediate Mack-type bootstrap is included %to be able 
to distinguish between the effects caused by the backward resampling proposed in Step 2 and by the different centering used in Step 5 of Section \ref{subsec_alternative_boot} on the finite sample performance. For this purpose, we introduce a novel centering term $\widehat R_{I,n}^{++}$ defined by
\begin{align}\label{centering_term_alternative++}
\widehat R_{I,n}^{++}=\sum^{I+n}_{i=0}C_{I-i,i}^+\prod^{I+n-1}_{j=i}\widehat f_{j,n}^*,
\end{align}
which deviates from $\widehat R_{I,n}^{+}$ in \eqref{centering_term_alternative} as it relies on $\widehat f_{j,n}^*$ %defined 
in \eqref{eqf_boot}, but is based on (parametrically generated)
\begin{align}
F_{i,j}^*|C_{i,j}^*,\mathcal{D}_{I,n}\sim \left(\widehat f_{j,n}, \frac{\widehat\sigma^2_{j,n}}{C_{i,j}^*}\right),
\end{align}
instead of $\widehat f_{j,n}^+$ defined in \eqref{eq_def_f+}. This choice of the centering term still resembles the decomposition in \eqref{eqpredroot_alternative_boot_decomp_appearance}, that shares the (sign) properties of \eqref{eqpredroot_decomp_first_appearance}, which is not the case for \eqref{eqpredroot_boot_decomp_appearance}. For all bootstraps, whenever a \emph{parametric} distribution is used to generate the \emph{upper} bootstrap loss triangle, we choose the same parametric distribution family used already for the lower triangle (to generate $R_{I,n}^*$ and $R_{I,n}^{+}$). However, as we do not know the correct parametric family of distributions of the $\F$'s, we make use of all three distribution families in (i)-(iii) for all three bootstrap approaches, respectively, %Finally, knowing the true parametric family of distributions of the $\F$'s, which allows to simulate the "true" predictive root of the reserve, we compare the simulation results 
to also investigate the effect of a misspecified parametric family of distributions to generate $R_{I,n}^*$ and $R_{I,n}^+$.

In the %supplementary material \cite{Mack_boot_supplement}, 
Appendix \ref{sec_add_simulation},
we provide also %corresponding 
simulation results that compare %just 
the distribution of the first (i.e.~the process uncertainty) parts of the bootstrap predictive roots $(R^*_{I,n}-\widehat{R}_{I,n})_1$ and $(R_{I,n}^+-\widehat{R}_{I,n}^+)_1$ conditional on $\mathcal{Q}_{I,n}^*=\mathcal{Q}_{I,n}$ or $\mathcal{Q}_{I,n}^+=\mathcal{Q}_{I,n}$ and $\mathcal{D}_{I,n}$, respectively, with the distribution of $(R_{I,n}-\widehat{R}_{I,n})_1$ conditional on $\mathcal{Q}_{I,n}$. Note that this distribution is straightforward to simulate. As expected, in view of Theorems \ref{Dist_R1_boot_new} and \ref{Dist_R1_alternative_boot_new}, we find %essentially 
no differences %between 
in the performances of both construction principles.

\subsection{Simulation results}\label{subsec_sim_together}\

First, we consider the bootstrap variances of the bootstrap predictive roots of the reserves obtained for the three Mack-type bootstraps under study. For both Setups a) and b), we find that the alternative Mack-type bootstrap variance is always 1-5 percentage points smaller than the bootstrap variances obtained for the other two approaches, which do not differ much. %(less than 1 percentage point). 
This result perfectly agrees to the findings of Theorem \ref{Dist_R2_boot_new}, where the (conditional) variance $\widetilde \Xi(\mathcal{Q}_{I,\infty})$, which is mimicked by the original Mack bootstrap and by the intermediate Mack-type bootstrap, is generally larger than the variance $\Xi(\mathcal{Q}_{I,\infty})$ found in Theorem \ref{Dist_R2}, which is mimicked by the alternative Mack bootstrap correctly according to Theorem \ref{Dist_R2_alternative_boot_new}.

Next, we consider the whole distributions of the bootstrap predictive roots $R_{I,n}^{*(m)}-\widehat{R}_{I,n}^{(m)}$, $R_{I,n}^{+(m)}-\widehat{R}_{I,n}^{+(m)}$ and $R_{I,n}^{+(m)}-\widehat{R}_{I,n}^{++(m)}$ conditional on $\mathcal{Q}_{I,n}^*=\mathcal{Q}_{I,n}$ or $\mathcal{Q}_{I,n}^+=\mathcal{Q}_{I,n}$ and $\mathcal{D}_{I,n}$, respectively, for $m=1,\dots, 500$. Using the Kolmogorov-Smirnov test of level $\alpha=5\%$ to test the null hypotheses
\begin{align*}
& H_0^*: \mathcal{L}\left(R_{I,n}^{*(m)}-\widehat R_{I,n}^{(m)}|(\mathcal{Q}_{I,n}^{*(m)}=\mathcal{Q}_{I,n}^{(m)},\mathcal{D}_{I,n}^{(m)})\right)=\mathcal{L}\left(R_{I,n}^{(m)}-\widehat R_{I,n}^{(m)}|\mathcal{Q}_{I,n}^{(m)}\right),	\\
& H_0^+: \mathcal{L}\left(R_{I,n}^{+(m)}-\widehat R_{I,n}^{+(m)}|(\mathcal{Q}_{I,n}^{+(m)}=\mathcal{Q}_{I,n}^{(m)},\mathcal{D}_{I,n}^{(m)})\right)=\mathcal{L}\left(R_{I,n}^{(m)}-\widehat R_{I,n}^{(m)}|\mathcal{Q}_{I,n}^{(m)}\right),	\\
& H_0^{++}: \mathcal{L}\left(R_{I,n}^{+(m)}-\widehat R_{I,n}^{++(m)}|(\mathcal{Q}_{I,n}^{+(m)}=\mathcal{Q}_{I,n}^{(m)},\mathcal{D}_{I,n}^{(m)})\right)=\mathcal{L}\left(R_{I,n}^{(m)}-\widehat R_{I,n}^{(m)}|\mathcal{Q}_{I,n}^{(m)}\right)
\end{align*}
%that the (conditional) distributions of $R_{I,n}^{*(m)}-\widehat R_{I,n}^{(m)}|(\mathcal{Q}_{I,n}^{*(m)}=\mathcal{Q}_{I,n}^{(m)},\mathcal{D}_{I,n}^{(m)})$, $R_{I,n}^{+(m)}-\widehat R_{I,n}^{+(m)}|(\mathcal{Q}_{I,n}^{+(m)}=\mathcal{Q}_{I,n}^{(m)},\mathcal{D}_{I,n}^{(m)})$ and $R_{I,n}^{+(m)}-\widehat R_{I,n}^{++(m)}|(\mathcal{Q}_{I,n}^{+(m)}=\mathcal{Q}_{I,n}^{(m)},\mathcal{D}_{I,n}^{(m)})$, respectively, are equal to the distribution $R_{I,n}^{(m)}-\widehat R_{I,n}^{(m)}|\mathcal{Q}_{I,n}^{(m)}$ conditional on $\mathcal{Q}_{I,n}$ 
for $m=1,\dots, 500$. The resulting percentages of \emph{failed rejections of the null hypotheses} for all three bootstrap approaches, for different $n$ and different families of distributions are summarized in %Tables \ref{tab_ks_predroot_a} and \ref{tab_ks_predroot_b} for setups a) and b), respectively.
Table \ref{tab_ks_predroot_a}. While the percentages increase for growing $n$, for all bootstraps, % and in both setups a) and b), 
the alternative Mack-type bootstrap consistently achieves percentages that are %always 
higher by 1-3 percentage points in comparison to %the percentages of 
to the two other bootstraps, which turn out to be quite similar throughout. %Moreover, the percentages obtained for setup a) are higher compared to setup b), and choosing the true distributional family for $\F^*$ appears to be less important than for setup b). In particular, when choosing a log-normal distribution instead of a truncated normal distribution or vice versa leads to the lowest percentages of failed rejections for setup b).

\begin{small}
\begin{table}[t]
\begin{tabular}{|c|c|lll|lll|lll|}
\hline
\multicolumn{2}{|c|}{\textbf{chosen distribution}}                      & \multicolumn{3}{c}{\cellcolor[HTML]{D9D9D9}\textbf{gamma}} & \multicolumn{3}{c}{\cellcolor[HTML]{D9D9D9}\textbf{log-normal}} & \multicolumn{3}{c|}{\cellcolor[HTML]{D9D9D9}\textbf{trunc. normal}} \\
\hline
\multicolumn{1}{|l|}{\textbf{true distribution}}                  & n  & oMB                & aMB                & iMB                & oMB                 & aMB                  & iMB                  & oMB                  & aMB                   & iMB                   \\\hline
\cellcolor[HTML]{D9D9D9}                                         & 0  & 0.21               & 0.22              & 0.21              & 0.30                & 0.33                & 0.29                & 0.29                 & 0.28                 & 0.21                 \\
\cellcolor[HTML]{D9D9D9}                                         & 10 & 0.38               & 0.49              & 0.38              & 0.47                & 0.48                & 0.43                & 0.37                 & 0.41                 & 0.38                 \\
\cellcolor[HTML]{D9D9D9}                                         & 20 & 0.47               & 0.56              & 0.47              & 0.51                & 0.56                & 0.51                & 0.54                 & 0.53                 & 0.49                 \\
\cellcolor[HTML]{D9D9D9}                                         & 30 & 0.58               & 0.64              & 0.58              & 0.56                & 0.61                & 0.55                & 0.60                 & 0.65                 & 0.59                 \\
\multirow{-5}{*}{\cellcolor[HTML]{D9D9D9}\textbf{gamma}}         & 40 & 0.66               & 0.70              & 0.66              & 0.61                & 0.66                & 0.60                & 0.72                 & 0.76                 & 0.70                 \\\hline
\cellcolor[HTML]{D9D9D9}                                         & 0  & 0.20               & 0.22              & 0.20              & 0.27                & 0.25                & 0.25                & 0.24                 & 0.24                 & 0.22                 \\
\cellcolor[HTML]{D9D9D9}                                         & 10 & 0.37               & 0.38              & 0.37              & 0.38                & 0.40                & 0.37                & 0.37                 & 0.41                 & 0.36                 \\
\cellcolor[HTML]{D9D9D9}                                         & 20 & 0.45               & 0.49              & 0.45              & 0.48                & 0.55                & 0.52                & 0.45                 & 0.51                 & 0.49                 \\
\cellcolor[HTML]{D9D9D9}                                         & 30 & 0.51               & 0.55              & 0.51              & 0.57                & 0.60                & 0.55                & 0.51                 & 0.56                 & 0.54                 \\
\multirow{-5}{*}{\cellcolor[HTML]{D9D9D9}\textbf{log-normal}}    & 40 & 0.57               & 0.62              & 0.57              & 0.60                & 0.63                & 0.60                & 0.63                 & 0.66                 & 0.64                 \\ \hline
\cellcolor[HTML]{D9D9D9}                                         & 0  & 0.13               & 0.16              & 0.13              & 0.30                & 0.33                & 0.29                & 0.24                 & 0.27                 & 0.20                 \\
\cellcolor[HTML]{D9D9D9}                                         & 10 & 0.42               & 0.45              & 0.42              & 0.44                & 0.45                & 0.43                & 0.38                 & 0.47                 & 0.42                 \\
\cellcolor[HTML]{D9D9D9}                                         & 20 & 0.54               & 0.58              & 0.54              & 0.53                & 0.55                & 0.52                & 0.58                 & 0.62                 & 0.59                 \\
\cellcolor[HTML]{D9D9D9}                                         & 30 & 0.57               & 0.63              & 0.57              & 0.57                & 0.60                & 0.57                & 0.67                 & 0.70                 & 0.68                 \\
\multirow{-5}{*}{\cellcolor[HTML]{D9D9D9}\textbf{trunc. normal}} & 40 & 0.60               & 0.66              & 0.60              & 0.57                & 0.62                & 0.59                & 0.75                 & 0.78                 & 0.76 \\    \hline           
\end{tabular}
\caption{Percentages of failed rejections for Kolmogorov-Smirnov tests of level $\alpha=5\%$ for the null hypotheses $H_0^*$, $H_0^+$ and $H_0^{++}$,
 %$\mathcal{L^*}((R_{I,n}^*-\widehat R_{I,n})|\mathcal{Q}_{I,n}^*=\mathcal{Q}_{I,n})=\mathcal{L}((R_{I,n}-\widehat R_{I,n})_1|\mathcal{Q}_{I,n})$ and $\mathcal{L^*}((R_{I,n}^+-\widehat R_{I,n}^+)|\mathcal{Q}_{I,n}^+=\mathcal{Q}_{I,n})=\mathcal{L}((R_{I,n}-\widehat R_{I,n})_1|\mathcal{Q}_{I,n})$, $\mathcal{L^*}((R_{I,n}^+-\widehat R_{I,n}^++)|\mathcal{Q}_{I,n}^+=\mathcal{Q}_{I,n})=\mathcal{L}((R_{I,n}-\widehat R_{I,n})_1|\mathcal{Q}_{I,n})$, 
respectively, for the original Mack bootstrap (oMB), the alternative Mack bootstrap (aMB) and the intermediate Mack bootstrap (iMB) for different parametric families of distributions of $\F^*$ for $i+j\geq I$, for $I=10$ and different $n$ in Setup a).}\label{tab_ks_predroot_a}
\end{table}
\end{small}

We consider also the average over %the $M=500$ 
all simulations of the squared mean of the deviation of the bootstrap distribution given $\mathcal{Q}_{I,n}^*=\mathcal{Q}_{I,n}$ or $\mathcal{Q}_{I,n}^+=\mathcal{Q}_{I,n}$, respectively, and $\mathcal{D}_{I,n}$ and its true distribution given $\mathcal{Q}_{I,n}$. Therefore, we calculate the mean squared error of each simulation for $b=1, \dots, 10,000$ and then consider the root of the overall mean of the mean squared error (RMMSE) over all $M=500$ simulations, that is,  
\begin{align}
RMMSE_{oMB}=\sqrt{\frac{1}{500}\sum^{500}_{m=1}\frac{1}{10,000}\sum^{10,000}_{b=1}\left(R_{I,n}^{*(b)(m)}-\widehat R_{I,n}^{(m)}-\left(R_{I,n}^{(b)(m)}-\widehat R_{I,n}^{(m)}\right)\right)^2,}
\end{align}
where $R_{I,n}^{*(b)(m)}-\widehat R_{I,n}^{(m)}$ represents the $b$th ordered Mack-type bootstrap predictive root and $R_{I,n}^{(b)(m)}-\widehat R_{I,n}^{(m)}$ is the $b$th ordered true simulated predictive root for the $m$th simulation for $m=1, \dots, 500$. Similarly, we calculate $RMMSE_{aMB}$ and $RMMSE_{iMB}$ for the alternative Mack bootstrap and for the intermediate Mack bootstrap, respectively.

The results obtained for all $RMMSE$s are summarized for all three bootstrap approaches in %Tables \ref{tab_mse_predroot_a} and \ref{tab_mse_predroot_b}, respectively, for setups a) and b). 
Table \ref{tab_mse_predroot_a}. For increasing $n$, the $RMMSE$s are decreasing for all bootstrap approaches in both setups, while the alternative Mack bootstrap has the smallest $RMMSE$ in most cases in comparison to the intermediate and the original Mack bootstraps. 

\begin{small}
\begin{table}[t]
\begin{tabular}{|c|c|lll|lll|lll|}
\hline
\multicolumn{2}{|c|}{\textbf{chosen distribution}}                      & \multicolumn{3}{c}{\cellcolor[HTML]{D9D9D9}\textbf{gamma}}                & \multicolumn{3}{c}{\cellcolor[HTML]{D9D9D9}\textbf{log-normal}}           & \multicolumn{3}{c|}{\cellcolor[HTML]{D9D9D9}\textbf{trunc. normal}}        \\
\hline
\multicolumn{1}{|l|}{\textbf{true distribution}}                   & n  & \multicolumn{1}{c}{oMB} & \multicolumn{1}{c}{aMB} & \multicolumn{1}{c|}{iMB} & \multicolumn{1}{c}{oMB} & \multicolumn{1}{c}{aMB} & \multicolumn{1}{c|}{iMB} & \multicolumn{1}{c}{oMB} & \multicolumn{1}{c}{aMB} & \multicolumn{1}{c|}{iMB} \\
\hline
\cellcolor[HTML]{D9D9D9}                                         & 0  & 99.720                  & 99.706                 & 99.600                 & 99.244                  & 99.202                 & 99.967                 & 98.210                  & 98.156                 & 98.822                 \\
\cellcolor[HTML]{D9D9D9}                                         & 10 & 94.120                  & 93.786                 & 94.320                 & 94.659                  & 94.376                 & 94.885                 & 97.713                  & 97.239                 & 97.854                 \\
\cellcolor[HTML]{D9D9D9}                                         & 20 & 92.800                  & 92.438                 & 92.760                 & 93.024                  & 92.738                 & 93.048                 & 95.153                  & 94.935                 & 95.296                 \\
\cellcolor[HTML]{D9D9D9}                                         & 30 & 86.660                  & 85.957                 & 86.120                 & 86.483                  & 85.935                 & 86.143                 & 88.281                  & 87.613                 & 87.790                 \\
\multirow{-5}{*}{\cellcolor[HTML]{D9D9D9}\textbf{gamma}}         & 40 & 81.910                  & 81.667                 & 81.510                 & 84.329                  & 81.215                 & 83.828                 & 84.832                  & 83.421                 & 84.990                 \\
\hline
\cellcolor[HTML]{D9D9D9}                                         & 0  & 99.070                  & 99.080                 & 99.873                 & 99.787                  & 97.971                 & 98.538                 & 98.197                  & 98.112                 & 98.623                 \\
\cellcolor[HTML]{D9D9D9}                                         & 10 & 93.790                  & 93.587                 & 94.040                 & 94.558                  & 94.389                 & 94.720                 & 97.669                  & 97.389                 & 97.867                 \\
\cellcolor[HTML]{D9D9D9}                                         & 20 & 92.150                  & 91.996                 & 92.314                 & 90.629                  & 89.843                 & 90.179                 & 95.680                  & 95.365                 & 95.657                 \\
\cellcolor[HTML]{D9D9D9}                                         & 30 & 86.390                  & 85.607                 & 85.805                 & 87.071                  & 84.159                 & 85.457                 & 88.669                  & 87.945                 & 88.147                 \\
\multirow{-5}{*}{\cellcolor[HTML]{D9D9D9}\textbf{log-normal}}    & 40 & 81.510                  & 81.226                 & 81.420                 & 82.191                  & 81.522                 & 82.955                 & 83.895                  & 83.665                 & 84.530                 \\
\hline
\cellcolor[HTML]{D9D9D9}                                         & 0  & 96.770                  & 93.993                 & 94.836                 & 97.694                  & 97.519                 & 98.403                 & 93.974                  & 93.450                 & 94.670                 \\
\cellcolor[HTML]{D9D9D9}                                         & 10 & 94.970                  & 91.815                 & 95.185                 & 93.988                  & 93.774                 & 94.298                 & 92.498                  & 92.216                 & 92.630                 \\
\cellcolor[HTML]{D9D9D9}                                         & 20 & 90.740                  & 89.767                 & 90.745                 & 91.329                  & 90.016                 & 90.045                 & 89.120                  & 88.406                 & 89.117                 \\
\cellcolor[HTML]{D9D9D9}                                         & 30 & 86.540                  & 85.786                 & 86.033                 & 86.193                  & 85.567                 & 85.784                 & 87.758                  & 86.896                 & 87.120                 \\
\multirow{-5}{*}{\cellcolor[HTML]{D9D9D9}\textbf{trunc. normal}} & 40 & 82.650                  & 82.277                 & 83.865                 & 81.090                  & 81.374                 & 82.839                 & 83.343                  & 82.455                 & 82.840                \\
\hline
\end{tabular}
\caption{Root of the overall mean of the mean squared error (RMMSE) ($\times 10^{-3}$) for different Mack-type bootstraps, different distributional assumptions and different $n$ and $I=10$ for Setup a)}\label{tab_mse_predroot_a}
\end{table}
\end{small}

\section{Conclusion}\label{sec_conclusion}

In this paper, we adopt the stochastic and asymptotic framework that was proposed by \cite{jentschasympThInf} to derive asymptotic theory in Mack's model, also for investigating the consistency properties of the Mack bootstrap proposal.
%In this paper, for investigating the consistency properties of the Mack bootstrap proposal, we adopt the stochastic and asymptotic framework proposed by \cite{jentschasympThInf} to derive asymptotic theory in Mack's model. 
For this purpose, the (conditional) asymptotic theory derived in \cite{jentschasympThInf} serves well as benchmark results for the Mack bootstrap approximations. By splitting the predictive root of the reserve into two additive parts corresponding to process and estimation uncertainty, our approach enables - for the first time - a rigorous investigation of the validity of the Mack bootstrap. We prove that the (conditional) distribution of the asymptotically dominating process uncertainty part is correctly mimicked by the Mack bootstrap if the parametric family of distributions of the individual development factors is correctly specified in Mack's bootstrap. % proposal. 
Otherwise, this will be generally not the case. In contrast, the corresponding (conditional) distribution of the estimation uncertainty part is generally not correctly captured by the bootstrap. Altogether, as the process uncertainty part dominates asymptotically, this proves asymptotic validity of the Mack bootstrap for the whole predictive root of the reserve. However, it also proves that asymptotic pertinence in the sense of \cite{pan2016bootstrap} does not hold.

To remedy this, we propose a more natural alternative Mack-type bootstrap, that uses a different centering and that is designed to capture correctly also the (conditional) distribution of the estimation uncertainty part by using a backward resampling approach. Under suitable assumptions, we demonstrate that the newly proposed alternative Mack-type bootstrap can be indeed asymptotically valid and pertinent.

Our findings are illustrated by simulations, which show that the alternative Mack-type bootstrap performs superior to the original Mack bootstrap in finite samples. An intermediate Mack-type bootstrap provides evidence that the backward resampling %appears to be critical and 
is mainly responsible for this improvement.
%\bigskip
\newpage
\bibliographystyle{apalike}  
\bibliography{references}

\begin{thebibliography}{}

\bibitem[Beutner et~al., 2021]{beutner2021}
Beutner, E., Heinemann, A., and Smeekes, S. (2021).
\newblock {A justification of conditional confidence intervals}.
\newblock {\em Electronic Journal of Statistics}, 15(1):2517 -- 2565.

\bibitem[Bj{\"o}rkwall et~al., 2009]{Bjoerkwall2009}
Bj{\"o}rkwall, S., H{\"o}ssjer, O., and Ohlsson, E. (2009).
\newblock Non-parametric and parametric bootstrap techniques for age-to-age
  development factor methods in stochastic claims reserving.
\newblock {\em Scandinavian Actuarial Journal}, 2009(4):306--331.

\bibitem[Bj{\"o}rkwall et~al., 2010]{bjorkwall2010bootstrapping}
Bj{\"o}rkwall, S., H{\"o}ssjer, O., and Ohlsson, E. (2010).
\newblock Bootstrapping the separation method in claims reserving.
\newblock {\em ASTIN Bulletin: The Journal of the IAA}, 40(2):845--869.

\bibitem[Brockwell and Davis, 1991]{brockwell1991time}
Brockwell, P.~J. and Davis, R.~A. (1991).
\newblock {\em Time Series: Theory and Methods}.
\newblock Springer, New York.

\bibitem[Bruce et~al., 2008]{GRIRO2008}
Bruce, N., Chen, C., Dunne, G., Hinder, I., McMurrough, T., Meyers, G., White,
  A., and Wright, T. (2008).
\newblock Best {E}stimates and {R}eserving {U}ncertainity.
\newblock Actuarial Profession General Insurance (GIRO) Convention. Available
  at:
  https://www.actuaries.org.uk/system/files/documents/pdf/bhprizegibson.pdf.

\bibitem[England, 2002]{england2002addendum}
England, P. (2002).
\newblock Addendum to "analytic and bootstrap estimates of prediction errors in
  claims reserving".
\newblock {\em Insurance: Mathematics and Economics}, 31(3):461--466.

\bibitem[England and Verrall, 1999]{england1999analytic}
England, P. and Verrall, R. (1999).
\newblock Analytic and bootstrap estimates of prediction errors in claims
  reserving.
\newblock {\em Insurance: Mathematics and Economics}, 25(3):281--293.

\bibitem[England and Verrall, 2006]{england2006predictive}
England, P.~D. and Verrall, R.~J. (2006).
\newblock Predictive distributions of outstanding liabilities in general
  insurance.
\newblock {\em Annals of Actuarial Science}, 1(2):221.

\bibitem[Gibson et~al., 2007]{GRIRO2007}
Gibson, L., Archer-Lock, P., Bruce, N., Collins, A., Dunne, G., Felisky, K.,
  Hamilton, A., Jewell, M., Lo, J., Locke, J., Marshall, D., Nicholson, E.,
  Thomas, L., Wilcox, S., Winer, J., and Wright, T. (2007).
\newblock Best {E}stimates and {R}eserving {U}ncertainity.
\newblock Actuarial Profession General Insurance (GIRO) Convention. Available
  at:
  https://www.actuaries.org.uk/system/files/documents/pdf/bhprizegibson.pdf.

\bibitem[Hartl, 2010]{hartl2010bootstrapping}
Hartl, T. (2010).
\newblock Bootstrapping generalized linear models for development triangles
  using deviance residuals.
\newblock In {\em CAS E--Forum Fall}.

\bibitem[Mack, 1993]{mack1993distribution}
Mack, T. (1993).
\newblock Distribution-free {C}alculation of the {S}tandard {E}rror of {C}hain
  {L}adder {R}eserve {E}stimates.
\newblock {\em ASTIN Bulletin: The Journal of the IAA}, 23(2):213--225.

\bibitem[Pan and Politis, 2016]{pan2016bootstrap}
Pan, L. and Politis, D.~N. (2016).
\newblock Bootstrap prediction intervals for linear, nonlinear and
  nonparametric autoregressions.
\newblock {\em Journal of Statistical Planning and Inference}, 177:1--27.

\bibitem[Paparoditis and Shang, 2021]{PaparoditisShang2021}
Paparoditis, E. and Shang, H.~L. (2021).
\newblock Bootstrap prediction bands for functional time series.
\newblock {\em Journal of the American Statistical Association}.

\bibitem[Peremans et~al., 2017]{peremans2017robust}
Peremans, K., Segaert, P., Van~Aelst, S., and Verdonck, T. (2017).
\newblock Robust bootstrap procedures for the chain-ladder method.
\newblock {\em Scandinavian Actuarial Journal}, 2017(10):870--897.

\bibitem[Peters et~al., 2010]{peters2010chain}
Peters, G.~W., W{\"u}thrich, M.~V., and Shevchenko, P.~V. (2010).
\newblock Chain ladder method: Bayesian bootstrap versus classical bootstrap.
\newblock {\em Insurance: Mathematics and Economics}, 47(1):36--51.

\bibitem[Pinheiro et~al., 2003]{pinheiro2003bootstrap}
Pinheiro, P.~J., Andrade~e Silva, J.~M., and de~Lourdes~Centeno, M. (2003).
\newblock Bootstrap methodology in claim reserving.
\newblock {\em Journal of Risk and Insurance}, 70(4):701--714.

\bibitem[Renshaw and Verrall, 1998]{renshaw1998stochastic}
Renshaw, A.~E. and Verrall, R.~J. (1998).
\newblock A {S}tochastic {M}odel {U}nderlying the {C}hain-{L}adder {T}echnique.
\newblock {\em British Actuarial Journal}, 4(4):903--923.

\bibitem[Steinmetz and Jentsch, 2022]{jentschasympThInf}
Steinmetz, J. and Jentsch, C. (2022).
\newblock Asymptotic theory for mack's model.
\newblock {\em Insurance: Mathematics and Economics}, 107:223--268.

\bibitem[Tee et~al., 2017]{tee2017comparison}
Tee, L., K{\"a}{\"a}rik, M., and Viin, R. (2017).
\newblock On comparison of stochastic reserving methods with bootstrapping.
\newblock {\em Risks}, 5(1):2.

\bibitem[Verdonck and Debruyne, 2011]{verdonck2011influence}
Verdonck, T. and Debruyne, M. (2011).
\newblock The influence of individual claims on the chain-ladder estimates:
  Analysis and diagnostic tool.
\newblock {\em Insurance: Mathematics and Economics}, 48(1):85--98.

\end{thebibliography}

\newpage
\appendix

\section{Auxiliary results for Section 4}

\subsection*{Mack bootstrap asymptotics for parameter estimators}\

The following theorem is the Mack bootstrap version of the (unconditional!) Theorem 3.1 in \cite{jentschasympThInf} adapted to the asymptotic framework of Section 2.1. 

\begin{thm}[Asymptotic normality of $\widehat f_{j,n}^*$ conditional on $\mathcal{D}_{I,n}$]\label{CLT_f_boot}
Suppose Assumptions 2.2, 2.3, 2.5 and 4.9 are satisfied and let $\widehat f_{j,n}^*$, $j=0,\ldots,I+n-1$ be defined as in (3.5) according to the Mack bootstrap scheme of Section 3.1. Then, as $n\rightarrow \infty$, the following holds:
\begin{itemize}
	\item[(i)] For each fixed $j\in\mathbb{N}_0=\{0,1,2,\ldots\}$, we have 
\begin{align*}
\sqrt{I+n-j}\left(\widehat f_{j,n}^*-\widehat f_{j,n}\right) \overset{d}{\longrightarrow} \mathcal{N}\left(0, \frac{\s}{\mu_j}\right)	\quad	\text{in probability},
\end{align*}
where ``$\overset{d}{\longrightarrow}$'' denotes convergence in distribution. 
	\item[(ii)] For each fixed $K\in\mathbb{N}_0$, let $\underline{\widehat f^*}_{K,n}=(\widehat f_{0,n}^*,\widehat f_{1,n}^*,\ldots,\widehat f_{K,n}^*)^\prime$ be the $(K+1)$-dimensional Mack bootstrap version of $\underline{\widehat f}_{K,n}=(\widehat f_{0,n},\widehat f_{1,n},\ldots,\widehat f_{K,n})^\prime$. Then, we have
\begin{align*}
J^{1/2}\left(\underline{\widehat f^*}_{K,n}-\underline{\widehat f}_{K,n}\right) \overset{d}{\longrightarrow} \mathcal{N}\left(0, \mathbf{\Sigma}_{K,\underline f}\right)	\quad	\text{in probability},
\end{align*}
where $J^{1/2}=diag\left(\sqrt{I+n+1-j\vphantom{I^2_1}},j=0,\ldots,K\right)$ is a diagonal $(K+1)\times (K+1)$ matrix of inflation factors and $\mathbf{\Sigma}_{K,\underline f}=J_g(\underline{\mu}_K)\mathbf{\Sigma}_{K,\underline{C}}J_g(\underline{\mu}_K)'=diag\left(\frac{\sigma_0^2}{\mu_0},\frac{\sigma_1^2}{\mu_1},\ldots,\frac{\sigma_K^2}{\mu_K}\right)$ is a diagonal $(K+1)\times (K+1)$ covariance matrix, where
\begin{align}\label{Sigma_Ck}
\mathbf{\Sigma}_{K,\underline{C}}=Cov(\underline{C}_{i,K})=\left(\begin{array}{c}\left(\prod_{k=\min(j_1,j_2)}^{\max(j_1,j_2)-1} f_k\right)\tau_{\min(j_1,j_2)}^2 \\ j_1,j_2=0,\ldots,K+1\end{array}\right)
\end{align}
is a $(K+1)\times (K+1)$ matrix,
\begin{align}\label{Jacobian}
J_g(\underline{x})=\left(\begin{array}{ccccc}-\frac{x_1}{x_0^2} & \frac{1}{x_0} & 0 & \cdots & 0	\\ 0 & -\frac{x_2}{x_1^2} & \frac{1}{x_1}	& \ddots & \vdots	\\	\vdots & \ddots & \ddots & \ddots	& 0	\\ 0 & \cdots & 0 & -\frac{x_{K+1}}{x_K^2} & \frac{1}{x_K}\end{array}\right)
\end{align}
is a $(K+1)\times (K+2)$ matrix, and $\underline{\mu}_K=(\mu_0,\ldots,\mu_{K+1})^\prime$ as derived in the proof of Theorem 3.1 in \cite{jentschasympThInf}. 
\end{itemize}
\end{thm}

As the \emph{unconditional} limiting distributions obtained in Theorem \ref{CLT_f_boot} above and in Theorem 3.1 in \cite{jentschasympThInf} coincide, the Mack bootstrap is \emph{unconditionally}, that is without conditioning on $\mathcal{Q}_{I,n}^*=\mathcal{Q}_{I,n}$, consistent for an arbitrary, but fixed number of estimators of development factors. That is, for each fixed $K\in\mathbb{N}_0$, we have
\begin{align*}
d_K\left(\mathcal{L}^*\left(J^{1/2}\left(\underline{\widehat f^*}_{K,n}-\underline{\widehat f}_{K,n}\right)\right),\mathcal{L}\left(J^{1/2}\left(\underline{\widehat f}_{K,n}-\underline{f}_K\right)\right)\right)=o_P(1),
\end{align*}
where $\underline{f}_K=(f_0,f_1,\ldots,f_K)^\prime$ and $d_K$ denotes the Kolmogorov distance between two probability distributions.

The following direct corollary is the Mack bootstrap version of Corollary 3.2 in \cite{jentschasympThInf} adapted to the asymptotic framework of Section 2.1.

\begin{cor}[Asymptotic normality for products of $\widehat f_{j,n}^*$'s conditional on $\mathcal{D}_{I,n}$]\label{CLTProdf}
Suppose the assumptions of Theorem \ref{CLT_f_boot} hold. Then, as $n\rightarrow \infty$, the following holds:
\begin{itemize}
		\item[(i)] For each fixed $K\in\mathbb{N}_0$ and $i=0,\ldots,K$, we have		
		\begin{align*}
		\sqrt{I+n+1}\left(\prod_{j=i}^{K}\widehat f_{j,n}^*-\prod_{j=i}^{K}\widehat f_{j,n}\right)\overset{d}{\longrightarrow} \mathcal{N}\left(0,\sum_{j=i}^{K}\frac{\s}{\mu_j}\prod^{K}_{l=i, l \neq j}f^2_l\right)	\quad	\text{in probability}.
		\end{align*}
		\item[(ii)] For each fixed $K\in\mathbb{N}_0$, we have also joint convergence, that is,
		\begin{align*}
		\sqrt{I}\begin{pmatrix}
		\prod^K_{j=i}\widehat f_{j,n}^*-\prod_{j=i}^{K}\widehat f_{j,n}	\\
		i=0,\ldots,K
		\end{pmatrix} \overset{d}{\longrightarrow} \mathcal{N}\left(0, \bm{\Sigma}_{K,\prod \f}\right)	\quad	\text{in probability},
		\end{align*}
where $\bm{\Sigma}_{K,\prod \f}=J_h(\underline{f}_K)\mathbf{\Sigma}_{K,\underline f}J_h(\underline{f}_K)'=(\bm{\Sigma}_{K,\prod \f}(i_1,i_2))_{i_1,i_2=0,\ldots,K}$ is a $(K+1)\times (K+1)$ covariance matrix with entries
\begin{align*}
\bm{\Sigma}_{K,\prod \f}(i_1,i_2) &= \sum^K_{j=\max(i_1,i_2)}\frac{\sigma^2_j}{\mu_j}\prod^{K}_{l=\max(i_1,i_2), l \neq j}f^2_l\prod^{\max(i_1,i_2)-1}_{m=\min(i_1,i_2)}f_m,
\end{align*}
for $i_1,i_2=0,\dots,K$. Here, $\mathbf{\Sigma}_{K,\underline f}$ is defined in Theorem \ref{CLT_f_boot}(ii) and 
\begin{align}\label{Jacobian2}
J_h(\underline{x})=\left(\begin{array}{cccc}\prod_{l=0,l\not=0}^K x_l & \prod_{l=0,l\not=1}^K x_l & \cdots & \prod_{l=0,l\not=K}^K x_l	\\0 & \prod_{l=1,l\not=1}^K x_l & \ddots & \vdots	\\	\vdots & \ddots & \ddots & \vdots	\\ 0 & \cdots	&  0	&	\prod_{l=K,l\not=K}^K x_l\end{array}\right)
\end{align}
as derived in the proof of Corollary 3.2 in \cite{jentschasympThInf}.
	\end{itemize}
\end{cor}

\subsection{Proof of Theorem \ref{CLT_f_boot}}
By construction of the Mack bootstrap estimators $\widehat f_{j,n}^*$, $j=0,\ldots,I+n-1$ according to (3.5, for each fixed $K\in\mathbb{N}_0$, the $K+1$ estimators $\widehat f_{0,n}^*,\widehat f_{1,n}^*,\ldots,\widehat f_{K,n}^*$ are independent conditional on $\mathcal{D}_{I,n}$. Hence, it is actually sufficient to prove part $(i)$. For any fixed $j$ and from (2.10) and (3.5), using $C_{i,j+1} = C_{i,j}F_{i,j}$, we get immediately
\begin{align*}
\sqrt{I+n-j}\left(\widehat f_{j,n}^*-\widehat f_{j,n}\right) &= \sqrt{I+n-j}\left(\frac{\sum_{i=-n}^{I-j-1}C_{i,j}F_{i,j}^*}{\sum_{k=-n}^{I-j-1}C_{k,j}} - \frac{\sum_{i=-n}^{I-j-1}C_{i,j}}{\sum_{k=-n}^{I-j-1}C_{k,j}}\widehat f_{j,n}\right)	\\
&= \sum_{i=-n}^{I-j-1}\frac{C_{i,j}\left(F_{i,j}^*-\widehat f_{j,n}\right)}{\frac{1}{\sqrt{I+n-j}}\sum_{k=-n}^{I-j-1}C_{k,j}}	\\
&=: \sum_{i=-n}^{I-j-1} Z_{i,n}.
\end{align*}
Noting that, for all $j$, $(Z_{i,n}, i=-n,\ldots,I-j-1,\ n\in\mathbb{N}_0)$ forms a triangular array of random variables that are independent conditional on $\mathcal{D}_{I,n}$, we can make use of a (conditional) Lyapunov CLT to prove asymptotic normality. First, for the bootstrap mean, using measurability of all $C_{i,j}$'s and of $\widehat f_{j,n}$ in $Z_{i,n}$ with respect to $\mathcal{D}_{I,n}$, we get
\begin{align*}
E^*\left(Z_{i,n}^*\right) = E^*\left(\frac{C_{i,j}\left(F_{i,j}^*-\widehat f_{j,n}\right)}{\frac{1}{\sqrt{I+n-j}}\sum_{k=-n}^{I-j-1}C_{k,j}}\right) = \frac{C_{i,j}}{{\frac{1}{\sqrt{I+n-j}}\sum_{k=-n}^{I-j-1}C_{k,j}}}\left(E^*\left(F_{i,j}^*\right)-\widehat f_{j,n}\right).
\end{align*}
Further, by the construction of Mack's bootstrap, for any fixed $j$ and $i=-n, \dots, I-j-1$, we have $E^*(r_{i,j}^*)=0$ such that
\begin{align*}
E^*(F_{i,j}^*) = E^*\left(\widehat f_{j,n}+\frac{\widehat \sigma_{j,n}}{\sqrt{C_{i,j}}}r_{i,j}^*\right) = \widehat f_{j,n}+\frac{\widehat \sigma_{j,n}}{\sqrt{C_{i,j}}}E^*\left(r_{i,j}^*\right) = \widehat f_{j,n}
\end{align*}
leading to $E^*(Z_{i,n}^*)=0$. Second, for the bootstrap variance, we get
\begin{align*}
Var^*\left(Z_{i,n}^*\right) = \frac{C_{i,j}^2}{\left(\frac{1}{\sqrt{I+n-j}}\sum_{k=-n}^{I-j-1}C_{k,j}\right)^2} Var^*\left(F_{i,j}^*\right)
\end{align*}
and, from the particular construction of Mack's bootstrap leading to $E^*(r_{ij}^*)=0$ and $E^*(r_{ij}^{*2})=1$, we obtain
\begin{align*}
Var^*(F_{i,j}^*) &= E^*\left(\left(\widehat f_{j,n}+\frac{\widehat \sigma_{j,n}}{\sqrt{C_{i,j}}}r_{i,j}^*\right)^2\right)-\widehat f_{j,n}^2	\\
&= \widehat f_{j,n}^2 +2\widehat f_{j,n}\frac{\widehat \sigma_{j,n}}{\sqrt{C_{i,j}}}E^*(r_{i,j}^*) + \left(\frac{\widehat \sigma_{j,n}}{\sqrt{C_{i,j}}}\right)^2 E^*(r_{i,j}^{*2})-\widehat f_{j,n}^2	\\
&= \frac{\widehat \sigma_{j,n}^2}{C_{i,j}}
\end{align*}
such that
\begin{align*}
Var^*\left(Z_{i,n}^*\right) = \frac{C_{i,j}\widehat \sigma_{j,n}^2}{\left(\frac{1}{\sqrt{I+n-j}}\sum_{k=-n}^{I-j-1}C_{k,j}\right)^2}  
\end{align*}
and, altogether,
\begin{align*}
Var^*\left(\sum_{i=-n}^{I-j-1}Z_{i,n}^*\right) &= \sum_{i=-n}^{I-j-1}\frac{C_{i,j}\widehat \sigma_{j,n}^2}{\left(\frac{1}{\sqrt{I+n-j}}\sum_{k=-n}^{I-j-1}C_{k,j}\right)^2} = \frac{\left(\frac{1}{I+n-j}\sum_{i=-n}^{I-j-1}C_{i,j}\right)\widehat \sigma_{j,n}^2}{\left(\frac{1}{I+n-j}\sum_{k=-n}^{I-j-1}C_{k,j}\right)^2}	\\
&= \frac{\widehat \sigma_{j,n}^2}{\frac{1}{I+n-j}\sum_{k=-n}^{I-j-1}C_{k,j}}.
\end{align*}
Letting $n\rightarrow \infty$, making use of Assumption 2.5, we get $\widehat\sigma_{j,n}^2\rightarrow \sigma_j^2$ by Theorem 3.5 in \cite{jentschasympThInf}, as well as
\begin{align*}
\frac{1}{I+n-j}\sum_{k=-n}^{I-j-1}C_{k,j}\overset{p}{\longrightarrow} \mu_j
\end{align*}
by a WLLN using that, for all $j$, $(C_{k,j},k\in\mathbb{Z},\ k\leq I-j-1)$ are iid by Assumption 2.1(iii) with (finite) mean $\mu_j$ and variance $\tau_j^2$ according to (2.19) and (2.20), respectively.

Finally, it remains to prove a Lyapunov condition to complete the proof. Choosing $\delta=2$ for the Lyapunov condition, for any $j$, it is sufficient to show that
\begin{align*}
\sum_{i=-n}^{I-K-1}E^*\left(\left(\frac{C_{i,j}\left(F_{i,j}^*-\widehat f_{j,n}\right)}{\frac{1}{\sqrt{I+n-K}}\sum_{k=-n}^{I-K-1}C_{k,j}}\right)^4\right)\overset{p}{\longrightarrow} 0.
\end{align*}
Due to measurability of all $C_{i,j}$'s with respect to $\mathcal{D}_{I,n}$, we get
\begin{align*}
& \sum_{i=-n}^{I-K-1}E^*\left(\left(\frac{C_{i,j}\left(F_{i,j}^*-\widehat f_{j,n}\right)}{\frac{1}{\sqrt{I+n-K}}\sum_{k=-n}^{I-K-1}C_{k,j}}\right)^4\right)	\\
&= \frac{1}{\left(\frac{1}{I+n-K}\sum_{k=-n}^{I-K-1}C_{k,j}\right)^4}\frac{1}{(I+n-K)^2}\sum_{i=-n}^{I-K-1}C_{i,j}^4E^*\left(\left(F_{i,j}^*-\widehat f_{j,n}\right)^4\right).
\end{align*}
Further, as $\frac{1}{I+n-K}\sum_{k=-n}^{I-K-1}C_{k,j}=O_P(1)$, it is sufficient to show that 
\begin{align*}
\frac{1}{I+n-K}\sum_{i=-n}^{I-K-1}C_{i,j}^4E^*\left(\left(F_{i,j}^*-\widehat f_{j,n}\right)^4\right)=O_P(1).
\end{align*}
For this purpose, we have to compute $E^*((F_{i,j}^*-\widehat f_{j,n})^4)$ next. By plugging-in for $F_{i,j}^*$, we get
\begin{align*}
E^*\left(\left(F_{i,j}^*-\widehat f_{j,n}\right)^4\right) %= E^*\left(\left(\widehat f_{j,n}+\frac{\widehat \sigma_{j,n}}{\sqrt{C_{i,j}}}r_{i,j}^*-\widehat f_{j,n}\right)^4\right) 
= E^*\left(\left(\frac{\widehat \sigma_{j,n}}{\sqrt{C_{i,j}}}r_{i,j}^*\right)^4\right) = \frac{\widehat \sigma_{j,n}^4}{C_{i,j}^2}E^*\left(r_{i,j}^{*4}\right)
\end{align*}
leading to
\begin{align*}
\frac{1}{I+n-K}\sum_{i=-n}^{I-K-1}C_{i,j}^4E^*\left(\left(F_{i,j}^*-\widehat f_{j,n}\right)^4\right) = \frac{\widehat \sigma_{j,n}^4}{I+n-K}\sum_{i=-n}^{I-K-1}C_{i,j}^2E^*\left(r_{i,j}^{*4}\right).
\end{align*}
Further, as $\frac{\widehat \sigma_{j,n}^4}{I+n-K}\sum_{i=-n}^{I-K-1}C_{i,j}^2=O_P(1)$, it remains to show that $E^*\left(r_{i,j}^{*4}\right)=O_P(1)$ holds as well. By construction, we have\footnote{Note that we implicitly assume that \emph{only} $\widehat\sigma_{I+n-1,n}^2$ is estimated as zero; see Section 3.1.} 
\begin{align*}
E^*\left(r_{i,j}^{*4}\right) = \frac{2}{(I+n+1)(I+n)-2}\sum_{s=0}^{I+n-2}\sum_{t=-n}^{I-s-1} \widetilde r_{t,s}^4.
\end{align*}
In the following, suppose for convenience that $\widetilde r_{t,s}=\widehat r_{t,s}$. However, the arguments for $\widetilde r_{t,s}$ including re-centering (and re-scaling) are essentially the same, but tedious and lengthy. In this case, by plugging-in for $\widetilde r_{t,s}$, we get 
\begin{align*}
E^*\left(r_{i,j}^{*4}\right) &= \frac{2}{(I+n+1)(I+n)-2}\sum_{s=0}^{I+n-2}\sum_{t=-n}^{I-s-1} \left(\frac{\sqrt{C_{t,s}}(F_{t,s}-\widehat f_{s,n})}{\widehat{\sigma}_{s,n}}\right)^4	\\
&= \frac{2}{(I+n+1)(I+n)-2}\sum_{s=0}^{I+n-2}\sum_{t=-n}^{I-s-1} \frac{C_{t,s}^2}{\widehat{\sigma}_{s,n}^4}\left(F_{t,s}^4-4F_{t,s}^3\widehat f_{s,n}+6F_{t,s}^2\widehat f_{s,n}^2-4F_{t,s}\widehat f_{s,n}^3+\widehat f_{s,n}^4\right)
\end{align*}
By Assumption 4.9, for $n\rightarrow \infty$, 
we have
\begin{align*}
\sup_{j=0,\ldots,I+n-1}\frac{\widehat f_{j,n}}{f_j}=O_P(1)	\quad	\text{and}	\quad	\sup_{j=0,\ldots,I+n-2}\frac{\sigma_j^2}{\widehat \sigma_{j,n}^2}=O_P(1).
\end{align*}
Hence, we can bound $E^*(r_{i,j}^{*4})$ above by
\begin{align*}
& O_P(1)\left(\frac{2}{(I+n+1)(I+n)-2}\sum_{s=0}^{I+n-2}\sum_{t=-n}^{I-s-1} \frac{C_{t,s}^2}{\sigma_s^4}\left(F_{t,s}^4-4F_{t,s}^3f_s+6F_{t,s}^2f_s^2-4F_{t,s}f_s^3+f_s^4\right)\right)	\\
=& O_P(1)\left(\frac{2}{(I+n+1)(I+n)-2}\sum_{s=0}^{I+n-2}\sum_{t=-n}^{I-s-1} \frac{C_{t,s}^2}{\sigma_s^4}\left(F_{t,s}-f_s\right)^4\right).
\end{align*}
Finally, the term in brackets on the last right-hand side is a sum consisting of non-negative summands, which is also $O_P(1)$ as its expectation is bounded because the $\kappa_j^{(4)}$'s defined in (2.21) are assumed to form a bounded sequence $((\kappa_j^{(4)}/\sigma_j^4), j\in\mathbb{N}_0)$ again according to Assumption 4.9.	\hfill	$\square$

\subsection{Proof of Corollary \ref{CLTProdf}}\label{proofCLTProdf}
The proof follows from an application of the delta method and Theorem \ref{CLT_f_boot} and is completely analogous to the proof of Corollary 3.2 in \cite{jentschasympThInf}.	\hfill $\square$

\bigskip

\section{Proofs of Section 4}

\subsection{Proof of Theorem 4.8}

As the (conditional) $L_2$-convergence result in Theorem 4.1 implies the (conditional) convergence in distribution in (4.6), for $n\rightarrow \infty$, it remains to show
\begin{align}\label{eq_goal_asymptotics}
(R_{I,n}^*-\widehat{R}_{I,n})_1|\left(\mathcal{Q}_{I,n}^*=\mathcal{Q}_{I,n}, \mathcal{D}_{I,n}\right) \overset{d}{\longrightarrow} \mathcal{G}_1|\mathcal{Q}_{I,\infty}
\end{align}
with $E^*((R_{I,n}^*-\widehat{R}_{I,n})_1|\mathcal{Q}_{I,n}^*=\mathcal{Q}_{I,n})\rightarrow 0$ and 
\begin{align}
Var^*\left((R_{I,n}^*-\widehat{R}_{I,n})_1|\mathcal{Q}_{I,n}^*=\mathcal{Q}_{I,n}\right)\rightarrow Var\left((R_{I,\infty}-\widehat{R}_{I,\infty})_1|\mathcal{Q}_{I,\infty}\right) 
\end{align}
in probability, respectively.
 
Nevertheless, the asymptotic theory for $(R_{I,n}^*-\widehat{R}_{I,n})_1$ conditional on $\mathcal{Q}_{I,n}^*=\mathcal{Q}_{I,n}$ and $\mathcal{D}_{I,n}$ is not straightforward as it is composed of sums and products consisting asymptotically of infinitely many summands and factors. Hence, we decompose $(R_{I,n}^*-\widehat{R}_{I,n})_1$ by truncating these sums and products to be able to apply Proposition 6.3.9 in \cite{brockwell1991time}. For this purpose, let $K\in \mathbb{N}_0$ be fixed and suppose $I,n\in\mathbb{N}_0$ are large enough such that $K<I+n-1$. Then, we have\\

\begin{align*}
(R_{I,n}^*-\widehat{R}_{I,n})_1 =& \sum^{I+n}_{i=0} C_{I-i,i}^*\left(\prod^{I+n-1}_{j=i}F_{I-i,j}^*-\prod_{j=i}^{I+n-1}\widehat f_{j,n}^*\right)
\\
=&\sum_{i=0}^K C_{I-i,i}^*\left(\prod^{K}_{j=i}F_{I-i,j}^*-\prod_{j=i}^{K}\widehat f_{j,n}^*\right)\\
&+\sum_{i=0}^K C_{I-i,i}^*\left(\prod^{K}_{j=i}F_{I-i,j}^*\left(\prod^{I+n-1}_{l=K+1}F_{I-i,l}^*-1\right)-\prod^{K}_{j=i}\widehat f_{j,n}^*\left(\prod^{I+n-1}_{l=K+1}\widehat f_{l,n}^*-1\right)\right)\\
&+ \sum^{I+n}_{i=K+1}C_{I-i,i}^*\left(\prod^{I+n-1}_{j=i}F_{I-i,j}^*-\prod^{I+n-1}_{j=i}\widehat f_{j,n}^*\right)\\
=:& A_{1,K,I,n}^*+A_{2,K,I,n}^*+A_{3,K,I,n}^*.
\end{align*}
Hence, to derive the claimed conditional limiting distribution, it suffices to show that, a) for all $K\in \mathbb{N}_0$, $A_{1,K,I,n}^*|(\mathcal{Q}_{I,n}^*=\mathcal{Q}_{I,n}, \mathcal{D}_{I,n})\overset{d}{\rightarrow} \mathcal{G}_{1,K}|\mathcal{Q}_{I,\infty}$ in probability as $n\rightarrow\infty$ for some (conditional) distribution $\mathcal{G}_{1,K}|\mathcal{Q}_{I,\infty}$, b) $\mathcal{G}_{1,K}|\mathcal{Q}_{I,\infty}\overset{d}{\rightarrow}\mathcal{G}_{1}|\mathcal{Q}_{I,\infty}$ as $K\rightarrow \infty$, and c) that, for all $\epsilon>0$, we have
\begin{align}
\lim\limits_{K\rightarrow \infty}\limsup\limits_{n\rightarrow\infty}P^*\left(|A_{2,K,I,n}^*|>\epsilon|\mathcal{Q}_{I,n}^*=\mathcal{Q}_{I,n}\right)=0	\quad	\text{and}	\quad	\lim\limits_{K\rightarrow \infty}\limsup\limits_{n\rightarrow\infty}P^*\left(|A_{3,K,I,n}^*|>\epsilon|\mathcal{Q}_{I,n}^*=\mathcal{Q}_{I,n}\right)=0.
\end{align}
We begin with showing part a). The parametric family of (conditional) distributions used to generate the $F_{i,j}|C_{i,j}$ and $F_{i,j}^*|C_{i,j}^*$ is continuous with respect to $C_{i,j}$, $f_j$, $\sigma_j^2$ and $C_{i,j}^*$, $\widehat f_{j,n}^*$, $\widehat \sigma_{j,n}^2$, respectively, by Assumption 4.6. Hence, as $\widehat f_{j,n}-f_j=O_P((I+n-1)^{-1/2})$, $\widehat f_{j,n}^*-\widehat f_{j,n}=O_{P^*}((I+n-1)^{-1/2})$ and $\widehat \sigma_{j,n}^2-\sigma_j^2=O_P((I+n-1)^{-1/2})$ holds for all fixed $j\in\mathbb{N}_0$, we can conclude that, for all fixed $K\in\mathbb{N}_0$ and as $n\rightarrow \infty$, that
\begin{align}
\sum_{i=0}^K C_{I-i,i}^*\left(\prod^{K}_{j=i}F_{I-i,j}^*-\prod_{j=i}^{K}\widehat f_{j,n}^*\right)|(\mathcal{Q}_{I,n}^*=\mathcal{Q}_{I,n}, \mathcal{D}_{I,n})\overset{d}{\longrightarrow}\sum_{i=0}^K C_{I-i,i}\left(\prod^{K}_{j=i}F_{I-i,j}-\prod_{j=i}^{K}f_j\right)|\mathcal{Q}_{I,\infty}%\sim \mathcal{G}_{1,K}|\mathcal{Q}_{I,\infty}
\end{align}
in probability, which proves $A_{1,K,I,n}^*|(\mathcal{Q}_{I,n}^*=\mathcal{Q}_{I,n}, \mathcal{D}_{I,n})\overset{d}{\rightarrow}\mathcal{G}_{1,K}|\mathcal{Q}_{I,\infty}$. For part b), by letting also $K\rightarrow \infty$, we get immediately 
\begin{align}
\sum_{i=0}^K C_{I-i,i}\left(\prod^{K}_{j=i}F_{I-i,j}-\prod_{j=i}^{K}f_j\right)|\mathcal{Q}_{I,\infty}\overset{d}{\longrightarrow}\sum_{i=0}^\infty C_{I-i,i}\left(\prod^\infty_{j=i}F_{I-i,j}-\prod_{j=i}^\infty f_j\right)|\mathcal{Q}_{I,\infty}\sim \mathcal{G}_1|\mathcal{Q}_{I,\infty},
\end{align}
which proves $\mathcal{G}_{1,K}|\mathcal{Q}_{I,\infty}\overset{d}{\rightarrow}\mathcal{G}_1|\mathcal{Q}_{I,\infty}$. Before we prove part c), let us also consider mean and variance of $A_{1,K,I,n}^*$ (conditional on $\mathcal{Q}_{I,n}^*=\mathcal{Q}_{I,n}$ and $\mathcal{D}_{I,n}$). For the mean, using measurability of $C_{I-i,i}$ with respect to $\mathcal{D}_{I,n}$ and the law of iterated expectations, we have%\footnote{CJ: Notation hier einfuehren!}
\begin{align*}
& E^*(A_{1,K,I,n}^*|\mathcal{Q}_{I,n}^*=\mathcal{Q}_{I,n})	\\
&= E^*(E^*(A_{1,K,I,n}^*|\mathcal{Q}_{I,n}^*=\mathcal{Q}_{I,n},\mathcal{F}_{I,n}^*)|\mathcal{Q}_{I,n}^*=\mathcal{Q}_{I,n})	\\
&= \sum_{i=0}^K C_{I-i,i} E^*\left(E^*\left(\left.\prod^{K}_{j=i}F_{I-i,j}^*-\prod_{j=i}^{K}\widehat f_{j,n}^*\right|\mathcal{Q}_{I,n}^*=\mathcal{Q}_{I,n},\mathcal{F}_{I,n}^*\right)|\mathcal{Q}_{I,n}^*=\mathcal{Q}_{I,n}\right)	\\
&=0
\end{align*}
due to
\begin{align}
E^*\left(\left.\prod^{K}_{j=i}F_{I-i,j}^*-\prod_{j=i}^{K}\widehat f_{j,n}^*\right|\mathcal{Q}_{I,n}^*=\mathcal{Q}_{I,n},\mathcal{F}_{I,n}^*\right) &= E^*\left(\prod^{K}_{j=i}F_{I-i,j}^*|\mathcal{Q}_{I,n}^*=\mathcal{Q}_{I,n},\mathcal{F}_{I,n}^*\right)-\prod_{j=i}^{K}\widehat f_{j,n}^*	\nonumber	\\
&=\prod_{j=i}^{K}\widehat f_{j,n}^*-\prod_{j=i}^{K}\widehat f_{j,n}^*=0	\label{double_cond_mean_A1KIn}
\end{align}
using similar arguments as used to show $E(\prod^{K}_{j=i}F_{I-i,j})=\prod^{K}_{j=i}f_j$. 
Similarly, using the law of total variance and \eqref{double_cond_mean_A1KIn}, we get for the variance
\begin{align}
Var^*(A_{1,K,I,n}^*|\mathcal{Q}_{I,n}^*=\mathcal{Q}_{I,n}) =& E^*\left(Var^*(A_{1,K,I,n}^*|\mathcal{Q}_{I,n}^*=\mathcal{Q}_{I,n},\mathcal{F}_{I,n}^*)|\mathcal{Q}_{I,n}^*=\mathcal{Q}_{I,n}\right)	\nonumber	\\
&+ Var^*\left(E^*(A_{1,K,I,n}^*|\mathcal{Q}_{I,n}^*=\mathcal{Q}_{I,n},\mathcal{F}_{I,n}^*)|\mathcal{Q}_{I,n}^*=\mathcal{Q}_{I,n}\right)	\nonumber	\\
=& E^*\left(\sum_{i=0}^K C_{I-i,i}\sum_{j=i}^K\left(\prod_{k=i}^{j-1}\widehat{f}_{k,n}^*\right)\widehat{\sigma}^2_{j,n}\left(\prod_{l=j+1}^K\widehat{f}^{*2}_{l,n}\right)|\mathcal{Q}_{I,n}^*=\mathcal{Q}_{I,n}\right)	\nonumber	\\
=& \sum_{i=0}^K C_{I-i,i}\sum_{j=i}^K\widehat{\sigma}^2_{j,n}E^*\left(\left(\prod_{k=i}^{j-1}\widehat{f}_{k,n}^*\right)\left(\prod_{l=j+1}^K\widehat{f}^{*2}_{l,n}\right)\right)	\label{cond_var_A1KIn}
\end{align}
due to the fact that $\widehat{f}_{k,n}^*$'s are independent of the condition $\mathcal{Q}_{I,n}^*=\mathcal{Q}_{I,n}$ and because of
\begin{align*}
Var^*(A_{1,K,I,n}^*|\mathcal{Q}_{I,n}^*=\mathcal{Q}_{I,n},\mathcal{F}_{I,n}^*)=\sum_{i=0}^K C_{I-i,i}\sum_{j=i}^K\left(\prod_{k=i}^{j-1}\widehat{f}_{k,n}^*\right)\widehat{\sigma}^2_{j,n}\left(\prod_{l=j+1}^K\widehat{f}^{*2}_{l,n}\right)
\end{align*}
obtained by similar arguments as used in the proof of \cite[Theorem 4.3]{jentschasympThInf} and using the measurability of $C_{I-i,i}$ and $\widehat{\sigma}^2_{j,n}$ with respect to $\mathcal{D}_{I,n}$. Now, using similar arguments as in \cite[Theorem 4.7]{jentschasympThInf} and exploiting the fact that the $\widehat{f}_{k,n}^*$'s are stochastically independent conditional on $\mathcal{D}_{I,n}$, for the expectation in \eqref{cond_var_A1KIn}, we get
\begin{align*}
E^*\left(\left(\prod_{k=i}^{j-1}\widehat{f}_{k,n}^*\right)\left(\prod_{l=j+1}^K\widehat{f}^{*2}_{l,n}\right)\right) &= \left(\prod_{k=i}^{j-1}E^*\left(\widehat{f}_{k,n}^*\right)\right)\left(\prod_{l=j+1}^KE^*\left(\widehat{f}^{*2}_{l,n}\right)\right)	\\
&= \left(\prod_{k=i}^{j-1}\widehat{f}_{k,n}\right)\left(\prod_{l=j+1}^K\left(\frac{\widehat \sigma_{l,n}^2}{\sum_{k=-n}^{I-l-1}C_{k,l}}+\widehat{f}^{2}_{l,n}\right)\right)		\\
&= \left(\prod_{k=i}^{j-1}\widehat{f}_{k,n}\right)\left(\prod_{l=j+1}^K\widehat{f}^{2}_{l,n}\right)+O_P\left(\frac{1}{I+n}\right)
\end{align*}
due to, for all $c\in\{0,\ldots,K\}$, we have
\begin{align}
E^*\left(\widehat f_{c,n}^{*2}|\mathcal{B}_{I,n}(c)\right) = \frac{\widehat \sigma_{c,n}^2}{\sum_{k=-n}^{I-c-1}C_{k,c}}+\widehat{f}_{c,n}^2,	\label{cond_boot_expectation_f_squared}
\end{align}
where $\mathcal{B}_{I,n}(k)=\left\{C_{i,j}|i=-n,\ldots,I,~j=0,\ldots,k,~i+j\leq I+n\right\}$ denotes all elements of $\mathcal{D}_{I,n}$ up to its $k$th column, and because of $\widehat \sigma_{l,n}^2\rightarrow \sigma_l^2$ in probability for all $l\in\{0,\ldots,K\}$ and
\begin{align*}
\frac{1}{\sum_{k=-n}^{I-l-1}C_{k,l}}\leq \frac{1}{(I+n-l)\epsilon^{l}}\leq \frac{1}{(I+n-K)\epsilon^{K}}=O\left(\frac{1}{I+n}\right)
\end{align*}
as $K$ is fixed. This leads to
\begin{align}
Var^*(A_{1,K,I,n}^*|\mathcal{Q}_{I,n}^*=\mathcal{Q}_{I,n}) =& \sum_{i=0}^K C_{I-i,i}\sum_{j=i}^K\left(\prod_{k=i}^{j-1}\widehat{f}_{k,n}\right)\widehat{\sigma}^2_{j,n}\left(\prod_{l=j+1}^K\widehat{f}^2_{l,n}\right)+O_P\left(\frac{1}{I+n}\right)	\\
\overset{p}{\longrightarrow}& \sum_{i=0}^K C_{I-i,i}\sum_{j=i}^K\left(\prod_{k=i}^{j-1}f_k\right)\sigma^2_j\left(\prod_{l=j+1}^Kf^2_l\right)
\end{align}
as $n\rightarrow \infty$ for all $K$ fixed, because $\widehat f_{j,n}-f_j=O_P((I+n-1)^{-1/2})$ and $\widehat \sigma_{j,n}^2-\sigma_j^2=O_P((I+n-1)^{-1/2})$ for all $j\in\mathbb{N}_0$. Finally, letting $K\rightarrow \infty$, we get
\begin{align}
\sum_{i=0}^K C_{I-i,i}\sum_{j=i}^K\left(\prod_{k=i}^{j-1}f_k\right)\sigma^2_j\left(\prod_{l=j+1}^Kf^2_l\right)\longrightarrow \sum_{i=0}^\infty C_{I-i,i}\sum_{j=i}^\infty\left(\prod_{k=i}^{j-1}f_k\right)\sigma^2_j\left(\prod_{l=j+1}^\infty f^2_l\right),
\end{align}
which equals $Var((R_{I,\infty}-\widehat{R}_{I,\infty})_1|\mathcal{Q}_{I,\infty})$. Hence, it remains to show part c) to complete the proof. We begin with showing part c) for $A_{2,K,I,n}^*$. By similar arguments used above, for the mean, we have $E^*(A_{2,K,I,n}^*|\mathcal{Q}_{I,n}^*=\mathcal{Q}_{I,n})=0$ due to $E^*(A_{2,K,I,n}^*|\mathcal{Q}_{I,n}^*=\mathcal{Q}_{I,n},\mathcal{F}_{I,n}^*)=0$ and, for the variance, we have $Var^*(A_{2,K,I,n}^*|\mathcal{Q}_{I,n}^*=\mathcal{Q}_{I,n})=E^*((A_{2,K,I,n}^*)^2|\mathcal{Q}_{I,n}^*=\mathcal{Q}_{I,n})=E^*(E^*((A_{2,K,I,n}^*)^2|\mathcal{Q}_{I,n}^*=\mathcal{Q}_{I,n},\mathcal{F}_{I,n}^*)|\mathcal{Q}_{I,n}^*=\mathcal{Q}_{I,n})$. For the inner expectation, using stochastic independence over accident years leading to stochastic independent summands of $A_{2,K,I,n}^*$ (conditional on $\mathcal{Q}_{I,n}^*=\mathcal{Q}_{I,n}$, $\mathcal{D}_{I,n}$ and $\mathcal{F}_{I,n}^*$), we get
\begin{align*}
&E^*((A_{2,K,I,n}^*)^2|\mathcal{Q}_{I,n}^*=\mathcal{Q}_{I,n},\mathcal{F}_{I,n}^*)	\\
&= \sum^{K}_{i=0}C_{I-i,i}^2E^*\left(\left.\left(\prod^{K}_{j=i}F^*_{I-i,j}\left(\prod^{I+n-1}_{l=K+1}F_{I-i,l}^*-1\right)-\prod^{K}_{j=i}\widehat f_{j,n}^*\left(\prod^{I+n-1}_{l=K+1}\widehat f_{l,n}^*-1\right)\right)^2\right|\mathcal{Q}_{I,n}^*=\mathcal{Q}_{I,n},\mathcal{F}^*_{I,n}\right)	\\
&= \sum^{K}_{i=0}C_{I-i,i}^2\left[E^*\left(\left.\left(\prod^{K}_{j=i}F_{I-i,j}^*\left(\prod^{I+n-1}_{l=K+1}F^*_{I-i,l}-1\right)\right)^2\right|\mathcal{Q}_{I,n}^*=\mathcal{Q}_{I,n},\mathcal{F}_{I,n}^*\right)-\left(\prod^{K}_{j=i}\widehat f_{j,n}^*\left(\prod^{I+n-1}_{l=K+1}\widehat f_{l,n}^*-1\right)\right)^2\right].
\end{align*}
For the term corresponding to the first term in brackets on the last right-hand side, we get
\begin{align}
& \sum^{K}_{i=0}C_{I-i,i}^2E^*\left(\left.\left(\prod^{K}_{j=i}F^*_{I-i,j}\left(\prod^{I+n-1}_{l=K+1}F^*_{I-i,l}-1\right)\right)^2\right|\mathcal{Q}_{I,n}^*=\mathcal{Q}_{I,n},\mathcal{F}_{I,n}^*\right)	\nonumber\\
=& \sum^{K}_{i=0}C_{I-i,i}^2E^*\left(\left.\prod^{I+n-1}_{j=i}F^{*2}_{I-i,j}-2\left(\prod^{K}_{j=i}F^{*2}_{I-i,j}\right)\left(\prod^{I+n-1}_{l=K+1}F^*_{I-i,l}\right)+\prod^{K}_{j=i}F_{I-i,j}^{*2}\right|\mathcal{Q}_{I,n}^*=\mathcal{Q}_{I,n},\mathcal{F}_{I,n}^*\right).	\label{A2KI_second_moment_decomposition}
\end{align}
Using linearity of expectations, for the first expectation on the last right-hand side of \eqref{A2KI_second_moment_decomposition}, due to $F^*_{i,j}=\frac{C^*_{i,j+1}}{C^*_{i,j}}$ and $C_{I-i,i}^*=C_{I-i,i}$, we get
\begin{align*}
&E^*\left(\prod^{I+n-1}_{j=i}F_{I-i,j}^{*2}\Big|\mathcal{Q}_{I,n}^*=\mathcal{Q}_{I,n},\mathcal{F}_{I,n}^*\right)\\
&= E^*\left(\left(\prod^{I+n-2}_{j=i}F_{I-i,j}^{*2}\right)E^*\left(F_{I-i,I+n-1}^{*2}|\mathcal{Q}_{I,n}^*=\mathcal{Q}_{I,n},C^*_{I-i,i},\ldots,C^*_{I-i,I+n-1},\mathcal{F}_{I,n}^*\right)\Big|\mathcal{Q}_{I,n}^*=\mathcal{Q}_{I,n},\mathcal{F}_{I,n}^*\right)	\\
&= E^*\left(\prod^{I+n-2}_{j=i}F_{I-i,j}^{*2}\left(\frac{\widehat\sigma_{I+n-1,n}^2}{C^*_{I-i,I+n-1}}+\widehat f_{I+n-1,n}^{*2}\right)|\mathcal{Q}_{I,n}^*=\mathcal{Q}_{I,n},\mathcal{F}_{I,n}^*\right)	\\
&= E^*\left(\prod^{I+n-2}_{j=i}F_{I-i,j}^{*2}\frac{1}{C^*_{I-i,I+n-1}}|\mathcal{Q}_{I,n}^*=\mathcal{Q}_{I,n},\mathcal{F}_{I,n}^*\right)\widehat \sigma_{I+n-1,n}^2+E^*\left(\prod^{I+n-2}_{j=i}F_{I-i,j}^{*2}\Big|\mathcal{Q}_{I,n}^*=\mathcal{Q}_{I,n},\mathcal{F}_{I,n}^*\right)\widehat f_{I+n-1,n}^{*2}	\\
&= E^*\left(\prod^{I+n-2}_{j=i}F^*_{I-i,j}|\mathcal{Q}_{I,n}^*=\mathcal{Q}_{I,n},\mathcal{F}_{I,n}^*\right)\frac{\widehat \sigma_{I+n-1,n}^2}{C_{I-i,i}}+E^*\left(\prod^{I+n-2}_{j=i}F_{I-i,j}^{*2}\Big|\mathcal{Q}_{I,n}^*=\mathcal{Q}_{I,n},\mathcal{F}^*_{I,n}\right)\widehat f_{I+n-1,n}^{2*}	\\
&= \frac{\prod^{I+n-2}_{j=i}\widehat f_{j,n}^*\widehat \sigma_{I+n-1,n}^2}{C_{I-i,i}}+E^*\left(\prod^{I+n-2}_{j=i}F_{I-i,j}^{*2}\Big|\mathcal{Q}_{I,n}^*=\mathcal{Q}_{I,n}, \mathcal{F}_{I,n}^*\right)\widehat f_{I+n-1,n}^{*2}.
\end{align*}
By recursively plugging-in, we get
\begin{align*}
E^*\left(\prod^{I+n-1}_{j=i}F_{I-i,j}^{*2}\Big|\mathcal{Q}_{I,n}^*=\mathcal{Q}_{I,n},\mathcal{F}_{I,n}^*\right) = \frac{1}{C_{I-i,i}}\sum_{k=i}^{I+n-1}\left(\prod_{j=i}^{k-1}\widehat f_{j,n}^*\right)\widehat \sigma_{k,n}^2\left(\prod_{h=k+1}^{I+n-1}\widehat f_{h,n}^{*2}\right)+\prod_{j=i}^{I+n-1} \widehat f_{j,n}^{*2}.
\end{align*}
Similarly, for the second expectation in \eqref{A2KI_second_moment_decomposition}, we get
\begin{align*}
& E^*\left(-2\left(\prod^{K}_{j=i}F_{I-i,j}^{*2}\right)\left(\prod^{I+n-1}_{l=K+1}F^*_{I-i,l}\right)\Big|\mathcal{Q}_{I,n}^*=\mathcal{Q}_{I,n},\mathcal{F}^*_{I,n}\right)	\\
&= -2\left(\frac{1}{C_{I-i,i}}\sum_{k=i}^{K}\left(\prod_{j=i}^{k-1}\widehat f^*_{j,n}\right)\widehat \sigma_{k,n}^2\left(\prod_{h=k+1}^{K}\widehat f_{h,n}^{*2}\right)+\prod_{j=i}^{K} \widehat f_{j,n}^{*2}\right)\left(\prod^{I+n-1}_{l=K+1}\widehat f_{l,n}^*\right)
\end{align*}
and for the third one, we have
\begin{align*}
E^*\left(\prod^{K}_{j=i}F_{I-i,j}^{*2}\Big|\mathcal{Q}_{I,n}^*=\mathcal{Q}_{I,n},\mathcal{F}_{I,n}^*\right) = \frac{1}{C_{I-i,i}}\sum_{k=i}^{K}\left(\prod_{j=i}^{k-1}\widehat f_{j,n}^*\right)\widehat \sigma_{k,n}^2\left(\prod_{h=k+1}^{K}\widehat f_{h,n}^{*2}\right)+\prod_{j=i}^{K}\widehat  f_{j,n}^{*2}.
\end{align*}
Altogether, for all $K<I+n-1$, this leads to
\begin{align*}
& E^*((A_{2,K,I,n}^{*})^2|\mathcal{Q}_{I,n}^*=\mathcal{Q}_{I,n},\mathcal{F}_{I,n}^*)	\\
&= \sum^{K}_{i=0}C_{I-i,i}\left[\sum_{k=i}^{I+n-1}\left(\prod_{j=i}^{k-1}\widehat f^*_{j,n}\right)\widehat \sigma_{k,n}^2\left(\prod_{h=k+1}^{I+n-1}\widehat f_{h,n}^{*2}\right)-2\left(\sum_{k=i}^{K}\left(\prod_{j=i}^{k-1}\widehat f_{j,n}^*\right)\widehat \sigma_{k,n}^2\left(\prod_{h=k+1}^{K}\widehat f_{h,n}^{*2}\right)\right)\left(\prod^{I+n-1}_{l=K+1}\widehat f_{l,n}^*\right)\right.	\\
& \quad	+\left.\sum_{k=i}^{K}\left(\prod_{j=i}^{k-1}\widehat f_{j,n}^*\right)\widehat \sigma_{k,n}^2\left(\prod_{h=k+1}^{K}\widehat f_{h,n}^{*2}\right)\right].
\end{align*}
Plugging-in and making use of the fact that the $\widehat f_{j,n}^*$'s are stochastically independent conditional on $\mathcal{D}_{I,n}$ and $\mathcal{Q}_{I,n}^*=\mathcal{Q}_{I,n}$, this leads to
\begin{align*}
& Var^*(A_{2,K,I,n}^*|\mathcal{Q}_{I,n}^*=\mathcal{Q}_{I,n}) \\
&= \sum^{K}_{i=0}C_{I-i,i}\left[\sum_{k=i}^{I+n-1}\widehat \sigma_{k,n}^2\left(\prod_{j=i}^{k-1}E^*\left(\widehat f^*_{j,n}|\mathcal{Q}_{I,n}^*=\mathcal{Q}_{I,n}\right)\right)\left(\prod_{h=k+1}^{I+n-1}E^*\left(\widehat f_{h,n}^{*2}|\mathcal{Q}_{I,n}^*=\mathcal{Q}_{I,n}\right)\right)\right.	\\
& \quad-2\sum_{k=i}^{K}\widehat \sigma_{k,n}^2\left(\prod_{j=i}^{k-1}E^*\left(\widehat f_{j,n}^*|\mathcal{Q}_{I,n}^*=\mathcal{Q}_{I,n}\right)\right)\left(\prod_{h=k+1}^{K}E^*\left(\widehat f_{h,n}^{*2}|\mathcal{Q}_{I,n}^*=\mathcal{Q}_{I,n}\right)\right)\left(\prod^{I+n-1}_{l=K+1}E^*\left(\widehat f_{l,n}^*|\mathcal{Q}_{I,n}^*=\mathcal{Q}_{I,n}\right)\right)	\\
& \qquad\qquad\qquad	+\left.\sum_{k=i}^{K}\widehat \sigma_{k,n}^2\left(\prod_{j=i}^{k-1}E^*\left(\widehat f_{j,n}^*\right)\right)\left(\prod_{h=k+1}^{K}E^*\left(\widehat f_{h,n}^{*2}\right)\right)\right]	\\
&= \sum^{K}_{i=0}C_{I-i,i}\left[\sum_{k=i}^{I+n-1}\widehat \sigma_{k,n}^2\left(\prod_{j=i}^{k-1}\widehat{f}_{j,n}\right)\left(\prod_{h=k+1}^{I+n-1}\left(\widehat{f}_{h,n}^2+\frac{\widehat \sigma_{h,n}^2}{\sum^{I-h-1}_{p=-n}C_{p,h}}\right)\right)\right.	\\
& \qquad\qquad\qquad-2\sum_{k=i}^{K}\widehat \sigma_{k,n}^2\left(\prod_{j=i}^{k-1}\widehat{f}_{j,n}\right)\left(\prod_{h=k+1}^{K}\left(\widehat{f}_{h,n}^2+\frac{\widehat \sigma_{h,n}^2}{\sum^{I-h-1}_{p=-n}C_{p,h}}\right)\right)\left(\prod^{I+n-1}_{l=K+1}\widehat{f}_{l,n}\right)	\\
& \qquad\qquad\qquad	+\left.\sum_{k=i}^{K}\widehat \sigma_{k,n}^2\left(\prod_{j=i}^{k-1}\widehat{f}_{j,n}\right)\left(\prod_{h=k+1}^{K}\left(\widehat{f}_{h,n}^2+\frac{\widehat \sigma_{h,n}^2}{\sum^{I-h-1}_{p=-n}C_{p,h}}\right)\right)\right]	\\
&= \sum^{K}_{i=0}C_{I-i,i}\left[\sum_{k=i}^{K}\left(\prod_{j=i}^{k-1}\widehat{f}_{j,n}\right)\widehat \sigma_{k,n}^2\left(\prod_{h=k+1}^{K}\left(\widehat{f}_{h,n}^2+\frac{\widehat \sigma_{h,n}^2}{\sum^{I-h-1}_{p=-n}C_{p,h}}\right)\right)\right.	\\
& \qquad\qquad\qquad \times\left(\left(\prod_{h=K+1}^{I+n-1}\left(\widehat{f}_{h,n}^2+\frac{\widehat \sigma_{h,n}^2}{\sum^{I-h-1}_{p=-n}C_{p,h}}\right)\right)-\left(\prod^{I+n-1}_{l=K+1}\widehat{f}_{l,n}\right)\right)	\\
& \qquad\qquad\qquad+\sum_{k=K+1}^{I+n-1}\left(\prod_{j=i}^{k-1}\widehat{f}_{j,n}\right)\widehat \sigma_{k,n}^2\left(\prod_{h=k+1}^{I+n-1}\left(\widehat{f}_{h,n}^2+\frac{\widehat \sigma_{h,n}^2}{\sum^{I-h-1}_{p=-n}C_{p,h}}\right)\right)	\\
& \qquad\qquad\qquad	+\left.\sum_{k=i}^{K}\left(\prod_{j=i}^{k-1}\widehat{f}_{j,n}\right)\widehat \sigma_{k,n}^2\left(\prod_{h=k+1}^{K}\left(\widehat{f}_{h,n}^2+\frac{\widehat \sigma_{h,n}^2}{\sum^{I-h-1}_{p=-n}C_{p,h}}\right)\right)\left(1-\prod^{I+n-1}_{l=K+1}\widehat{f}_{l,n}\right)\right]
\end{align*}
obtained by re-arranging terms and due to $E^*(\widehat f^*_{j,n}|\mathcal{Q}_{I,n}^*=\mathcal{Q}_{I,n})=E^*(\widehat f^*_{j,n})=\widehat f_{j,n}$ and
\begin{align*}
E^*\left(\widehat f_{j,n}^{*2}|\mathcal{Q}_{I,n}^*=\mathcal{Q}_{I,n}\right)=E^*\left(\widehat f_{j,n}^{*2}\right)=\widehat{f}_{j,n}^2+\frac{\widehat \sigma_{j,n}^2}{\sum^{I-j-1}_{p=-n}C_{p,j}}
\end{align*}
for all $j$. Next, to argue that $Var^*(A_{2,K,I,n}^*|\mathcal{Q}_{I,n}^*=\mathcal{Q}_{I,n})\geq 0$ vanishes in probability for $K\rightarrow \infty$ and $n\rightarrow \infty$ afterwards, it suffices to show that its unconditional expectation is bounded for $K\rightarrow \infty$ and that its bound converges to zero as $n\rightarrow \infty$. We get
\begin{align*}
& E\left(Var^*(A_{2,K,I,n}^*|\mathcal{Q}_{I,n}^*=\mathcal{Q}_{I,n})\right)	\\
=& \sum^{K}_{i=0}E\left[C_{I-i,i}\sum_{k=i}^{K}\left(\prod_{j=i}^{k-1}\widehat{f}_{j,n}\right)\widehat \sigma_{k,n}^2\left(\prod_{h=k+1}^{K}\left(\widehat{f}_{h,n}^2+\frac{\widehat \sigma_{h,n}^2}{\sum^{I-h-1}_{p=-n}C_{p,h}}\right)\right)\right.	\\
& \qquad\qquad\qquad \left.\times\left(\prod_{h=K+1}^{I+n-1}\left(\widehat{f}_{h,n}^2+\frac{\widehat \sigma_{h,n}^2}{\sum^{I-h-1}_{p=-n}C_{p,h}}\right)-\prod^{I+n-1}_{l=K+1}\widehat{f}_{l,n}\right)\right]	\\
&+\sum^{K}_{i=0}E\left[C_{I-i,i}\sum_{k=K+1}^{I+n-1}\left(\prod_{j=i}^{k-1}\widehat{f}_{j,n}\right)\widehat \sigma_{k,n}^2\left(\prod_{h=k+1}^{I+n-1}\left(\widehat{f}_{h,n}^2+\frac{\widehat \sigma_{h,n}^2}{\sum^{I-h-1}_{p=-n}C_{p,h}}\right)\right)\right]	\\
&+\sum^{K}_{i=0}E\left[C_{I-i,i}\sum_{k=i}^{K}\left(\prod_{j=i}^{k-1}\widehat{f}_{j,n}\right)\widehat \sigma_{k,n}^2\left(\prod_{h=k+1}^{K}\left(\widehat{f}_{h,n}^2+\frac{\widehat \sigma_{h,n}^2}{\sum^{I-h-1}_{p=-n}C_{p,h}}\right)\right)\left(1-\prod^{I+n-1}_{l=K+1}\widehat{f}_{l,n}\right)\right]	\\
=& \sum^{K}_{i=0}E\left[C_{I-i,i}\sum_{k=i}^{K}\left(\prod_{j=i}^{k-1}\widehat{f}_{j,n}\right)\widehat \sigma_{k,n}^2\left(\prod_{h=k+1}^{K}\left(\widehat{f}_{h,n}^2+\frac{\widehat \sigma_{h,n}^2}{\sum^{I-h-1}_{p=-n}C_{p,h}}\right)\right)\right.	\\
& \qquad\qquad\qquad \left.\times\left(\prod_{h=K+1}^{I+n-1}\left(\widehat{f}_{h,n}^2+\frac{\widehat \sigma_{h,n}^2}{\sum^{I-h-1}_{p=-n}C_{p,h}}\right)-\prod^{I+n-1}_{l=K+1}\widehat{f}_{l,n}\right)\right]	\\
&+\sum^{K}_{i=0}E\left[C_{I-i,i}\sum_{k=K+1}^{I+n-1}\left(\prod_{j=i}^{k-1}\widehat{f}_{j,n}\right)\widehat \sigma_{k,n}^2\left(\prod_{h=k+1}^{I+n-1}\left(\widehat{f}_{h,n}^2+\frac{\widehat \sigma_{h,n}^2}{\sum^{I-h-1}_{p=-n}C_{p,h}}\right)\right)\right]	\\
&+\sum^{K}_{i=0}E\left[C_{I-i,i}\sum_{k=i}^{K}\left(\prod_{j=i}^{k-1}\widehat{f}_{j,n}\right)\widehat \sigma_{k,n}^2\left(\prod_{h=k+1}^{K}\left(\widehat{f}_{h,n}^2+\frac{\widehat \sigma_{h,n}^2}{\sum^{I-h-1}_{p=-n}C_{p,h}}\right)\right)\left(1-\prod^{I+n-1}_{l=K+1}\widehat{f}_{l,n}\right)\right].
\end{align*}
Now, let us consider the three terms on the last right-hand side separately. Using $\sum^{I-h-1}_{p=-n}C_{p,h}\geq(I+n-h)\epsilon^h\geq\epsilon^h$, the first one can be bounded by
\begin{align*}
& \sum^{K}_{i=0}E\left[C_{I-i,i}\sum_{k=i}^{K}\left(\prod_{j=i}^{k-1}\widehat{f}_{j,n}\right)\widehat \sigma_{k,n}^2\left(\prod_{h=k+1}^{K}\left(\widehat{f}_{h,n}^2+\frac{\widehat \sigma_{h,n}^2}{\epsilon^h}\right)\right)\left(\prod_{h=K+1}^{I+n-1}\left(\widehat{f}_{h,n}^2+\frac{\widehat \sigma_{h,n}^2}{\epsilon^h}\right)-\prod^{I+n-1}_{l=K+1}\widehat{f}_{l,n}\right)\right]	\\
=& \sum^{K}_{i=0}E\left[C_{I-i,i}\sum_{k=i}^{K}\left(\prod_{j=i}^{k-1}\widehat{f}_{j,n}\right)\widehat \sigma_{k,n}^2\left(\prod_{h=k+1}^{I+n-1}\left(\widehat{f}_{h,n}^2+\frac{\widehat \sigma_{h,n}^2}{\epsilon^h}\right)\right)\right]	\\
&- \sum^{K}_{i=0}E\left[C_{I-i,i}\sum_{k=i}^{K}\left(\prod_{j=i}^{k-1}\widehat{f}_{j,n}\right)\widehat \sigma_{k,n}^2\left(\prod_{h=k+1}^{K}\left(\widehat{f}_{h,n}^2+\frac{\widehat \sigma_{h,n}^2}{\epsilon^h}\right)\right)\left(\prod^{I+n-1}_{l=K+1}\widehat{f}_{l,n}\right)\right].	
\end{align*}
Next, using the law of iterated expectations and
\begin{align*}
E\left[\widehat{f}_{c,n}^2+\frac{\widehat \sigma_{c,n}^2}{\epsilon^{c}}|\mathcal{B}_{I,n}(c)\right]=\frac{\sigma_c^2}{\sum_{k=-n}^{I-c-1}C_{k,c}}+f_c^2+\frac{\sigma_{c}^2}{\epsilon^{c}}\leq f_c^2+\frac{2\sigma_{c}^2}{\epsilon^{c}}
\end{align*}
for all $c\in\{0,\ldots,I+n-1\}$, where $\mathcal{B}_{I,n}(k)=\{C_{i,j}|i=-n, \dots, I, j=0, \dots, k, i+j\leq I\}$, the first term on the right-hand side above becomes
\begin{align*}
\sum^{K}_{i=0}E\left[C_{I-i,i}\sum_{k=i}^{K}\left(\prod_{j=i}^{k-1}\widehat{f}_{j,n}\right)\widehat \sigma_{k,n}^2\left(\prod_{h=k+1}^K\left(\widehat{f}_{h,n}^2+\frac{\widehat \sigma_{h,n}^2}{\epsilon^h}\right)\right)\right]\left(\prod_{l=K+1}^{I+n-1}\left(f_l^2+\frac{2\sigma_l^2}{\epsilon^{l}}\right)\right)
\end{align*}
and, similarly, for the second term, we obtain
\begin{align*}
- \sum^{K}_{i=0}E\left[C_{I-i,i}\sum_{k=i}^{K}\left(\prod_{j=i}^{k-1}\widehat{f}_{j,n}\right)\widehat \sigma_{k,n}^2\left(\prod_{h=k+1}^{K}\left(\widehat{f}_{h,n}^2+\frac{\widehat \sigma_{h,n}^2}{\epsilon^h}\right)\right)\right]\left(\prod^{I+n-1}_{l=K+1}f_l\right).
\end{align*}
Together, this term becomes
\begin{align*}
\sum^{K}_{i=0}E\left[C_{I-i,i}\sum_{k=i}^{K}\left(\prod_{j=i}^{k-1}\widehat{f}_{j,n}\right)\widehat \sigma_{k,n}^2\left(\prod_{h=k+1}^K\left(\widehat{f}_{h,n}^2+\frac{\widehat \sigma_{h,n}^2}{\epsilon^h}\right)\right)\right]\left(\prod_{l=K+1}^{I+n-1}\left(f_l^2+\frac{2\sigma_l^2}{\epsilon^{l}}\right)-\prod^{I+n-1}_{l=K+1}f_l\right),
\end{align*}
which, using similar arguments as above, can be bounded by
\begin{align*}
\left[\sum^{K}_{i=0}\mu_i\sum_{k=i}^{K}\left(\prod_{j=i}^{k-1}f_j\right)\sigma_k^2\left(\prod_{h=k+1}^K\left(f_h^2+\frac{2\sigma_h^2}{\epsilon^h}\right)\right)\right]\left(\prod_{l=K+1}^{I+n-1}\left(f_l^2+\frac{2\sigma_l^2}{\epsilon^{l}}\right)-\prod^{I+n-1}_{l=K+1}f_l\right).
\end{align*}
Now, letting $n\rightarrow\infty$, we get the following upper bound
\begin{align*}
\mu_\infty^2\left(\prod_{h=0}^\infty\left(f_h^2+\frac{2\sigma_h^2}{\epsilon^h}\right)\right)\sum^{K}_{i=0}\sum_{k=i}^{K}\sigma_k^2\left(\prod_{l=K+1}^\infty\left(f_l^2+\frac{2\sigma_l^2}{\epsilon^{l}}\right)-\prod^\infty_{l=K+1}f_l\right)<\infty,
\end{align*}
which is finite using $\prod^{\infty}_{j=0} x_j<\infty$ if and only if $\sum^{\infty}_{j=0} (x_j-1)<\infty$ for $x_j\geq 1$ for all $j$, and as we have
\begin{align*}
\sum^{\infty}_{h=0} \left(f_{h}^2+\frac{2\sigma_{h}^2}{\epsilon^{h}}-1\right) =& \sum^{\infty}_{h=0} (f_{h}^2-1)+\sum^{\infty}_{h=0}\frac{2\sigma_{h}^2}{\epsilon^{h}} = \sum^{\infty}_{h=0} (f_{h}-1)(f_{h}+1)+2\sum^{\infty}_{h=0}\frac{\sigma_{h}^2}{\epsilon^{h}}	\\
\leq& \sup_{h\in\mathbb{N}}(f_{h}+1)\sum^{\infty}_{h=0} (f_{h}-1)+2\sum^{\infty}_{h=0}\frac{\sigma_{h}^2}{\epsilon^h}<\infty
\end{align*}
by Assumptions 2.4 and 4.2. Now, letting also $K\rightarrow \infty$, the term $\sum^{K}_{i=0}\sum_{k=i}^{K}\sigma_k^2$ also remains bounded due to
\begin{align*}
\sum^{K}_{i=0}\sum_{k=i}^{K}\sigma_k^2 = \sum^{K}_{j=0}\sigma_{j}^2\sum_{i=0}^{j} = \sum^{K}_{j=0}(j+1)\sigma_{j}^2 \leq \sum^\infty_{j=0}(j+1)\sigma_j^2<\infty.
\end{align*}
Finally, as $f_l\rightarrow 1$ and $\sigma_l^2/\epsilon^l\rightarrow 0$ for $l\rightarrow \infty$, we get $\prod_{l=K+1}^\infty(f_l^2+\frac{2\sigma_l^2}{\epsilon^{l}})\rightarrow 1$ and $\prod^\infty_{l=K+1}f_l\rightarrow 1$ for $K\rightarrow \infty$ leading to
\begin{align*}
\prod_{l=K+1}^\infty\left(f_l^2+\frac{2\sigma_l^2}{\epsilon^{l}}\right)-\prod^\infty_{l=K+1}f_l\rightarrow 0.
\end{align*}
Similarly, using the same arguments, the second term in the representation of $E(Var^*(A_{2,K,I,n}^*|\mathcal{Q}_{I,n}^*=\mathcal{Q}_{I,n}))$ above can be bounded by
\begin{align*}
\sum^{K}_{i=0}\mu_i\sum_{k=K+1}^{I+n-1}\left(\prod_{j=i}^{k-1}f_j\right)\sigma_k^2\left(\prod_{h=k+1}^{I+n-1}\left(f_h^2+\frac{2\sigma_h^2}{\epsilon^h}\right)\right),
\end{align*}
which, for $n\rightarrow \infty$, can be bounded by
\begin{align*}
\mu_\infty^2\left(\prod_{h=0}^\infty\left(f_h^2+\frac{2\sigma_h^2}{\epsilon^h}\right)\right)\sum^{K}_{i=0}\sum_{k=K+1}^\infty \sigma_k^2 \leq \mu_\infty^2\left(\prod_{h=0}^\infty\left(f_h^2+\frac{2\sigma_h^2}{\epsilon^h}\right)\right)\sum_{k=K+1}^\infty (k+1)\sigma_k^2<\infty,
\end{align*}
which vanishes for $K\rightarrow \infty$.

Finally, for the third term in the representation of $E(Var^*(A_{2,K,I,n}^*|\mathcal{Q}_{I,n}^*=\mathcal{Q}_{I,n}))$, we get
\begin{align*}
& \sum^{K}_{i=0}E\left[C_{I-i,i}\sum_{k=i}^{K}\left(\prod_{j=i}^{k-1}\widehat{f}_{j,n}\right)\widehat \sigma_{k,n}^2\left(\prod_{h=k+1}^{K}\left(\widehat{f}_{h,n}^2+\frac{\widehat \sigma_{h,n}^2}{\sum^{I-h-1}_{p=-n}C_{p,h}}\right)\right)\left(1-\prod^{I+n-1}_{l=K+1}\widehat{f}_{l,n}\right)\right]	\\
=&  \left(\sum^{K}_{i=0}E\left[C_{I-i,i}\sum_{k=i}^{K}\left(\prod_{j=i}^{k-1}\widehat{f}_{j,n}\right)\widehat \sigma_{k,n}^2\left(\prod_{h=k+1}^{K}\left(\widehat{f}_{h,n}^2+\frac{\widehat \sigma_{h,n}^2}{\sum^{I-h-1}_{p=-n}C_{p,h}}\right)\right)\right]\right)\left(1-\prod^{I+n-1}_{l=K+1}f_l\right)
\end{align*}
While, for $n\rightarrow \infty$, the second factor $1-\prod^{I+n-1}_{l=K+1}f_l$ can be bounded by $1-\prod^\infty_{l=K+1}f_l$, which converges to $0$ due to $\prod^\infty_{l=K+1}f_l\rightarrow 1$ for $K\rightarrow \infty$, the first factor above can be bounded by 
\begin{align*}
\sum^{K}_{i=0}\mu_i\sum_{k=i}^{K}\left(\prod_{j=i}^{k-1}f_j\right)\sigma_k^2\left(\prod_{h=k+1}^{K}\left(f_h^2+\frac{\sigma_h^2}{\epsilon^h}\right)\right),
\end{align*}
which, for $n\rightarrow \infty$, can be bounded further by
\begin{align*}
\mu_\infty^2\left(\prod_{h=0}^\infty\left(f_h^2+\frac{\widehat \sigma_h^2}{\epsilon^h}\right)\right)\sum^{K}_{j=0}(j+1)\sigma_{j}^2,
\end{align*}
which is bounded as $\sum^{K}_{j=0}(j+1)\sigma_{j}^2\rightarrow \sum^\infty_{j=0}(j+1)\sigma_{j}^2<\infty$ for $K\rightarrow \infty$. This completes the first part of c) for term $A_{2,K,I,n}^*$. Continuing with $A_{3,K,I,n}^*$ to prove also the second part of c), we have $E^*(A^*_{3,K,I,n}|\mathcal{Q}_{I,n}^*=\mathcal{Q}_{I,n})=0$ using the law of iterated expectations. By using similar calculations as for $A_{2,K,I,n}^*$, we get
\begin{align*}
E^*((A^*_{3,K,I,n})^2|\mathcal{Q}_{I,n}^*=\mathcal{Q}_{I,n},\mathcal{F}_{I,n}^*) &= \sum^{I+n}_{i=K+1}C_{I-i,i}\sum^{I+n-1}_{j=i}\left(\prod^{j-1}_{h=i}\widehat f_{h,n}^*\right)\widehat \sigma_{j,n}^2\left(\prod^{I+n-1}_{l=j+1}\widehat f_{l,n}^{*2}\right)
\end{align*}
and
\begin{align*}
Var^*(A^*_{3,K,I,n}|\mathcal{Q}_{I,n}^*=\mathcal{Q}_{I,n}) &= E^*\left(E^*((A^*_{3,K,I,n})^2|\mathcal{Q}_{I,n}^*=\mathcal{Q}_{I,n},\mathcal{F}_{I,n}^*)|\mathcal{Q}_{I,n}^*=\mathcal{Q}_{I,n}\right)	\\
&= \sum^{I+n}_{i=K+1}C_{I-i,i}\sum^{I+n-1}_{j=i}\left(\prod^{j-1}_{h=i}\widehat f_{h,n}\right)\widehat \sigma_{j,n}^2\left(\prod^{I+n-1}_{l=j+1}\left(\widehat f_{l,n}^{2}+\frac{\widehat \sigma_{l,n}^2}{\sum^{I-j+1}_{k=-n}C_{k,l}}\right)\right).
\end{align*}
Hence, to show that $Var^*(A_{3,K,I,n}^*)\geq 0$ vanishes in probability, we prove that $E(Var^*(A_{3,K,I,n}^*|\mathcal{Q}_{I,n}^*=\mathcal{Q}_{I,n}))$ is bounded for $n\rightarrow \infty$ and that its bound converges to zero as $K\rightarrow \infty$. By plugging-in and using similar arguments as above for $A_{2,K,I,n}^*$, we get
\begin{align*}
E(Var^*(A_{3,K,I,n}^*|\mathcal{Q}_{I,n}^*=\mathcal{Q}_{I,n})) &\leq \sum^{I+n}_{i=K+1}\mu_{i}\sum^{I+n-1}_{j=i}\left(\prod^{j-1}_{h=i} f_h\right) \sigma_j^2\left(\prod^{I+n-1}_{l=j+1}\left( f_l^{2}+\frac{2 \sigma_l^2}{\epsilon^l}\right)\right)\\
&\leq \mu_\infty^2\left(\prod^\infty_{l=0}\left( f_l^{2}+\frac{2 \sigma_l^2}{\epsilon^l}\right)\right)\sum^{I+n}_{i=K+1}\sum^{I+n-1}_{j=i} \sigma_j^2,
\end{align*}
which is bounded for $n\rightarrow \infty$ due to
\begin{align*}
\sum^{I+n}_{i=K+1}\sum^{I+n-1}_{j=i} \sigma_j^2 \leq \sum^{I+n-1}_{j=K+1}(j+1)\sigma_j^2\rightarrow \sum^\infty_{j=K+1}(j+1)\sigma_j^2<\infty
\end{align*}
and this bound vanishes for $K\rightarrow \infty$. $\hfill\Box$

\bigskip

\subsection{Proof of Theorem 4.10}

Similar to the proof of Theorem 4.8 for the conditional limiting behavior of $(R_{I,n}^*-\widehat{R}_{I,n})_1$ and to the proof of Theorem 4.7 in \cite{jentschasympThInf} for the (unconditional!) limiting behavior of $(R_{I,n}-\widehat{R}_{I,n})_2$, we decompose $(R_{I,n}^*-\widehat{R}_{I,n})_2$ by truncating the sums and products to be able to apply Proposition 6.3.9 in \cite{brockwell1991time}. For this purpose, let $K\in \mathbb{N}_0$ be fixed and suppose $I,n\in\mathbb{N}_0$ are large enough such that $K<I+n-1$. Then, after inflating $(R_{I,n}^*-\widehat{R}_{I,n})_2$ with $\sqrt{I+n+1}$, we get
\begin{align*}
& \sqrt{I+n+1} \left(R_{I,n}^*-\widehat{R}_{I,n}\right)_2	\\
=& \sqrt{I+n+1}\sum^{I+n}_{i=0}C_{I-i,i}^*\left(\prod_{j=i}^{I+n-1}\widehat{f}_{j,n}^*-\prod_{j=i}^{I+n-1}\widehat{f}_{j,n}\right)\\
=&\sqrt{I+n+1}\sum_{i=0}^K C_{I-i,i}^*\left(\prod_{j=i}^{K}\widehat{f}_{j,n}^*-\prod_{j=i}^{K}\widehat{f}_{j,n}\right)\\
&+\sqrt{I+n+1}\sum_{i=0}^K C_{I-i,i}^*\left(\prod^{K}_{j=i}\widehat{f}_{j,n}^*\left(\prod^{I+n-1}_{l=K+1}\widehat{f}_{l,n}^*-1\right)-\prod^{K}_{j=i}\widehat{f}_{j,n}\left(\prod^{I+n-1}_{l=K+1}\widehat f_{l,n}-1\right)\right)\\
&+ \sqrt{I+n+1}\sum^{I+n}_{i=K+1}C_{I-i,i}^*\left(\prod^{I+n-1}_{j=i}\widehat{f}_{j,n}^*-\prod^{I+n-1}_{j=i}\widehat{f}_{j,n}\right)\\
=:& B_{1,K,I,n}^*+B_{2,K,I,n}^*+B_{3,K,I,n}^*,
\end{align*}
Hence, to derive the claimed conditional limiting distribution, it suffices to show that, a) for all $K\in \mathbb{N}_0$, $B_{1,K,I,n}^*|(\mathcal{Q}_{I,n}^*=\mathcal{Q}_{I,n}, \mathcal{D}_{I,n})\overset{d}{\rightarrow} \mathcal{\widetilde G}_{2,K}|\mathcal{Q}_{I,\infty}$
%\footnote{CJ: Hier jetzt mit Platzhalter! Besser Notation ohne Tilde moeglich? Ist das einfach $\mathcal{G}_{2,K}|\mathcal{Q}_{I,\infty}=\mathbf{q}_{I,\infty}$ bzw. $\mathcal{G}_2|\mathcal{Q}_{I,\infty}=\mathbf{q}_{I,\infty}$ aus dem ersten Paper!? $\mathbf{q}_{I,n}$ vs. $\mathbf{q}_{I,\infty}$} 
in probability as $n\rightarrow\infty$ for some (conditional) distribution $\mathcal{\widetilde G}_{2,K}|\mathcal{Q}_{I,\infty}$, b) $\mathcal{\widetilde G}_{2,K}|\mathcal{Q}_{I,\infty}\overset{d}{\rightarrow}\mathcal{G}_{2}|\mathcal{Q}_{I,\infty}$ as $K\rightarrow \infty$, and c) that, for all $\epsilon>0$, we have
\begin{align*}
\lim\limits_{K\rightarrow \infty}\limsup\limits_{n\rightarrow\infty}P^*\left(|B_{2,K,I,n}^*|>\epsilon|\mathcal{Q}_{I,n}^*=\mathcal{Q}_{I,n}\right)=0	\quad	\text{and}	\quad	\lim\limits_{K\rightarrow \infty}\limsup\limits_{n\rightarrow\infty}P^*\left(|A_{3,K,I,n}^*|>\epsilon|\mathcal{Q}_{I,n}^*=\mathcal{Q}_{I,n}\right)=0.
\end{align*}
We begin with part a). That is, for each fixed $K\in\mathbb{N}_0$, we consider 
\begin{align}\label{B1KI*}
B_{1,K,I,n}^*=\sqrt{I+n+1}\sum_{i=0}^KC_{I-i,i}^*\left(\prod_{j=i}^{K}\widehat{f}_{j,n}^*-\prod_{j=i}^{K}\widehat{f}_{j,n}\right),
\end{align}
where 
\begin{align*}%\label{eqf2}
\widehat{f}_{j,n}=\frac{\sum_{k=-n}^{I-j-1}C_{k,j+1}}{\sum_{k=-n}^{I-j-1}C_{k,j}}	\quad	\text{and}	\quad	\widehat{f}_{j,n}^*=\frac{\sum_{k=-n}^{I-j-1}C_{k,j}F_{k,j}^*}{\sum_{k=-n}^{I-j-1}C_{k,j}}.
\end{align*}
In contrast to the situation in the proof of Theorem 4.7 in \cite{jentschasympThInf}, where all $\widehat{f}_{j,n}$'s are indeed affected by conditioning on $\mathcal{Q}_{I,n}$, here, conditional on $\mathcal{D}_{I,n}$, all $\widehat{f}_{j,n}^*$'s are \emph{independent} of the condition $\mathcal{Q}_{I,n}^*=\mathcal{Q}_{I,n}$. Hence, for $n\rightarrow \infty$, the (unconditional!) asymptotic bootstrap theory derived in Theorem \ref{CLT_f_boot} and Corollary \ref{CLTProdf} lead to
\begin{align*}
B_{1,K,I,n}^*|(\mathcal{Q}_{I,n}^*=\mathcal{Q}_{I,n},\mathcal{D}_{I,n}) \overset{d}{\longrightarrow}\left\langle \mathcal{Q}_{I,K-I},\mathbf{Y}_K \right\rangle|\mathcal{Q}_{I,\infty}
\end{align*}
in probability, where $\mathcal{Q}_{I,K-I}=\{C_{I-i,i}|i=0,\ldots,I+(K-I)=K\}$, $\left\langle \cdot,\cdot \right\rangle$ denotes the Euclidean inner product in $\mathbb{R}^{K+1}$, and $\mathbf{Y}_K=(Y_i,i=0,\ldots,K)$ is a $(K+1)$-dimensional multivariate normally distributed random variable with $\mathbf{Y}_K\sim\mathcal{N}\left(0,\mathbf{\Sigma}_{K, \prod f_j}\right)$ with $\mathbf{\Sigma}_{K, \prod f_j}$ defined in Corollary \ref{CLTProdf}.

Further, letting $K\rightarrow \infty$, we get $\left\langle \mathcal{Q}_{I,K-I},\mathbf{Y}_K \right\rangle|\mathcal{Q}_{I,\infty}\overset{d}{\rightarrow} \left\langle \mathcal{Q}_{I,\infty},\mathbf{Y}_\infty \right\rangle|\mathcal{Q}_{I,\infty}$, where $\mathcal{Q}_{I,\infty}=\{C_{I-i,i}|i\in\mathbb{N}_0\}$, and $\mathbf{Y}_\infty=(Y_i,i\in\mathbb{N}_0)$ denotes a centered Gaussian process with covariance
\begin{align}
Cov(Y_{i_1},Y_{i_2}) = \lim_{K\rightarrow \infty} \mathbf{\Sigma}_{K, \prod f_j}(i_1,i_2) = \sum^\infty_{j=\max(i_1,i_2)}\frac{\sigma^2_j}{\mu_j}\prod^{\infty}_{l=\max(i_1,i_2), l \neq j}f^2_l\prod^{\max(i_1,i_2)-1}_{m=\min(i_1,i_2)}f_m
\end{align}
for $i_1,i_2\in\mathbb{N}_0$. Moreover, as $\mathcal{Q}_{I,\infty}$ and $\mathbf{Y}_\infty$ are stochastically independent, conditional on $\mathcal{Q}_{I,\infty}$, the variance of $\left\langle \mathcal{Q}_{I,\infty},\mathbf{Y}_\infty \right\rangle$ computes to
\begin{align*}
&Var\left(\left\langle \mathcal{Q}_{I,\infty},\mathbf{Y}_\infty \right\rangle|\mathcal{Q}_{I,\infty}\right) =
\sum^\infty_{i=0}Var(C_{I-i,i}Y_{i}|\mathcal{Q}_{I,\infty})+\underset{i_1\neq i_2}{\sum^\infty_{i_1,i_2=0}}Cov(C_{I-i_1,i_1}Y_{i_1},C_{I-i_2, i_2}Y_{i_2}|\mathcal{Q}_{I,\infty})\\
&= \sum^\infty_{i=0}C_{I-i, i}^2Var(Y_{i})+\underset{i_1\neq i_2}{\sum^\infty_{i_1,i_2=0}}C_{I-i_1,i_1}C_{I-i_2,i_2}Cov(Y_{i_1},Y_{i_2})\\
&= \sum^\infty_{i=0}C_{I-i, i}^2
	\sum^\infty_{j=i}\frac{\sigma^2_j}{\mu_j}\underset{l \neq j}{\prod^{\infty}_{l=i}}f^2_l	 +\underset{i_1\neq i_2}{\sum^\infty_{i_1,i_2=0}}C_{I-i_1,i_1}C_{I-i_2,i_2}\sum^\infty_{j=\max(i_1,i_2)}\frac{\sigma^2_j}{\mu_j}\underset{l \neq j}{\prod^{\infty}_{l=\max(i_1,i_2)}}f^2_l\prod^{\max(i_1,i_2)-1}_{m=\min(i_1,i_2)}f_m\\
&	=\sum_{j=0}^{\infty}\frac{\s}{\mu_jf_j^2}\left(\sum^j_{i=0}C_{I-i, i}^2
\prod^\infty_{k=i}f_k^2	\right)
	+2\sum_{j=1}^{\infty}\frac{\s}{\mu_j}\left(\sum^j_{i_1=1}C_{I-i_1,i_1}\sum^{i_1-1}_{i_2=0}C_{I-i_2,i_2}\prod^{i_1-1}_{l=i_2}f_l\right)\prod^\infty_{k=j}f_k^2
	\\
&	=\sum_{j=0}^{\infty}\frac{\s}{\mu_jf_j^2}\left(\frac{j+1}{j+1}\sum^j_{i=0}C_{I-i, i}^2
\prod^\infty_{k=i}f_k^2	\right)
	+2\sum_{j=1}^{\infty}\frac{\s}{\mu_j}\left(\frac{j}{j}\sum^j_{i_1=1}C_{I-i_1,i_1}\frac{i_1}{i_1}\sum^{i_1-1}_{i_2=0}C_{I-i_2,i_2}\prod^{i_1-1}_{l=i_2}f_l\right)\prod^\infty_{k=j}f_k^2
	\\
&	\leq\mu_\infty^2\sum_{j=0}^{\infty}\frac{\s}{\mu_jf_j^2}(j+1)\left(\frac{1}{j+1}\sum^j_{i=0}C_{I-i, i}^2\right)
	+2\mu_\infty^2\sum_{j=1}^{\infty}\frac{\s}{\mu_j}(j+1)^2\left(\frac{1}{j}\sum^j_{i_1=1}C_{I-i_1,i_1}\frac{1}{i_1}\sum^{i_1-1}_{i_2=0}C_{I-i_2,i_2}\right)
	\\
& = O_P(1)
\end{align*}
due to $\sum_{j=0}^\infty (j+1)^2\sigma_j^2<\infty$ by Assumption 2.4.

We continue with showing part c) for $B_{2,K,I,n}^*$. Using similar arguments as above, we have to consider 
\begin{align*}
& B_{2,K,I,n}^*|(\mathcal{Q}_{I,n}^*=\mathcal{Q}_{I,n},\mathcal{D}_{I,n})	\\
=& \sqrt{I+n+1}\sum_{i=0}^K C_{I-i,i}^*\left(\prod^{K}_{j=i}\widehat{f}_{j,n}^*\left(\prod^{I+n-1}_{l=K+1}\widehat{f}_{l,n}^*-1\right)-\prod^{K}_{j=i}\widehat{f}_{j,n}\left(\prod^{I+n-1}_{l=K+1}\widehat f_{l,n}-1\right)\right)|(\mathcal{Q}_{I,n}^*=\mathcal{Q}_{I,n},\mathcal{D}_{I,n}).
\end{align*}
Using the unbiasedness of $\widehat{f}_{j,n}^*$ conditional on $\mathcal{Q}_{I,n}^*=\mathcal{Q}_{I,n},\mathcal{D}_{I,n}$ for $\widehat{f}_{j,n}$, that is,
\begin{align*}
E^*(\widehat{f}_{j,n}^*|\mathcal{Q}_{I,n}^*=\mathcal{Q}_{I,n})=E^*(\widehat{f}_{j,n}^*)=\widehat{f}_{j,n}
\end{align*}
for all $j$, and the independence of the $\widehat{f}_{j,n}^*$'s conditional on $\mathcal{Q}_{I,n}^*=\mathcal{Q}_{I,n}$ and $\mathcal{D}_{I,n}$, we have $E^*(B_{2,K,I,n}^*|\mathcal{Q}_{I,n}^*=\mathcal{Q}_{I,n})=0$ by construction. Hence, it remains to show that $Var^*(B_{2,K,I,n}^*|\mathcal{Q}_{I,n}^*=\mathcal{Q}_{I,n})$ is bounded in probability for $n\rightarrow \infty$ and its bound vanishes for $K\rightarrow \infty$ afterwards. 
Now, to compute the bootstrap variance $Var^*(B_{2,K,I,n}^*|\mathcal{Q}_{I,n}^*=\mathcal{Q}_{I,n})$, for any fixed $K\in\mathbb{N}_0$ and $n\in\mathbb{N}_0$ large enough such that $K<I+n-1$, we get 
\begin{align*}%\label{B2KI_var}
& Var^*(B_{2,K,I,n}^*|\mathcal{Q}_{I,n}^*=\mathcal{Q}_{I,n})	\\
&= (I+n+1) \sum_{i_1,i_2=0}^K C_{I-i_1,i_1}C_{I-i_2,i_2}Cov^*\left(\prod^{I+n-1}_{j_1=i_1}\widehat f_{j_1,n}^*-\prod^{K}_{j_1=i_1}\widehat f_{j_1,n}^*,\prod^{I+n-1}_{j_2=i_2}\widehat f_{j_2,n}^*-\prod^{K}_{j_2=i_2}\widehat f_{j_2,n}^*|\mathcal{Q}_{I,n}^*=\mathcal{Q}_{I,n}\right)	\\
&= (I+n+1) \sum_{i_1,i_2=0}^K C_{I-i_1,i_1}C_{I-i_2,i_2}Cov^*\left(\prod^{I+n-1}_{j_1=i_1}\widehat f_{j_1,n}^*-\prod^{K}_{j_1=i_1}\widehat f_{j_1,n}^*,\prod^{I+n-1}_{j_2=i_2}\widehat f_{j_2,n}^*-\prod^{K}_{j_2=i_2}\widehat f_{j_2,n}^*\right)	\\
&\leq 2(I+n+1) \sum_{i_2=0}^K\sum_{i_1=0}^{i_2}C_{I-i_1,i_1}C_{I-i_2,i_2} Cov^*\left(\prod^{I+n-1}_{j_1=i_1}\widehat f_{j_1,n}^*-\prod^{K}_{j_1=i_1}\widehat f_{j_1,n}^*,\prod^{I+n-1}_{j_2=i_2}\widehat f_{j_2,n}^*-\prod^{K}_{j_2=i_2}\widehat f_{j_2,n}^*\right),	%\label{B2KI_var_2}
\end{align*}
using that, conditional on $\mathcal{D}_{I,n}$, the $\widehat f_{j,n}^*$'s are independent of the condition $\mathcal{Q}_{I,n}^*=\mathcal{Q}_{I,n}$.

To calculate the covariance on the last right-hand side, for $i_1\leq i_2$, first, we consider the mixed moment
\begin{align*}
& E^*\left(\left(\prod^{I+n-1}_{j_1=i_1}\widehat f_{j_1,n}^*-\prod^{K}_{j_1=i_1}\widehat f_{j_1,n}^*\right)\left(\prod^{I+n-1}_{j_2=i_2}\widehat f_{j_2,n}^*-\prod^{K}_{j_2=i_2}\widehat f_{j_2,n}^*\right)\right)	\\
&= E^*\left(\left(\prod^{i_2-1}_{j_1=i_1}\widehat f_{j_1,n}^*\right)\left(\prod^{I+n-1}_{j_2=i_2}\widehat f_{j_2,n}^{*2}\right)\right)-2E^*\left(\left(\prod^{i_2-1}_{j_1=i_1}\widehat f_{j_1,n}^*\right)\left(\prod^{K}_{j_2=i_2}\widehat f_{j_2,n}^{*2}\right)
\left(\prod^{I+n-1}_{j_3=K+1}\widehat f_{j_3,n}^*\right)\right)	\\
& \quad+E^*\left(\left(\prod^{i_2-1}_{j_1=i_1}\widehat f_{j_1,n}^*\right)\left(\prod^{K}_{j_2=i_2}\widehat f_{j_2,n}^{*2}\right)\right)	\\
&= \left(\prod^{i_2-1}_{j_1=i_1}\widehat f_{j_1,n}\right)\left(\prod^{I+n-1}_{j_2=i_2}\left(\widehat f_{j_2,n}^{2}+\frac{\widehat\sigma_{j_2,n}^2}{\sum^{I-{j_2}-1}_{k=-n}C_{k,j_2}}\right)\right) \\
&-2\left(\prod^{i_2-1}_{j_1=i_1}\widehat f_{j_1,n}\right)\left(\prod^{K}_{j_2=i_2}\left(\widehat f_{j_2,n}^{2}+\frac{\widehat\sigma_{j_2,n}^2}{\sum^{I-{j_2}-1}_{k=-n}C_{k,j_2}}\right)\right)
\left(\prod^{I+n-1}_{j_3=K+1}\widehat f_{j_3,n}\right)\\
 &\quad +\left(\prod^{i_2-1}_{j_1=i_1}\widehat f_{j_1,n}\right)\left(\prod^{K}_{j_2=i_2}\left(\widehat f_{j_2,n}^{2}+\frac{\widehat\sigma_{j_2,n}^2}{\sum^{I-{j_2}-1}_{k=-n}C_{k,j_2}}\right)\right),
\end{align*}
since $\widehat{f}_{j,n}^*$ and $\widehat{f}_{k,n}^*$ are independent for $j\neq k$ and $j,k\in\{0,\dots, I+n-1\}$ conditional on $\mathcal{D}_{I,n}$. Similarly, we have
\begin{align*}
E^*\left(\prod^{I+n-1}_{j_1=i_1}\widehat f_{j_1,n}^*-\prod^{K}_{j_1=i_1}\widehat f_{j_1,n}^*\right)=\prod^{I+n-1}_{j_1=i_1}\widehat f_{j_1}-\prod^{K}_{j_1=i_1}\widehat f_{j_1}
\end{align*}
leading to
\begin{align*}
& Cov^*\left(\prod^{I+n-1}_{j_1=i_1}\widehat f_{j_1,n}^*-\prod^{K}_{j_1=i_1}\widehat f_{j_1,n}^*,\prod^{I+n-1}_{j_2=i_2}\widehat f_{j_2,n}^*-\prod^{K}_{j_2=i_2}\widehat f_{j_2,n}^*\right)	\\
=& \left(\prod^{i_2-1}_{j_1=i_1}\widehat f_{j_1,n}\right)\left(\prod^{I+n-1}_{j_2=i_2}\left(\widehat f_{j_2,n}^{2}+\frac{\widehat\sigma_{j_2,n}^2}{\sum^{I-{j_2}-1}_{k=-n}C_{k,j_2}}\right)-\prod^{I+n-1}_{j_2=i_2}\widehat f_{j_2,n}^{2}\right) \\
&-2\left(\prod^{i_2-1}_{j_1=i_1}\widehat f_{j_1,n}\right)\left(\prod^{K}_{j_2=i_2}\left(\widehat f_{j_2,n}^{2}+\frac{\widehat\sigma_{j_2,n}^2}{\sum^{I-{j_2}-1}_{k=-n}C_{k,j_2}}\right)-\prod^{K}_{j_2=i_2}\widehat f_{j_2,n}^{2}\right)
\left(\prod^{I+n-1}_{j_3=K+1}\widehat f_{j_3,n}\right)\\
& +\left(\prod^{i_2-1}_{j_1=i_1}\widehat f_{j_1,n}\right)\left(\prod^{K}_{j_2=i_2}\left(\widehat f_{j_2,n}^{2}+\frac{\widehat\sigma_{j_2,n}^2}{\sum^{I-{j_2}-1}_{k=-n}C_{k,j_2}}\right)-\prod^{K}_{j_2=i_2}\widehat f_{j_2,n}^{2}\right)	\\
=& \left(\prod^{i_2-1}_{j_1=i_1}\widehat f_{j_1,n}\right)\left[\sum^{I+n-1}_{j_3=i_2} \frac{\widehat\sigma_{j_3,n}^2}{\sum^{I-{j_3}-1}_{k=-n}C_{k,j_3}}\left(\prod^{j_3-1}_{j_4=i_2}\widehat{f}_{j_4,n}^2 \right)\left(\prod^{I+n-1}_{j_2=j_3+1}\left(\widehat f_{j_2,n}^{2}+\frac{\widehat\sigma_{j_2,n}^2}{\sum^{I-{j_2}-1}_{k=-n}C_{k,j_2}}\right)\right)\right.\\
&-2\left.\sum^{K}_{j_3=i_2} \frac{\widehat\sigma_{j_3,n}^2}{\sum^{I-{j_3}-1}_{k=-n}C_{k,j_3}}\left(\prod^{j_3-1}_{j_4=i_2}\widehat{f}_{j_4,n}^2 \right)\left(\prod^{K}_{j_2= j_3+1}\left(\widehat f_{j_2,n}^{2}+\frac{\widehat\sigma_{j_2,n}^2}{\sum^{I-{j_2}-1}_{k=-n}C_{k,j_2}}\right)\right)\left(\prod^{I+n-1}_{j_4=K+1}\widehat{f}_{j_4,n}\right)\right. \\
&+\left.\sum^{K}_{j_3=i_2} \frac{\widehat\sigma_{j_3,n}^2}{\sum^{I-{j_3}-1}_{k=-n}C_{k,j_3}}\left(\prod^{j_3-1}_{j_4=i_2}\widehat{f}_{j_4,n}^2 \right)\left(\prod^{K}_{j_2= j_3+1}\left(\widehat f_{j_2,n}^{2}+\frac{\widehat\sigma_{j_2,n}^2}{\sum^{I-{j_2}-1}_{k=-n}C_{k,j_2}}\right)\right)\right].
\end{align*}
By rearranging the terms in brackets on the last right-hand side above, it becomes%\footnote{CJ: Checken! Stimmt das!?}
\begin{align*}
& \sum^{I+n-1}_{j_3=K+1} \frac{\widehat\sigma_{j_3,n}^2}{\sum^{I-{j_3}-1}_{k=-n}C_{k,j_3}}\left(\prod^{j_3-1}_{j_4=i_2}\widehat{f}_{j_4,n}^2 \right)\left(\prod^{I+n-1}_{j_2=j_3+1}\left(\widehat f_{j_2,n}^{2}+\frac{\widehat\sigma_{j_2,n}^2}{\sum^{I-{j_2}-1}_{k=-n}C_{k,j_2}}\right)\right)\\
&+\sum^{K}_{j_3=i_2} \frac{\widehat\sigma_{j_3,n}^2}{\sum^{I-{j_3}-1}_{k=-n}C_{k,j_3}}\left(\prod^{j_3-1}_{j_4=i_2}\widehat{f}_{j_4,n}^2 \right)\left(\prod^{I+n-1}_{j_2=j_3+1}\left(\widehat f_{j_2,n}^{2}+\frac{\widehat\sigma_{j_2,n}^2}{\sum^{I-{j_2}-1}_{k=-n}C_{k,j_2}}\right)\right)\left(1-\prod^{I+n-1}_{j_4=K+1}\widehat{f}_{j_4,n}\right)\\
&+\sum^{K}_{j_3=i_2} \frac{\widehat\sigma_{j_3,n}^2}{\sum^{I-{j_3}-1}_{k=-n}C_{k,j_3}}\left(\prod^{j_3-1}_{j_4=i_2}\widehat{f}_{j_4,n}^2 \right)\left[\prod^{I+n-1}_{j_2=j_3+1}\left(\widehat f_{j_2,n}^{2}+\frac{\widehat\sigma_{j_2,n}^2}{\sum^{I-{j_2}-1}_{k=-n}C_{k,j_2}}\right)-\prod^{K}_{j_2= j_3+1}\left(\widehat f_{j_2,n}^{2}+\frac{\widehat\sigma_{j_2,n}^2}{\sum^{I-{j_2}-1}_{k=-n}C_{k,j_2}}\right)\right]	\\
& \quad\quad\times\left(\prod^{I+n-1}_{j_4=K+1}\widehat{f}_{j_4,n}\right)\\
&+\sum^{K}_{j_3=i_2} \frac{\widehat\sigma_{j_3,n}^2}{\sum^{I-{j_3}-1}_{k=-n}C_{k,j_3}}\left(\prod^{j_3-1}_{j_4=i_2}\widehat{f}_{j_4,n}^2 \right)\left(\prod^{K}_{j_2= j_3+1}\left(\widehat f_{j_2,n}^{2}+\frac{\widehat\sigma_{j_2,n}^2}{\sum^{I-{j_2}-1}_{k=-n}C_{k,j_2}}\right)\right)\left(1-\prod^{I+n-1}_{j_4=K+1}\widehat{f}_{j_4,n}\right).
\end{align*}
Now, following the same steps as in the proof of Theorem 4.7 in \citet{jentschasympThInf}, we can compute the unconditional expectation of the above. Using $C_{i,j}>\epsilon^j$, $E(\widehat f_{c,n}|\mathcal{B}_{I,n}(c))=f_c$, $E(\widehat \sigma^2_{c,n}|\mathcal{B}_{I,n}(c))=\sigma_c^2$ as well as
\begin{align}\label{cond_second_moment_fhat}
E\left(\widehat f_{c,n}^2|\mathcal{B}_{I,n}(c)\right)=\frac{\sigma_c^2}{\sum_{k=-n}^{I-c-1}C_{k,c}}+f_c^2\leq \frac{\sigma_c^2}{(I+n-c)\epsilon^c}+f_c^2 \leq \frac{\sigma_c^2}{\epsilon^c}+f_c^2
\end{align}
for all $c\in\{0,\ldots,I+n-1\}$, 
where $\mathcal{B}_{I,n}(k)=\left\{C_{i,j}|i=-n,\ldots,I,~j=0,\ldots,k,~i+j\leq I+n\right\}$, we can argue that $Var^*(B_{2,K,I,n}^*|\mathcal{Q}_{I,n}^*=\mathcal{Q}_{I,n})\geq 0$ vanishes in probability for $n\rightarrow \infty$ and $K\rightarrow \infty$ afterwards, by showing that its unconditional expectation is bounded for $n\rightarrow \infty$ and that its bound converges to zero as $K\rightarrow \infty$.

Using that $E(Cov^*(\prod^{I+n-1}_{j_1=i_1}\widehat f_{j_1,n}^*-\prod^{K}_{j_1=i_1}\widehat f_{j_1,n}^*,\prod^{I+n-1}_{j_2=i_2}\widehat f_{j_2,n}^*-\prod^{K}_{j_2=i_2}\widehat f_{j_2,n}^*))$ can be bounded by
\begin{align*}
\left(\prod^{i_2-1}_{j_1=i_1}f_{j_1}\right)&\left\{\sum^{I+n-1}_{j_3=K+1} \frac{\sigma_{j_3}^2}{(I+n-j_3)\epsilon^{j_3}}\left(\prod^{j_3-1}_{j_4=i_2}\left(f_{j_4}^2+\frac{\sigma_{j_4}^2}{\epsilon^{j_4}}\right) \right)\left(\prod^{I+n-1}_{j_2=j_3+1}\left(f_{j_2}^{2}+\frac{2\sigma_{j_2}^2}{\epsilon^{j_2}}\right)\right)\right.\\
&+\sum^{K}_{j_3=i_2} \frac{\sigma_{j_3}^2}{(I+n-j_3)\epsilon^{j_3}}\left(\prod^{j_3-1}_{j_4=i_2}\left(f_{j_4}^2+\frac{\sigma_{j_4}^2}{\epsilon^{j_4}}\right) \right)\left(\prod^{I+n-1}_{j_2=j_3+1}\left(f_{j_2}^{2}+\frac{2\sigma_{j_2,n}^2}{\epsilon^{j_2}}\right)\right)\left(\prod^{I+n-1}_{j_4=K+1}f_{j_4}-1\right)\\
&+\sum^{K}_{j_3=i_2} \frac{\sigma_{j_3}^2}{(I+n-j_3)\epsilon^{j_3}}\left(\prod^{j_3-1}_{j_4=i_2}\left(f_{j_4}^2+\frac{\sigma_{j_4}^2}{\epsilon^{j_4}}\right) \right)	\\
& \qquad\times \left[\prod^{K}_{j_2= j_3+1}\left(f_{j_2}^{2}+\frac{2\sigma_{j_2,n}^2}{\epsilon^{j_2}}\right)\left(\prod^{I+n-1}_{j_2=K+1}\left(f_{j_2}^{2}+\frac{2\sigma_{j_2}^2}{\epsilon^{j_2}}\right)-1\right)\right]	\left(\prod^{I+n-1}_{j_4=K+1}f_{j_4}\right)\\
&+\left.\sum^{K}_{j_3=i_2} \frac{\sigma_{j_3}^2}{(I+n-j_3)\epsilon^{j_3}}\left(\prod^{j_3-1}_{j_4=i_2}\left(f_{j_4}^2+\frac{\sigma_{j_4}^2}{\epsilon^{j_4}}\right) \right)\left(\prod^{K}_{j_2= j_3+1}\left(f_{j_2}^{2}+\frac{2\sigma_{j_2}^2}{\epsilon^{j_2}}\right)\right)\left(\prod^{I+n-1}_{j_4=K+1}f_{j_4}-1\right)\right\},
\end{align*}
we can bound also $E(Var^*(B_{2,K,I,n}^*|\mathcal{Q}_{I,n}^*=\mathcal{Q}_{I,n}))$ from above. Precisely, putting everything together, we get
\begin{align*}
& E(Var^*(B_{2,K,I,n}^*|\mathcal{Q}_{I,n}^*=\mathcal{Q}_{I,n}))	\\
&\leq 2 \sum_{i_2=0}^K\sum_{i_1=0}^{i_2}\mu_{i_1}\mu_{i_2}\left(\prod^{i_2-1}_{j_1=i_1}f_{j_1}\right)	\\
& \quad\times\left\{\sum^{I+n-1}_{j_3=K+1} \frac{(I+n+1)\sigma_{j_3}^2}{(I+n-j_3)\epsilon^{j_3}}\left(\prod^{j_3-1}_{j_4=i_2}\left(f_{j_4}^2+\frac{\sigma_{j_4}^2}{\epsilon^{j_4}}\right) \right)\left(\prod^{I+n-1}_{j_2=j_3+1}\left(f_{j_2}^{2}+\frac{2\sigma_{j_2}^2}{\epsilon^{j_2}}\right)\right)\right.\\
&\qquad+\sum^{K}_{j_3=i_2} \frac{(I+n+1)\sigma_{j_3}^2}{(I+n-j_3)\epsilon^{j_3}}\left(\prod^{j_3-1}_{j_4=i_2}\left(f_{j_4}^2+\frac{\sigma_{j_4}^2}{\epsilon^{j_4}}\right) \right)\left(\prod^{I+n-1}_{j_2=j_3+1}\left(f_{j_2}^{2}+\frac{2\sigma_{j_2,n}^2}{\epsilon^{j_2}}\right)\right)\left(\prod^{I+n-1}_{j_4=K+1}f_{j_4}-1\right)\\
&\qquad+\sum^{K}_{j_3=i_2} \frac{(I+n+1)\sigma_{j_3}^2}{(I+n-j_3)\epsilon^{j_3}}\left(\prod^{j_3-1}_{j_4=i_2}\left(f_{j_4}^2+\frac{\sigma_{j_4}^2}{\epsilon^{j_4}}\right) \right)	\\
& \qquad\qquad\times \left[\prod^{K}_{j_2= j_3+1}\left(f_{j_2}^{2}+\frac{2\sigma_{j_2,n}^2}{\epsilon^{j_2}}\right)\left(\prod^{I+n-1}_{j_2=K+1}\left(f_{j_2}^{2}+\frac{2\sigma_{j_2}^2}{\epsilon^{j_2}}\right)-1\right)\right]	\left(\prod^{I+n-1}_{j_4=K+1}f_{j_4}\right)\\
&\qquad+\left.\sum^{K}_{j_3=i_2} \frac{(I+n+1)\sigma_{j_3}^2}{(I+n-j_3)\epsilon^{j_3}}\left(\prod^{j_3-1}_{j_4=i_2}\left(f_{j_4}^2+\frac{\sigma_{j_4}^2}{\epsilon^{j_4}}\right) \right)\left(\prod^{K}_{j_2= j_3+1}\left(f_{j_2}^{2}+\frac{2\sigma_{j_2}^2}{\epsilon^{j_2}}\right)\right)\left(\prod^{I+n-1}_{j_4=K+1}f_{j_4}-1\right)\right\}
\end{align*}
and the leading term of the last right-hand side becomes
\begin{align*}
& 2 \sum_{i_2=0}^K\sum_{i_1=0}^{i_2}\mu_{i_1}\mu_{i_2}\left(\prod^{i_2-1}_{j_1=i_1}f_{j_1}\right)	\\
& \quad\times\left\{\sum^{I+n-1}_{j_3=K+1} \frac{\sigma_{j_3}^2}{\epsilon^{j_3}}\left(\prod^{j_3-1}_{j_4=i_2}\left(f_{j_4}^2+\frac{\sigma_{j_4}^2}{\epsilon^{j_4}}\right) \right)\left(\prod^{I+n-1}_{j_2=j_3+1}\left(f_{j_2}^{2}+\frac{2\sigma_{j_2}^2}{\epsilon^{j_2}}\right)\right)\right.\\
&\qquad+\sum^{K}_{j_3=i_2} \frac{\sigma_{j_3}^2}{\epsilon^{j_3}}\left(\prod^{j_3-1}_{j_4=i_2}\left(f_{j_4}^2+\frac{\sigma_{j_4}^2}{\epsilon^{j_4}}\right) \right)\left(\prod^{I+n-1}_{j_2=j_3+1}\left(f_{j_2}^{2}+\frac{2\sigma_{j_2,n}^2}{\epsilon^{j_2}}\right)\right)\left(\prod^{I+n-1}_{j_4=K+1}f_{j_4}-1\right)\\
&\qquad+\sum^{K}_{j_3=i_2} \frac{\sigma_{j_3}^2}{\epsilon^{j_3}}\left(\prod^{j_3-1}_{j_4=i_2}\left(f_{j_4}^2+\frac{\sigma_{j_4}^2}{\epsilon^{j_4}}\right) \right)	\left[\prod^{K}_{j_2= j_3+1}\left(f_{j_2}^{2}+\frac{2\sigma_{j_2,n}^2}{\epsilon^{j_2}}\right)\left(\prod^{I+n-1}_{j_2=K+1}\left(f_{j_2}^{2}+\frac{2\sigma_{j_2}^2}{\epsilon^{j_2}}\right)-1\right)\right]	\left(\prod^{I+n-1}_{j_4=K+1}f_{j_4}\right)\\
&\qquad+\left.\sum^{K}_{j_3=i_2} \frac{\sigma_{j_3}^2}{\epsilon^{j_3}}\left(\prod^{j_3-1}_{j_4=i_2}\left(f_{j_4}^2+\frac{\sigma_{j_4}^2}{\epsilon^{j_4}}\right) \right)\left(\prod^{K}_{j_2= j_3+1}\left(f_{j_2}^{2}+\frac{2\sigma_{j_2}^2}{\epsilon^{j_2}}\right)\right)\left(\prod^{I+n-1}_{j_4=K+1}f_{j_4}-1\right)\right\},
\end{align*}
which can be bounded further by
\begin{align*}
& 2 \mu_\infty\left(\prod^\infty_{j_4=0}\left(f_{j_4}^2+\frac{2\sigma_{j_4}^2}{\epsilon^{j_4}}\right)\right)\sum_{i_2=0}^K\sum_{i_1=0}^{i_2}\mu_{i_1}\mu_{i_2}\left\{\sum^\infty_{j_3=K+1} \frac{\sigma_{j_3}^2}{\epsilon^{j_3}}+\sum^{K}_{j_3=i_2} \frac{\sigma_{j_3}^2}{\epsilon^{j_3}}\left(\prod^\infty_{j_4=K+1}f_{j_4}-1\right)\right.\\
&\qquad\left.+\sum^{K}_{j_3=i_2} \frac{\sigma_{j_3}^2}{\epsilon^{j_3}}	\left[\left(\prod^\infty_{j_2=K+1}\left(f_{j_2}^{2}+\frac{2\sigma_{j_2}^2}{\epsilon^{j_2}}\right)-1\right)\right]	\left(\prod^\infty_{j_4=K+1}f_{j_4}\right)+\sum^{K}_{j_3=i_2} \frac{\sigma_{j_3}^2}{\epsilon^{j_3}}\left(\prod^\infty_{j_4=K+1}f_{j_4}-1\right)\right\}.
\end{align*}
Now, considering the four terms in brackets separately, for the first one, we can argue that it vanishes asymptotically due to
\begin{align*}
\sum_{i_2=0}^K\sum_{i_1=0}^{i_2}\mu_{i_1}\mu_{i_2} \sum^\infty_{j_3=K+1} \frac{\sigma_{j_3}^2}{\epsilon^{j_3}}
&\leq \left(\frac{1}{K+1}\sum_{i_2=0}^K\mu_{i_2}\right)\left(\frac{1}{i_2+1}\sum_{i_1=0}^{i_2}\mu_{i_1}\right) \sum^\infty_{j_3=K+1} (j_3+1)^2\frac{\sigma_{j_3}^2}{\epsilon^{j_3}}\rightarrow 0
\end{align*}
for $K\rightarrow\infty$, because the sequence $(\frac{1}{j+1}\sum^j_{i=0}\mu_{i},j\in\mathbb{N}_0)$ is converging and, consequently, also bounded, and due to $\sum_{j=0}^\infty (j+1)^2\frac{\sigma_j^2}{\epsilon^j}<\infty$ by Assumption 4.2. Similarly, using that $\prod^\infty_{j=K+1}f_j\rightarrow 1$ and $\prod^\infty_{j=K+1}(f_j^{2}+\frac{2\sigma_j^2}{\epsilon^j})\rightarrow 1$ for $K\rightarrow \infty$, we can also show that the other three terms vanish asymptotically. This completes the first part of c) for $B_{2, K, I,n}^*$.

Similarly, for showing part c) for $B_{3,K,I,n}^*$, we have to consider
\begin{align*}
B_{3,K,I,n}^*|(\mathcal{Q}_{I,n}^*=\mathcal{Q}_{I,n},\mathcal{D}_{I,n}) = \sqrt{I+n+1}\sum^{I+n}_{i=K+1}C_{I-i,i}^*\left(\prod^{I+n-1}_{j=i}\widehat{f}_{j,n}^*-\prod^{I+n-1}_{j=i}\widehat{f}_{j,n}\right)|(\mathcal{Q}_{I,n}^*=\mathcal{Q}_{I,n},\mathcal{D}_{I,n}).
\end{align*}
By the same arguments as used above for $B_{2,K,I,n}^*$, we get $E^*(B_{3,K,I,n}^*|\mathcal{Q}_{I,n}^*=\mathcal{Q}_{I,n})=0$ and for any fixed $K\in\mathbb{N}_0$ and $n\in\mathbb{N}_0$ large enough such that $K<I+n-1$, we have
\begin{align*}%\label{B3KI_var}
& Var^*(B_{3,K,I,n}^*|\mathcal{Q}_{I,n}^*=\mathcal{Q}_{I,n})	\\
=& (I+n+1) \sum_{i_1,i_2=K+1}^{I+n-1} C_{I-i_1,i_1}C_{I-i_2,i_2}Cov^*\left(\prod^{I+n-1}_{j_1=i_1}\widehat f^*_{j_1,n},\prod^{I+n-1}_{j_2=i_2}\widehat f^*_{j_2,n}|\mathcal{Q}_{I,n}^*=\mathcal{Q}_{I,n}\right)	\\
\leq& 2(I+n+1) \sum_{i_2=K+1}^{I+n-1}\sum_{i_1=K+1}^{i_2} C_{I-i_1,i_1}C_{I-i_2,i_2} Cov^*\left(\prod^{I+n-1}_{j_1=i_1}\widehat f^*_{j_1,n},\prod^{I+n-1}_{j_2=i_2}\widehat f^*_{j_2,n}\right).%\label{B3KI_var_2}
\end{align*}
To calculate the covariance on the last right-hand side, for $i_1\leq i_2$, we consider the mixed moment
\begin{align*} E^*\left(\left(\prod^{i_2-1}_{j_1=i_1}\widehat f_{j_1,n}^*\right)\left(\prod^{I+n-1}_{j_2=i_2}\widehat f_{j_2,n}^{*2}\right)\right),
\end{align*}
which is just the first term of the mixed moment of the covariance calculated for $B_{2,K,I,n}^*$. By using similar calculations to get $E^*(\widehat{f}_{c,n}^{*2})$ (for $B^*_{2,K,I,n}$), we obtain 
\begin{align*}
&Cov^*\left(\prod^{I+n-1}_{j_1=i_1}\widehat f^*_{j_1,n},\prod^{I+n-1}_{j_2=i_2}\widehat f^*_{j_2,n}\right)\\
&=\left[E^*\left(\left(\prod^{i_2-1}_{j_1=i_1}\widehat f_{j_1,n}^*\right)\left(\prod^{I+n-1}_{j_2=i_2}\widehat f_{j_2,n}^{*2}\right)\right)-E^*\left(\prod^{I+n-1}_{j_1=i_1}\widehat f^*_{j_1,n}\right)E^*\left(\prod^{I+n-1}_{j_2=i_2}\widehat f^*_{j_2,n}\right)\right]\\
&=\left[\left(\prod^{i_2-1}_{j_1=i_1}E^*\left(\widehat f_{j_1,n}^*\right)\right)\left(\prod^{I+n-1}_{j_2=i_2}E^*\left(\widehat f_{j_2,n}^{*2}\right)\right)-\left(\prod^{I+n-1}_{j_1=i_1}\widehat f_{j_1,n}\right)\left(\prod^{I+n-1}_{j_2=i_2}\widehat f_{j_2,n}\right)\right]\\
&=\left[\left(\prod^{i_2-1}_{j_1=i_1}\widehat f_{j_1,n}\right)\left(\prod^{I+n-1}_{j_2=i_2}\left(\widehat f_{j_2,n}^{2}+\frac{\widehat \sigma_{j_2,n}^2}{\sum_{k=-n}^{I-j_2-1}C_{k,j_2}}\right)\right)-\left(\prod^{I+n-1}_{j_1=i_1}\widehat f_{j_1,n}\right)\left(\prod^{I+n-1}_{j_2=i_2}\widehat f_{j_2,n}\right)\right]	\\
&=\left(\prod^{i_2-1}_{j_1=i_1}\widehat f_{j_1,n}\right)\left[\prod^{I+n-1}_{j_2=i_2}\left(\widehat f_{j_2,n}^{2}+\frac{\widehat \sigma_{j_2,n}^2}{\sum_{k=-n}^{I-j_2-1}C_{k,j_2}}\right)-\prod^{I+n-1}_{j_2=i_2}\widehat f_{j_2,n}\right]	\\
&= \sum_{j_4=i_2}^{I+n-1} \frac{\widehat \sigma_{j_4,n}^2}{\sum_{k=-n}^{I-j_4-1}C_{k,j_4}}\left(\prod^{i_2-1}_{j_1=i_1}\widehat f_{j_1,n}\right)\left(\prod_{j_2=i_2}^{j_4-1} \widehat f_{j_2,n}^{2}\right)\left(\prod_{j_3=j_4+1}^{I+n-1} \left(\widehat f_{j_3,n}^2+\frac{\widehat \sigma_{j_3,n}^2}{\sum_{k=-n}^{I-j_3-1}C_{k,j_3}}\right)\right)	\\
&\leq \sum_{j_4=i_2}^{I+n-1} \frac{(I+n+1)\widehat \sigma_{j_4,n}^2}{(I+n-j_4)\epsilon^{j_4}}\left(\prod^{i_2-1}_{j_1=i_1}\widehat f_{j_1,n}\right)\left(\prod_{j_2=i_2}^{j_4-1} \widehat f_{j_2,n}^{2}\right)\left(\prod_{j_3=j_4+1}^{I+n-1} \left(\widehat f_{j_3,n}^2+\frac{\widehat \sigma_{j_3,n}^2}{\epsilon^{j_3}}\right)\right).
\end{align*}
Noting that all involved summands and factors are non-negative, taking expectations of the last right-hand side and using the law of iterative expectations and $C_{i,j}>\epsilon^j$, $E(\widehat f_{c,n}|\mathcal{B}_{I,n}(c))=f_c$, $E(\widehat \sigma^2_{c,n}|\mathcal{B}_{I,n}(c))=\sigma_c^2$ as well as \eqref{cond_second_moment_fhat}, we get
\begin{align*}
\sum_{j_4=i_2}^{I+n-1} \frac{\sigma_{j_4}^2}{(I+n-j_4)\epsilon^{j_4}}\left(\prod^{i_2-1}_{j_1=i_1}f_{j_1}\right)\left(\prod_{j_2=i_2}^{j_4-1} \left(f_{j_2}^{2}+\frac{\sigma_{j_2}^2}{\epsilon^{j_2}}\right)\right)\left(\prod_{j_3=j_4+1}^{I+n-1} \left(f_{j_3}^2+\frac{2\sigma_{j_3}^2}{\epsilon^{j_3}}\right)\right)
\end{align*}
such that the leading term of $E(Var^*(B_{3,K,I,n}^*|\mathcal{Q}_{I,n}^*=\mathcal{Q}_{I,n}))$ becomes
\begin{align*}
& 2\sum_{i_2=K+1}^{I+n-1}\sum_{i_1=K+1}^{i_2} \mu_{i_1}\mu_{i_2}\sum_{j_4=i_2}^{I+n-1} \frac{\sigma_{j_4}^2}{\epsilon^{j_4}}\left(\prod^{i_2-1}_{j_1=i_1}f_{j_1}\right)\left(\prod_{j_2=i_2}^{j_4-1} \left(f_{j_2}^{2}+\frac{\sigma_{j_2}^2}{\epsilon^{j_2}}\right)\right)\left(\prod_{j_3=j_4+1}^{I+n-1} \left(f_{j_3}^2+\frac{2\sigma_{j_3}^2}{\epsilon^{j_3}}\right)\right)		\\
\leq& 2\mu_\infty\left(\prod_{j_3=0}^\infty \left(f_{j_3}^2+\frac{2\sigma_{j_3}^2}{\epsilon^{j_3}}\right)\right)\sum_{i_2=K+1}^{I+n-1}\sum_{i_1=K+1}^{i_2} \mu_{i_1}\mu_{i_2}\sum_{j_4=i_2}^{I+n-1} \frac{\sigma_{j_4}^2}{\epsilon^{j_4}}	\\
\leq& 2\mu_\infty\left(\prod_{j_3=0}^\infty \left(f_{j_3}^2+\frac{2\sigma_{j_3}^2}{\epsilon^{j_3}}\right)\right)\sum_{i_2=K+1}^{I+n-1}\sum_{i_1=K+1}^{i_2} \mu_{i_1}\mu_{i_2}\sum_{j_4=i_2}^{I+n-1} \frac{\sigma_{j_4}^2}{\epsilon^{j_4}}.
\end{align*}
%\footnote{CJ: Kann man hier nicht die $\mu_i$s gleich loswerden? Noch zurueckaendern!}
For the triple sum on the last right-hand side, we get
\begin{align*}
\sum_{i_2=K+1}^{I+n-1}\sum_{i_1=K+1}^{i_2} \mu_{i_1}\mu_{i_2}\sum_{j_4=i_2}^{I+n-1} \frac{\sigma_{j_4}^2}{\epsilon^{j_4}} =& \sum_{i_2=1}^{I+n-K-1}\sum_{i_1=K+1}^{i_2+K} \mu_{i_1}\mu_{i_2+K}\sum_{j_4=i_2+K}^{I+n-1} \frac{\sigma_{j_4}^2}{\epsilon^{j_4}}	\\
=& \sum_{i_2=1}^{I+n-K-1}i_2\left(\frac{1}{i_2}\sum_{i_1=K+1}^{i_2+K} \mu_{i_1}\right)\mu_{i_2+K}\sum_{j_4=i_2+K}^{I+n-1} \frac{\sigma_{j_4}^2}{\epsilon^{j_4}}	\\
\leq& const. \sum_{i_2=1}^{I+n-K-1}i_2 \mu_{i_2+K}\sum_{j_4=i_2+K}^{I+n-1} \frac{\sigma_{j_4}^2}{\epsilon^{j_4}}	\\
=& const. \sum_{j=1+K}^{I+n-1}\frac{\sigma_j^2}{\epsilon^j}\sum_{l=1}^{j-K}l\mu_{l+K}.
\end{align*}
Further, the sequence $(\mu_{i},i\in\mathbb{N}_0)$ shares the properties of $(C_{I-i,i},i\in\mathbb{N}_0)$ in a deterministic sense  
such that $\sum_{l=1}^{j-K}l\mu_{l+K}\leq const. l^2$.
Consequently, we have
\begin{align*}
\sum_{j=1+K}^{I+n-1}\frac{\sigma_j^2}{\epsilon^j}\sum_{l=1}^{j-K}l\mu_{l+K}\leq const. \sum_{j=1+K}^\infty j^2\frac{\sigma_j^2}{\epsilon^j} \rightarrow 0
\end{align*}
as $K\rightarrow \infty$ by Assumption 4.2.	$\hfill\square$

\subsection{Proof of Theorem 4.11}
The proof is analogous to the proof of Theorem 4.12 and Corollary 4.13 in \cite{jentschasympThInf}. The claimed uncorrelatedness of $(R^*_{I,n}-\widehat R_{I,n})_1$ and $(R^*_{I,n}-\widehat R_{I,n})_2$ conditional on $\mathcal{Q}_{I,n}^*=\mathcal{Q}_{I,n}$ and $\mathcal{D}_{I,n}$, follows from

\thickmuskip=0.01mu
\medmuskip=0.01mu
\thickmuskip=0.01mu
\begin{align*}
	& Cov^*\left((R_{I,n}^*-\widehat R_{I,n}^*)_1,\sqrt{I+n+1}(R_{I,n}^*-\widehat R_{I,n}^*)_2|\mathcal{Q}_{I,n}^*=\mathcal{Q}_{I,n}\right)	\\
	=&Cov^*\left(\sum_{i=0}^{I+n}C_{I-i,i}\left(\prod^{I+n-1}_{j=i}F_{I-i,j}^*-\prod_{j=i}^{I+n-1}{\widehat f}_{j,n}^*\right),\sqrt{I+n+1}\sum^{I+n}_{i=0}C_{I-i,i}\left(\prod_{j=i}^{I+n-1}{\widehat f}_{j,n}^*-\prod_{j=i}^{I+n-1}\widehat{f}_{j,n}\right)|\mathcal{Q}_{I,n}^*=\mathcal{Q}_{I,n}\right)	\\
	=& \sqrt{I+n+1}\sum_{i_1,i_2=0}^{I+n}E^*\left(C_{I-i_1,i_1}\left(\prod^{I+n-1}_{j_1=i_1}F_{I-i_1,j_1}^*-\prod_{j_1=i_1}^{I+n-1}{\widehat f}^*_{j_1,n}\right)C_{I-i_2,i_2}\left(\prod_{j_2=i_2}^{I+n-1}{\widehat f}^*_{j_2}-\prod_{j_2=i_2}^{I+n-1}\widehat{f}_{j_2,n}\right)|\mathcal{Q}_{I,n}^*=\mathcal{Q}_{I,n}\right)	\\
	=& 0
\end{align*}
since for all $i_1,i_2=0, \dots , I+n$, we have 
\thickmuskip=0.01mu
\medmuskip=0.01mu
\thickmuskip=0.01mu
\begin{align*}	
&E^*\left(C_{I-i_1,i_1}C_{I-i_2,i_2}\left(\prod_{j_2=i_2}^{I+n-1}{\widehat f^*}_{j_2,n}-\prod_{j_2=i_2}^{I+n-1}\widehat{f}_{j_2,n}\right)\mspace{-5mu} E^*\left(\prod^{I+n-1}_{j_1=i_1}F^*_{I-i_1,j_1}-\prod_{j_1=i_1}^{I+n-1}{\widehat f}^*_{j_1,n} \bigg|\mathcal{F}^*_{I,n}, \mathcal{Q}_{I,n}^*=\mathcal{Q}_{I,n}\right)|\mathcal{Q}_{I,n}^*=\mathcal{Q}_{I,n}\right)	\\
&=0,
\end{align*}
because the inner conditional expectation on the last right-hand side is zero. \hfill $\square$	

\bigskip

\section{Proofs of Section 6}

\subsection{Proof of Theorem 6.1}
Following the technique of proof in Theorem 4.8 and using $\widehat f_{j,n}-f_j=O_P((I+n-1)^{-1/2})$, $\widehat f_{j,n}^*-\widehat f_{j,n}=O_{P^*}((I+n-1)^{-1/2})$ and $\widehat \sigma_{j,n}^2-\sigma_j^2=O_P((I+n-1)^{-1/2})$ leads to the same limiting result also for the process uncertainty part $(R_{I,n}^+-\widehat{R}_{I,n}^+)_1$ of the alternative Mack bootstrap. \hfill $\square$

\subsection{Proof of Theorem 6.3}

Following the technique of proof in Theorem 4.10 and exploiting the limiting properties from Assumption 6.2, we get the claimed asymptotic results.	\hfill $\square$

\subsection{Proof of Theorem 6.4}
Based on the results established in Theorems 6.1 and 6.3, the arguments are completely analogous to those used in the proof of Theorem 4.11.		\hfill $\square$ 

\bigskip	  
		
\section{Conditional versions of the CLTs from \cite{jentschasympThInf}}

For the sake of completeness, in Theorem \ref{CLT_f_cond} and Corollary \ref{CLTProdf_cond} below, we summarize the results from Theorem C.1(ii,iv) and Corollary C.2(ii,iv) in \cite{jentschasympThInf}.

\begin{thm}[Asymptotic normality of $\fh$ conditionally on $\mathcal{Q}_{I,n}$; Theorem C.1(ii,iv)  in \cite{jentschasympThInf}]\label{CLT_f_cond}
	Suppose Assumptions 2.2, 2.3, 2.4, 4.2 and 4.3 are satisfied. Then, as $n\rightarrow \infty$, the following holds:
	\begin{itemize}
		\item[(i)] For each fixed $K\in\mathbb{N}_0$, let $\underline{f}_K=(f_{0},f_{1},\ldots,f_{K})^\prime$ and define
		\begin{align*}
		\underline{f}_{K,n}(\mathcal{Q}_{I,n})=(f_{0,n}(\mathcal{Q}_{I,n}),f_{1,n}(\mathcal{Q}_{I,n}),\ldots,f_{K,n}(\mathcal{Q}_{I,n}))^\prime.
		\end{align*}
		Then, unconditionally, we have
		\begin{align*}
		J_n^{1/2}\left(\underline{f}_{K,n}(\mathcal{Q}_{I,n})-\underline{f}_K\right) \overset{d}{\longrightarrow} \mathcal{N}\left(0, \bm{\Sigma}_{K,\underline{f}}^{(1)}\right),
		\end{align*}
		where $J_n^{1/2}=diag\left(\sqrt{I+n+1-j\vphantom{I^2}},j=0,\ldots,K\right)$ is a diagonal $(K+1)\times (K+1)$ matrix of inflation factors and the variance-covariance matrix
		\begin{align*}
		\bm{\Sigma}_{K,\underline{f}}^{(1)}=J_g\left(\underline{\mu}_K\right)\mathbf{\Sigma}_{K,\underline{C}}^{(1)}J_g\left(\underline{\mu}_K\right)^\prime,
		\end{align*}
		where $\mathbf{\Sigma}_{K,\underline{C}}^{(1)}$ is defined in \eqref{Sigma_Ck1}, has entries 
		\begin{align*}
		&\bm{\Sigma}_{K,\underline{f}}^{(1)}(j_1,j_2)	\\
		=&\frac{f_{j_1}f_{j_2}E\left(E(C_{i,j_1}|C_{i,\infty})E(E(C_{i,j_2}|C_{i,\infty}))\right)+E\left(E(C_{i,j_1+1}|C_{i,\infty})E(E(C_{i,j_2+1}|C_{i,\infty}))\right)}{\mu_{j_1}\mu_{j_2}}\\
		&+\frac{-f_{j_2}E\left(E(C_{i,j_{1}+1}|C_{i,\infty})E(C_{i,j_2}|C_{i,\infty})\right)-f_{j_1}E\left(E(C_{i,j_{1}}|C_{i,\infty})E(C_{i,j_2+1}|C_{i,\infty})\right)}{\mu_{j_1}\mu_{j_2}}
		\end{align*}
		for $j_1, j_2=0, \dots, K$.
		\item[(ii)] For each fixed $K\in\mathbb{N}_0$, let $\underline{\widehat f}_{K,n}=(\widehat f_{0,n},\widehat f_{1,n},\ldots,\widehat f_{K,n})^\prime$. Then, conditionally on $\mathcal{Q}_{I,n}$, we have
		\begin{align*}
		J_n^{1/2}\left(\underline{\widehat f}_{K,n}-\underline{f}_{K,n}(\mathcal{Q}_{I,n})\right)|\mathcal{Q}_{I,n} \overset{d}{\longrightarrow} \mathcal{N}\left(0,\bm{\Sigma}_{K,\underline{f}}^{(2)}\right),
		\end{align*}
		where the variance-covariance matrix
		\begin{align*}
		\bm{\Sigma}_{K,\underline{f}}^{(2)}=J_g\left(\underline{\mu}_K\right)\mathbf{\Sigma}_{K,\underline{C}}^{(2)}J_g\left(\underline{\mu}_K\right)^\prime,
		\end{align*}
		where $\mathbf{\Sigma}_{K,\underline{C}}^{(2)}$ is defined in \eqref{Sigma_Ck2}, has entries $\bm{\Sigma}_{K,\underline{f}}^{(2)}(j,j)=\sigma_{f_j,2}^{2}=\frac{\sigma^2_j}{\mu_j}-\sigma_{f_j,1}^{2}$ for $j=0, \dots, K$ and $\bm{\Sigma}_{K,\underline{f}}^{(2)}(j_1,j_2)=-\bm{\Sigma}_{K,\underline{f}}^{(1)}(j_1,j_2)$ for $j_1,j_2=0,\dots, K$, $j_1\not= j_2$.
	\end{itemize}
\end{thm}
We obtain 
\begin{align}\label{Sigma_Ck1}
\mathbf{\Sigma}_{K,\underline{f}}^{(1)}=J_g\left(\underline{\mu}_K\right)\mathbf{\Sigma}_{K,\underline{C}}^{(1)}J_g\left(\underline{\mu}_K\right)^\prime
\end{align}
and
\begin{align}\label{Sigma_Ck2}
\mathbf{\Sigma}_{K,\underline{C}}^{(2)}=E(Var(\underline{C}_{i,K}|C_{i,\infty})),
\end{align}
where $C_{i,\infty} = C_{i,0}\prod_{k=0}^\infty F_{i,k}$. Note that, due to the law of total variance, we have
\begin{align}
\mathbf{\Sigma}_{K,\underline{C}}=\mathbf{\Sigma}_{K,\underline{C}}^{(1)}+\mathbf{\Sigma}_{K,\underline{C}}^{(2)}	\quad	\text{and}	\quad	\mathbf{\Sigma}_{K,\underline{f}}=\mathbf{\Sigma}_{K,\underline{f}}^{(1)}+\mathbf{\Sigma}_{K,\underline{f}}^{(2)},
\end{align}
where $\mathbf{\Sigma}_{K,\underline f}=diag\left(\frac{\sigma_0^2}{\mu_0},\frac{\sigma_1^2}{\mu_1},\ldots,\frac{\sigma_K^2}{\mu_K}\right)$. 

\begin{cor}[Asymptotic normality for products of $\widehat f_{j,n}$'s conditionally on $\mathcal{Q}_{I,n}$; Corollary C.2(ii,iv) in \cite{jentschasympThInf}]\label{CLTProdf_cond}
	Suppose the assumptions of Theorem \ref{CLT_f_cond} hold. Then, as $n\rightarrow \infty$, the following holds:
	\begin{itemize}
		\item[(i)] For each fixed $K\in\mathbb{N}_0$, unconditionally, we have also joint convergence, that is,
		\begin{align*}
		\sqrt{I+n+1}\begin{pmatrix}
		\prod^K_{j=i}f_{j,n}(\mathcal{Q}_{I,n})-\prod^K_{j=i}f_{j,n}	\\
		i=0,\ldots,K
		\end{pmatrix} \overset{d}{\longrightarrow} \mathcal{N}\left(0, \mathbf{\Sigma}_{K,\prod f_j}^{(1)}\right)
		\end{align*}
		where $\mathbf{\Sigma}_{K,\prod f_j}^{(1)}=J_h(\underline{f}_K)\mathbf{\Sigma}_{K,\underline{f}}^{(1)}J_h(\underline{f}_K)^\prime$
		with $J_h(\cdot)$ as defined in \eqref{Jacobian2}.
		\item[(ii)] For each fixed $K\in\mathbb{N}_0$, conditionally on $\mathcal{Q}_{I,n}$, we have also joint convergence, that is,
		\begin{align*}
		\sqrt{I+n+1}\begin{pmatrix}
		\prod^K_{j=i}\widehat f_{j,n}-\prod^K_{j=i}f_{j,n}(\mathcal{Q}_{I,n})	\\
		i=0,\ldots,K
		\end{pmatrix}|\mathcal{Q}_{I,n} \overset{d}{\longrightarrow} \mathcal{N}\left(0,\mathbf{\Sigma}_{K,\prod f_j}^{(2)}\right)
		\end{align*}
		where $\mathbf{\Sigma}_{K,\prod f_j}^{(2)}=\mathbf{\Sigma}_{K,\prod f_j}-\mathbf{\Sigma}_{K,\prod f_j}^{(1)}$, where $\mathbf{\Sigma}_{K,\prod f_j}$ is defined in Corollary 3.2(ii) in \cite{jentschasympThInf}.
	\end{itemize}
\end{cor}

%\newpage

\bigskip

\section{Additional Simulation Results}\label{sec_add_simulation}

Note that the first parts of the alternative Mack bootstrap predictive root of the reserve and the intermediate Mack bootstrap predictive root of the reserve
are equal. Hence, the following findings hold for both approaches.

 Moreover, in both setups a) and b) for the different distributional assumptions, we applied the Kolmogorov-Smirnov test of level $\alpha=5\%$ to test $(R_{I,n}^*-\widehat R_{I,n}^*)_1^{(m)}$ given $\mathcal{Q}_{I,n}^{(m)*}=\mathcal{Q}_{I,n}^{(m)}$ and $\mathcal{D}_{I,n}^{(m)}$ for $m=1,\dots, 500$ is normally distributed with zero mean and variance as in 4.5.

For setup a), it fails to reject the null hypothesis of a Gaussian distribution for about 92\% out of $M=500$ samples, if the gamma distribution is used, for about 87\% in the case of a log-normal, and for about 95\% for a truncated normal distribution to generate the lower bootstrap triangle. The picture is essentially the same for all $n\in\{0,10,20,30,40\}$. In comparison, for setup b), the test does always reject the null for the gamma and for log-normal distribution, but only in about 28\% out of $M=500$  for the truncated normal distribution. The results are pretty similar for all $n\in\{0,10,20,30,40\}$.

These findings can be explained by a property of the gamma and the log-normal distribution. Both tend to 'lose' their skewness and excess of kurtosis for $\frac{C^*_{i,j}}{\widehat \sigma_{j,n}^2}$ growing large in this parameter setting. Hence, as the range for the entries of the first column in setup a) is $[120\times 10^6,350\times 10^6]$ with $[120\times 10^4,350\times 10^4]$ for setup b), we observe more skewness and more excess kurtosis in b) in comparison to a). In particular, this demonstrates that the distribution of the (asymptotically dominating) process uncertainty terms $(R_{I,n}^*-\widehat{R}_{I,n})_1|(\mathcal{Q}_{I,n}^*=\mathcal{Q}_{I,n}, \mathcal{D}_{I,n})$ and $(R_{I,n}^+-\widehat{R}_{I,n}^+)_1| (\mathcal{Q}_{I,n}^+=\mathcal{Q}_{I,n}, \mathcal{D}_{I,n})$, respectively, generally does depend on the distribution (family) of the individual development factors also for large (effective) number of accident years $I+n+1$.

\begin{figure}
	\includegraphics[width=3.5cm,	height=3.5cm]{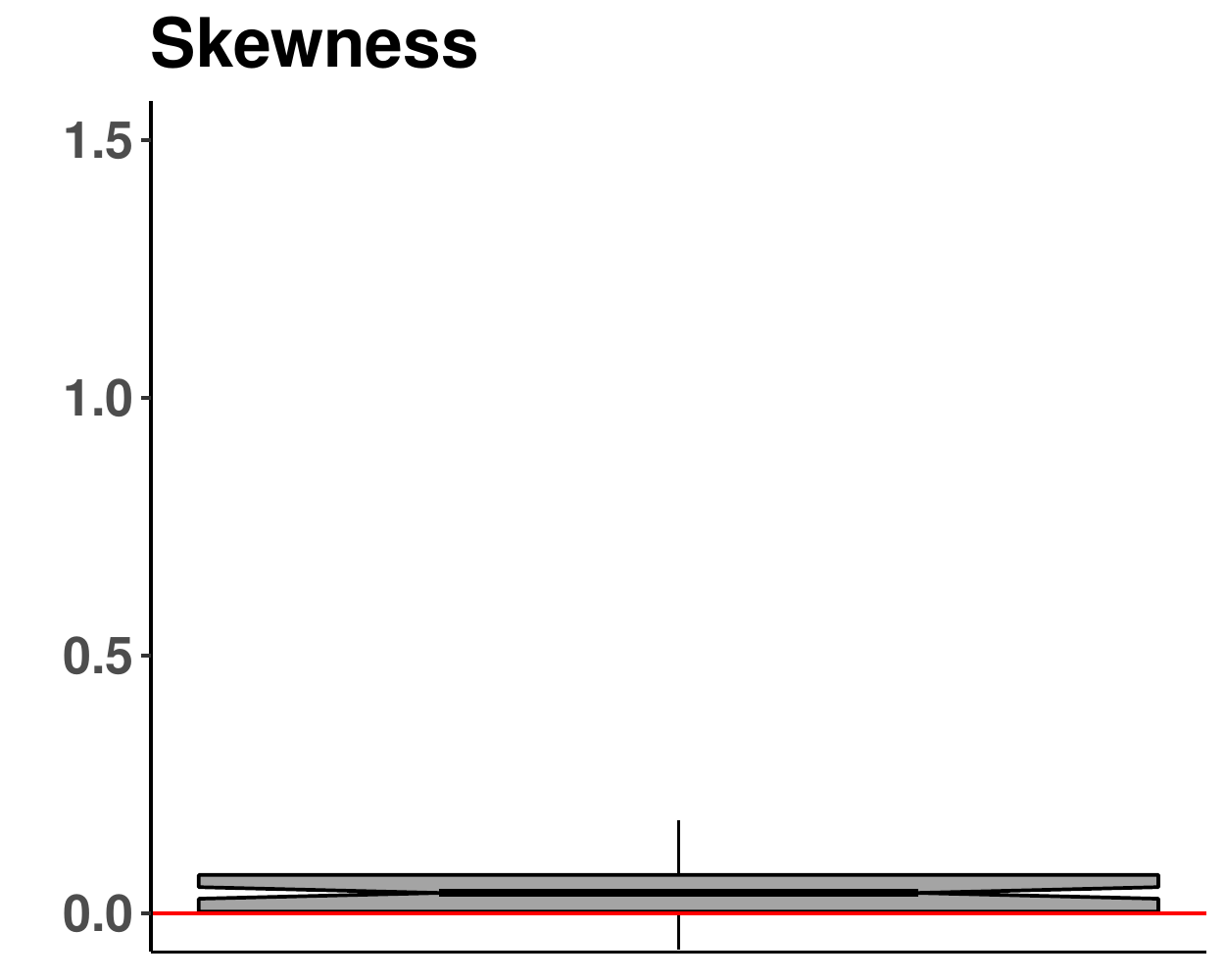}
	\includegraphics[width=3.5cm,	height=3.5cm]{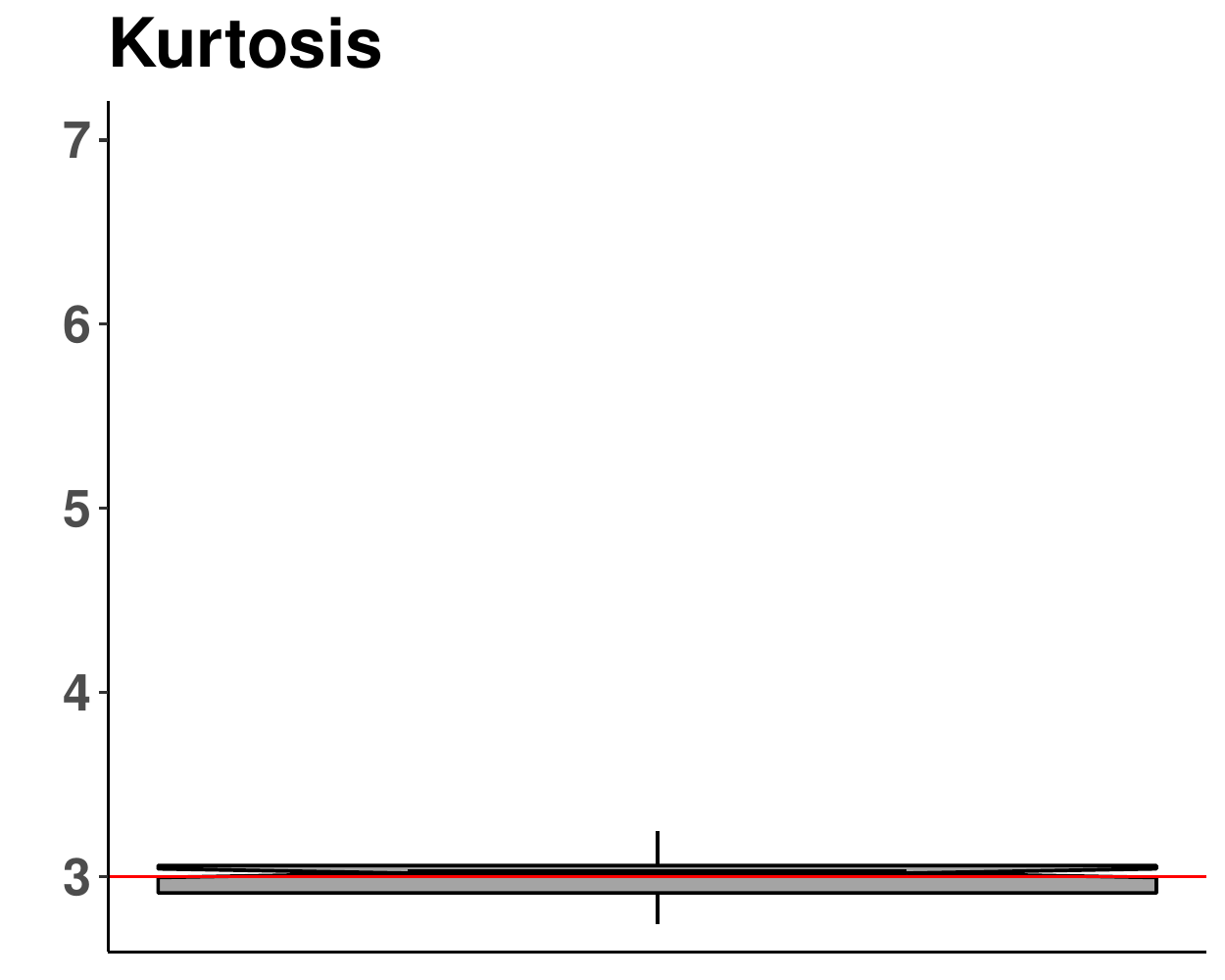}
	\includegraphics[width=5cm,	height=3.5cm]{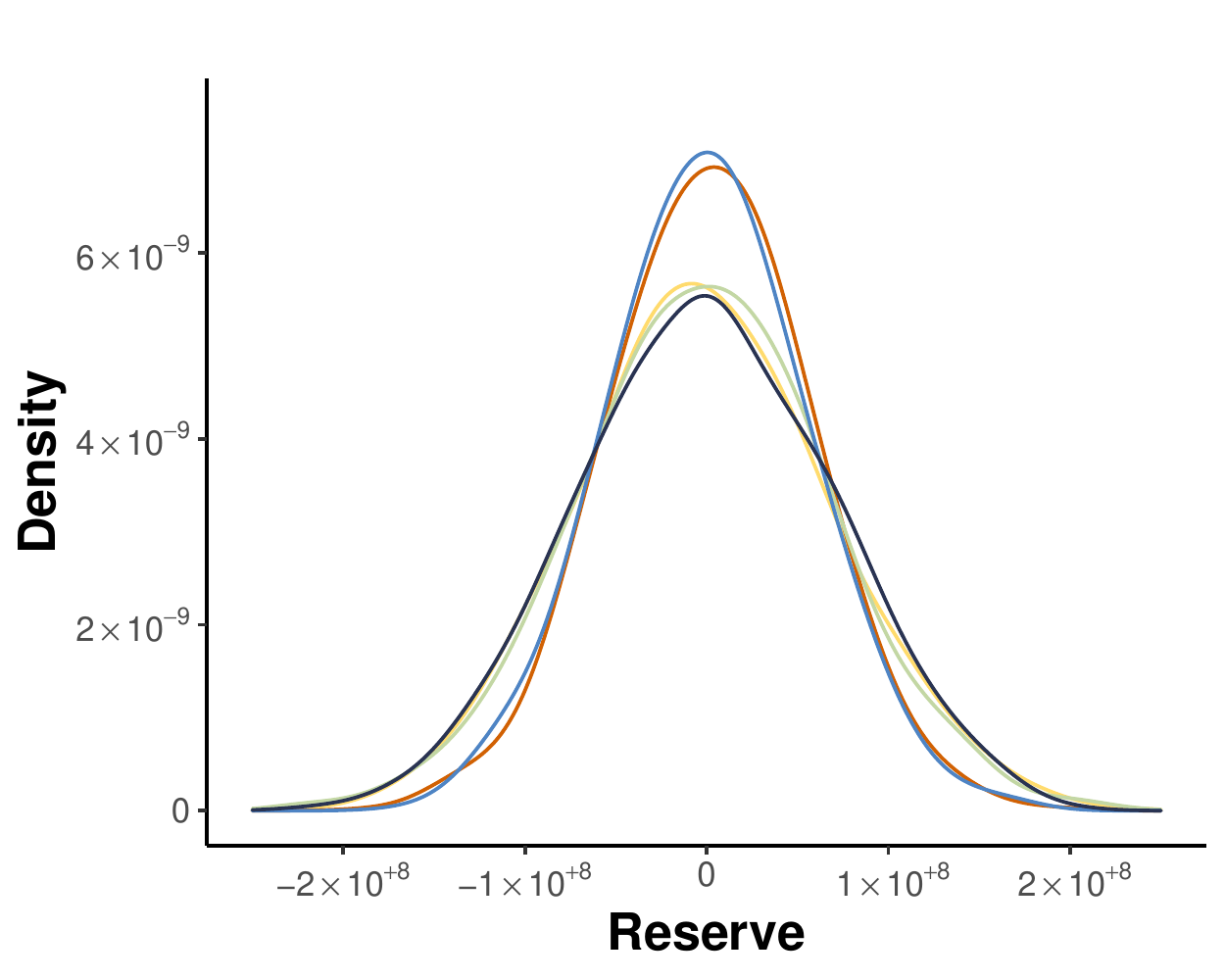}\\ 	
	\includegraphics[width=3.5cm,	height=3.5cm]{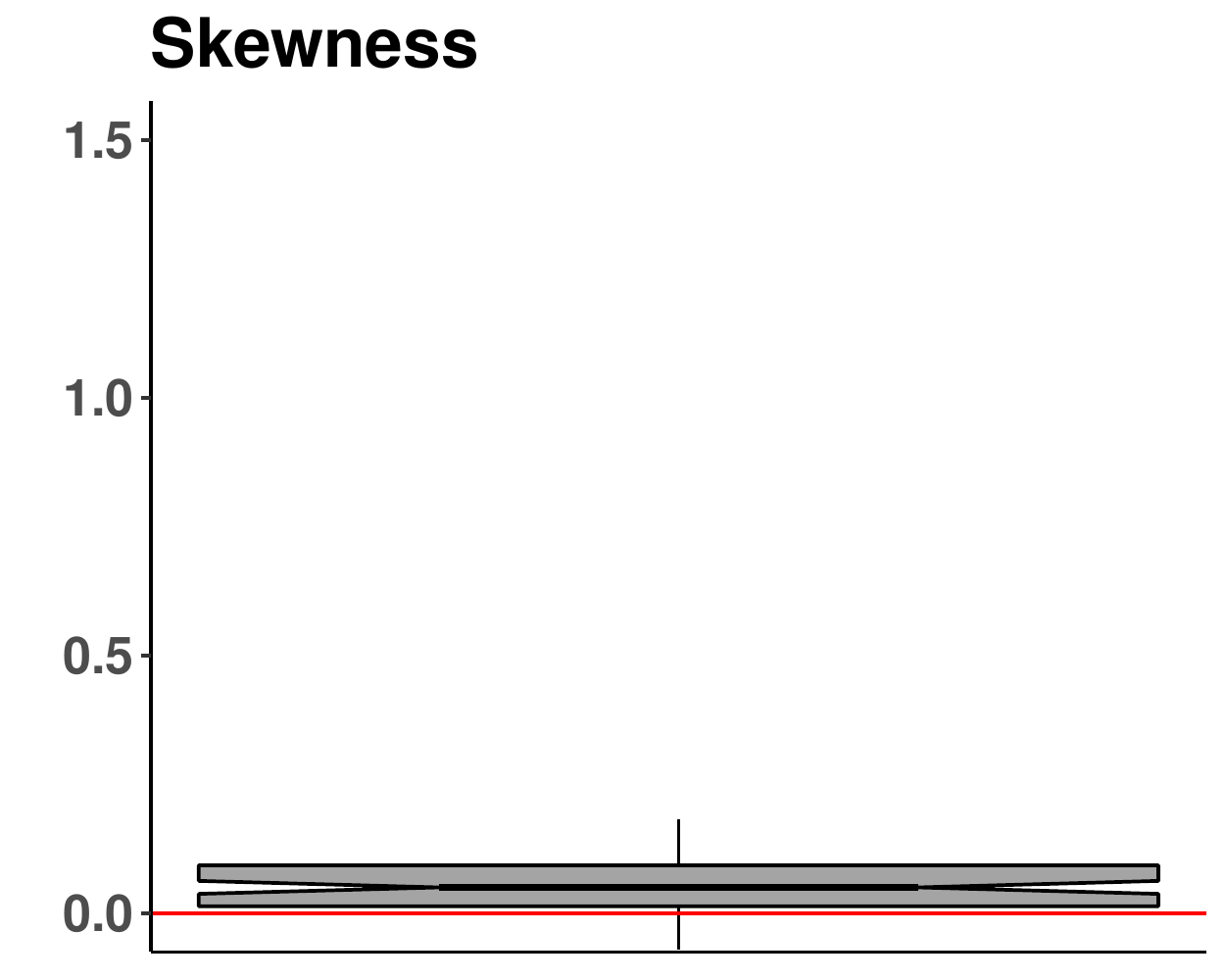}
	\includegraphics[width=3.5cm,	height=3.5cm]{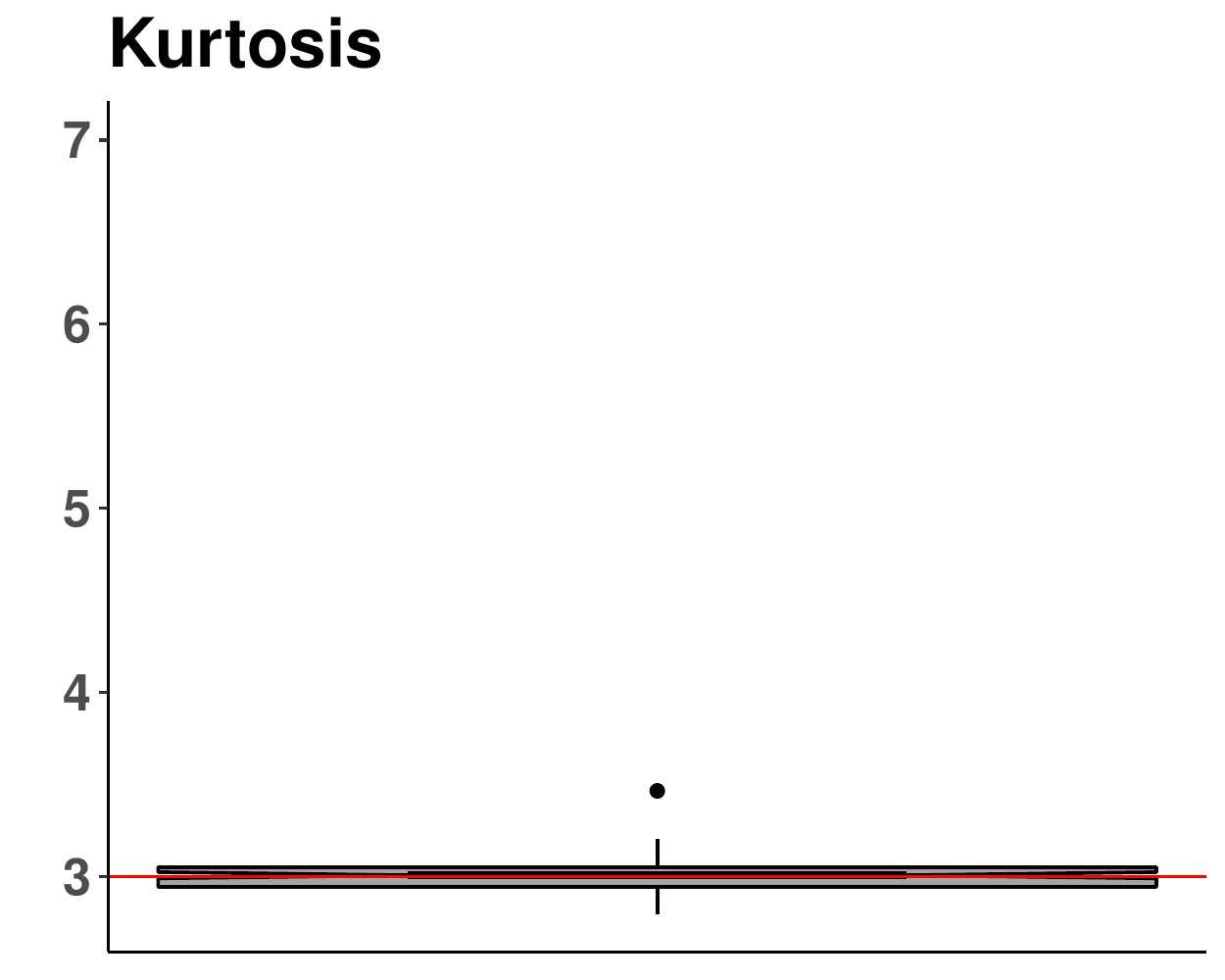}
	\includegraphics[width=5cm,	height=3.5cm]{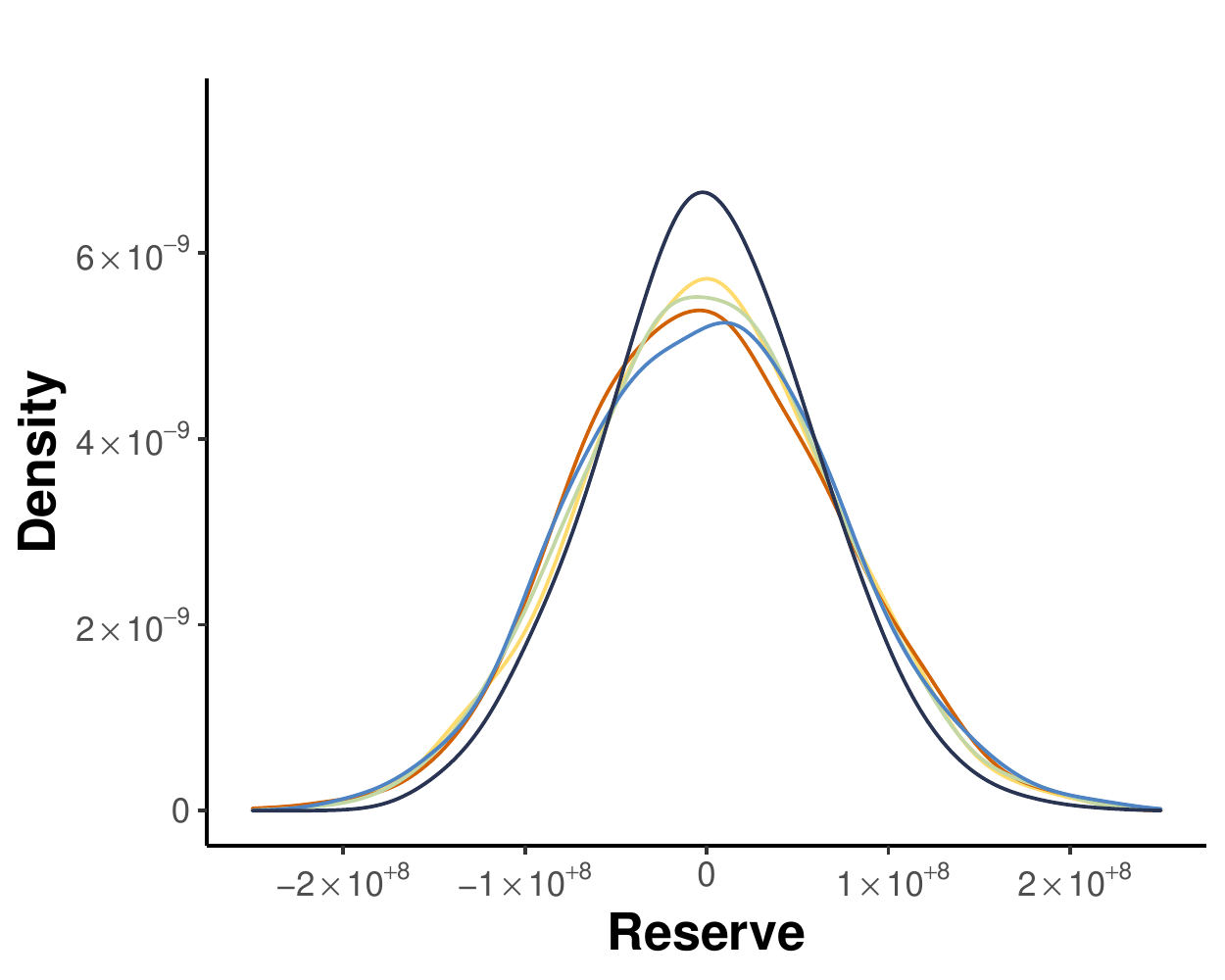}\\
	\includegraphics[width=3.5cm,	height=3.5cm]{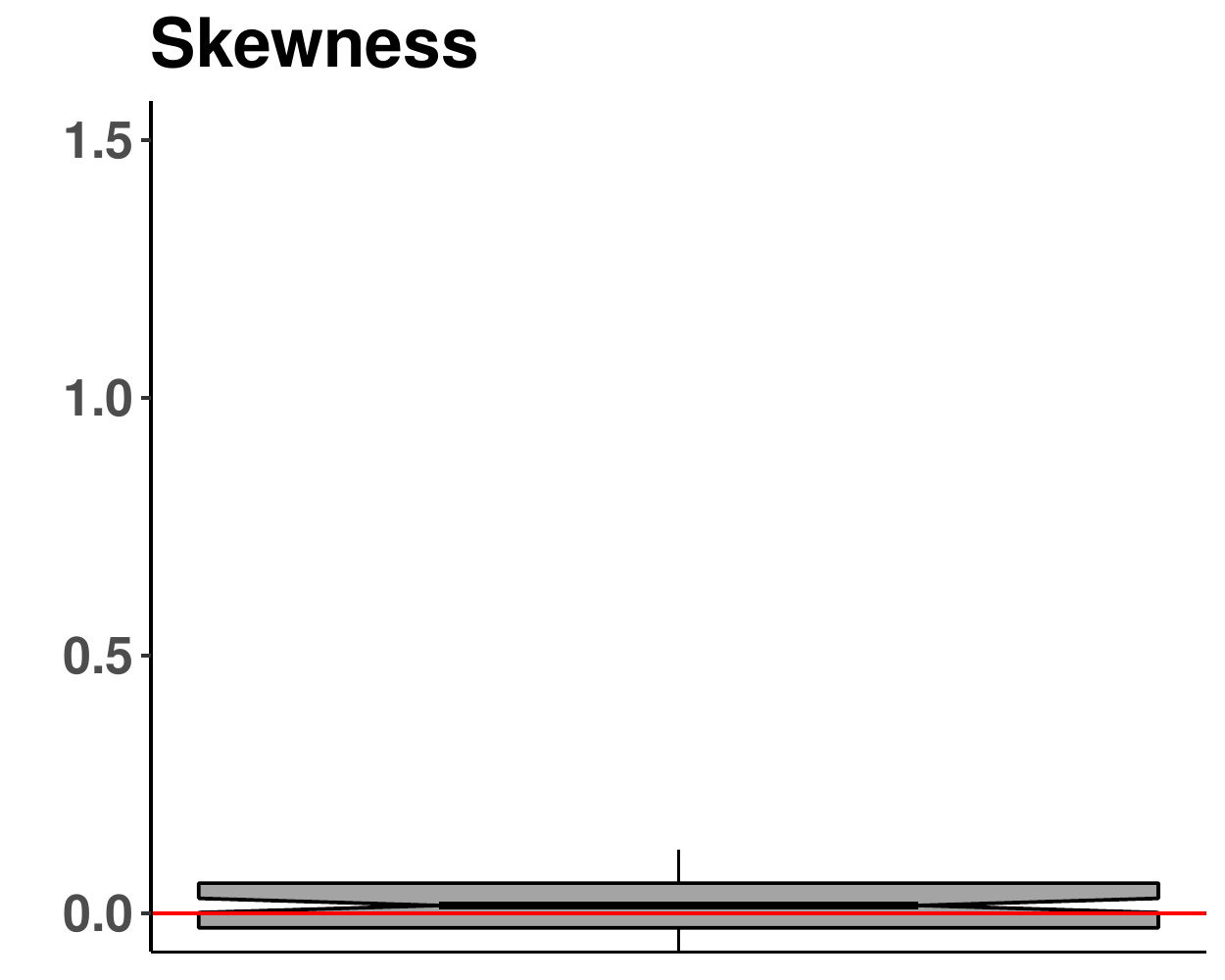}
	\includegraphics[width=3.5cm,	height=3.5cm]{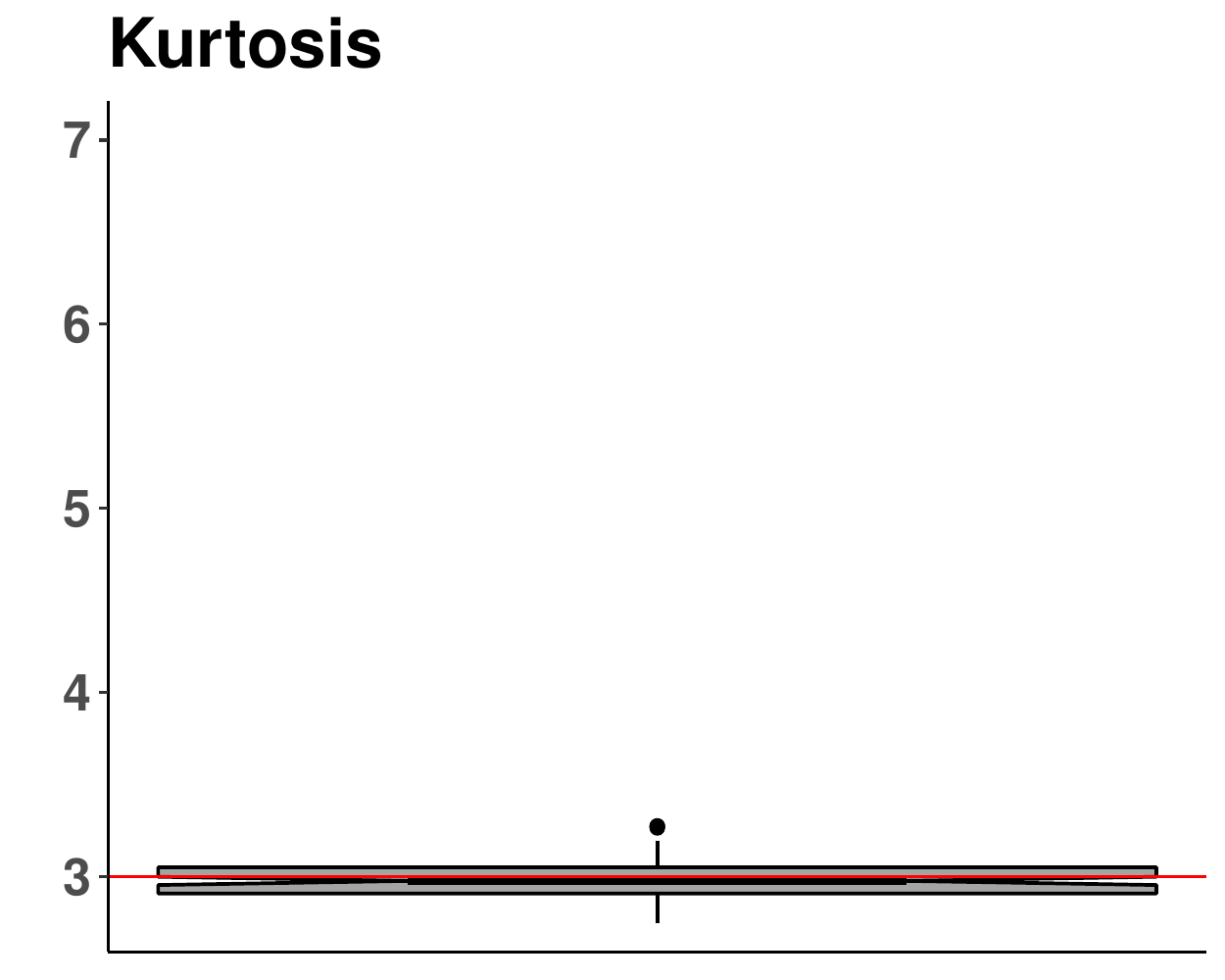}
	\includegraphics[width=5cm,	height=3.5cm]{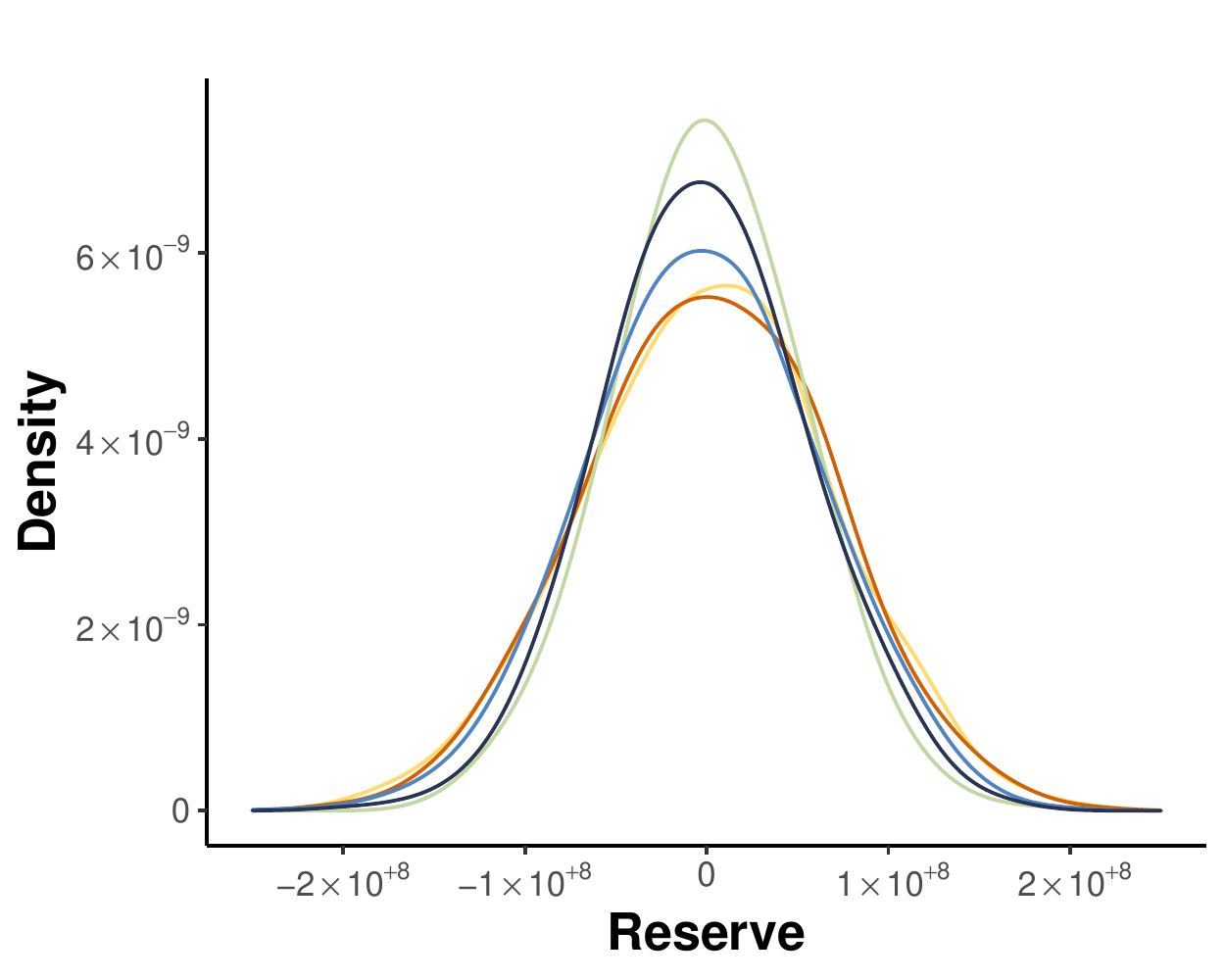}	
	\caption{Boxplots of skewness and kurtosis as well as five arbitrarily selected density plots for the simulated Mack type bootstrap conditional distribution of $(R_{I,n}^*-\widehat{R}_{I,n})_1$ given $\mathcal{Q}_{I,n}^{*}=\mathcal{Q}_{I,n}$ and $\mathcal{D}_{I,n}$ for $n=10$ and $I=10$ for the setup of a), where $F_{i,j}^*$ follows a (conditional) gamma (top), log-normal (center) and truncated normal distribution (bottom).}\label{ProcessUncertainityBS_30a}
\end{figure}

\begin{figure}
	\includegraphics[width=3.5cm,	height=3.5cm]{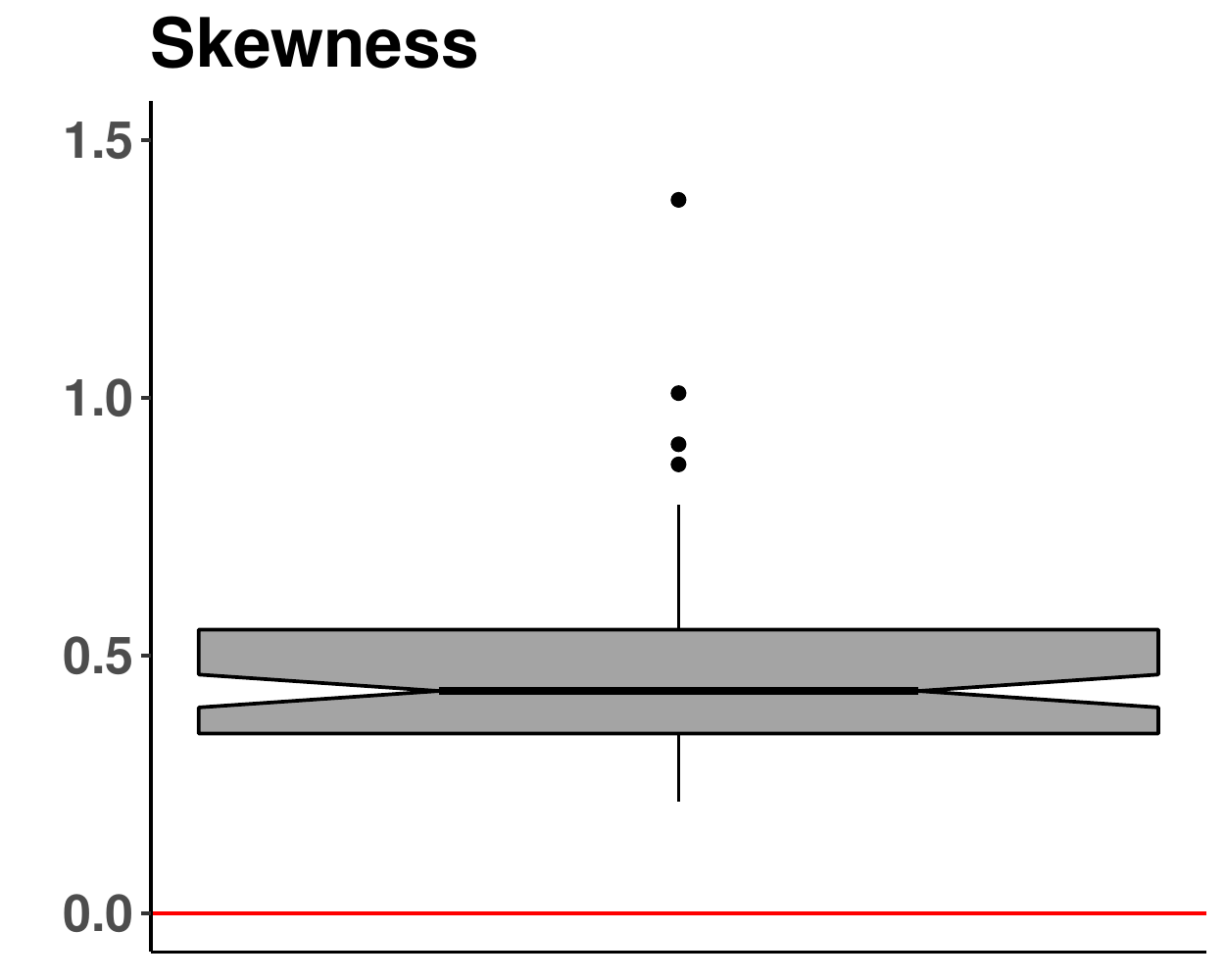}
	\includegraphics[width=3.5cm,	height=3.5cm]{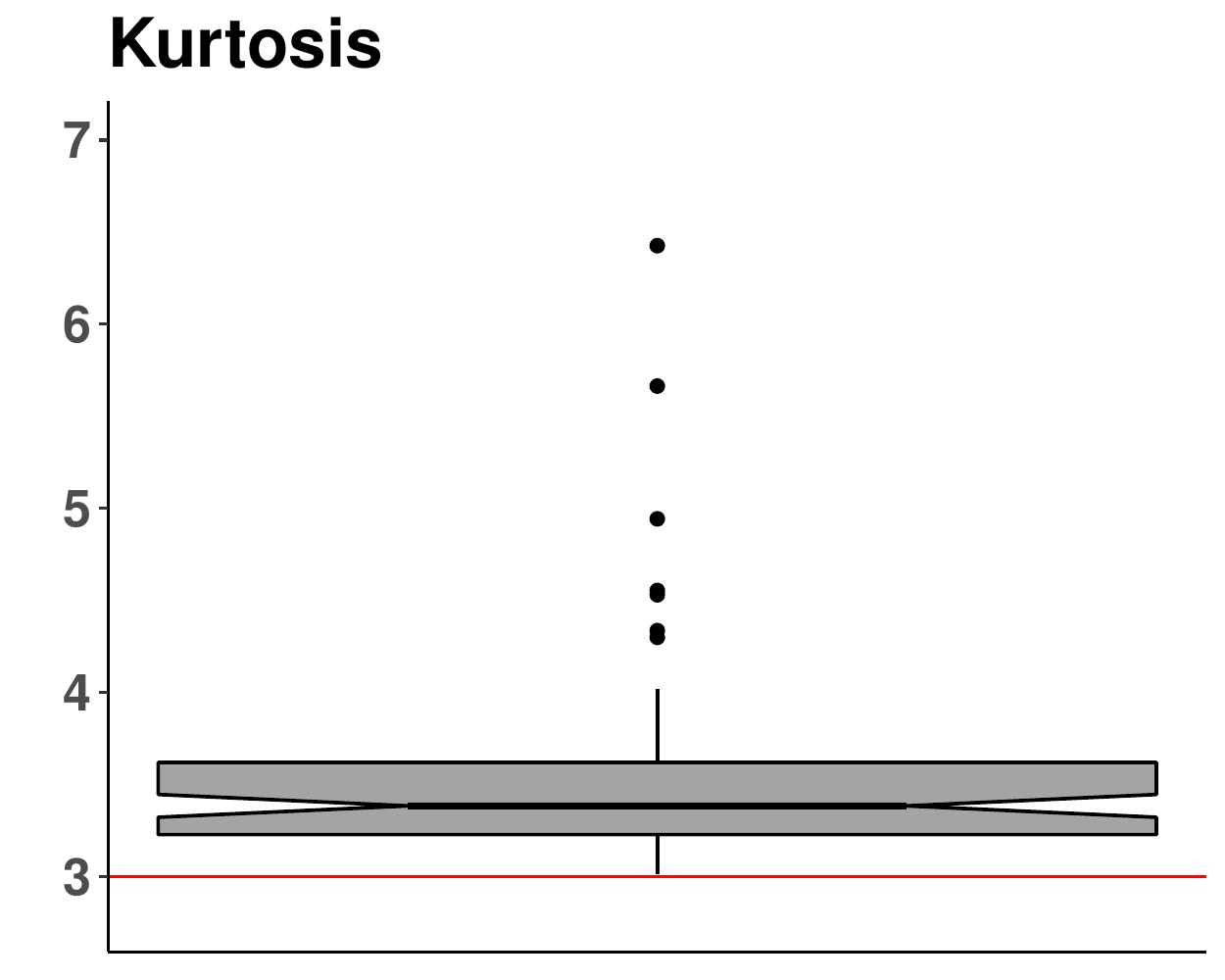}
	\includegraphics[width=5cm,	height=3.5cm]{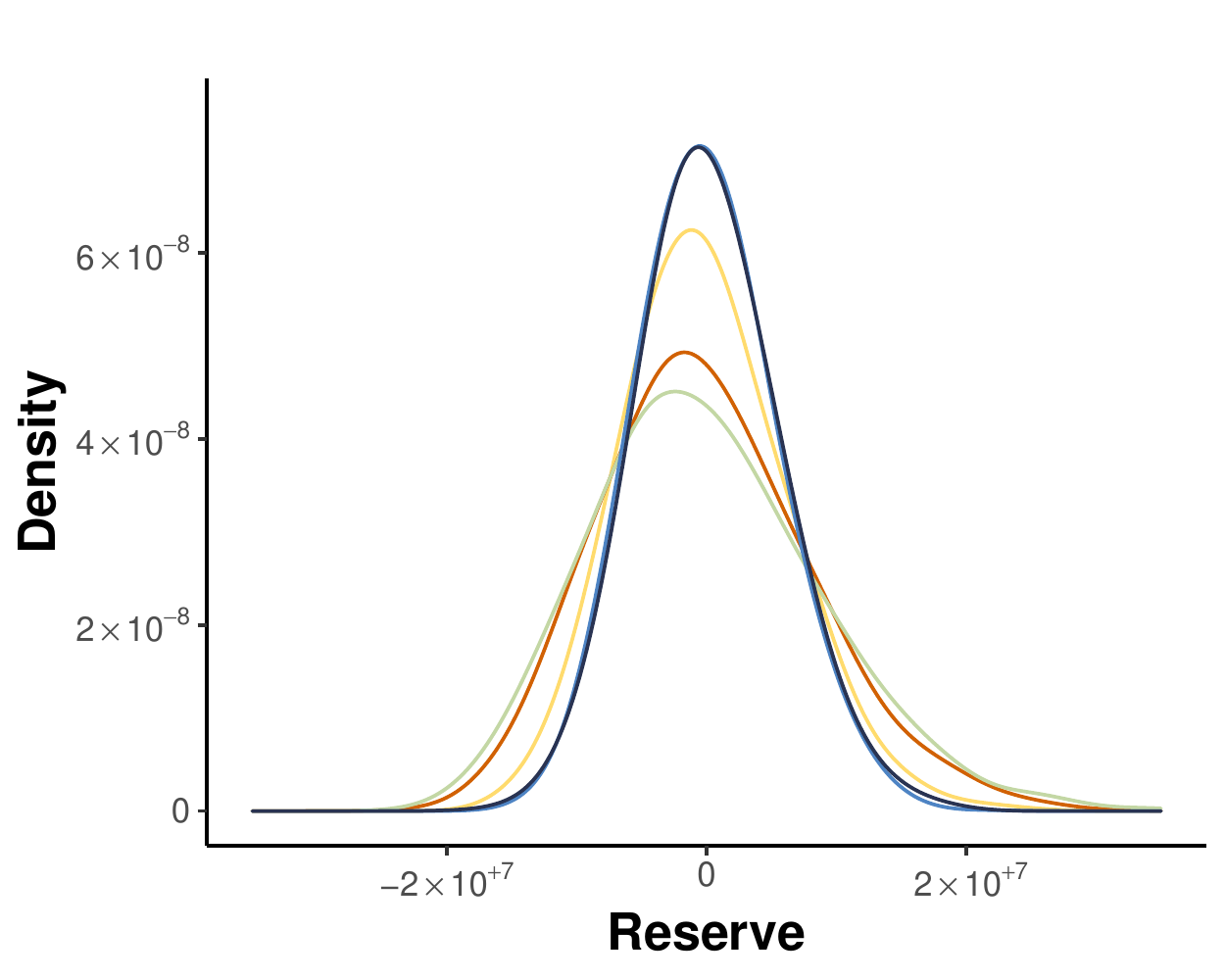}\\ 	
	\includegraphics[width=3.5cm,	height=3.5cm]{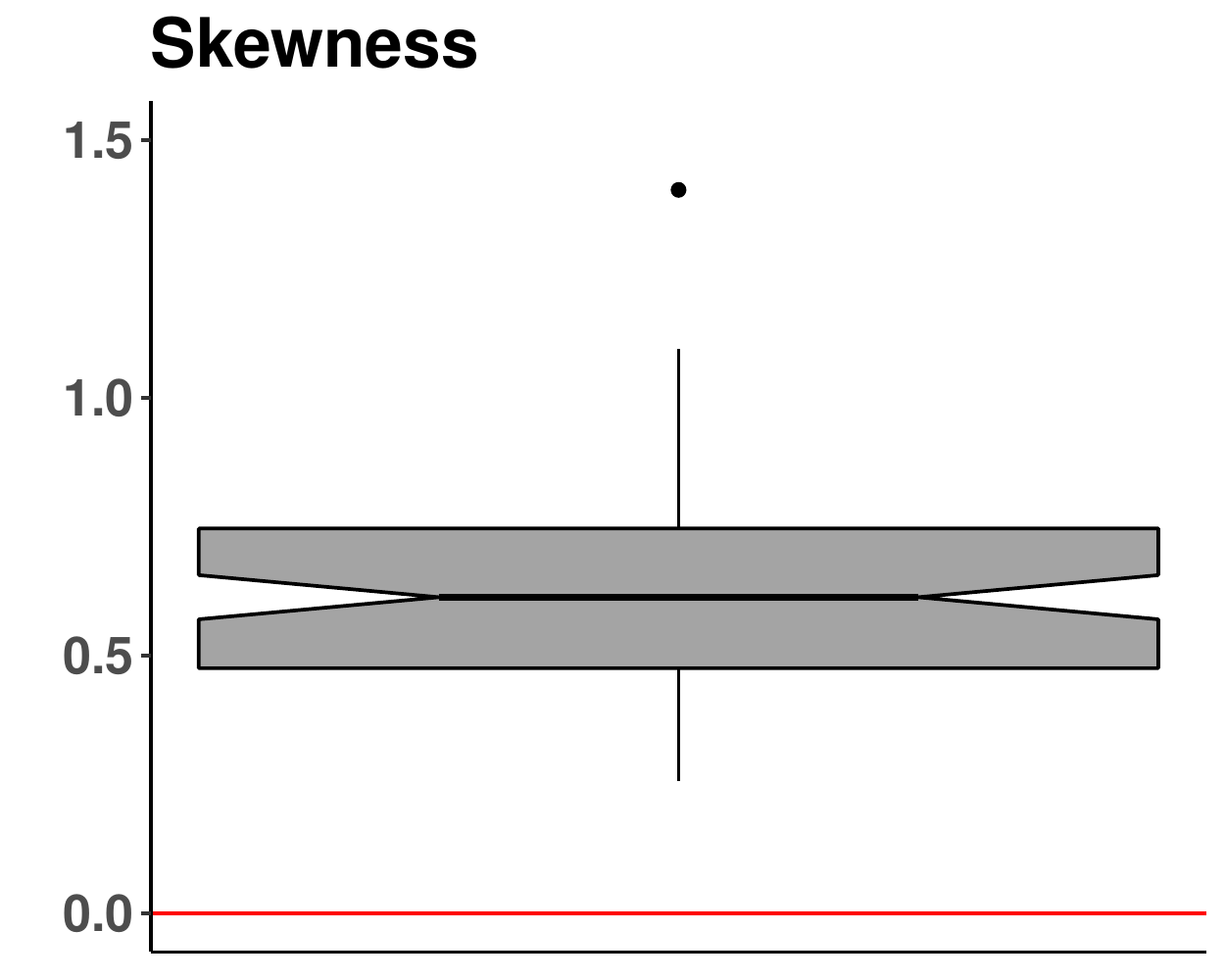}
	\includegraphics[width=3.5cm,	height=3.5cm]{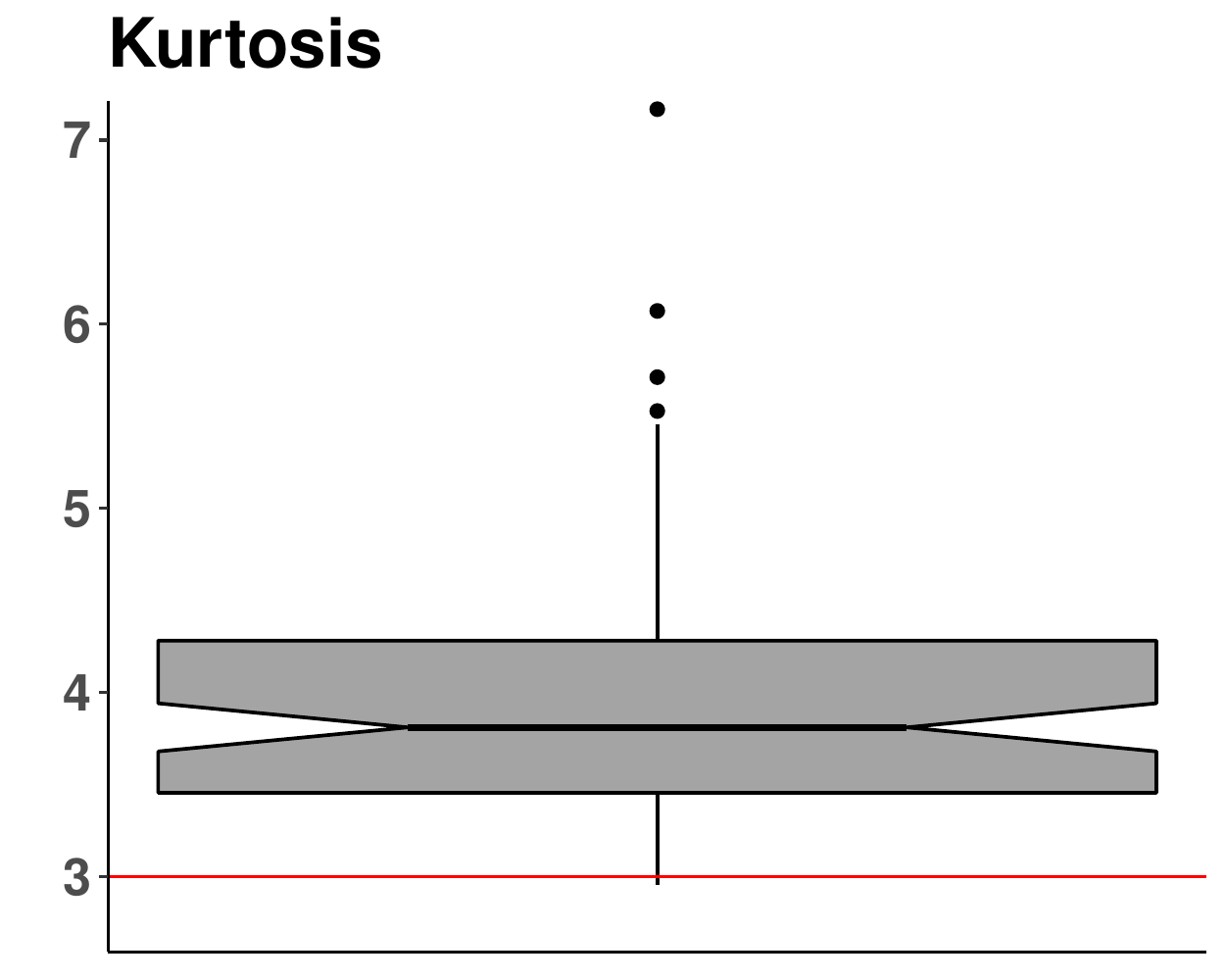}
	\includegraphics[width=5cm,	height=3.5cm]{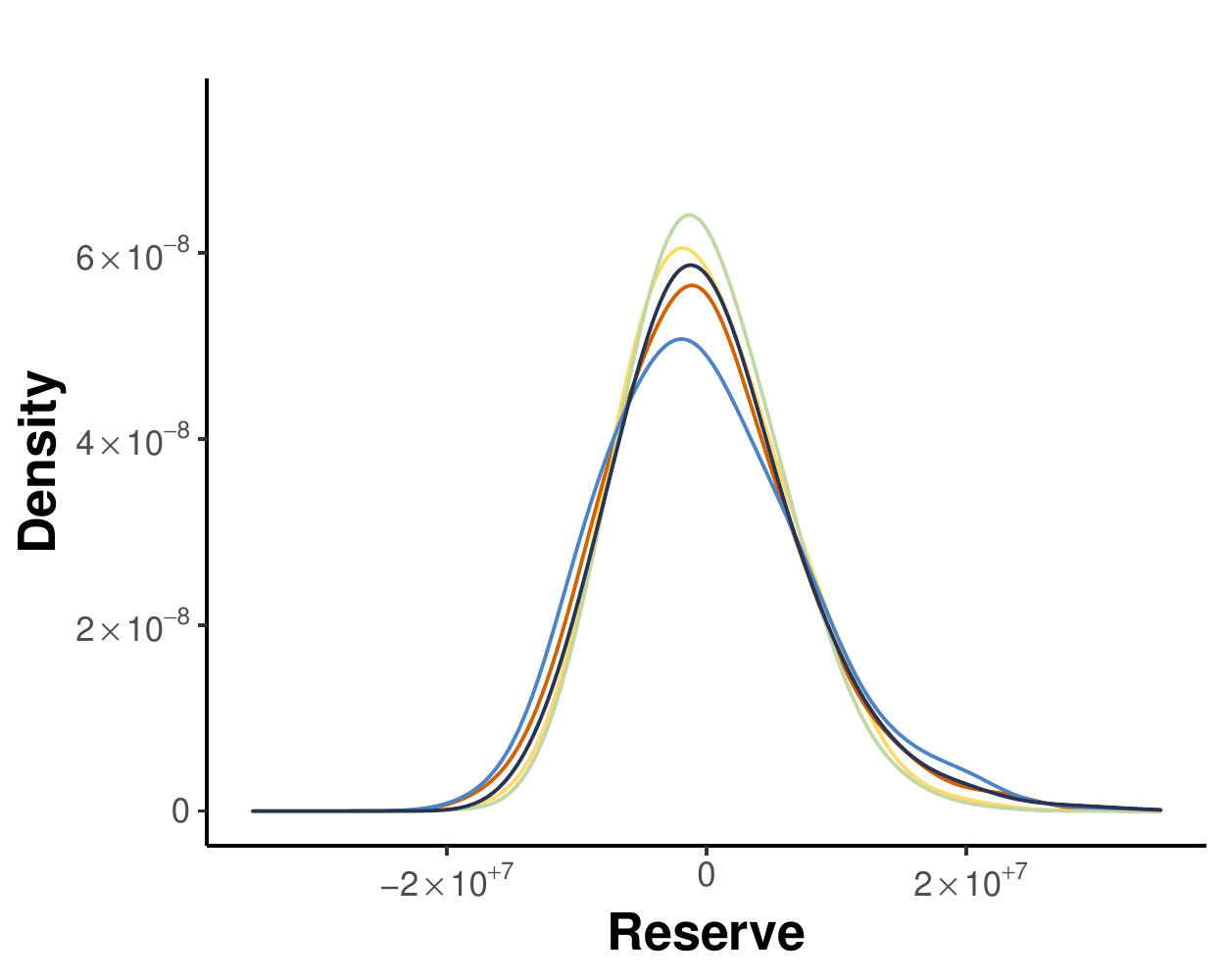}\\
	\includegraphics[width=3.5cm,	height=3.5cm]{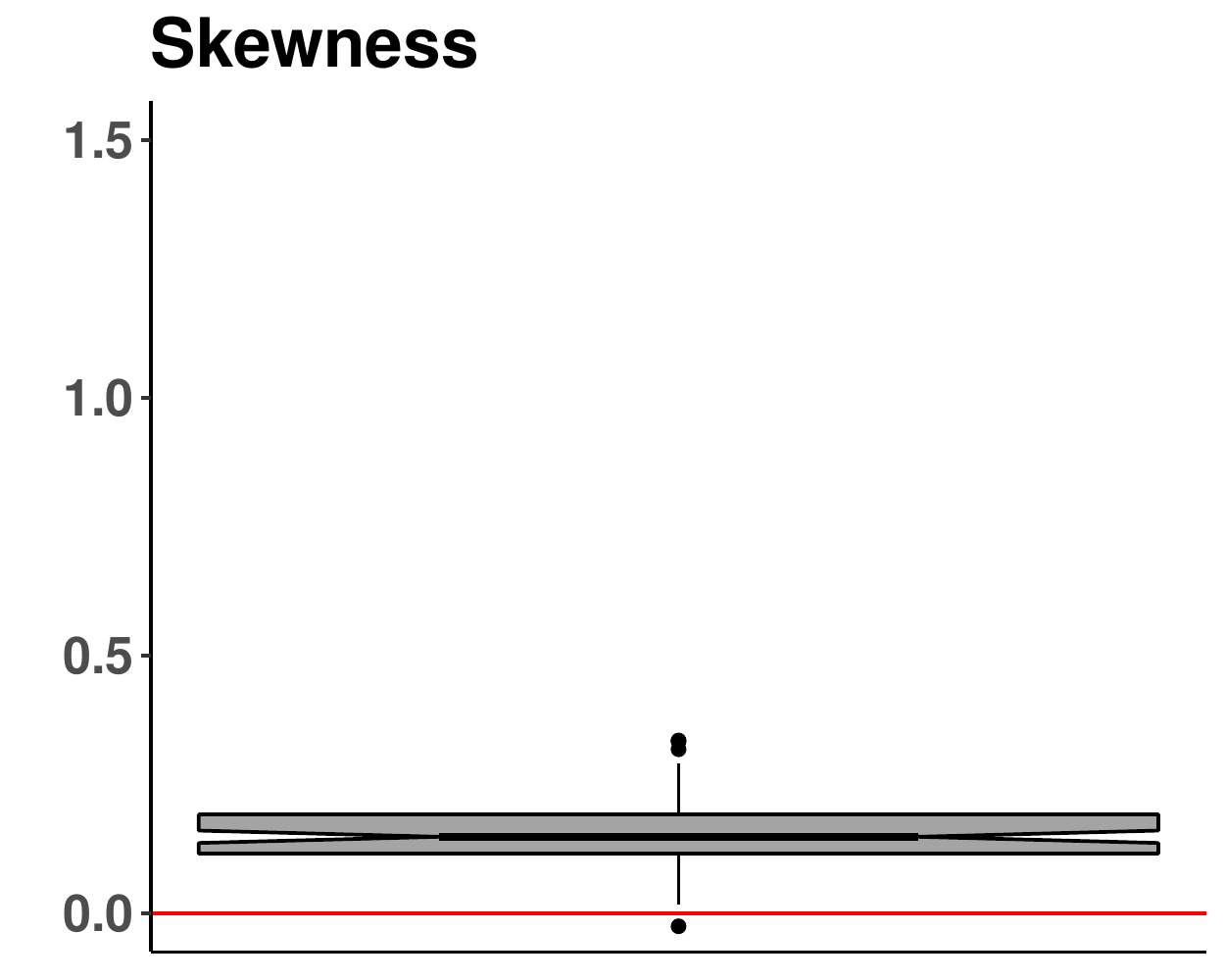}
	\includegraphics[width=3.5cm,	height=3.5cm]{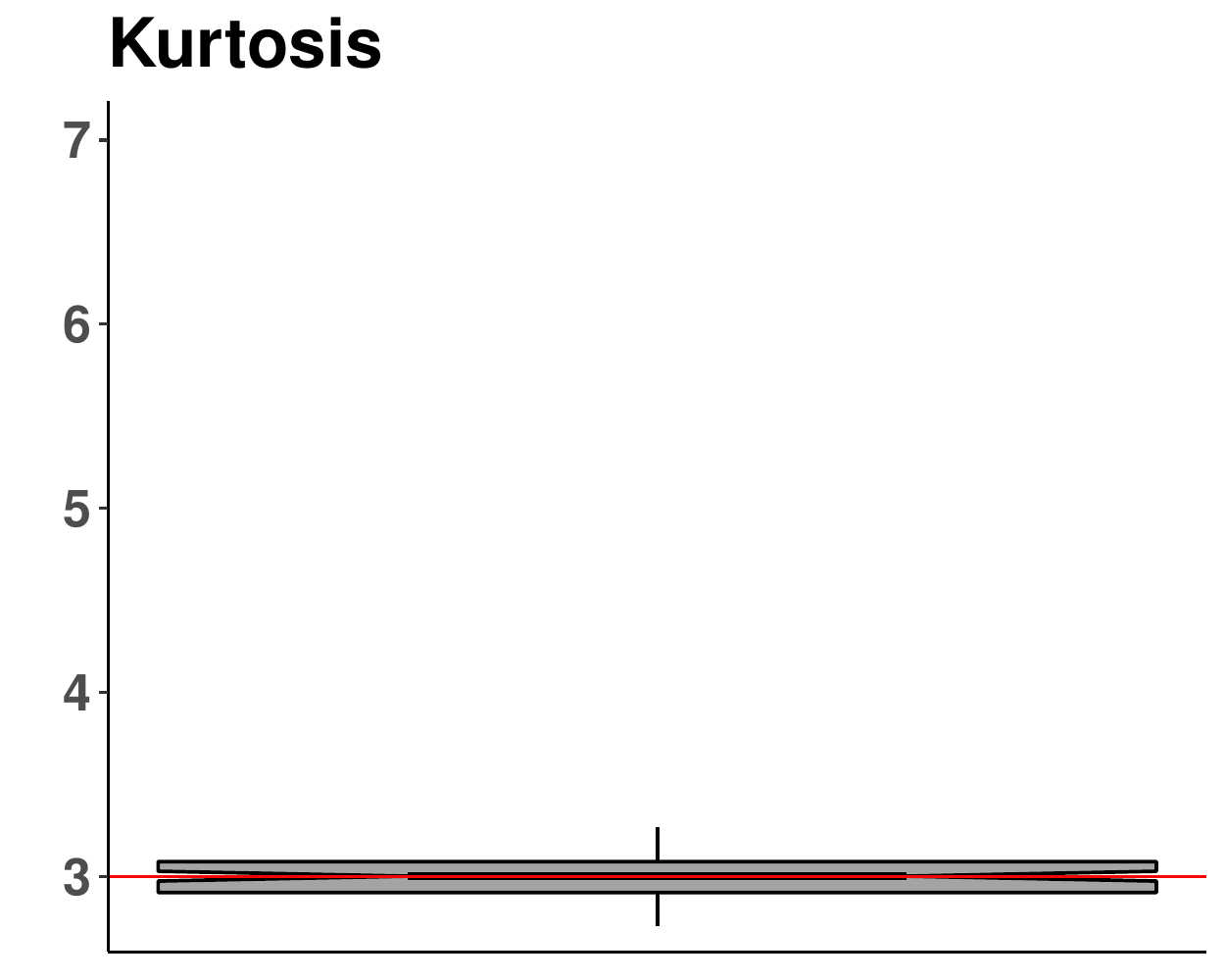}
	\includegraphics[width=5cm,	height=3.5cm]{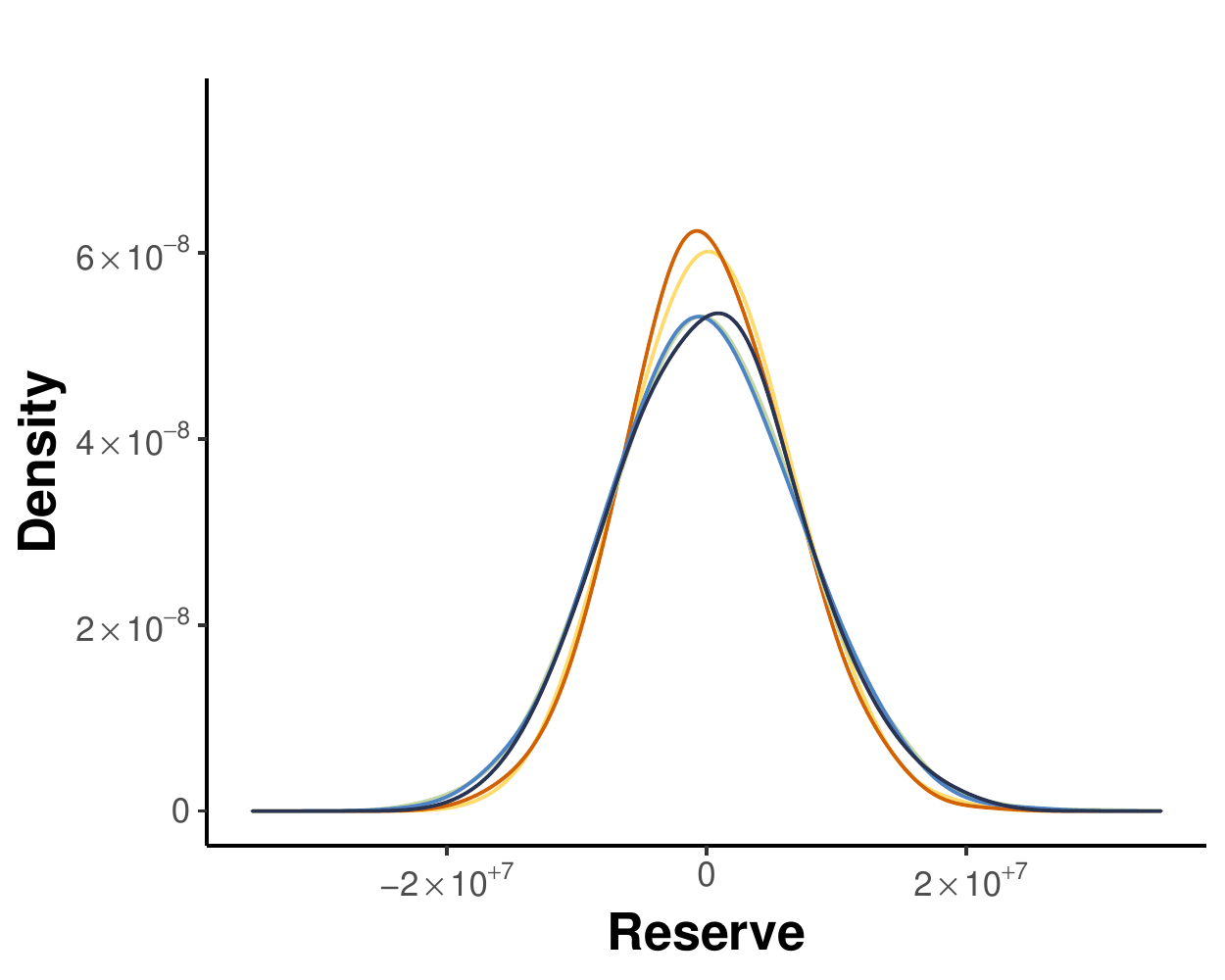}	
	\caption{Boxplots of skewness and kurtosis as well as five arbitrarily selected density plots for the simulated Mack type bootstrap conditional distribution of $(R_{I,n}^*-\widehat{R}_{I,n})_1$ given $\mathcal{Q}_{I,n}^{*}=\mathcal{Q}_{I,n}$ and $\mathcal{D}_{I,n}$ for $n=10$ and $I=10$ for the setup of b), where $F_{i,j}^*$ follows a (conditional) gamma (top), log-normal (center) and truncated normal distribution (bottom).}\label{ProcessUncertainityBS_30b}
\end{figure}

As a summary, we show boxplots of skewness and kurtosis as well as arbitrarily chosen density plots for both settings a) and b) in Figures \ref{ProcessUncertainityBS_30a} and \ref{ProcessUncertainityBS_30b} for $I=10$ and $n=10$ and for all three different distribution assumptions in (i), (ii), (iii) generated by the original Mack bootstrap. The results do not change for the alternative Mack bootstrap. 

Next we compare the bootstrap distribution of  $(R_{I,n}^*-\widehat R_{I,n})^{(m)}_1|(\mathcal{Q}_{I,n}^{(m)*}=\mathcal{Q}_{I,n}^{(m)}, \mathcal{D}_{I,n}^{(m)})$ and  $(R_{I,n}^+-\widehat R_{I,n}^+)^{(m)}_1| (\mathcal{Q}_{I,n}^{(m)*}=\mathcal{Q}_{I,n}^{(m)}, \mathcal{D}_{I,n}^{(m)})$, respectively, to the distribution $(R_{I,n}-\widehat R_{I,n})_1^{(m)}|\mathcal{Q}_{I,n}^{(m)}$ obtained by Monte Carlo Simulation for $m=1, \dots, 500$. We apply the Kolmogorov-Smirnov test of level $\alpha=5\%$.

Tables \ref{tab_KS_part1_pred_ab} and \ref{KS_part1_pred_a_betest_app} summarize the results for setup a) and b), respectively, for the original Mack and alternative Mack bootstrap. 
The results of the original Mack and the alternative Mack bootstrap do not differ.

In general, for increasing $n$ the percentages of fail to reject the null hypothesis increase. If we choose the true underlying distribution, we fail to reject the null hypothesis more frequently than if we choose the wrong distribution. For setup b) it is more important to choose the true underlying distribution compared to setup b). If the underlying distribution of the individual development factors is skewed, the chosen distribution for $\F^*$ and $\F^+$, respectively, for the lower triangle should be skewed. For example, if we choose a gamma distribution for $\F^*$ instead of a log-normal distribution as true distributional family of $\F^*$,  the percentage to fail to reject the null hypothesis is higher compared to if we choose a truncated normal distribution, e.g., for $n=40$, we get that 69\% out of $M=500$ fail to reject the null hypothesis assuming a gamma distribution compared to 31\% using a truncated normal distribution (cf. Table \ref{KS_part1_pred_a_betest_app}). 
Also, if the true underlying distribution is a truncated normal distribution and we choose a gamma, then 50\% out of $M=500$ fail to reject the null hypothesis or to assume a log-normal distribution, then 39\% and if we choose the true underlying distribution, then 84\% for $n=40$ (cf. Table \ref{KS_part1_pred_a_betest_app}).

For setup a) the effect of choosing the wrong distribution is not as high as for b). We can explain this with the property of the gamma and the log-normal distribution. Both tend to 'lose' their skewness and excess of kurtosis for $\frac{C^*_{i,j}}{\widehat \sigma_{j,n}^2}$ growing large in this parameter setting.

Tables \ref{tab_ks_predroot_b} and \ref{tab_mse_predroot_b} here contain simulations results according to Tables 3 and 4 in Section 7, but for Setup b).

\begin{small}
\begin{table}[h]
\begin{tabular}{|c|c|p{1cm}p{1cm}|p{1cm}p{1cm}|p{1cm}p{1cm}|}
\hline
\multicolumn{2}{|c|}{\textbf{chosen   distribution } }                                    & \multicolumn{2}{c}{\cellcolor[HTML]{D9D9D9}\textbf{gamma} }       & \multicolumn{2}{c}{\cellcolor[HTML]{D9D9D9}\textbf{log-normal}}   & \multicolumn{2}{c|}{\cellcolor[HTML]{D9D9D9}\textbf{trunc. normal}} \\
\hline
\multicolumn{1}{|l|}{\textbf{true   distribution}  }               & n  & oMB & aMB & oMB & aMB & oMB  & aMB \\
\hline
\cellcolor[HTML]{D9D9D9}                                & 0  & 0.57                        & 0.55                       & 0.45                        & 0.50                       & 0.49                         & 0.51                       \\
\cellcolor[HTML]{D9D9D9}                                & 10 & 0.66                        & 0.69                       & 0.65                        & 0.70                       & 0.58                         & 0.71                       \\
\cellcolor[HTML]{D9D9D9}                                & 20 & 0.73                        & 0.72                       & 0.72                        & 0.79                       & 0.68                         & 0.80                       \\
\cellcolor[HTML]{D9D9D9}                                & 30 & 0.75                        & 0.73                       & 0.75                        & 0.83                       & 0.72                         & 0.81                       \\
\multirow{-5}{*}{\cellcolor[HTML]{D9D9D9}\textbf{gamma}}         & 40 & 0.80                        & 0.76                       & 0.79                        & 0.87                       & 0.79                         & 0.83                       \\
\hline
\cellcolor[HTML]{D9D9D9}                                & 0  & 0.44                        & 0.47                       & 0.57                        & 0.56                       & 0.45                         & 0.47                       \\
\cellcolor[HTML]{D9D9D9}                                & 10 & 0.60                        & 0.61                       & 0.69                        & 0.68                       & 0.60                         & 0.62                       \\
\cellcolor[HTML]{D9D9D9}                                & 20 & 0.69                        & 0.73                       & 0.78                        & 0.77                       & 0.70                         & 0.64                       \\
\cellcolor[HTML]{D9D9D9}                                & 30 & 0.70                         & 0.77                       & 0.83                        & 0.80                       & 0.78                         & 0.70                       \\
\multirow{-5}{*}{\cellcolor[HTML]{D9D9D9}\textbf{log-normal}}    & 40 & 0.81                        & 0.80                        & 0.89                        & 0.85                       & 0.80                         & 0.73                       \\
\hline
\cellcolor[HTML]{D9D9D9}                                & 0  & 0.46                        & 0.49                       & 0.45                        & 0.52                       & 0.50                          & 0.48                       \\
\cellcolor[HTML]{D9D9D9}                                & 10 & 0.62                        & 0.63                       & 0.71                        & 0.67                       & 0.59                         & 0.57                       \\
\cellcolor[HTML]{D9D9D9}                                & 20 & 0.67                        & 0.68                       & 0.78                        & 0.75                       & 0.71                         & 0.71                       \\
\cellcolor[HTML]{D9D9D9}                                & 30 & 0.72                        & 0.71                       & 0.80                        & 0.80                       & 0.75                         & 0.73                       \\
\multirow{-5}{*}{\cellcolor[HTML]{D9D9D9}\textbf{trunc. normal}} & 40 & 0.76                        & 0.73                       & 0.82                        & 0.82                       & 0.84                         & 0.85            \\
\hline
\end{tabular}
		\caption{Process Uncertainty: Percentages of failed rejections for Kolmogorov-Smirnov tests of level $\alpha=5\%$ for the null hypotheses $H_0^*$, $H_0^+$ and $H_0^{++}$,
		%$\mathcal{L^*}((R_{I,n}^*-\widehat R_{I,n})_1|\mathcal{Q}_{I,n}^*=\mathcal{Q}_{I,n})=\mathcal{L}((R_{I,n}-\widehat R_{I,n})_1|\mathcal{Q}_{I,n})$ and $\mathcal{L^*}((R_{I,n}^+-\widehat R_{I,n}^+)_1|\mathcal{Q}_{I,n}^+=\mathcal{Q}_{I,n})=\mathcal{L}((R_{I,n}-\widehat R_{I,n})_1|\mathcal{Q}_{I,n})$, 
	respectively, for the original Mack bootstrap (oMB) and the alternative Mack bootstrap (aMB) for different parametric families of distributions of $\F^*$ for $i+j\geq I$, for $I=10$ and different $n$ in Setup a).}	\label{tab_KS_part1_pred_ab}%
	\end{table}
\end{small}

\begin{small}
\begin{table}[h]
\begin{tabular}{|c|c|p{1cm}p{1cm}|p{1cm}p{1cm}|p{1cm}p{1cm}|}
\hline
\multicolumn{2}{|c|}{\textbf{chosen   distribution } }                                    & \multicolumn{2}{c}{\cellcolor[HTML]{D9D9D9}\textbf{gamma} }       & \multicolumn{2}{c}{\cellcolor[HTML]{D9D9D9}\textbf{log-normal}}   & \multicolumn{2}{c|}{\cellcolor[HTML]{D9D9D9}\textbf{trunc. normal}} \\
\hline
\multicolumn{1}{|l|}{\textbf{true   distribution}  }               & n  & oMB & aMB & oMB & aMB & oMB  & aMB \\
\hline
\cellcolor[HTML]{D9D9D9}                                & 0  & 0.52                        & 0.48                       & 0.40                        & 0.35                       & 0.34                         & 0.34                       \\
\cellcolor[HTML]{D9D9D9}                                & 10 & 0.66                        & 0.63                       & 0.44                        & 0.53                       & 0.49                         & 0.47                       \\
\cellcolor[HTML]{D9D9D9}                                & 20 & 0.77                        & 0.72                       & 0.66                        & 0.66                       & 0.58                         & 0.51                       \\
\cellcolor[HTML]{D9D9D9}                                & 30 & 0.80                        & 0.76                       & 0.70                        & 0.68                       & 0.60                         & 0.53                       \\
\multirow{-5}{*}{\cellcolor[HTML]{D9D9D9}\textbf{gamma}}         & 40 & 0.83                        & 0.80                       & 0.71                        & 0.75                       & 0.61                         & 0.57                       \\
 \hline
\cellcolor[HTML]{D9D9D9}                                & 0  & 0.32                        & 0.33                       & 0.44                        & 0.41                       & 0.18                         & 0.18                       \\
\cellcolor[HTML]{D9D9D9}                                & 10 & 0.50                        & 0.54                       & 0.55                        & 0.50                       & 0.21                         & 0.22                       \\
\cellcolor[HTML]{D9D9D9}                                & 20 & 0.60                        & 0.55                       & 0.65                        & 0.63                       & 0.23                         & 0.25                       \\
\cellcolor[HTML]{D9D9D9}                                & 30 & 0.61                        & 0.60                       & 0.68                        & 0.72                       & 0.29                         & 0.28                       \\
\multirow{-5}{*}{\cellcolor[HTML]{D9D9D9}\textbf{log-normal}}    & 40 & 0.69                        & 0.67                       & 0.78                        & 0.75                       & 0.31                         & 0.30                       \\
\hline
\cellcolor[HTML]{D9D9D9}                                & 0  & 0.30                        & 0.33                       & 0.21                        & 0.26                       & 0.49                         & 0.51                       \\
\cellcolor[HTML]{D9D9D9}                                & 10 & 0.40                        & 0.45                       & 0.29                        & 0.33                       & 0.64                         & 0.60                       \\
\cellcolor[HTML]{D9D9D9}                                & 20 & 0.43                        & 0.49                       & 0.39                        & 0.35                       & 0.78                         & 0.78                       \\
\cellcolor[HTML]{D9D9D9}                                & 30 & 0.45                        & 0.54                       & 0.40                        & 0.38                       & 0.80                         & 0.85                       \\
\multirow{-5}{*}{\cellcolor[HTML]{D9D9D9}\textbf{trunc. normal}} & 40 & 0.50                        & 0.59                       & 0.42                        & 0.45                       & 0.84                         & 0.87       \\
\hline
\end{tabular}
		\caption{Process Uncertainty: Percentages of failed rejections for Kolmogorov-Smirnov tests of level $\alpha=5\%$ for the null hypotheses $H_0^*$, $H_0^+$ and $H_0^{++}$,
		%$\mathcal{L^*}((R_{I,n}^*-\widehat R_{I,n})_1|\mathcal{Q}_{I,n}^*=\mathcal{Q}_{I,n})=\mathcal{L}((R_{I,n}-\widehat R_{I,n})_1|\mathcal{Q}_{I,n})$ and $\mathcal{L^*}((R_{I,n}^+-\widehat R_{I,n}^+)_1|\mathcal{Q}_{I,n}^+=\mathcal{Q}_{I,n})=\mathcal{L}((R_{I,n}-\widehat R_{I,n})_1|\mathcal{Q}_{I,n})$, 
respectively, for the original Mack bootstrap (oMB) and the alternative Mack bootstrap (aMB) for different parametric families of distributions of $\F^*$ for $i+j\geq I$, for $I=10$ and different $n$ in setup b).}\label{KS_part1_pred_a_betest_app}%
	\end{table}
\end{small}

%\tr{ACHTUNG: Hier kommen die ehemaligen Table 4 und 6}

\begin{small}
\begin{table}[t]
\begin{tabular}{|c|c|lll|lll|lll|}
\hline
\multicolumn{2}{|c|}{\textbf{chosen distribution}}                      & \multicolumn{3}{c}{\cellcolor[HTML]{D9D9D9}\textbf{gamma}} & \multicolumn{3}{c}{\cellcolor[HTML]{D9D9D9}\textbf{log-normal}} & \multicolumn{3}{c|}{\cellcolor[HTML]{D9D9D9}\textbf{trunc. normal}} \\ \hline
\multicolumn{1}{|l|}{\textbf{true distribution}}                  & n  & oMB                & aMB                & iMB                & oMB                 & aMB                  & iMB                  & oMB                  & aMB                   & iMB                   \\\hline
\cellcolor[HTML]{D9D9D9}                                         & 0  & 0.10               & 0.11              & 0.09              & 0.11                & 0.14                & 0.12                & 0.09                 & 0.10                 & 0.10                 \\
\cellcolor[HTML]{D9D9D9}                                         & 10 & 0.19               & 0.21              & 0.18              & 0.23                & 0.25                & 0.23                & 0.16                 & 0.17                 & 0.17                 \\
\cellcolor[HTML]{D9D9D9}                                         & 20 & 0.28               & 0.35              & 0.32              & 0.36                & 0.40                & 0.37                & 0.21                 & 0.22                 & 0.22                 \\
\cellcolor[HTML]{D9D9D9}                                         & 30 & 0.38               & 0.42              & 0.40              & 0.40                & 0.44                & 0.39                & 0.28                 & 0.28                 & 0.28                 \\
\multirow{-5}{*}{\cellcolor[HTML]{D9D9D9}\textbf{gamma}}         & 40 & 0.52               & 0.56              & 0.51              & 0.51                & 0.53                & 0.50                & 0.41                 & 0.41                 & 0.41                 \\\hline
\cellcolor[HTML]{D9D9D9}                                         & 0  & 0.07               & 0.10              & 0.08              & 0.09                & 0.10                & 0.08                & 0.07                 & 0.08                 & 0.08                 \\
\cellcolor[HTML]{D9D9D9}                                         & 10 & 0.16               & 0.19              & 0.17              & 0.20                & 0.22                & 0.20                & 0.16                 & 0.15                 & 0.15                 \\
\cellcolor[HTML]{D9D9D9}                                         & 20 & 0.31               & 0.33              & 0.30              & 0.22                & 0.26                & 0.23                & 0.22                 & 0.22                 & 0.22                 \\
\cellcolor[HTML]{D9D9D9}                                         & 30 & 0.34               & 0.38              & 0.33              & 0.30                & 0.33                & 0.30                & 0.27                 & 0.27                 & 0.27                 \\
\multirow{-5}{*}{\cellcolor[HTML]{D9D9D9}\textbf{log-normal}}    & 40 & 0.39               & 0.42              & 0.39              & 0.45                & 0.47                & 0.44                & 0.29                 & 0.28                 & 0.28                 \\\hline
\cellcolor[HTML]{D9D9D9}                                         & 0  & 0.09               & 0.11              & 0.08              & 0.15                & 0.17                & 0.16                & 0.24                 & 0.23                 & 0.23                 \\
\cellcolor[HTML]{D9D9D9}                                         & 10 & 0.16               & 0.18              & 0.15              & 0.21                & 0.24                & 0.20                & 0.34                 & 0.34                 & 0.34                 \\
\cellcolor[HTML]{D9D9D9}                                         & 20 & 0.28               & 0.30              & 0.27              & 0.26                & 0.30                & 0.27                & 0.41                 & 0.42                 & 0.42                 \\
\cellcolor[HTML]{D9D9D9}                                         & 30 & 0.36               & 0.39              & 0.35              & 0.31                & 0.35                & 0.32                & 0.56                 & 0.55                 & 0.55                 \\
\multirow{-5}{*}{\cellcolor[HTML]{D9D9D9}\textbf{trunc. normal}} & 40 & 0.43               & 0.45              & 0.42              & 0.36                & 0.39                & 0.36                & 0.61                 & 0.62                 & 0.60      \\\hline          
\end{tabular}
\caption{Percentages of failed rejections for Kolmogorov-Smirnov tests of level $\alpha=5\%$ for the null hypotheses $H_0^*$, $H_0^+$ and $H_0^{++}$,
 %$\mathcal{L^*}((R_{I,n}^*-\widehat R_{I,n})|\mathcal{Q}_{I,n}^*=\mathcal{Q}_{I,n})=\mathcal{L}((R_{I,n}-\widehat R_{I,n})_1|\mathcal{Q}_{I,n})$ and $\mathcal{L^*}((R_{I,n}^+-\widehat R_{I,n}^+)|\mathcal{Q}_{I,n}^+=\mathcal{Q}_{I,n})=\mathcal{L}((R_{I,n}-\widehat R_{I,n})_1|\mathcal{Q}_{I,n})$, $\mathcal{L^*}((R_{I,n}^+-\widehat R_{I,n}^++)|\mathcal{Q}_{I,n}^+=\mathcal{Q}_{I,n})=\mathcal{L}((R_{I,n}-\widehat R_{I,n})_1|\mathcal{Q}_{I,n})$, 
respectively, for the original Mack bootstrap (oMB), the alternative Mack bootstrap (aMB) and the intermediate Mack bootstrap (iMB) for different parametric families of distributions of $\F^*$ for $i+j\geq I$, for $I=10$ and different $n$ in Setup b)}\label{tab_ks_predroot_b}
\end{table}
\end{small}

\begin{small}
\begin{table}[t]
\begin{tabular}{|c|c|lll|lll|lll|}
\hline
\multicolumn{2}{|c|}{\textbf{chosen distribution}}                      & \multicolumn{3}{c}{\cellcolor[HTML]{D9D9D9}\textbf{gamma}}                & \multicolumn{3}{c}{\cellcolor[HTML]{D9D9D9}\textbf{log-normal}}           & \multicolumn{3}{c|}{\cellcolor[HTML]{D9D9D9}\textbf{trunc. normal}}        \\
\hline 
\multicolumn{1}{|l|}{\textbf{true distribution}}                   & n  & \multicolumn{1}{c}{oMB} & \multicolumn{1}{c}{aMB} & \multicolumn{1}{c|}{iMB} & \multicolumn{1}{c}{oMB} & \multicolumn{1}{c}{aMB} & \multicolumn{1}{c|}{iMB} & \multicolumn{1}{c}{oMB} & \multicolumn{1}{c}{aMB} & \multicolumn{1}{c|}{iMB} \\
\hline
\cellcolor[HTML]{D9D9D9}                                         & 0  & 9.881                   & 9.874                  & 9.961                  & 9.925                   & 9.850                  & 9.919                  & 9.841                   & 9.822                  & 9.849                  \\
\cellcolor[HTML]{D9D9D9}                                         & 10 & 9.644                   & 9.582                  & 9.680                  & 9.722                   & 9.650                  & 9.742                  & 9.414                   & 9.328                  & 9.442                  \\
\cellcolor[HTML]{D9D9D9}                                         & 20 & 9.479                   & 9.317                  & 9.476                  & 9.459                   & 9.362                  & 9.477                  & 9.254                   & 9.118                  & 9.265                  \\
\cellcolor[HTML]{D9D9D9}                                         & 30 & 8.799                   & 8.623                  & 8.754                  & 8.757                   & 8.598                  & 8.692                  & 9.194                   & 8.871                  & 8.648                  \\
\multirow{-5}{*}{\cellcolor[HTML]{D9D9D9}\textbf{gamma}}         & 40 & 8.452                   & 8.449                  & 8.544                  & 8.513                   & 8.431                  & 8.534                  & 8.706                   & 8.556                  & 8.589                  \\
\hline
\cellcolor[HTML]{D9D9D9}                                         & 0  & 9.990                   & 9.914                  & 9.982                  & 9.983                   & 9.893                  & 9.961                  & 9.868                   & 9.857                  & 9.827                  \\
\cellcolor[HTML]{D9D9D9}                                         & 10 & 9.831                   & 9.752                  & 9.803                  & 9.652                   & 9.635                  & 9.682                  & 9.626                   & 9.537                  & 9.636                  \\
\cellcolor[HTML]{D9D9D9}                                         & 20 & 9.457                   & 9.399                  & 9.468                  & 9.342                   & 9.243                  & 9.354                  & 9.339                   & 9.292                  & 9.355                  \\
\cellcolor[HTML]{D9D9D9}                                         & 30 & 8.959                   & 8.933                  & 8.998                  & 8.771                   & 8.764                  & 8.775                  & 8.712                   & 8.654                  & 8.658                  \\
\multirow{-5}{*}{\cellcolor[HTML]{D9D9D9}\textbf{log-normal}}    & 40 & 8.656                   & 8.584                  & 8.643                  & 8.348                   & 8.325                  & 8.346                  & 8.249                   & 8.184                  & 8.199                  \\
\hline
\cellcolor[HTML]{D9D9D9}                                         & 0  & 9.830                   & 9.789                  & 9.843                  & 9.894                   & 9.876                  & 9.997                  & 9.881                   & 9.874                  & 9.910                  \\
\cellcolor[HTML]{D9D9D9}                                         & 10 & 9.520                   & 9.513                  & 9.524                  & 9.676                   & 9.583                  & 9.517                  & 9.543                   & 9.474                  & 9.575                  \\
\cellcolor[HTML]{D9D9D9}                                         & 20 & 9.167                   & 9.156                  & 9.234                  & 9.234                   & 9.227                  & 9.234                  & 9.344                   & 9.323                  & 9.345                  \\
\cellcolor[HTML]{D9D9D9}                                         & 30 & 8.745                   & 8.692                  & 8.698                  & 8.672                   & 8.660                  & 8.687                  & 8.683                   & 8.630                  & 9.143                  \\
\multirow{-5}{*}{\cellcolor[HTML]{D9D9D9}\textbf{trunc. normal}} & 40 & 8.388                   & 8.356                  & 8.498                  & 8.414                   & 8.376                  & 8.497                  & 8.388                   & 8.321                  & 8.367   \\              
\hline
\end{tabular}\caption{Root of the overall mean of the mean squared error (RMMSE) ($\times 10^{-3}$) for different Mack-type bootstraps, different distributional assumptions and different $n$ and $I=10$ for Setup b)}
\label{tab_mse_predroot_b}
\end{table}
\end{small}

\end{document}